%% file: iwacomp.tex
\documentclass[10pt]{amsart}
\usepackage{amsmath,amsfonts,amscd,amssymb}

\let\oldlabel\label

\def\comment{}
\def\endcomment{}

\csname@addtoreset\endcsname{equation}{section}
\def\theequation{\thesection.\arabic{equation}}
\csname@addtoreset\endcsname{subsection}{equation}
\let\oldsubsection\subsection
\def\subsection{\stepcounter{equation}\oldsubsection}

\theoremstyle{plain}    
\newtheorem{theorem}[equation]{Theorem}
\newtheorem{corollary}[equation]{Corollary}
\newtheorem{lemma}[equation]{Lemma}
\newtheorem{proposition}[equation]{Proposition}
\theoremstyle{definition}
\newtheorem{remark}[equation]{Remark}
\newtheorem{hypothesis}[equation]{Hypothesis}
\newtheorem{example}[equation]{Example}
\newtheorem{definition}[equation]{Definition}
\newtheorem{conjecture}[equation]{Conjecture}

\newcounter{tabcounter}
\def\tabname#1#2{\def\thetabcounter{#1-#2}{\bf Table #1-#2}%
   \refstepcounter{tabcounter}\oldlabel{#1-#2}}

\input cyracc.def       

\def\sha{\text{\tencyr\cyracc{Sh}}}
\font\eightcyr=wncyr8

\def\smallsha{\text{\eightcyr\cyracc{Sh}}}

\def\vQ{{\mathbb Q}/{\mathbf }}
\def\vQ{{\mathbb Q}}
\def\F{{\mathbb F}}
\def\vQl{{\mathbb Q}_l}
\def\vQp{{{\mathbb Q}_p}}
\def\vQb{\bar{\mathbb Q}}
\def\vQlb{{\bar{\mathbb Q}_l}}
\def\vQpb{{\bar{\mathbb Q}_p}}
\def\Finfty{{F_\infty}}
\def\vQGL{\Q(E[5^\infty])}
\def\vZ{{\mathbb Z}}
\def\vZp{{\mathbb Z}_p}

\def\vR{{\mathbb R}}
\def\vC{{\mathbb C}}

\def\Lvmp{L_{\scriptscriptstyle v\!\not\>|\!mp}}
\def\Lvten{L_{\scriptscriptstyle v\!\not\>| 10}}
\let\Q\vQ
\let\Z\vZ

\def\vt{\otimes}

\def\cC{{\mathcal C}}
\def\cD{{\mathcal D}}
\def\vO{{\mathcal O}}
\let\cO\vO
\def\vA{{\mathcal A}}
\let\BigO O
\def\cL{{\mathcal L}}

\def\gna{({\textstyle\frac{n\pi^d}{\sqrt{N}}})}

\def\ufootnote#1{\insert\footins{\noindent\footnotesize{{\hskip
  1.5em}\llap{${}^{\vphantom a}$}#1}}}

\def\newmathop#1{\expandafter\gdef\csname #1\endcsname{\mathop{\rm #1}\nolimits}}
\def\newvmathop#1{\expandafter\gdef\csname v#1\endcsname{\mathop{\rm #1}\nolimits}}

\let\vdim\dim
\let\vdet\det

\newvmathop{Hom}
\newvmathop{ord}
\newvmathop{Re}
\newvmathop{GL}
\newvmathop{Frob}
\newvmathop{PGL}
\newvmathop{Disc}
\newvmathop{Res}
\newvmathop{Tr}
\def\vNorm{N} 
\newvmathop{if}
\newvmathop{iff}
\newvmathop{and}
\newvmathop{or}
\newvmathop{for}
\newvmathop{sign}
\newvmathop{rk}
\let\vrank\vrk
\newvmathop{corank}
\newvmathop{coker}
\newvmathop{codim}
\newvmathop{Ind}
\newvmathop{Tor}
\newvmathop{Res}
\newvmathop{tw}
\newmathop{Gal}
\newmathop{Sel}
\newvmathop{Gal}
\newvmathop{Sel}
\newvmathop{Ak}

\def\vTRIV{{\mathbf 1}}
\def\vmod{\quad {\rm mod}\,}

\let\lar\longrightarrow
\let\iso\cong
\let\tensor\otimes
\def\notdiv{\hbox{$\not|\,$}}

\def\XE#1{{X(E/#1^{cyc})}}
\def\XEna#1{{X(E/#1)}}
\def\Selpi{\Sel_{p^{\infty}}}
\def\Selfivei{\Sel_{5^{\infty}}}
\def\Selthreei{\Sel_{3^{\infty}}}
\def\Selpn{\Sel_{p^n}}
\def\E{E_{21}}
\def\tildeE{\tilde{E}_{21}}

\def\KL{F}

\def\n#1{n(#1)}

\def\fM{{\mathfrak M}}

\def\refAppB{Appendix B}

\def\bigk{L}

\begin{document}

\hskip-10cm{{\tiny.}}

\vskip-1cm

\title{Computations in non-commutative Iwasawa theory}
\author{T. Dokchitser, V. Dokchitser}
\date{13 September 2005}
\ufootnote{2000 {\it Mathematics Subject Classification\/}: Primary 11G05; Secondary 11G40, 11R23}
\address{Tim and Vladimir Dokchitser\vskip 0mm
Department of Pure Mathematics and Mathematical Statistics\vskip 0mm
Centre for Mathematical Sciences\vskip 0mm
University of Cambridge\vskip 0mm
Wilberforce Road\vskip 0mm
Cambridge CB3 0WB\vskip 0mm
United Kingdom}
\email{t.dokchitser@dpmms.cam.ac.uk, v.dokchitser@dpmms.cam.ac.uk}


\begin{abstract}
We study special values of $L$-functions of elliptic curves over $\Q$
twisted by Artin representations that factor through a false Tate
curve extension $\Q(\mu_{p^\infty},\sqrt[p^\infty]m)/\Q$.
In this setting,
we explain how to compute $L$-functions and the corresponding
Iwasawa-theoretic invariants of non-abelian twists of elliptic curves.
Our results provide both theoretical and computational
evidence for the main conjecture of non-commutative Iwasawa theory.
\end{abstract}

\maketitle

\tableofcontents

\section{Introduction}
\label{Sintro}

\endcomment


Let $E/\Q$ be an elliptic curve and let $F_\infty=\bigcup_n F_n$
be an infinite Galois extension of $\Q$. The type of questions that we are
concerned with is how do the Mordell-Weil group and the Selmer group
of $E$ change in the subfields of $F_\infty$, and
the conjectural relations to the special values of twisted
$L$-functions $L(E,\tau,s)|_{s=1}$ for Artin representations $\tau$ that
factor through $\Gal(F_\infty/\Q)$.

Although for an arbitrary algebraic extension of $\Q$ the questions
might be hopelessly hard, a good deal is known for the $p$-adic
extension of $\Q$ obtained by adjoining all $p$-power roots of unity,
$\Q(\mu_{p^\infty})=\bigcup_n\Q(\mu_{p^n})$. This is the
subject of cyclotomic Iwasawa theory of elliptic curves.
%
%
We will consider
the false Tate curve extension of $\Q$,
given by
$$
  F_n=\Q(\mu_{p^n},\sqrt[p^n]m), \qquad
  \Finfty=\bigcup_{n=1}^\infty F_n
$$
for some fixed odd prime $p$ and $p$-power free integer $m$.
The motivation is that this is the simplest non-abelian $p$-adic Lie
extension to which the conjectures of non-commutative Iwasawa theory apply.

Assume that $E$ has good ordinary reduction at $p$.
For Artin representations $\tau$ of $\Gal(\Finfty/\Q)$,
the authors of \cite{CFKSV,FK} have proposed precise 
modifications $\cL_E(\tau)$ of the $L$-values $L(E,\tau,1)$, that are
supposed to be interpolated by a non-abelian $p$-adic $L$-function
of $E$. For self-dual Artin representations $\tau$ with values in
$\vGL_n(\vQp)$, their Main Conjecture asserts that
\begin{equation}\label{conMC}
  \vord_p\cL_E(\tau)=\vord_p\chi_{na}(E,\tau),
\end{equation}
where $\chi_{na}(E,\tau)$ is the non-abelian Euler characteristic of
the dual of the $p^\infty$-Selmer group of $E/\Finfty$ twisted by $\tau$
(see \S\ref{ssnonabeuler}).
Moreover, the conjectures imply that if $\tau$ is congruent to $\tau'$
modulo $p$, then so are the $p$-adic numbers
$\cL_E(\tau)$ and $\cL_E(\tau')$.

The purpose of this paper is to provide some evidence for these conjectures.
We will be concerned with two particular 
representations $\sigma$ and $\rho$ of $G$ which are congruent modulo $p$.
First, under suitable hypotheses (see Theorem \ref{main}) we show that
\begin{equation}\label{TheThm}
  \vord_p\chi_{na}(E,\sigma)=0 \quad\iff\quad\vord_p\chi_{na}(E,\rho)=0,
\end{equation}
as predicted by the conjectures. We also show that the corresponding
statement for $\vord_p\cL_E(\sigma)$ and $\vord_p\cL_E(\rho)$ implied by
the Main Conjecture also follows from the conjecture of Birch and
Swinnerton-Dyer (see \S\ref{ssBSD}). Finally, we test the full congruence
between $\cL_E(\sigma)$ and $\cL_E(\rho)$ numerically for varying
$E$, $p$ and $m$.

Let us begin with a description of the Artin representations of
$G=\Gal(\Finfty/\Q)$ (see, e.g., \cite{vd} for details).
For $n\ge1$, let $\rho_n$ denote the representation of $G$ obtained
by inducing any character of exact order
$p^n$ of $\Gal(F_n/\Q(\mu_{p^n}))$ to $\Gal(F_n/\Q)$.
Then $\rho_n$ is irreducible, and every irreducible Artin representation
of $G$ is of the form $\chi$ or $\rho_n\chi$ for some $n$, where
$\chi$ is a 1-dimensional character of $\Gal(\Q(\mu_{p^\infty})/\Q)$.
In particular, let $\sigma_n$ denote the sum of all 1-dimensional characters
of $\Gal(\Q(\mu_{p^n})/\Q)$. For $n\ge 1$, both
$\sigma_n$ and $\rho_n$ are defined over $\Q$ and, moreover, the two
representations are congruent modulo $p$ (meaning that their reductions
modulo $p$ have isomorphic semi-simplifications). For simplicity, we
will denote $\sigma=\sigma_1$ and $\rho=\rho_1$.

Now, for any elliptic curve $E/\Q$ and any Artin representation $\tau$,
the twisted $L$-function $L(E,\tau,s)$ is an analytic function
for $\vRe s>3/2$ that is conjectured to have an analytic continuation to $\vC$.
In fact, for Artin representations of $G$, this conjecture follows
from modularity of elliptic curves over $\Q$ 
and cyclic base change from the theory of automorphic
forms (\cite{vd}, Theorem 14).

Deligne's conjecture \cite{DelV} then states that the quotient
of $L(E,\tau,1)$ by the period
${\Omega_+(E)^{d^+(\tau)}\Omega_-(E)^{d^-(\tau)}}$ is an algebraic number
(see the list of notation \S\ref{intronot}).
This is also true for Artin representations of $G$ by the
results of \cite{BD}. There it is also shown that
further modifying the value by the global $\epsilon$-factor of $\tau$,
$$
  L^*(E,\tau) = \frac{L(E,\tau,1)}{\Omega_+(E)^{d^+(\tau)}\Omega_-(E)^{d^-(\tau)}}\>\epsilon(\tau^{-1})
$$
lies in the field of definition of $\tau$. In particular,
$L^*(E,\sigma_n)$ and $L^*(E,\rho_n)$ are rational numbers for $n\ge 1$.




Assume that $E/\Q$ has good ordinary reduction at $p$. For
the false Tate curve extension $\Finfty/\Q$,
the following modification of the $L$-value $L(E,\tau,1)$ is proposed
in \cite{CFKSV}:
\begin{equation}\label{conPL}
  \cL_E(\tau) = \frac{\Lvmp (E,\tau,1)}{\Omega_+(E)^{d^+(\tau)}\Omega_-(E)^{d^-(\tau)}}
                     \epsilon_p(\tau) \cdot\frac{P_p(\tau^*,u^{-1})}{P_p(\tau,w^{-1})}\cdot u^{-\n{\tau}}\>.
\end{equation}
Here $\Lvmp (E,\tau,s)$ is the $L$-function $L(E,\tau,s)$ with local factors
at primes dividing $m$ and $p$ removed. (These are the primes with infinite
ramification in $\Finfty/\Q$.) See \S\ref{intronot} for the definitions
of the other quantities in \eqref{conPL}.
Then $\cL_E(\tau)\in\vQpb$ and,
furthermore, one has the following conjectures
(see \cite{CFKSV}, 5.6-5.10 and the penultimate paragraph
of the paper). The definition of the non-abelian Euler characteristic
$\chi_{na}(E,\tau)$ will be given in \S\ref{ssnonabeuler}.

\begin{conjecture}\label{maincon}
Let $E$ be an elliptic curve over $\Q$ with good ordinary reduction at $p$
and let $\tau: G\to\vGL_n(\vQp)$ be a self-dual Artin representation.
Then $L(E,\tau,1)\ne 0$ if and only if $\chi_{na}(E,\tau)$ is finite.
In this case, equation \eqref{conMC} holds.
\end{conjecture}
Another consequence of the conjectures is that 
\begin{equation}\label{TheCong}
  \cL_E(\sigma_n) \equiv \cL_E(\rho_n) \>\>\mod p\>.
\end{equation}
(This follows from the fact that the representations are congruent
modulo $p$ and the conjectural integrality of the non-abelian
$p$-adic $L$-function, see Remark \ref{integrality}.)
We can verify \eqref{TheCong} numerically in many cases
(see the Tables in Appendix B) for $n=1$ and $p=3, 5$ and $7$.
For $p\!>\!7$ or $n\!>\!1$, our approach is too computationally demanding.

In fact, this type of congruence appears in the unpublished work of
Balister \cite{Bal}, which was the original motivation for our study
of these $L$-values. We also note that for $\sigma$ and $\rho$
(i.e. $n=1$) and $p=3$,
Bouganis \cite{Boug} has proved the congruence \eqref{TheCong} in some cases.

There is a remarkable recent result due to Kato \cite{KatK} that
the conjectured existence and integrality of the $p$-adic
$L$-function implies
congruences to a higher power $p^k$ for certain
Artin representations of $G$.
It would be very
interesting to test these higher
congruences, at least numerically. Unfortunately, the conductors
of the relevant twists of $E$ are extremely large.

For $p=3$, if $\cL_E(\sigma)$ and $\cL_E(\rho)$ are $p$-adic units, then
the congruence between them follows from the Birch--Swinnerton-Dyer conjecture.
The point is that $\cL_E(\sigma)$ and $\cL_E(\rho)$ can be related to the
Birch--Swinnerton-Dyer quotient for $E/\Q$, $E/\Q(\mu_p)$ and $E/\Q(\sqrt[p]m)$.
This gives a formula for the quotient $\cL_E(\sigma)/\cL_E(\rho)$, where
the only difficult terms are the orders of the Tate-Shafarevich groups.
As these are squares, for $p=3$ they do not affect the congruence modulo 3.
For $p\ge 5$ this argument gives the congruence modulo squares in $\F_p$,
but the full congruence appears to be stronger than
what can be deduced from the Birch--Swinnerton-Dyer conjecture.

The structure of the paper is as follows.

In \S\ref{Stwistedl} and \S\ref{Scyciwec} we fix the notation and review
some basic results concerning
systems of $l$-adic representations associated to twists of elliptic
curves, Selmer groups and the cyclotomic Iwasawa theory of elliptic curves.
Also, in \S\ref{SSdivcharelt} we prove an auxiliary result
concerning characteristic elements that we use later on
(Proposition \ref{divcharelt}, Corollary \ref{chardivcor}).

Then we turn to false Tate curve extensions and non-commutative
Iwasawa theory of elliptic curves in \S\ref{sec-nonab}.
We prove the relation \eqref{TheThm} in Theorem \ref{main}.
From here, the Main Conjecture
implies the corresponding relation between the values
$\cL_E(\sigma)$ and $\cL_E(\rho)$ of the $p$-adic $L$-function, that
we show to be compatible with the Birch--Swinnerton-Dyer conjecture
(Proposition \ref{bsd-pred}).

Next, we turn to numerical verification of the congruence
\eqref{TheCong} for $n=1$. In \S\ref{Smethod} we explain how to compute
$L(E,\sigma,1)$, $L(E,\rho,1)$ and various related arithmetic invariants.
We give an example for the elliptic curve 21A4 over $\Q(\mu_5,\sqrt[5]2)$
illustrating both our results and computations in \S\ref{Sexample}, and
tabulate all other examples in \refAppB. In \S\ref{gl2ex} we also
illustrate how similar computations can be carried out when the false
Tate curve extension is replaced by a $\vGL_2$ extension,
by taking the curves 11A2 and 11A3 over their 5-division fields.

Appendix A by J. Coates and R. Sujatha deals with theoretical questions
related to regularity and
Heegner like phenomena for the extension $\Finfty/\Q$.
In particular, they generalise \eqref{TheThm} to a larger class of curves
and $\rho_n$ in place of $\rho$ (Corollary \ref{CorA19}).

\subsection{Notation}
\label{intronot}

By a number field we will always mean a finite extension of $\Q$.
We fix, once and for all, embeddings of $\vQb$ into $\vQlb$ for all $l$.
Throughout the paper we keep the following notation:
\smallskip

\begin{tabular}{ll}
$p$                  & odd prime. \cr
$\mu_{p^n}$          & group of $p^{n}$-th roots of unity. \cr
$\mu_{p^{\infty}}$   & union of $\mu_{p^n}$ for $n\ge 1$. \cr
$k_v$                & completion of a field $k$ at a prime $v$,\cr
$\F_v$               & residue field at a prime $v$,\cr
$\tau^*$             & contragredient representation of an Artin representation $\tau$,\cr
$d^\pm(\tau)$        & dimensions of the $\pm 1$ eigenspaces of complex conjugation on $\tau$,\cr
$\epsilon(\tau)$     & the global $\epsilon$-factor of $\tau$,\cr
$\epsilon_p(\tau)$   & the local $\epsilon$-factor of $\tau$ at $p$ (see \S\ref{Sepsilon}),\cr
$\n{\tau}$           & the $p$-valuation of the conductor $N(\tau)$ (see \S\ref{ssCond}),\cr
$|S|$                & cardinality of $S$,\cr
$A[p]$               & $p$-torsion of an abelian group $A$,\cr
$A[p^\infty]$        & $p$-primary component of an abelian group $A$.\cr
\end{tabular}

\smallskip
\noindent
Notation for the false Tate curve extension and its representations:
\smallskip

\begin{tabular}{ll}
$m$      & an integer $\ge 2$ such that $n^p\notdiv m$ for $n>1$,\cr
$F$      & $=F_1=\vQ(\mu_p,\sqrt[p]{m})$,\cr
$K$      & $=\vQ(\mu_p)$,\cr
$\sigma$ & $=\sigma_1$, the regular representation of $\vGal(K/\vQ)$,\cr
$\rho$   & $=\rho_1$, the unique $(p-1)$-dimensional irreducible representation\cr
         & of $\vGal(F/\vQ)$, as above. \cr
\end{tabular}

\smallskip
\noindent
Notation relating to elliptic curves and Iwasawa theory:
\smallskip

\begin{tabular}{ll}
$E$                 & elliptic curve over $\vQ$ with good ordinary reduction at $p$, \cr
$P_q(-,T)$           & local polynomial of the $L$-series $L(-,s)$ at $q$ (see \S\ref{SSlfun}),\cr
$u, w$               & $p$-adic numbers, chosen so that $u$ is a $p$-adic unit, defined by \cr
                     & $P_p(E,T)=1\!-\!a_p T\!+\!pT^2 = (1\!-\!uT)(1\!-\!wT)$. \cr
\end{tabular}

\noindent
We denote by $\Omega_\pm(E)$
the periods of $E$, defined by integrating the N\'eron differential
of a global minimal Weierstrass equation over the generators
of the $\pm1$-eigenspaces $H_1(E(\vC),\Z)^\pm$ of complex conjugation.
As usual, $\Omega_-$ is chosen to lie in $i\vR$.

\subsection*{Acknowledgements}
We would like to thank J. Coates for numerous discussions and comments,
and for his constant encouragement.
We also thank J. Coates, T.~Fukaya, K. Kato, R. Sujatha and
O. Venjakob who inspired the creation of this paper by referring to it
in \cite{CFKSV} and \cite{KatK}.
The first author expresses his gratitude to the University of Edinburgh,
where part of this research was carried out.

\comment

\section{Twisted $L$-functions}
\label{Stwistedl}


We briefly recall the definition of $L$-functions of elliptic curves
by Artin representations, and the invariants attached to them
which we will need.
We refer to \cite{SerF} and \cite{SerZ} for
$\zeta$- and $L$-functions of varieties,
and to \cite{DelV}, \cite{DelC} \S\S3-4 and \cite{TatN} \S4
for their twists.

\subsection{Systems of $l$-adic representations}

Let $E/\Q$ be an arbitrary elliptic curve and
$\tau: \Gal(\vQb/\Q)\to\vGL_n(\vQb)$ be an Artin representation. Both $E$
and $\tau$ determine (a compatible system of) $l$-adic representations
for primes $l$ of $\Q$. In case of $E$, the $l$-adic representation is
$M_l(E)=H^1_{et}(E,\Z_l)\tensor_{\Z_l}\bar\Q_l$ or, equivalently, the dual of the
$l$-adic Tate module $T_l(E)$ with scalars extended to~$\vQlb$.
The $l$-adic representation that corresponds to $\tau$
is~$M_l(\tau)=\tau\tensor\vQlb$. Now we can
construct a system of representations
$$
  M_l(E,\tau) = M_l(E)\tensor_{\vQlb} M_l(\tau) \>.
$$

\subsection{$L$-functions}
\label{SSlfun}

To a system of $l$-adic representations $M=\{M_l\}_l$ we can associate an $L$-function
$L(M,s)$ as follows. For a prime $q$ of $\Q$, the local polynomials of $L(M,s)$ are
\begin{equation}
\label{locpol}
  P_q(M,T) = \det\bigl(1-\vFrob_q^{-1} T\bigm|M_l^{I_q}\bigr)
\end{equation}
for any prime $l\ne q$. We define the local $L$-factor
$$
  L_q(M,s) = P_q(M,q^{-s})^{-1}
$$
and the global $L$-function (or $L$-series)
$$
  L(M,s) = \prod_q L_q(M,s) \>.
$$
We write
$$
  L(E,s) = L(M(E),s), \quad L(\tau,s) = L(M(\tau),s), \quad
  L(E,\tau,s) = L(M(E,\tau),s)\>.
$$
The $L$-series $L(\tau,s)$ converges to an analytic function on the
half-plane $\vRe s>1$, and it is classical that it has
a meromorphic continuation to the entire complex plane.
The $L$-series $L(E,s)$ and $L(E,\tau,s)$
converge to analytic functions on $\vRe s>3/2$.
It follows from the modularity of elliptic curves \cite{Wil, TW, BCDT}
that $L(E,s)$ possesses an analytic continuation to $\vC$.
We will be concerned with $L(E,\tau,s)$ when $\tau$ factors
through a false Tate curve extension, in which case there is also an
analytic continuation (\cite{vd}, Theorem 14).

Finally, recall that $L$-functions are multiplicative,
$$
  L(E,\tau_1\oplus\tau_2,s) = L(E,\tau_1,s) L(E,\tau_2,s) \>.
$$

\subsection{Conductors}
\label{ssCond}

In the formula \eqref{locpol} defining
the local polynomials of $L(M,s)$, the inertia group $I_q$ acts trivially
on $M_l$ for all but finitely many primes $q$. Therefore for almost all primes
$q$ (the {\em good\/} primes of $M$) we have $\deg P_q(M,T)=\dim M_l$.
The primes $q$ with $\deg P_q(M,T)<\dim M_l$ are called {\em bad\/} primes
of $M$.

For a prime $q$ the local conductor $N_q(M)$ is given by (see \cite{SerF}, \S2)
$$
  N_q(M) = q^{t_q+\delta_q}, \qquad t_q=\vcodim_{M_l}M_l^{I_q}, \quad \delta_q\ge 0.
$$
The term $\delta_q$ is defined in terms of the representation of the wild
inertia subgroup on $M_l$. It is zero if and only if this subgroup acts
trivially, i.e. $M_l$ is tamely ramified.
Finally, the (global) conductor is $N(M)=\prod_q N_q(M)$.
Once again, we write
$$
  N(E) = N(M(E)), \quad N(\tau) = N(M(\tau)), \quad
  N(E,\tau) = N(M(E,\tau))\>.
$$

\subsection{Functional equation}

The twisted $L$-functions $L(E,\tau,s)$ conjecturally satisfy a functional
equation of the following form (see, e.g, \cite{DelV} 5.2, \cite{TatN} 4.5).
Let
$$
  \hat L(E,\tau,s) = \Bigl(\frac{N(E,\tau)}{\pi^{2\dim\tau}}\Bigr)^{s/2}
     \Gamma\Bigl(\frac s2\Bigr)^{\dim\tau} \Gamma\Bigl(\frac{s+1}2\Bigr)^{\dim\tau} L(E,\tau,s) \>.
$$
Then, conjecturally,
\begin{equation}\label{funeq}
  \hat L(E,\tau,s) = w(E,\tau) \hat L(E,\tau^*,2-s)
\end{equation}
with $w(E,\tau)$ an algebraic number of complex absolute value 1. If
$\tau\iso\tau^*$, then $w(E,\tau)=\pm 1$ and we call it the sign in the
functional equation.

\subsection{$L$-functions of elliptic curves over number fields}
\label{SSLfEcNf}

Recall that for an elliptic curve $E$ a number field $k$, the $L$-function
$L(E/k,s)$ is given by the Euler product
$$
  L(E/k,s) = \prod_v L_v(E/k,s) = \prod_v P_v(E/k,\vNorm_{k/\Q}(v)^{-s})^{-1} \>,
$$
where $v$ runs over primes of $k$ and
$$
  P_v(E/k,T) = \det\bigl(1-\vFrob_v^{-1} T\bigm|M_l(E)^{I_v}\bigr), \qquad v\notdiv l\>.
$$

The polynomials $P_v(E/k,T)$ depend on the reduction type of $E$
over the local field $k_v$ and are given explicitly by:
\begin{equation}\label{pvclass}
P_v(E/k,T) = \left\{
   \begin{array}{ll}
     1 \!-\! a_v T \!+\! \vNorm_{k/\Q}(v) T^2, & {\rm good\ reduction,} \cr
     1 - T, & {\rm split\ multiplicative\ reduction,} \cr
     1 + T, & {\rm non\hbox{-}split\ multiplicative\ reduction,} \cr
     1,     & {\rm additive\ reduction.} \cr
   \end{array}
 \right.
\end{equation}
In the case of good reduction, $P_v(E/k,1)=|\tilde E(\F_v)|$,
the number of points on the reduction of $E$ at $v$.

Suppose the elliptic curve $E$ is defined over $\Q$ and let $k/\Q$ be a
finite extension. Define $R_k$ to be the representation of $\Gal(\vQb/\Q)$
induced from the trivial representation of $\Gal(\vQb/k)$.
In particular, if $k/\Q$ is Galois, $R_k$ is the regular representation
of $\Gal(k/\Q)$ considered as a representation of $\Gal(\vQb/\Q)$. Also,
for $k$ arbitrary, $R_k$ is the permutation representation of $\Gal(\vQb/\Q)$
on the set of embeddings of $k$ into $\vQb$.
One can show that the local factors satisfy Artin formalism,
$$
  L_q(E, R_k, s) = \prod_{v|q} L_v(E/k, s) \>.
$$
In particular, taking the product over all primes,
$$
  L(E, R_k, s) = L(E/k, s) \>.
$$
This allows us to write $L(E/k, s)$ as a product of twisted $L$-functions
of $E/\Q$ that correspond to irreducible pieces of $R_k$.
If $k/\Q$ is abelian, $L(E/k, s)$ becomes a product of $L$-functions of
$E$ twisted by Dirichlet characters.

\section{Cyclotomic Iwasawa theory of elliptic curves}
\label{Scyciwec}

In this section we introduce the main tools from cyclotomic Iwasawa theory
of elliptic curves. We use the following notation

\begin{tabular}{ll}
$p$            & odd prime.\cr
$k$            & a number field. \cr
$k^{cyc}$      & the maximal pro-$p$ extension of $k$ in $k(\mu_{p^{\infty}})$. \cr
$\Gamma_k$     & the Galois group $\Gal(k^{cyc}/k)$.  \cr
$\Lambda_k$    & Iwasawa algebra of $\Gamma_k$.
\end{tabular}

\noindent
Recall that $\Lambda_k=\varprojlim_H
\vZp[\Gamma_k/H]$ where the limit is taken over the open subgroups $H$
of $\Gamma_k$. So it is simply the completion of the group ring
$\vZp[\Gamma_k]$ in the profinite topology.
The group $\Gamma_k$ is isomorphic to $\vZp$ and there is an isomorphism
$$
  \Lambda_k\iso\vZp[[T]]
$$
obtained by sending a fixed topological generator of $\Gamma_k$ to $1+T$.

We recall the definition of the $p^n$- and $p^\infty$-Selmer groups
of an elliptic curve $E$ over an arbitrary algebraic extension of $\Q$.
(These definitions do not require that $p$ is odd or that $E$ is
good ordinary at $p$.)

If $k$ is a number field, we have an exact sequence
$$
  0 \lar \Selpn(E/k) \lar H^1(k,E[p^n]) \lar \prod_v H^1(k_v,E(\bar{k_v}))[p^n]
$$
where $v$ runs over the places of $k$ (\cite{Sil1}, \S X.4). For arbitrary
$\bigk\subset\vQb$, let
$$
  \Selpn(E/\bigk) = \varinjlim_{k\subset \bigk} \Selpn(E/k)
$$
where $k$ runs over finite extensions of $\vQ$ and the maps that define the
inductive system are induced by the restriction maps on cohomology.
Finally, define the $p^\infty$-Selmer group by
$$
  \Selpi(E/\bigk) = \varinjlim_{n\ge 1} \Selpn(E/\bigk) \>.
$$
Over a number field, the group $\Selpi(E/k)$
is a part of the fundamental exact sequence
\begin{equation}
\label{selexactseq}
  0 \lar E(k)\tensor_{\Z}(\Q_p/\Z_p) \lar \Selpi(E/k) \lar \sha(E/k)[p^\infty] \lar 0\>,
\end{equation}
where $\sha(E/k)$ is the Tate-Shafarevich group of $E/k$.


The Selmer group $\Selpi(E/k^{cyc})$ is a discrete $\Lambda_k$-module.
Let $\XE{k}$ be its Pontryagin dual,
$$
  \XE{k} = \vHom(\Selpi(E/k^{cyc}),\vQp/\vZp) \>,
$$
where $\vQp/\vZp$ is a $\Lambda_k$-module with trivial $\Gamma_k$-action.
This is a compact $\Lambda_k$-module and
it is an elementary fact in Iwasawa theory that $\XE{k}$ is
finitely generated over $\Lambda_k$ (see, e.g., \cite{Man} Thm. 4.5.(a)
or \cite{CS}).

\subsection{Torsion modules}

A $\Lambda_k$-module is said to be {\em torsion\/} if for each element
of the module there is some non-zero element in $\Lambda_k$ annihilating it.

We recall an important result due to Mazur and its corollary which we will
use later on:

\begin{theorem}
\label{mazcon}
{\rm(Mazur's control theorem, \cite{Maz}, \cite{Gre} Theorem 4.1).}
Let $k$ be a number field and $E/k$ and elliptic
curve with good ordinary reduction at all primes above $p$. Assume that
$k_\infty=\bigcup_n k_n$ is a Galois extension of $k$ with Galois
group $\vZ_p$. Then the natural maps
$$
  \Selpi(E/k_n) \lar \Selpi(E/k_\infty)^{\Gal(k_\infty/k_n)}
$$
have finite kernels and cokernels. Their orders are bounded as $n\to\infty$.
\end{theorem}

\begin{theorem}\label{mazcon2}
{\rm(\cite{Maz}, \cite{Gre} Corollary 4.9)}
Let $k$ be a number field and let $E/k$ be an elliptic curve
with good ordinary reduction at all primes above $p$.
If $\Selpi(E/k)$ is finite, then $\XE{k}$ is $\Lambda_k$-torsion.
\end{theorem}

Mazur has conjectured that $\XE{k}$ is always $\Lambda_k$-torsion when
$E$ has good ordinary reduction at all primes $v$ of $k$ dividing $p$.
The best result to date in this direction is the following
deep theorem due to Kato.

\begin{theorem}\label{katotorsion}
{\rm(Kato \cite{KatP})}
If $k/\vQ$ is abelian and $E/\vQ$ has good
ordinary reduction at $p$, then $\XE{k}$ is $\Lambda_k$-torsion.
\end{theorem}

\subsection{Structure theory}
A map of $\Lambda_k$-modules is a {\em pseudo-isomorphism}
(denoted $X\iso_{ps}Y$) if it has finite kernel and cokernel.
There is a well-known structure theory for finitely generated $\Lambda_k$-torsion
modules up to pseudo-isomorphism
(see e.g. \cite{Bour} Ch. 7 or \cite{Was} Theorem 13.12).
For such a module $X$, one has
\begin{equation}\label{eq:structure}
  X \iso_{ps} \bigoplus_i \Lambda_k/{p^{\mu_i}} \oplus \bigoplus_j \Lambda_k/{f_j^{m_j}},
\end{equation}
Here the direct sums are finite, $\mu_i,m_j\ge 1$ are integers and
$f_j$ are elements of $\Lambda_k$. The $f_j$ are defined up to units
in $\Lambda_k$ and they can
be chosen to map to irreducible distinguished polynomials
under the isomorphism $\Lambda_k\iso \vZp[[T]]$.
Recall that {\em distinguished} means monic with all other coefficients
divisible by $p$.
With such a choice of a canonical form for the $f_j$, the above
decomposition is unique up to order. The product
$$
  f_{X} =\prod_i p^{\mu_i} \prod_j f_j^{m_j}
$$
is called a characteristic element for $X$. It is uniquely defined up to
a $\Lambda_k$-unit and the characteristic ideal $(f_X)\subset\Lambda_k$
is well-defined. Characteristic elements are
multiplicative in short exact sequences.


Now let $E/\vQ$ be an elliptic curve with good ordinary reduction at $p$.
Provided that $\XE{k}$ is a $\Lambda_k$-torsion module
(e.g. if $k/\Q$ is abelian), we can decompose $\XE{k}$ as in
\eqref{eq:structure} and define the following invariants:

\begin{tabular}{ll}
$f_{E/k}$        & $=f_{\XE{k}}=\prod_i p^{\mu_i} \prod_j f_j^{m_j},$ the characteristic element.\cr
$\mu_{E/k}$      & $=\sum_{i} \mu_i$, the $\mu$-invariant. \cr
$\lambda_{E/k}$  & $=\deg f_{E/k}=\dim_{\vQp}(\XE{k}\tensor_{\vZp}\vQp)$, the $\lambda$-invariant.\cr
\end{tabular}

\begin{remark}
We note that $\mu_{E/k}=0$ if and only if $\XE{k}$ is a finitely generated
$\Z_p$-module, and $\lambda_{E/k}=\mu_{E/k}=0$ if and only if $\XE{k}$
is finite.
\end{remark}

\begin{remark}
\label{poweroft}
Assuming $\XE{k}$ is a $\Lambda_k$-torsion module, it is easy to see that
$f_{E/k}$ is divisible by $T$ if and only if the coinvariant
space $(\XE{k})_{\Gamma_k}$ is infinite.
By Mazur's control theorem (\ref{mazcon}),
this is true if and only if $\Selpi(E/k)$ is infinite.
In particular, $\Selpi(E/k)$ is finite if and only if $f_{E/k}$ has a
non-zero constant term.
\end{remark}

We also mention the following result.

\begin{theorem}{\rm(Matsuno \cite{Mat})}\label{thM}
Let $p$ be an odd prime.
Let $k$ be a totally imaginary algebraic number field. Let $E/k$ be an
elliptic curve which has good reduction at all primes above $p$. If
$\XE{k}$ is $\Lambda_k$-torsion and $\mu_{E/k}=0$, then $\XE{k}$
is $p$-torsion-free.
In particular, if $\lambda_{E/k}$ is also zero, then $\XE{k}=0$.
\end{theorem}

\subsection{Euler characteristics}
\label{ssCHI}

Let $E$ be an elliptic curve over a number field $k$ and assume that
$\XE{k}$ is $\Lambda_k$-torsion.
When $\Selpi(E/k)$ is finite, the homology groups $H_i(\Gamma_k,\XE{k})$
are finite ($i=0,1$), and the $\Gamma_k$-Euler characteristic is defined by
$$
  \chi_{cyc}(E/k) = \frac{|H_0(\Gamma_k,\XE{k})|}{|H_1(\Gamma_k,\XE{k})|} \>.
$$
(When $\Selpi(E/k)$ is infinite, the groups
are no longer finite, but one can nevertheless make sense
of the quotient and define a ``generalised Euler characteristic'',
see e.g. \cite{CSS}).

Of importance to us are the following two formulae.
First, $\chi_{cyc}(E/k)$ is a power of $p$, and we have the basic Euler
characteristic formula \cite{P-R, Schn}
\begin{equation}\label{eq:chicyc}
  \vord_p\chi_{cyc}(E/k)=
  \vord_p\>\frac{|\sha(E/k)[p^\infty]|\prod_{v|p}|\tilde{E}(\F_v)|^2\prod_v c_v}{|E(k)|^2}\>.
\end{equation}

Second,
under the isomorphism $\Lambda\iso\vZp[[T]]$, 
the constant term of the characteristic element recovers $\chi_{cyc}(E/k)$,
\begin{equation}\label{eq:chichar}
  \vord_p\chi_{cyc}(E/k)=\vord_p f_{E/k}(0)\>.
\end{equation}

\begin{lemma}
\label{chicyc0}
Let $E$ be an elliptic curve over a number field $k$ such that $E$ has
good ordinary reduction at all primes above $p$.
Suppose that $\Selpi(E/k)$ is finite. Then
$\XE{k}$ is finite if and only if $\chi_{cyc}(E/k)=1$.
\end{lemma}

\begin{proof}
By Theorem \ref{mazcon2}, $\XE{k}$ is $\Lambda_k$-torsion.
By \eqref{eq:chichar} we have $f_j\neq T$ for every $j$.
The $f_j$ are distinguished irreducible polynomials, so that every $f_j\neq T$
has a constant term divisible by $p$. Lemma follows from \eqref{eq:chichar}.
\end{proof}

\subsection{$\lambda$- and $\mu$-invariants in $p$-extensions}

The following theorem, due to Hachimori and Matsuno, concerns
the change in the $\lambda$-invariant of $\XE{k}$ when $k$ is replaced by
a finite $p$-extension of $k$. This will be very useful when studying
the behaviour of $E$ in a false Tate curve tower.

\begin{theorem}\label{thHM} {\rm(\cite{HM}, Theorem 3.1)}
Let $k'/k$ be a Galois extension of number fields of $p$-power degree.
Let $E/k$ be an elliptic curve, good ordinary at primes above $p$
and such that whenever $E$ has additive
reduction at a prime $v$, the reduction stays additive at every prime of
$k'^{cyc}$ above $v$.
If $X(E/k)$ is $\Lambda_k$-torsion and $\mu_{E/k}=0$, then $X(E/k')$
is $\Lambda_{k'}$-torsion, $\mu_{E/k'}=0$ and
\begin{equation}\label{hacmat}
 \lambda_{E/k'} = [k'^{cyc}:k^{cyc}] \lambda_{E/k}
    +   \sum_{v\in V_1} (e(v)-1)
    + 2 \sum_{v\in V_2} (e(v)-1) \>.
\end{equation}
Here $e(v)$ is the ramification index of $v$ in $k'^{cyc}/k^{cyc}$,
$V_1$ is the set of primes $v\notdiv p$ of $k'^{cyc}$
where $E$ has split multiplicative reduction
and $V_2$ is the set of primes $v\notdiv p$ of $k'^{cyc}$
where $E$ has good reduction and $E(k'^{cyc}_v)[p]\ne 0$.
\end{theorem}

\begin{remark}\label{HMpge5}
The condition that the reduction stays additive is automatically satisfied
if $p\ge 5$ (\cite{ST}, p. 498), and is vacuous if $E/k$ is semistable.
\end{remark}

It is useful to reformulate some of the conditions in the theorem
as follows:

\begin{lemma}
\label{lemhm}
Let $k'/k$ be a Galois extension of number fields of $p$-power degree.
Let $E/k$ be an elliptic curve such that whenever $E$ has additive
reduction at a prime of $k$, the reduction stays additive at every prime of
$k'$ above it. Let $v\notdiv p$ be a prime of $k$.
\begin{enumerate}
\item\label{ququ1}
$v$ ramifies in $k'$ if and only if
any prime above $v$ in $k^{cyc}$ ramifies in $k'^{cyc}$.
\item\label{ququ2}
$E$ has split multiplicative reduction at $v$ if and only if
$E$ has split multiplicative reduction at any prime of $k'^{cyc}$ above $v$.
\item\label{ququ3}
If $E/k$ has good reduction at $v\notdiv p$ and $w$
is a prime of $k'^{cyc}$ above $v$, then
$E(k'^{cyc}_w)[p]=0$ if and only if $\tilde E(\F_v)[p]=0$.
\end{enumerate}
\end{lemma}

\begin{proof}
First observe that the condition that the reduction stays additive
at every prime of $k'$ above a given prime of $k$ ensures that it also stays
additive over every prime of $k'^{cyc}$ above it.
This follows from the fact that $p$ is the only prime that ramifies
in the $p$-cyclotomic extension of a number field.

(\ref{ququ1}).
This also follows from the fact that $p$ is the only prime that ramifies
in the $p$-cyclotomic extension of a number field.

(\ref{ququ2}).
If $E/k$ has additive reduction, it stays additive in $k'^{cyc}$
by assumption. If the reduction is non-split multiplicative, it stays
the same because the degree of the extension is (pro-)odd.

(\ref{ququ3}).
We show the asserted equivalence in two steps.
First, for a prime $w$ of $k'^{cyc}$ above $v$ apply
the multiplication by $p$ map to the exact sequence
$$
  0 \to \hat E(m_w) \to E(k'^{cyc}_w) \to \tilde E(\F_w) \to 0
$$
where $m_w$ is the maximal ideal and $\F_w$ is the residue field of the
local field $k'^{cyc}_w$.
We get a piece of the kernel-cokernel exact sequence
$$
  \hat E(m_w)[p] \to E(k'^{cyc}_w)[p] \to \tilde E(\F_w)[p] \to \hat E(m_w)/p\hat E(m_w) \>.
$$
Since $[p]$ on $\hat E(m_w)$ is an isomorphism \cite{Sil1}, IV.2.3,
it follows that $E(k'^{cyc}_w)[p] \iso \tilde E(\F_w)[p]$.

It remains to show that $E(\F_v)[p]=0$ if and only if
$E(\F_w)[p]=0$. This follows from Nakayama's lemma: $G=\vGal(\F_w/\F_v)$
is pro-$p$, so for any discrete $p$-primary $G$-module $A=0$ if and only if
$A^G=0$.
\end{proof}

\begin{corollary}
\label{corhm}
Under the assumptions of Theorem \ref{thHM}, we have that
$\lambda_{E/k'}=\lambda_{E/k}$ if
and only if either $k'\subset k^{cyc}$ or the following conditions are
satisfied
\begin{enumerate}
\item $\lambda_{E/k}=0$.
\item There are no primes of split multiplicative reduction of $E/k$
that ramify in $k'/k$.
\item There are no primes $v\notdiv p$ of good reduction of $E/k$
that ramify in $k'/k$ and such that there is a non-trivial point of order p
on the reduced curve $\tilde E/\F_v$.
\end{enumerate}
\end{corollary}

\begin{proof}
This follows from Theorem \ref{thHM} and Lemma \ref{lemhm}.
\end{proof}

\subsection{Divisibility of characteristic elements}
\label{SSdivcharelt}

The main result of this section (Proposition \ref{divcharelt})
concerns the behaviour of the characteristic element $f_{E/k}$
when changing the base field.

\begin{lemma}
\label{h1etale}
Let $k'/k$ be a finite Galois extension with Galois group $\Delta$
and let $E$ be an elliptic curve over $k$. Then the natural restriction map
$$
  H^1(k,E[p^\infty]) \>\>{\buildrel{\rm Res}\over\lar}\>\>H^1(k',E[p^\infty])^\Delta
$$
has finite $p$-power kernel and cokernel. If either $p\notdiv|\Delta|$
or $E/k'$ has trivial $p$-torsion, then {\rm Res\/} is an isomorphism.
\end{lemma}

\begin{proof}
Denote $M=E[p^\infty](k')$. 
The map {\rm Res} is a part of the inflation-restriction sequence
$$
  0
  \lar
  H^1(\Delta,M)
  \lar
  H^1(k,E[p^\infty])
  \>\>{\buildrel{\rm Res}\over\lar}\>\>
  H^1(k',E[p^\infty])^\Delta
  \lar
  H^2(\Delta,M) \>.
$$
The groups $H^i(\Delta,M)$ are annihilated by $|\Delta|$.
Furthermore, they are cofinitely generated continuous
$\vZp$-modules (since $M$ is), and are therefore finite.
If either $p\notdiv|\Delta|$ or $M=0$, these groups are trivial
and $\vRes$ is an isomorphism.
\end{proof}

\begin{proposition}
\label{divcharelt}
Let $k'/k$ be a finite Galois extension of number fields
such that $k'\cap k^{cyc}=k$.
Let $E/k$ be an elliptic curve such that $\XE{k}$ is $\Lambda_k$-torsion
and $\XE{k'}$ is $\Lambda_{k'}$-torsion. Fix an isomorphism
$\Lambda_k\iso\vZ_p[[T]]$. It induces naturally an isomorphism
$\Lambda_{k'}\iso\vZ_p[[T]]$ and, with this identification,
$f_{E/k}$ divides $f_{E/k'}$.
\end{proposition}

\begin{proof}
Denote $\Delta=\vGal(k'/k)$.
By Lemma \ref{h1etale} and the definition of the Selmer group,
we have an exact sequence of $\Gamma_k$-modules
$$
  0 \lar A \lar \Selpi(k^{cyc}) \lar \Selpi(k'^{cyc})^\Delta \lar B \lar 0
$$
with $A$ finite. As Pontryagin dual is exact,
$$
  0 \lar \hat B \lar \XE{k'}/I \lar \XE{k} \lar \hat A \lar 0
$$
for some ideal $I$, so $\XE{k}$ is pseudo-isomorphic to a quotient of
$\XE{k'}/I$ as a $\Lambda_k$-module.

We have two structures of a $\Lambda_k$-module on $\XE{k'}/I$.

First, we have that $\Lambda_{k'}\iso\Lambda_k$ as follows: if
$\sigma\in\Gamma_{k'}=\Gal(k'^{cyc}/k')$, we can consider it as an element of
$\Gal(k'^{cyc}/k)$ and we take it modulo $\Delta=\Gal(k'^{cyc}/k^{cyc})$.
This gives an element of $\Gal(k^{cyc}/k)=\Gamma_k$. This gives an isomorphism
$\Gamma_{k'}\iso\Gamma_k$ since $\Gal(k'^{cyc}/k)$ is a direct product of
$\Delta$ and $\Gamma_{k'}$. Hence $\Lambda_{k'}\iso\Lambda_k$,
inducing the asserted isomorphism $\Lambda_{k'}\iso\vZ_p[[T]]$.
This also allows us to consider $\XE{k'}/I$ as a $\Lambda_k$-module.

Second, $\Sel(k')$ is a $\Gal(k'^{cyc}/k)$-module and $\Sel(k')^\Delta$ is
a $\Gal(k^{cyc}/k)$-module. This gives a $\Lambda_k$-module structure on
its Pontryagin dual $\XE{k'}/I$.

By inspection, these $\Lambda_k$-module structures are the same.
Hence the $\Lambda_k$-characteristic element of $\XE{k'}/I$ can be
identified with its $\Lambda_{k'}$-characteristic element.
By multiplicativity in short exact sequences, the former is divisible by
$f_{E/k}$ while the latter divides $f_{E/k'}$.
This completes the proof.

%
%
%
%
\end{proof}

\begin{corollary}
\label{chardivcor}
Let $k'/k$ be a finite Galois extension of number fields
such that $k'\cap k^{cyc}=k$. Let $E/k$ be an elliptic curve such
that $\XE{k}$ is $\Lambda_k$-torsion and $\XE{k'}$ is $\Lambda_{k'}$-torsion.
Then $\lambda_{E/k'}\ge\lambda_{E/k}$ and $\mu_{E/k'}\ge\mu_{E/k}$.
Moreover, assuming that $\Selpi(E/k)$ is finite, we have
$$
  \begin{array}{l}
    \lambda_{E/k'}=\lambda_{E/k},\cr
    \mu_{E/k'}=\mu_{E/k}
  \end{array}
     \quad\iff\quad
  \begin{array}{l}
    |\Selpi(E/k')|<\infty,\cr
    \chi_{cyc}(E/k')=\chi_{cyc}(E/k)
  \end{array} \>.
$$
\end{corollary}

\begin{proof}
The first statement follows directly from the proposition.
For the second statement, observe that
the proposition implies that $\lambda_{E/k'}=\lambda_{E/k}$
and $\mu_{E/k'}=\mu_{E/k}$ if and only if $f_{E/k'}=f_{E/k}$.
Now the implication ``$\Leftarrow$''
follows from formula \eqref{eq:chichar}. For the reverse implication,
$|\Selpi(E/k')|<\infty$ is a consequence of Remark \ref{poweroft}, and the
equality of Euler characteristics follows from formula \eqref{eq:chichar}.
\end{proof}


\section{False Tate curve extensions}
\label{sec-nonab}


We now turn to non-commutative Iwasawa theory in the setting of a false
Tate curve extension $\Finfty/\vQ$,
with $\Finfty=\vQ(\mu_{p^\infty}, \sqrt[p^\infty]m)$.
As always, $p$ is an odd prime, $m\ge 2$ a $p$-power free integer, and
$E/\vQ$ is an elliptic curve with good ordinary reduction at $p$.
We write

\begin{tabular}{ll}
$G_k$                  & $\Gal(\Finfty/k)$, for a number field $k\subset \Finfty$, \cr
$\Lambda(G_k)$         & Iwasawa algebra of $G_k$, that is $\varprojlim_{H\le G_k,\rm\ open} \vZp[G_k/H]$, \cr
$\XEna{\Finfty}$ & Pontryagin dual of $\Selpi(E/\Finfty)$, as in the cyclotomic case,\cr
$P_1^{(k)}$            & for $k\subset \Finfty$, this is the set of primes $w$ of $k$ with $w|m$, $w\notdiv p$, where\cr
                       & $E$ has split multiplicative reduction,\cr
$P_2^{(k)}$            & for $k\subset \Finfty$, this is the set of primes $w$ of $k$ with $w|m$, $w\notdiv p$, where\cr
                       & $E$ has good reduction and the reduced curve has a point of order $p$.\cr
\end{tabular}

\subsection{Non-abelian Euler characteristic}
\label{ssnonabeuler}

Let $k\subset\Finfty$ be a number field and $G_k=\Gal(\Finfty/k)$.
Similarly to the cyclotomic case, one can define a
non-abelian $G_k$-Euler characteristic of a compact $G_k$-module $X$ by
$$
  \chi(G_k,X) = \prod_{i=0}^2 |H_i(G_k,X)|^{(-1)^i} \>,
$$
provided the above groups are finite. For an elliptic curve $E/k$
we also write
$$
  \chi_{na}(E/k)=\chi(G_k,\XEna{\Finfty})
$$
and call it the non-abelian Euler characteristic of $E/k$.
One can also define it using the $G_k$-homology groups of $\Selpi(E/\Finfty)$,
as for instance in \cite{HV}. The equivalence of these definitions is
a simple duality argument (see, e.g, \cite{Hows} \S1.1).

For a $p$-adic representation $\tau: G_k\to\vGL_n(\vZ_p)$ with
finite image, write
$$
  \chi_{na}(E,\tau) = \chi_{na}(G_k,\XEna{\Finfty}\vt_{\vZp}\tau).
$$
In fact, one can define the Euler characteristic for all Artin twists
and not just those defined over $\vZp$ (see \cite{CFKSV}), but we will
not need this.

\subsection{A formula of Hachimori--Venjakob}

There is an explicit formula for $\chi_{na}(E/k)$ in terms of $\chi_{cyc}(E/k)$
due to Hachimori and Venjakob (\cite{HV}, thm 4.1
and the note following it). 

\begin{theorem}\label{thHV}
Assume $p\ge 5$. Let $k\subset \Finfty$ be a number field that
contains $K=\vQ(\mu_p)$. Suppose $E/k$ has good ordinary reduction at all
primes in $k$ above $p$, and that $\Selpi(E/k)$ is finite.
Then $\chi_{cyc}(E/k)$ and $\chi_{na}(E/k)$ are finite and satisfy
\begin{equation}\label{eqHV}
  \chi_{na}(E/k) = \chi_{cyc}(E/k) \prod_{v\in P_1^{(k)} \cup P_2^{(k)}} |L_v(E/k,1)|_p,
\end{equation}
\end{theorem}

\begin{remark}
In fact Hachimori and Venjakob have a slightly different
definition of $P_2^{(k)}$. For a prime $v\notdiv p$ where $E$ has
good reduction, their condition is that $E(k_v)[p]\ne 0$, where $k_v$
is the completion of $k$ at $v$. It is equivalent to our requirement
that $\tilde E(\F_v)[p]=0$ for the reduced curve,
since multiplication by $p$ is an automorphism of the formal group
of $E$ at $v$ (see \cite{Sil1}, VII.2.1, VII.2.2 and IV.2.3).
\end{remark}

\begin{remark}
\label{hvterms}
Let $v\notdiv p$ be a prime of $k$. Denote the residue field of $v$ by $\F_v$
and its number of elements by $Nv$.
Consider the local polynomial $P_v(p^{-s})=L_v(E,s)^{-1}$ and
let us write down how the valuation of $P_v(1/Nv)=L_v(E,1)^{-1}$
depends on the reduction type of $E$ at $v$:
\begin{enumerate}
  \item If $E$ has split multiplicative reduction at $v$, then
        $P_v(T)=1-T$ and $Nv\equiv 1$ modulo $p$, so
        $P_v(1/Nv)$ has positive $p$-valuation.
  \item If $E$ has non-split multiplicative reduction at $v$, then
        $P_v(T)=1+T$ and $Nv\equiv 1$ modulo $p$, so
        $P_v(1/Nv)$ is a $p$-adic unit.
  \item If $E$ has additive reduction at $v$, then $P_v(T)=1$.
  \item Finally, if $E$ has good reduction at $v$, we have
        $$\begin{array}{rl}
          P_v(1/Nv) = & 1-a_v\frac{1}{Nv}+Nv\frac{1}{(Nv)^2}\cr
            = & (Nv)^{-1} (Nv - a_v + 1) = (Nv)^{-1} \#\tilde E(\F_v)\>.
        \end{array}$$
        So $P_v(1/Nv)$ has positive $p$-valuation precisely when $E(\F_v)$
        has non-trivial $p$-torsion.
\end{enumerate}
Thus $P_1^{(k)}$ and $P_2^{(k)}$ contain precisely those primes
not dividing $p$ for which $|L_v(E/k,1)|_p$ is non-trivial.
(Note that this also applies when $p=3$).
In particular, every term in the product in \eqref{eqHV} gives a
non-trivial contribution and $\chi_{na}(E/k) = \chi_{cyc}(E/k)$ if
and only if both $P_1^{(k)}$ and $P_2^{(k)}$ are empty.
\end{remark}

\subsection{Artin Formalism}

Euler characteristics satisfy the ``Artin formalism''.
Recall that we write $F=\vQ(\mu_p,\sqrt[p]{m})$ and $K=\vQ(\mu_p)$,
and that we have defined in \S\ref{Sintro} two representations,
$\sigma$ and $\rho$, of $\vGal(F/\vQ)$.
Both $\sigma$ and $\rho$ can be realised over the integers, and
the regular representation of $\Gal(\KL/\vQ)$ is $\sigma\oplus\rho^{p-1}$.
The following result
is a direct consequence of \cite{CFKSV} Theorem 3.10, and the multiplicative
properties of Euler characteristics.

\begin{proposition}\label{prArtinChi}
Let $E/\Q$ be an elliptic curve with good ordinary reduction at $p$. Then
we have the following equalities
\begin{equation}\label{eqAF1}
   \chi_{na}(E/K) = \chi_{na}(E,\sigma)\>,
\end{equation}
\begin{equation}\label{eqAF2}
   \chi_{na}(E/\KL) = \chi_{na}(E,\rho)^{p-1} \cdot \chi_{na}(E,\sigma) \>.
\end{equation}
In particular, in the above equations the left-hand side is defined
if and only if the right-hand side is.
\end{proposition}


\section{Compatibility of the conjectures}
\label{Stesting}

In this section we establish the relation \eqref{TheThm} between
$\chi_{na}(E,\rho)$ and $\chi_{na}(E,\sigma)$.
This is the content of our main theorem (Theorem \ref{main});
see Theorem \ref{Th6} (Appendix) for a more general argument.

Theorem \ref{main} allows us to formulate explicit consequences of
Conjecture \ref{maincon} for the $L$-values $L(E,\sigma,1)$
and $L(E,\rho,1)$, namely Conjectures \ref{lpred} and \ref{lpred3}.
These are consistent with the congruence \eqref{TheCong}, and
we also show that they follow from the Birch--Swinnerton-Dyer conjecture.


\subsection{Main theorem}

Recall that $p$ is an odd prime, $m\ge 2$ a $p$-power free integer,
$K=\vQ(\mu_p)$, 
and $\sigma$ and $\rho$ are the
$(p-1)$-dimensional representations of $\Gal(\KL/\Q)$ defined
in \S\ref{Sintro}.

\begin{theorem}\label{main}
Let $E$ be an elliptic curve over $\vQ$ with good ordinary reduction at~$p\ge 5$.
Assume that $\mu_{E/K}=0$ and that $\Selpi(E/K)$ is finite.
Then $\chi_{na}(E,\sigma)$ is defined and,
moreover,
$$
  \chi_{na}(E,\sigma) = 1,
$$
if and only if
$$
  |\Selpi(E/\KL)|<\infty \quad \vand \quad \chi_{na}(E,\rho) = 1.
$$
\end{theorem}

\begin{proof}
The theorem follows from the four lemmas below.
To be able to use Iwasawa theory, we first show that $\XE{K}$ and $\XE{\KL}$
are torsion and that other relevant technical conditions are satisfied (Lemma \ref{mlem1}).
Then we prove the asserted equivalence in a sequence of steps
(Lemmas \ref{mlem2}, \ref{mlem3}, \ref{mlem4}):
$$
\begin{array}{l}
  \chi_{na}(E,\sigma) = 1\cr
  \iff\lambda_{E/\KL} = \lambda_{E/K}\cr
  \iff|\Selpi(E/\KL)|<\infty\>\>\vand\>\>\vord_p \chi_{cyc}(E/\KL) = \vord_p \chi_{cyc}(E/K)\cr
  \iff|\Selpi(E/\KL)|<\infty\>\>\vand\>\>\chi_{na}(E,\rho)=1\>.\cr
\end{array}
$$
\end{proof}

\begin{lemma}
\label{mlem1}
Under the assumptions of the theorem,
$\XE{K}$ is $\Lambda_K$-torsion, $\XE{\KL}$ is $\Lambda_{\KL}$-torsion,
$\mu_{E/\KL}=0$,
$\chi_{na}(E,\sigma)$ is defined and $\chi_{na}(E,\rho)$
is defined provided that $|\Selpi(E/\KL)|<\infty$.
\end{lemma}
\begin{proof}
By Theorem \ref{mazcon2}, $\XE{K}$ is
$\Lambda_K$-torsion. By Hachimori--Matsuno (Theorem \ref{thHM}),
$\XE{\KL}$ is $\Lambda_{\KL}$-torsion and $\mu_{E/\KL}=0$.
By Hachimori--Venjakob (Theorem \ref{thHV}),
$\chi_{na}(E/K)$ is defined and, provided that $|\Selpi(E/\KL)|$ is finite,
$\chi_{na}(E/\KL)$ is also defined. Now, Artin formalism
for non-abelian Euler characteristics applied to $\sigma$ and $\rho$
(equations \eqref{eqAF1} and \eqref{eqAF2}) proves the claim
for $\chi_{na}(E,\sigma)$ and $\chi_{na}(E,\rho)$.
\end{proof}

\begin{lemma}
\label{mlem2}
Under the assumptions of the theorem,
$$
  \chi_{na}(E,\sigma) = 1
$$
if and only if
$$
  \lambda_{E/\KL} = \lambda_{E/K}.
$$
\end{lemma}

\begin{proof}
By assumption, $\XE{K}$ is $\Lambda_{K}$-torsion and,
by Lemma \ref{mlem1}, $\XE{\KL}$ is $\Lambda_{\KL}$-torsion.
Thus $\lambda_{E/K}$ and $\lambda_{E/\KL}$ are defined.

Artin formalism (Proposition \ref{prArtinChi}, formula (\ref{eqAF1}))
together with the theorem of Hachimori and Venjakob (Theorem \ref{thHV}),
imply that
$$
\chi_{na}(E,\sigma) = \chi_{na}(E/K) =
\chi_{cyc}(E/K) \prod_{v\in P_1^{(K)} \cup P_2^{(K)}} |L_v(E/K,1)|_p.
$$
Now, each term $|L_v(E/K,1)|_p$ in the above product is divisible by $p$
(Remark \ref{hvterms}).
Moreover,
$\vord_p\chi_{cyc}(E/K)\ge 0$ (see formula (\ref{eq:chichar})),
and, as $\mu_{E/K}=0$ by hypothesis, $\vord_p\chi_{cyc}(E/K)=0$ if and
only if $\lambda_{E/K}=0$ (Lemma \ref{chicyc0}).

Therefore, $\vord_p \chi_{na}(E,\sigma)=0$ if and only if
$\lambda_{E/K}=0$ and the sets $P_1^{(K)}$ and $P_2^{(K)}$ are both empty.
By the theorem of Hachimori and Matsuno
(Theorem \ref{thHM} and Lemma \ref{lemhm}), this is
equivalent to $\lambda_{E/K}=\lambda_{E/\KL}$.
\end{proof}

\begin{lemma}
\label{mlem3}
Under the assumptions of the theorem,
$$
  \lambda_{E/\KL} = \lambda_{E/K}
$$
if and only if
$$
  |\Selpi(E/\KL)|<\infty \quad \vand \quad \vord_p \chi_{cyc}(E/\KL) = \vord_p \chi_{cyc}(E/K).
$$
\end{lemma}

\begin{proof}
This is a direct consequence of Lemma \ref{mlem1} and Corollary \ref{chardivcor}.
\end{proof}

\begin{lemma}
\label{mlem4}
Under the assumptions of the theorem and supposing, furthermore,
that $\Selpi(E/\KL)$ is finite, we have
$$
 \vord_p \chi_{cyc}(E/\KL) = \vord_p \chi_{cyc}(E/K).
$$
if and only if
$$
 \chi_{na}(E,\rho)=1.
$$
\end{lemma}

\begin{proof}
By the Artin formalism for non-abelian Euler characteristics
(Proposition \ref{prArtinChi}, formula (\ref{eqAF2})),
$$
\chi_{na}(E,\rho)^{p-1} =  \frac{\chi_{na}(E/\KL)}{\chi_{na}(E/K)}.
$$
Applying the formula of Hachimori and Venjakob for
non-abelian Euler characteristics (Theorem \ref{thHV}) for
both $E/\KL$ and $E/K$ gives
$$
  \chi_{na}(E,\rho)^{p-1} = \frac{\chi_{cyc}(E/\KL)}{\chi_{cyc}(E/K)} \cdot
  \frac{\prod_{w\in P_1^{(\KL)} \cup P_2^{(\KL)}} |L_w(E/\KL,1)|_p}{\prod_{v\in P_1^{(K)} \cup P_2^{(K)}} |L_v(E/K,1)|_p}.
$$
The primes $v$ in $P_1^{(K)}$ (respectively $P_2^{(K)}$) are in one-to-one
correspondence with the primes $w$ in $P_1^{(\KL)}$ (respectively $P_2^{(\KL)}$),
since all such primes ramify totally in the extension $\KL/K$.
Moreover, if $w|v$ is such a pair, both the reduction type of $E$
and the residue fields are the same at $v$ and at $w$.
So $|L_v(E/K,1)|_p=|L_w(E/\KL,1)|_p$ and the products in the above formula
cancel out,
$$
  \chi_{na}(E,\rho)^{p-1} = \frac{\chi_{cyc}(E/\KL)}{\chi_{cyc}(E/K)} \>.
$$
The lemma follows.
\end{proof}

The proof of Theorem \ref{main} is now complete.

\medskip

The difficult part of the theorem is the implication
$\chi_{na}(E,\rho)=1 \implies \chi_{na}(E,\sigma)=1$. For the converse, if
$\chi_{na}(E,\sigma)=1$, then it is not hard to show that the whole
module $\XEna{\Finfty}$ is trivial, so that $\chi_{na}=1$ for any twist:

\begin{theorem}\label{premain}
Let $E$ be an elliptic curve over $\vQ$ with good ordinary reduction at~$p\ge 5$.
Assume that $\Selpi(E/K)$ is finite and that $\chi_{na}(E,\sigma)$ is
defined and equals $1$. Then $\XEna{\Finfty}=0$. In particular,
$\chi_{na}(E,\tau)=1$ for every Artin representation $\tau$ that factors
through $\Finfty$.
\end{theorem}
\begin{proof}
By Theorem \ref{mazcon2}, $\XE{K}$ is $\Lambda_{K}$-torsion.
The formula of Hachimori and Venjakob
(see \ref{thHV}) for $\chi_{na}(E/K)=\chi_{na}(E,\sigma)$ shows
that $P_1^{(K)}$ and $P_2^{(K)}$ are both empty and
that $\chi_{cyc}(E/K)=1$. By Lemma \ref{chicyc0},
$\lambda_{E/K} = \mu_{E/K} =0$.

Now suppose that $k\subset \Finfty$ is a finite Galois extension of $K$.
Note that $P_1^{(k)}$ and $P_2^{(k)}$ remain empty. The theorem of
Hachimori and Matsuno (Theorem \ref{thHM}) shows that $\XE{k}$ is
$\Lambda_k$-torsion and that $\mu_{E/k}=\lambda_{E/k}=0$. It follows
that $\XE{k}$ is finite and hence, by the theorem of Matsuno
(Theorem \ref{thM}), $\XE{k}=0$.
Taking the limit over intermediate fields $K\subset k\subset\Finfty$,
we get $\XEna{\Finfty}=0$.
\end{proof}

\subsection{Conjectures for the $L$-values}

In view of our Theorem \ref{main}, the Main Conjecture
predicts the following behaviour of $L$-values.
Note that the conjectures below are a special case of
the congruence \eqref{TheCong}.

\begin{conjecture} {\rm(consequence of Conjecture \ref{maincon})}
\label{lpred}
Let $E$ be an elliptic curve over $\vQ$ with good ordinary reduction at~$p\ge 5$.
Assume that $\mu_{E/K}=0$ and that $\Selpi(E/K)$ is finite. Then
$$
  \vord_p(\cL_E(\rho)) = 0 \quad\iff\quad \vord_p(\cL_E(\sigma)) = 0 \>.
$$
\end{conjecture}

Theorem \ref{main} and the above conjecture
are stated under the assumption that $p\ge 5$.
The primary obstacle in generalising the theorem to $p=3$ is that
the result of Hachimori--Venjakob (Theorem \ref{thHV}), which gives an
explicit formula for $\chi_{na}$, assumes that $p\ge 5$.

We believe that Theorem \ref{thHV} also holds for $p=3$
under the additional assumption that whenever $E$ has additive reduction at
a prime of $K$, the reduction stays additive in $\KL$ (as in \ref{thHM}).
If that were the case, it would imply that our Theorem \ref{main}
is also valid for $p=3$ under the same hypothesis.
Then the Main Conjecture would imply the following:

\begin{conjecture}
\label{lpred3} {\rm(Version of \ref{lpred} for $p=3$.)}
Let $E$ be an elliptic curve over $\vQ$ with good ordinary reduction at~$p=3$.
Assume that whenever $E$ has additive reduction at a prime of $K$, the reduction stays additive in $\KL$.
Suppose furthermore that $\mu_{E/K}=0$ and that $\Selpi(E/K)$ is finite. Then
$$
  \vord_p(\cL_E(\rho)) = 0 \quad\iff\quad \vord_p(\cL_E(\sigma)) = 0 \>.
$$
\end{conjecture}

\subsection{Relation with the Birch--Swinnerton-Dyer conjecture}
\label{ssBSD}

Let $E/\Q$ be an elliptic curve and $k/\Q$ a number field.
Recall that a global Weierstrass minimal model of $E$ over $\Q$ does not
necessarily stay minimal over $k$. To correctly state the
Birch--Swinnerton-Dyer formula for $E/k$ using the periods of $E/\Q$,
we define the fractional ideal $\vA_{E/k}$ of $k$ by
$$
  H^0({\mathcal E}_k,\Omega^1) = \vA_{E/k}\> \omega_{E}
$$
where $\Omega^1$ is the sheaf of invariant differentials on the N\'eron model
${\mathcal E}_k$ of $E$ over $k$ and $\omega_{E}$ is the N\'eron
differential of $E$ over $\vQ$.
The ideal $\vA_{E/k}$ is $v$-adically trivial whenever the N\'eron model
of $E$ over $\vQl$ (with $v|l$) remains a N\'eron model of $E$ over $k_v$.
In particular this holds if $E/\vQl$ has good or multiplicative
reduction, so $\vA_{E/k}$ is the unit ideal if $E/\Q$ is semistable.

\begin{conjecture}
\label{bsd-l}
{\rm(Consequence of the Birch--Swinnerton-Dyer conjecture for $E/K$ and $E/\KL$, \cite{TatC}).}
Let $p$ be an odd prime and $E/\vQ$ an elliptic curve with $E(K)$ finite. Then
$\sha(E/K)$ and $\sha(E/\KL)$ are finite,
$$
  \frac{L(E,\sigma,1)\>\epsilon(\sigma) }{{\Omega_+(E)^{(p-1)/2}(2\Omega_-(E))^{(p-1)/2}}}\>
    =\frac{\vNorm_{K/\Q}(\vA_{E/K})\sha(E/K)\prod_v c_v}{|E(K)|^2}
$$
and
$$
  \frac{L(E,\rho,1)\>\epsilon(\rho) }{{\Omega_+(E)^{(p-1)/2}(2\Omega_-(E))^{(p-1)/2}}}
    = \sqrt[p-1]{\frac{|E(K)|^{2}\vNorm_{\KL/\Q}(\vA_{E/\KL})\sha(E/\KL)\prod_w c_w
    }{{|E(\KL)|^{2}\vNorm_{K/\Q}(\vA_{E/K})}\>{\sha(E/K)\prod_v c_v}}}
    \>.
$$
The second formula may read $0=0$, in which case it states
that $L(E,\rho,1)=0$ if and only if $E(\KL)$ is infinite.
\end{conjecture}

\begin{proposition}
\label{bsd-pred}
Let $p$ be an odd prime and $E/\vQ$ an elliptic curve with good
ordinary reduction at $p$ and $E(K)$ finite.
\begin{itemize}
\item[(a)]
Assume the Birch--Swinnerton-Dyer conjecture for $E/K$. Then
$$
  \vord_p\cL_E(\sigma) = \vord_p \biggl(\chi_{cyc}(E/K) \prod_{v\in P_1^{(K)} \cup P_2^{(K)}} L_v(E/K,1)^{-1}\biggr) \>.
$$
\item[(b)]
Assume the Birch--Swinnerton-Dyer conjecture for $E/K$ and $E/\KL$.
Then $\cL_E(\rho)=0$ if $\Selpi(E/\KL)$ is infinite, and otherwise
$$
  \vord_p\cL_E(\rho) =
  \frac1{p-1}\>
  \vord_p \biggl(\frac{\chi_{cyc}(E/\KL)}{\chi_{cyc}(E/K)} \cdot
  \frac{\prod_{w\in P_1^{(\KL)} \cup P_2^{(\KL)}} L_w(E/\KL,1)^{-1}}{\prod_{v\in P_1^{(K)} \cup P_2^{(K)}} L_v(E/K,1)^{-1}}\biggr)\>.
$$
\item[(c)]
Consequently, the conjectures \ref{lpred} and \ref{lpred3} follow from
the Birch--Swinnerton-Dyer conjecture.
\end{itemize}
\end{proposition}

\begin{proof}
(a) First we substitute the formula for $L(E,\sigma,1)$ from \ref{bsd-l}
into the formula \eqref{conPL} defining $\cL_E(\sigma)$,
$$
  \cL_E(\sigma) =
    \frac{2^{\frac{p-1}2}\vNorm_{K/\Q}(\vA_{E/K})\sha(E/K)\prod_v c_v}{|E(K)|^2}
    \,\cdot\,
    \prod_{l|mp} L_l(E,\sigma,1)^{-1}
    \,\cdot\,
    \frac{\epsilon_p(\sigma)}{\epsilon(\sigma)}
    \,\cdot\,
    \frac{P_p(\sigma,u^{-1})}{P_p(\sigma,w^{-1})}\cdot u^{-\n{\sigma}}
$$
The terms $2^{\frac{p-1}2}$, $\frac{\epsilon_p(\sigma)}{\epsilon(\sigma)}$,
$u^{-\n{\sigma}}$ and $\vNorm_{K/\Q}(\vA_{E/K})$ are $p$-adic units.
The formula \eqref{eq:chicyc} yields
$$
  \vord_p\cL_E(\sigma) = \vord_p\biggl(
    \frac{\chi_{cyc}(E/K)}{|\tilde E(\F_p)|^2}
    \,\cdot\,
    \prod_{l|mp} L_l(E,\sigma,1)^{-1}
    \,\cdot\,
    \frac{P_p(\sigma,u^{-1})}{P_p(\sigma,w^{-1})}
  \biggr) \>.
$$
By Remark \ref{hvterms},
$$
  \vord_p \prod_{l|m,l\ne p} L_l(E,\sigma,1)
    =
  \vord_p \prod_{v\in P_1^{(K)} \cup P_2^{(K)}} L_v(E/K,1) \>.
$$
Decomposing $\sigma$ into 1-dimensionals, all of which except $\vTRIV$ are
ramified at $p$, we see that $P_p(\sigma,T)=1-T$.
Hence, by Lemma \ref{lemuzhas},
$$
  \vord_p
\Bigl(
  \frac{P_p(\sigma,u^{-1})}{P_p(\sigma,w^{-1})}
  \cdot
  \frac{1}{|\tilde E(\F_p)|^2}
  \cdot
  L_p(E,\sigma,1)^{-1}
\Bigr)
    =
  0 \>.
$$
The assertion follows.

(b) To begin with, if $\Selpi(E/\KL)$ is infinite, then
finiteness of $\sha(E/\KL)$ (\ref{bsd-l}) implies that $E(\KL)$ is infinite,
which in turn (\ref{bsd-l} again) implies $L(E,\rho,1)=0$.

Now we proceed as in the proof of (a), first using the Birch--Swinnerton-Dyer
formula from Conjecture \ref{bsd-l} to eliminate $L(E,\rho,1)=\sqrt[p-1]{L(E/\KL,1)/L(E/K,1)}$
from \eqref{conPL}. Then we substitute $\chi_{cyc}$ from \eqref{eq:chicyc}
and deduce that
$$
  \vord_p\cL_E(\rho) = \vord_p\Biggl(
    \sqrt[p-1]{\frac{\chi_{cyc}(E/\KL)}{\chi_{cyc}(E/K)}}
    \,\cdot\,
    |E(\F_p)|^{-2\delta}
    \,\cdot\,
    \sqrt[p-1]{\frac{\prod_{w|mp} L_w(E/\KL,1)^{-1}}{\prod_{v|mp} L_v(E/K,1)^{-1}}}
    \,\cdot\,
$$
$$
    \,\cdot\,
    \sqrt[p-1]{\frac{\vNorm_{\KL/\Q}(\vA_{E/\KL})}{\vNorm_{K/\Q}(\vA_{E/K})}}
    \,\cdot\,
    \frac{\epsilon_p(\rho)}{\sqrt[p-1]{\Delta_{\KL}^{\scriptscriptstyle 1/2}/\Delta_{K}^{\scriptscriptstyle 1/2}}}
    \,\cdot\,
    \frac{P_p(\rho,u^{-1})}{P_p(\rho,w^{-1})}
    \,\cdot\,
    u^{-\n{\rho}}
    \,\cdot\,
    2^{\frac{p-1}2}
    \Biggr)
$$
where $\delta=0$ or 1 depending on whether the prime above $p$ ramifies or
splits in $\KL/K$.

Again $2^{\frac{p-1}2}$, $u^{-\n{\rho}}$, ${\vNorm_{\KL/\Q}(\vA_{E/\KL})}$
and $\vNorm_{K/\Q}(\vA_{E/K})$ are $p$-adic units.
By the conductor-discriminant formula and the fact that
the conductor of $\rho$ is $\epsilon_p(\rho)^2$ up to a unit,
so is the term $\epsilon_p(\rho)/\sqrt[p-1]\cdots$.
By Remark \ref{hvterms},
$$
   \vord_p \frac{\prod_{w|m, w\notdiv p} L_w(E/\KL,1)}{\prod_{v|m, v\notdiv p} L_v(E/K,1)}
     =
   \vord_p \frac{\prod_{w\in P_1^{(\KL)} \cup P_2^{(\KL)}}
      L_w(E/\KL,1)}{\prod_{v\in P_1^{(K)} \cup P_2^{(K)}} L_v(E/K,1)} \>.
$$
It is not hard to see that $P_p(\rho,T)=(1-T)^\delta$ with $\delta$ as above.
Then, by Lemma \ref{lemuzhas},
$$
  \vord_p
\Bigl(
  \frac{P_p(\rho,u^{-1})}{P_p(\rho,w^{-1})}
  \cdot
  \frac{1}{|\tilde E(\F_p)|^{2\delta}}
  \cdot
  L_p(E,\rho,1)^{-1}
\Bigr)
    =
  0 \>.
$$
The assertion follows.

(c) It remains to prove that Conjectures \ref{lpred}, \ref{lpred3} follow
from (a) and (b). For $p\ge 5$ we know (Theorem \ref{thHV})
that the right-hand sides in the formulae in (a) and (b) are
the $p$-adic valuations of $\chi_{na}(E,\sigma)$ and
$\chi_{na}(E,\rho)$. Then our main theorem (Theorem \ref{main})
implies Conjecture \ref{lpred}.

The reason that this does not apply when $p=3$ is that
we do not have Hachimori--Venjakob's formula (Theorem \ref{thHV}) in this case.
However, except for this formula and Artin formalism, our proof of
main theorem (Theorem \ref{main}) does not actually use anything else
about $\chi_{na}$. Thus, if we define
$$
  \chi'_{na}(E,\sigma)=\chi_{cyc}(E/K) \prod_{v\in P_1^{(K)} \cup P_2^{(K)}} L_v(E/K,1)^{-1}
$$
and
$$
  \chi'_{na}(E,\rho)=\biggl(
  \frac{\chi_{cyc}(E/\KL)}{\chi_{cyc}(E/K)} \cdot
  \frac{\prod_{w\in P_1^{(\KL)} \cup P_2^{(\KL)}} L_w(E/\KL,1)^{-1}}{\prod_{v\in P_1^{(K)} \cup P_2^{(K)}} L_v(E/K,1)^{-1}}
  \biggr)^{\frac1{p-1}} \>,
$$
the proof goes through word for word with $\chi'_{na}$ in place of $\chi_{na}$.
The condition on additive reduction in Conjecture \ref{lpred3} is needed for
the formula of Hachimori--Matsuno (Theorem \ref{thHM}) when $p=3$.
\end{proof}

\begin{lemma}
\label{lemuzhas}
Let $p$ be an odd prime and $E/\vQ$ an elliptic curve with good
ordinary reduction at $p$.
Let $\tau$ be an Artin representation with $P_p(\tau,T)=(1-T)^\delta$.
Then
$$
  |\tilde E(\F_p)|^{2\delta} \equiv
  \frac{P_p(\tau,u^{-1})}{P_p(\tau,w^{-1})}\>
  P_p(E,\tau,1/p) \vmod\>\vZp^* \>.
$$
\end{lemma}

\begin{proof}
As $E$ has good reduction at $p$, we have $P_p(E,\tau,T)=P_p(E,T)^\delta$.
By multiplicativity, it suffices to prove the formula when $\delta=1$.
Thus $P_p(\tau,T)=1-T$,
$$
\begin{array}{l}
  \displaystyle
  P_p(E,\tau,1/p) = P_p(E,1/p) =
    (1-a_pT+pT^2)\bigm|_{T=1/p}={|\tilde E(\F_p)|}/{p}\>,\cr
  \displaystyle
  P_p(\tau,w^{-1}) = 1-w^{-1} = 1-u/p \equiv 1/p \vmod\>\vZp^* \>,\cr
  \displaystyle
  P_p(\tau,u^{-1}) = 1-u^{-1} \equiv 1-u \equiv 1-(u+w)
    \equiv 1-a_p+p \equiv |\tilde E(\F_p)| \vmod\>\vZp^* \>.
\end{array}
$$
\end{proof}

\section{Computations}
\label{Smethod}

The main purpose of the remainder of the paper
is to provide numerical evidence for the congruence \eqref{TheCong}.
In view of Theorem \ref{main}, this also supports
Conjectures \ref{lpred} and \ref{lpred3}, and thus
the main conjecture of \cite{CFKSV}.

To do this, we pick a small odd prime $p$ and an elliptic curve $E/\Q$
with good ordinary reduction at $p$, finite $\Selpi(E/K)$ and $\mu_{E/K} = 0$.
Recall that we are working with $K=\Q(\mu_p)$,
$\KL=\Q(\mu_p,\sqrt[p]m)$ and two $(p\!-\!1)$-dimensional Artin representations
$\sigma$ and $\rho$ which factor through $\Gal(\KL/\Q)$.
The congruence reads
\begin{equation}
\label{TheCong1}
  \cL_E(\sigma)\equiv\cL_E(\rho)\mod p\>.
\end{equation}

For varying $m$, we compute $L(E,\sigma,1)$ and $L(E,\rho,1)$ numerically,
and deduce the values $\cL_E(\sigma)$ and $\cL_E(\rho)$.
The tables with our results can be found in \refAppB.
In all cases, the computations agree with the congruence.

In this section, we first explain how to compute the two $L$-values
$L(E,\sigma,1)$ and $L(E,\rho,1)$. Next, we show how we determine the local
epsilon factors of $\sigma$ and $\rho$ at $p$.
The other modifying factors in the expressions for $\cL_E$
are fairly straightforward to determine. Finally, we
make some remarks concerning the data in our tables.

\subsection{Computing $L$-values: analytic side}\label{SSComputeL}

To begin with, the two $L$-functions $L(E,\sigma,s)$ and $L(E,\rho,s)$
have an analytic continuation and satisfy the functional equation \eqref{funeq}.
This can be shown using cyclic base change from the theory of automorphic
forms (see \cite{vd}, Theorem 14).

There is a standard procedure for computing values of so-called
motivic $L$-functions. This method applies to $L$-functions of
twists of elliptic curves, and we will briefly outline it in our case.
See \cite{Tol}, \cite{Coh} 10.3 and \cite{td} for details.

Let $E$ be an elliptic curve over $\vQ$, and let $\tau$ be a $d$-dimensional
Artin representation, such that $\tau\iso\tau^{*}$. For brevity, denote
$$
  N=N(E,\tau) \>.
$$
Assume furthermore
that the twisted $L$-function $L(E,\tau,s)$ has an analytic continuation
to $\vC$ and satisfies the functional equation
\begin{equation}\label{funeq2}
  \hat L(E,\tau,s) = \pm \hat L(E,\tau,2-s) \>,
\end{equation}
where
$$
  \hat L(E,\tau,s) = \Bigl(\frac{N}{\pi^{2d}}\Bigr)^{s/2}
     \Gamma\Bigl(\frac s2\Bigr)^{d} \Gamma\Bigl(\frac{s+1}2\Bigr)^{d} L(E,\tau,s) \>.
$$
Expanding the Euler product defining the $L$-function, we see that it can
be written as a Dirichlet series,
$$
  L(E,\tau,s) = \sum_{n=1}^{\infty} \frac{a_n}{n^s}, \qquad \vRe(s)>3/2 \>.
$$

Consider the inverse Mellin transform $\phi(t)$ of the full $\Gamma$-factor,
defined by
$$
  \Gamma\Bigl(\frac s2\Bigr)^{d} \Gamma\Bigl(\frac{s+1}2\Bigr)^{d} =
  \int_{0}^{\infty} \phi(t) t^{s-1} dt \>.
$$
This is a real-valued function for $t>0$, which decays exponentially
as $t\to \infty$. One can show that the functional equation \eqref{funeq2}
is equivalent to
\begin{equation}\label{thetafuneq}
  t^2 \sum_{n=1}^{\infty} a_n \phi\left(n\frac{t\pi^d}{\sqrt N}\right) =
  \pm \sum_{n=1}^{\infty} a_n \phi\left(n\frac{\pi^d}{t\sqrt N}\right) \>.
\end{equation}

We remark that the above formula allows one to test the functional
equation numerically. Provided that $\phi(t)$ can be computed,
evaluate the two sides of \eqref{thetafuneq}, say for $t=1$, and
check that the two (exponentially convergent) sums yield the same result.

By considering the Mellin transform of
$\sum_{n=1}^{\infty} a_n \phi(n{t\pi^d}/{\sqrt N})$,
splitting the integral into two,
and using \eqref{thetafuneq}, one arrives at a formula for
$\hat L(E,\tau,s)$,
\begin{equation}\label{LviaG}
  \hat L(E,\tau,s) =
    \sum_{n=1}^\infty a_nG_s\gna
    \pm\sum_{n=1}^\infty a_nG_{2-s}\gna \>,
\end{equation}
where
$$
  G_s(t) = t^{-s}\>\int_t^\infty \phi(x)\> x^s \> \frac{dx}x \>,\qquad t>0 \>
$$
is the incomplete Mellin transform of $\phi(t)$.
Again, for fixed $s$, the function $G_s(t)$ decays
exponentially with $t$, so that \eqref{LviaG} gives an exponentially
converging series for $L(E,\tau,s)$ for arbitrary $s\in \vC$.
Note that, in particular, unlike the initial Dirichlet series,
this expression is valid not only when $\vRe(s)>3/2$.

The remaining issue is that of being able to efficiently compute
the function $G_s(t)$ for $s\in \vC$ and $t>0$. (In fact, we are only
interested in the critical value $L(E,\tau,1)$, i.e. only when $s=1$.)
There are various ways to compute $G_s(t)$. In principle, there
is a method due to Tollis \cite{Tol}, based on earlier work of Lavrik
\cite{Lav}, which is applicable in our situation and gives precise
and explicit error bounds. Let us note here that $G_s(t)$ is
totally independent of the elliptic curve $E$, and depends only
on the dimension $d$ of $\tau$. Thus this is purely a problem
in transcendental function theory.

However, although $G_s(t)$ depends only on $d$, the conductor $N$
influences the rate of convergence of the series \eqref{LviaG}.
To obtain an approximation to the value $L(E,\tau,1)$ we need
roughly $\sqrt N$ terms in the series. In our examples, we often
deal with very large $N$ (e.g. for $E=19A3$, $p=7$ and $m=2$,
we have $N(E,\rho)=2^{12} 7^{14} 19^{6} \approx 10^{23}$).
Thus the issue of computing $G_s(t)$ efficiently becomes significant.
We used empirical algorithms described in \cite{td} and implemented
as a PARI package ComputeL \cite{ComputeL} to compute the $L$-values.

\begin{remark}\label{ldirect}
If $\tau$ decomposes as $\tau\iso\bigoplus_i\tau_i$, then
$$
  L(E,\tau,s) = \prod_i L(E,\tau_i,s) \>.
$$
It is more efficient to compute the individual terms
in the right-hand side of the equation, as these have smaller conductors.
For instance, this applies to $\tau=\sigma$, which decomposes as a sum
of 1-dimensional representations. We should mention that these
need not satisfy $\tau_i\iso\tau_i^*$. However, it is easy to modify
the above method to deal with such representations as well (see
\cite{td}, Remark 2.7).
\end{remark}

\subsection{Computing $L$-values: arithmetic side}
\label{SSComArith}

We now specialise to our setting. As usual, we write $K=\vQ(\mu_p)$.
Recall that the regular representation of
$\vGal(\KL/\vQ)$ decomposes as $\sigma\oplus\rho^{p-1}$, where
$\sigma$ is a sum of $1$-dimensional representations, and $\rho$ is
irreducible. Both $\sigma$ and $\rho$ have dimension $p-1$, and
$\sigma^{*}\iso \sigma$ and $\rho^{*}\iso\rho$.

We aim to compute $L(E,\sigma,s)$ and $L(E,\rho,s)$ at $s=1$.
To apply the methods described in \S\ref{SSComputeL}, we need
to determine arithmetic invariants of the two $L$-functions:
the Dirichlet coefficients, the conductor and the sign in the
functional equation. We refer the reader to \S\ref{Stwistedl}
for the definitions of these quantities and the
basic results and notation.

\subsubsection{Dirichlet coefficients}

In order to compute the Dirichlet coefficients, we determine
the local factors of $L(E,\sigma,s)$ and $L(E,\rho,s)$ at all primes.
The two $L$-functions are related
to $L$-functions of $E$ over number fields as follows:
$$
  L(E,\sigma,s) = L(E/K,s),\qquad  L(E,\rho,s)=L(E/\vQ(\sqrt[p]{m}),s)/L(E/\vQ,s) \>.
$$
Similar formulae hold for the local $L$-factors.
If $q$ is a prime where $E/\vQ$ has good reduction, write the
local polynomial of $L(E/\vQ,s)$ at $q$ as
$$
  P_q(E,T) = 1 - a_q T + qT^2 = (1-\alpha_q T)(1-\beta_q T) \>.
$$
For any number field $k$,
$$
  L_{q}(E,R_k,s) = \prod_{v|q} L_v(E/k,s)
  = \prod_{v|q} (1 - (\alpha_q^{n_v} + \beta_q^{n_v})q^{-s} + q^{2-n_{v}s})^{-1}\>,
$$
where $n_v=[\F_v:\F_q]$.
Hence we can compute these local factors if we know how to determine
$a_q$ for rational primes $q$, and how primes decompose in $k$.
In particular, this gives us a method for computing the local factors
of $L(E,\sigma,s)$ and $L(E,\rho,s)$ for all primes of good
reduction of $E$.

It remains to deal with the finitely many factors $L_v(E/k,s)$, for
primes $v$ of $k$ dividing $N(E)$. Tate's algorithm
(\cite{TatA},\cite{Sil2} \S IV.9)
determines the reduction type of $E$ at $v$ and a local
minimal Weierstrass equation. The local factor is
then given by \eqref{pvclass}. If bad reduction of $E/\vQ_q$
becomes good for $E/k_v$, we can compute $a_v$ by brute force,
counting points on $E$ over the residue field of $v$.

\subsubsection{Conductors}

Next we consider the conductors $N(E,\sigma)$ and $N(E,\rho)$.
For every prime $q$ we need to determine the local conductor
of the $l$-adic representations $M_l(E)\vt M_l(\sigma)$ and
$M_l(E)\vt M_l(\rho)$. The conductors of $E$, $\sigma$ and $\rho$
themselves are computed using Tate's algorithm for $E$, and the
conductor-discriminant formula,
$$
  N(\sigma) = |\Delta_K|, \qquad
  N(\rho) = \frac{|\Delta_{\vQ(\sqrt[p]{m})}|}{|\Delta_{\Q}|} = |\Delta_{\vQ(\sqrt[p]{m})}| \>.
$$

We now describe the local conductors of $M_l(E)\vt M_l(\sigma)$ and
$M_l(E)\vt M_l(\rho)$. As the local conductor at $q$ depends only on
the action of the inertia subgroup $I_q$, it follows that
for a prime $q$ where either $M_l(E)$ or $M_l(\rho)$ is good,
\begin{equation}\label{tenscond}
  N_q(E,\rho) =  \left\{ \begin{array}{ll}
                         N_q(E)^{\dim\rho}, & \hbox{$\rho$ good at $q$,}  \cr
                         N_q(\rho)^2, & \hbox{$E$ has good reduction at $q$}.
                         \end{array}
                         \right.
\end{equation}
The same formula holds for $\sigma$.

It remains to deal with the case when $q$ is a bad prime for both $E$
and $\rho$. (This does not happen with $\sigma$, as the only bad prime for
$\sigma$ is $p$, and $E$ has good reduction at $p$.) The local
conductor has a tame and a wild part,
$N_q(E,\rho)=q^{t_q(E,\rho) + \delta_q(E,\rho)}$. The tame
contribution is, by definition,
$$
 t_q(E,\rho) = 2\dim\rho - \deg(P_q(E,\rho,T)).
$$

To compute the wild contribution, first observe that the only
prime where $\rho$ can be wildly ramified is $q=p$. Indeed,
primes $q\neq p$ are unramified in $K/\Q$, so the action of the
inertia subgroup $I_q$ on $\rho$ factors through $\Gal(\KL/K)$. The
latter group has no elements of order $q$, so wild inertia acts
trivially.

As $E$ has good reduction at $p$, we now see that there are no primes
$q$ where both $M_l(E)$ and $M_l(\rho)$ are wildly ramified. We
can thus use an analogue to formula \eqref{tenscond} for the wild
inertia subgroup,
\begin{equation*}
  \delta_q(E,\rho) =  \left\{ \begin{array}{ll}
                         (\dim\rho)\delta_q(E), & \hbox{$M_l(\rho)$ tamely ramified at $q$,}  \cr
                         2\delta_q(\rho), & \hbox{$M_l(E)$ tamely ramified at $q$}.
                         \end{array}
                         \right.
\end{equation*}

\subsubsection{Sign}

Finally, to determine the signs in the functional equations
of $L(E,\sigma,s)$ and $L(E,\rho,s)$, we use the following result:

\begin{theorem}\label{gensign} {\rm(\cite{vd}, Theorem 1)}
Let $E$ be an elliptic curve over $\vQ$. Let $\tau$ be an Artin
representation with $\tau \iso \tau^{*}$. Set $P_M$ (respectively, $P_A$)
to be the set of primes where $E$ has multiplicative (respectively, additive)
reduction. Suppose that no prime in $P_A$ is bad for $\tau$.

If the sign in the functional equation for $L(E,s)$ is $w_E$,
then the sign in the functional equation for $L(E,\tau,s)$ is
$$
 w_E^{\vdim \tau} (-1)^{d^-(\tau)} \prod_{p \in P_M} s_p^{\vdim\tau-\vdim\tau^{I_p}} \vdet(\vFrob^{-1}_p|\tau^{I_p})  \prod_{p \in P_A} \vdet(\vFrob^{-1}_p|\tau)^{N_p(E)},
$$
where $s_p$ is $-1$ if $E$ has split multiplicative reduction at $p$, and
$+1$ if the reduction is non-split.
\end{theorem}

The theorem allows to determine the sign for $L(E,\sigma,s)$. It also
applies to $L(E,\rho,s)$, except for the case when $E$ has additive
reduction at a prime dividing $m$, since such a prime is also bad for $\rho$.
In such cases, we compute the sign numerically from \eqref{thetafuneq}.

\subsection{The local $\epsilon$-factors}
\label{Sepsilon}

In this section we briefly sketch how the local epsilon factors at $p$
for $\rho$ and $\sigma$ can be computed. For the precise definitions
and properties we refer the reader to \cite{TatN}. The necessary
results in class field theory are contained in the articles of Serre and
Tate in \cite{CF}.

To define the local epsilon factor of an Artin representation over $\Q$,
it is first necessary to pick a Haar measure $\mu$ on $\vQp$ and an
additive character $s$ given by a homomorphism
from $(\vQp,+)$ to $\vC^*$. We take $\mu$ to be the canonical measure
determined by $\mu(\vZ_p)=1$ and
$$
  s(a p^{-n})=e^{2\pi i a/p^n} , \qquad a\in\Z_p\>.
$$
In particular $s(x)=1$ iff $x\in\Z_p$. We then write
$$
  \epsilon_p (\tau)=\epsilon_p (\tau,\mu,s) \>.
$$

The local epsilon factor $\epsilon_p (\sigma)$ can be computed by writing
$\sigma=\bigoplus_i\chi_i$ as the sum of all 1-dimensional representations
of $\Gal(K/\Q)$,
$$
  \epsilon_p (\sigma) = \prod_{i=1}^{p-1} \epsilon_p (\chi_i) \>.
$$
The local epsilon factor of the trivial representation
$\vTRIV$ is 1 and, for $\chi_i\ne\vTRIV$, it is given by
$$
  \epsilon_p (\chi_i) = \sum_{j=1}^{p-1} (\chi_i\circ\theta_p)(j/p)\cdot e^{2\pi i j/p} \>.
$$
Here $\theta_p$ is the local reciprocity map at $p$ given by
class field theory.

Next, we compute the local epsilon factor $\epsilon_p (\rho)$.
Both $\sigma$ and $\rho$ are induced from 1-dimensional
representations of $\Gal(\KL/K)$,
$$
  \sigma = \vInd \vTRIV, \qquad \rho = \vInd \psi\>,
$$
with any non-trivial 1-dimensional $\psi$. The inductive property of local
epsilon factors yields
$$
  \frac{\epsilon_p (\rho)}{\epsilon_p (\sigma)} = \frac{\epsilon_v (\chi)}{\epsilon_v ({\mathbf 1})} \>.
$$
Here $\epsilon_v $ is the local epsilon factor at the prime $v$ of $K$ above $p$
defined with any measure and the additive character $s\circ\vTr_{K_v/\Q_p}$.
Choose a uniformiser $\pi$ of $K_v$. Denote $\n{\psi}=\vord_v N(\psi)$,
which can be computed from the conductor-discriminant formula in $\KL/K$.
Then it is not hard to see that
$$
  \frac{\epsilon_v (\chi)}{\epsilon_v ({\mathbf 1})} =
  \sum_{x\in T} (\psi\circ\theta_v)(x)\cdot e^{2\pi i \vTr_{K_v/\Q_p}(x)} \>,
$$
where $\theta_v$ is the local reciprocity map at $v$ and $T$ is a set of
additive representatives of $\pi^{2-p-\n{\psi}}\vO_{K_v}^*$ modulo
$\pi^{2-p}\vO_{K_v}$, for instance
$$
  T = \Bigl\{
         \sum_{i=2-p-\n{\psi}}^{1-p} a_i\pi^i
      \Bigm|
         0\le a_i\le p-1,\>a_{2-p-\n{\psi}}\ne 0
      \Bigr\} \>.
$$

It remains to explain how compute $\psi\circ\theta_v(x)$ for $x\in T$.
First, we can approximate $x$ by $x'\in K$ such that
$\psi\circ\theta_v(x)=\psi\circ\theta_v(x')$ because $\psi\circ\theta$
is trivial on $1+\pi^{\n{\psi}}\vO_{K_v}$. Now, by the product formula
in global class field theory,
$$
  \prod_w \theta_w (x') = 1\>,
$$
where the product is taken over all (Archimedean and non-Archimedean) places
of $K$. Whenever $w$ is unramified in $\KL/K$,
$$
  \theta_w(x') = \vFrob_w^{-\vord_w(x')}\>,
$$
which is easy to compute from the explicit action of the Galois group
on the $p$-th roots of $m$. To make sure that the ramified primes do not
contribute to the product, use the Chinese remainder theorem for ideals
to choose $x'$ with
$$
  x'\equiv x\mod v^{\n{\psi}}, \qquad x'\equiv 1\mod \frac{\n{\psi}}{v^{\n{\psi}}} \>.
$$
The above congruences are taken modulo ideals, so $v$ is treated as a
prime ideal of $\vO_K$. Note that the conditions at the infinite places
are automatically satisfied as $K$ has no real embeddings.
Now $\theta_v(x)=\theta_v(x')=\prod_{w\ne v}\theta_w(x')^{-1}$, which can be explicitly computed.



\subsection{The finiteness of $\Selpi(E/K)$ and the $\mu$-invariant}

Conjectures \ref{lpred} and \ref{lpred3} require that the curve in
question has finite $\Selpi(E/K)$ and $\mu_{E/K}=0$.
In view of Theorem \ref{Th6} (Appendix), this restriction should
in fact not be necessary. We now indicate what we can say about
these conditions for the curves in our tables in \refAppB.

For $p=3$ we perform $3$-descent for $E/\vQ$ and for the quadratic
twist of $E$ by $-3$, using Magma \cite{Magma}.
For those curves in our tables for which 3-descent is implemented,
we find that $\Sel_3(E/K)$ is generated by the 3-torsion points, so the
Mordell-Weil rank of $E/K$ is zero and $\sha(E/K)[3]=0$.
It follows from \eqref{selexactseq} that $\Selthreei(E/K)=0$.
For $p=5$ and the curve $X_1(11)$ (Table \ref{5-11A3}),
Fisher \cite{FisD} has done $5$-descent over $\vQ(\mu_5)$,
and he has shown that $\Selfivei(E/K)=0$.
For all other curves, we appeal to the Birch--Swinnerton-Dyer
conjecture and deduce from the non-vanishing of the $L$-functions that
$\Selpi(E/K)$ is finite.

The check whether $\mu_{E/K}=0$ we first compute the Euler characteristic
$\chi_{cyc}(E/K)$. For this we need the order of $\sha(E/K)[p^\infty]$,
that we either know to be 1 or compute from $L(E/K,1)$ using
the second part of the Birch--Swinnerton-Dyer conjecture,
$$
  \chi_{cyc}(E/K) = L^*(E/K) \cdot |\tilde E(\F_p)|^2\quad\mod\Z_p^* \>.
$$
(Note that $|\tilde E(\F_p)|=P_p(E,1)$ can be read off from the tables.)

If $\chi_{cyc}(E/K)=1$, then by Lemma \ref{chicyc0} the module
$\XE{K}$ is finite and so $\lambda_{E/K}=\mu_{E/K}=0$.
When $\chi_{cyc}(E/K)$ is non-trivial, the question is whether
the power of $p$ comes from the $\lambda$-invariant.
If we can find a point of infinite order of $E$ in
$\vQ(\mu_{p^{\infty}})=K^{cyc}$, it guarantees that the $\lambda$-invariant
is non-trivial. Thus if $\chi_{cyc}(E/K)=p$, this ensures that
$\mu_{E/K}=0$. We can do this for $p=3$ in case of the curves
20A1, 92A1 and 116C2 where we have found explicit points of infinite
order over $\Q(\mu_9)$. For the curve 128B2 we also found such a point but,
as $\chi_{cyc}(E/K)=p^2$, we only know that $\lambda_{E/K}>0$
and $\mu_{E/K}\le 1$.

\subsection{Reliability of computations}

As explained earlier, the most computationally demanding task
is that of determining the value of $L(E,\rho,1)$. The time to
compute it numerically to a given precision is roughly proportional
to the square root of the conductor $N(E,\rho)$.
For $p=3$, we restrict ourselves to those $m<2000$ for which $N(E,\rho)$
is small enough that we can evaluate $L(E,\rho,1)$ to at least 6 digits
precision with $2\!\cdot\!10^6$ Dirichlet coefficients of $L(E,\rho,s)$.
For $p=5$ we require at least 2 digits precision with $1\!\cdot\!10^8$
coefficients.
In many cases the precision is much higher, up to 40 decimal digits.
The error bounds for the method \cite{td} that we use
have not been proved in general, but the ones for the computations
of $L(E,\sigma,1)$ for all $p$ and of $L(E,\rho,1)$ for $p=3$
can be established.

Recall that the $L$-functions $L(E,\sigma,s)$ and $L(E,\rho,s)$ are defined
on all of $\vC$ (and, in particular, at $s=1$) and the modified $L$-values
$L^*(E/\Q)$, $L^*(E/K)$ and $L^*(E/\Q(\sqrt[p]m))$
are rational (see \S\ref{Sintro}).
If we had an upper bound on their denominators, we could then use our
numerical approximations to find the actual values of
$L^*(E/\Q)$, $L^*(E/K)$ and $L^*(E/\Q(\sqrt[p]m))$ (or, equivalently, of
$L(E,1)$, $L(E,\sigma,1)$ and $L(E,\rho,1)$).
We do not have such upper bounds, so we used the conjectural ones predicted
by the Birch--Swinnerton-Dyer conjecture: the denominator of $L^*(E/k)$
is at most $|E(k)_{tors}|^2$.

As a sanity check, we compute the analytic order of $\sha$ of
$E/\Q$, $E/K$ and $E/\Q(\sqrt[p]m)$. To begin with, this number is always
an integer to the correct precision. Next, this integer
(possibly 0) is always a square. Finally, we use 2-descent over the three
fields to compute the order of the corresponding 2-Selmer groups and check
that they are consistent with our conjectural orders of $\sha$. Moreover,
over $\Q$ and, for $p=3$, over $K$ we use 3-descent for a similar comparison
between the 3-Selmer groups and $\sha[3]$. The descents are carried out
using Magma~\cite{Magma}.

We have also computed 4 examples for $p=7$ (Tables \ref{7-17A1}, \ref{7-19A3}),
but the conductors in these cases are too large to make the computations
reliable. Nevertheless, the $L$-values lead to plausible orders of $\sha$.

\subsection{Observations}
\label{Sobs}

For all curves in our tables, the congruence \eqref{TheCong1}
holds.
%
Note that for some curves, the congruence reads $0\equiv 0\mod p$ for
most $m$. In view of Conjecture \ref{maincon}, this corresponds to
$\chi_{na}(E,\sigma)$ being non-unit for such $m$.
This can happen for various reasons (cf. Theorem \ref{thHV}):

\begin{enumerate}
\item If $E(K)[p]\ne 0$, then every prime $v$ of $K$ with
$v\notdiv p$, $v|m$ and where $E$ has good reduction
lies in the set $P_2^{(K)}$ of Theorem \ref{thHV}
and thus contributes to $\chi_{na}(E,\sigma)$.
See, for instance, Tables \ref{3-20A3}, \ref{5-11A3}.
\item If $E(K)[p]=0$ but $\tilde E(\F_p)[p]\ne 0$, then $\chi_{na}(E,\sigma)$ is
non-unit for every $m$.
See, for instance, Tables \ref{3-128B2}, \ref{3-275B1}.
\item If $E(K)[p]=0$ and $p$ divides one of the Tamagawa numbers $c_v$ for
a prime $v$ of $K$, then again $\chi_{na}(E,\sigma)$ is
non-unit for every $m$.
See, for instance, Table \ref{3-116C2}, \ref{3-260A1}.
\end{enumerate}

\noindent
It is interesting to note that in cases (2) and (3) there appear to be
congruences modulo a higher power $p^{b+1}$. In all our examples, this
extra power $b$ is exactly the power of $p$ in the product
of $\prod_v c_v$ and $|\tilde E(\F_p)[p]|$. Note, however, that
(at least in case (3)) the congruence only holds when both terms
have valuation exactly $b$. It would be interesting to have an explanation
of this.



\section{Numerical examples}
\label{Sexample}

To illustrate the theory in the previous sections, let
us describe an example in detail. We take the elliptic curve $21A4$,
in the notation of Cremona \cite{Cre}, of conductor $21$.
We will twist it by Artin representations coming from
$\Q(\mu_5,\sqrt[5]{2})/\Q$. Thus
$$
\begin{array}{l}
  \E: y^2 + xy = x^3 + x    \qquad (21A4) \cr   
  p=5 \cr
  m=2 \cr
  K=\Q(\mu_5), \qquad \KL=\Q(\mu_5,\sqrt[5]{2}).
\end{array}
$$
In this section, we set $L=\Q(\sqrt[5]{2})$.
This example is the first row in Table \ref{5-21A4}.

\subsection{Reduction types}
\label{e21red}

First, consider $\E$ over $\Q$. Its standard invariants are
$$
  j(\E) = \frac{47^3}{3^2\cdot 7}, \qquad \Delta(\E)=-63=-3^2\cdot 7 \qquad N(\E)=21=3\cdot7.
$$
The curve has multiplicative reduction at $q=3$ and $q=7$ and good reduction
otherwise. The Kodaira symbols for the reduction types and the local Tamagawa
numbers are
$$
\begin{array}{lll}
  3:\quad & I_2\>\hbox{split}     & c_3=2, \cr
  7:\quad & I_1\>\hbox{non-split} & c_7=1.
\end{array}
$$

We will now look at $\E$ over the number fields $K$, $L$ and $\KL$.
The corresponding local information can be computed using Tate's algorithm.
However, for semistable curves (like $\E$), it can also be determined from the
reduction behaviour over $\Q$ as follows.

Recall that split multiplicative reduction stays split multiplicative
in any extension. Non-split multiplicative reduction becomes either
split or non-split, depending on the parity of the degree of the
residue field extension. Recall also from Tate's algorithm that
in case of multiplicative reduction,
$$
  c_q = \left\{
   \begin{array}{ll}
     \vord_q(\Delta(E))  , & \hbox{split}, \cr
     2                   , & \hbox{non-split, $\vord_q(\Delta(E))$ even}, \cr
     1                   , & \hbox{non-split, $\vord_q(\Delta(E))$ odd}. \cr
   \end{array}
 \right.
$$

Now we look at $\E/K$. The curve has good reduction
at all primes of $K$ not dividing 3 or 7. Both 3 and 7 are inert in $K/\Q$,
so there are unique primes $v_3|3$ and $v_7|7$ with
residue field extensions $[\F_{v_3}:\F_3]=4$ and $[\F_{v_7}:\F_7]=4$.
So for $\E/K$ we have
$$
\begin{array}{lll}
  v_3:\quad & I_2\>\hbox{split} & c_{v_3}=2, \cr
  v_7:\quad & I_1\>\hbox{split} & c_{v_7}=1.
\end{array}
$$

Next, consider $\E/L$. The primes 3 and 7 decompose in $L$ as
$3=z_3^{(1)}z_3^{(2)}$ and $7=z_7^{(1)}z_7^{(2)}$ with
$[\F_{z_3^{(1)}}:\F_3]=1$, $[\F_{z_7^{(2)}}:\F_7]=4$ and
$[\F_{z_3^{(1)}}:\F_3]=1$, $[\F_{z_7^{(2)}}:\F_7]=4$.
This can be seen by factoring $x^5-2$ modulo 3 and 7. So for $\E/L$ we have
$$
\begin{array}{lll}
  z_3^{(1)}, z_3^{(2)} :\quad & I_2\>\hbox{split} & c_{z_3^{(i)}}=2, \cr
  z_7^{(1)}:\quad & I_1\>\hbox{non-split} & c_{z_7^{(1)}}=1,\cr
  z_7^{(2)}:\quad & I_1\>\hbox{split} & c_{z_7^{(2)}}=1.
\end{array}
$$

Finally, take $\E/\KL$. The primes $v_3$ and $v_7$ of $K$ split in $\KL/K$
as can be seen from the decomposition of 3 and 7 in $L$. Hence, in $\KL/\Q$,
$$
  3 = w_3^{(1)}w_3^{(2)}w_3^{(3)}w_3^{(4)}w_3^{(5)}, \qquad
  7 = w_7^{(1)}w_7^{(2)}w_7^{(3)}w_7^{(4)}w_7^{(5)} \>,
$$
with all residue degrees equal to 4. For $\E/\KL$ we have
$$
\begin{array}{lll}
  w_3^{(i)}:\quad & I_2\>\hbox{split} & c_{w_3^{(i)}}=2, \cr
  w_7^{(i)}:\quad & I_1\>\hbox{split} & c_{w_7^{(i)}}=1.
\end{array}
$$

\subsection{Groups of Mordell, Weil, Selmer, Tate and Shafarevich}

The curve $\E/\Q$ has Mordell-Weil rank $0$,
and $\E(\Q)\iso \Z/4\Z$, with the point $(1,1)$ as a generator.
Using 2-descent over $\Q$, $K$ and $L$ (e.g, using Magma \cite{Magma}),
we find
$$
\begin{array}{l}
  \Sel_2(\E/\Q)\iso\Z/2\Z, \cr
  \Sel_2(\E/K)\iso\Z/2\Z, \cr
  \Sel_2(\E/L)\,\iso\Z/2\Z\>.
\end{array}
$$
Thus, $\E(k)/2\E(k)\to\Sel_2(\E/k)$ is an isomorphism for $k=\Q, K, L$,
mapping the 4-torsion point $(1,1)$ to the generator. In particular,
$\E$ has no points of infinite order and $\sha(\E)[2]$ is trivial over
the three fields. Note that this implies that $\E(\KL)$ has Mordell-Weil
rank 0 as well.

\subsection{Artin representations}
\label{e21rep}

Recall that we are interested in the two representations $\sigma$
and $\rho$ of the Galois group $\Gal(\KL/\vQ)$. Write $\zeta=e^{2\pi i/5}$.

The group $\Gal(\KL/\vQ)$ has order $20$ and acts faithfully on
the 5 roots of $x^5-2$. The subgroup $\Gal(\KL/L)$ fixes $\sqrt[5]{2}$
and $\Gal(\KL/K)$ acts as a 5-cycle. We fix their generators by the requirement
$$
 \begin{array}{ll}
 g(\sqrt[5]{2})=\sqrt[5]{2},      &  g(\zeta\sqrt[5]{2})=\zeta^2\sqrt[5]{2}, \cr
 h(\sqrt[5]{2})=\zeta\sqrt[5]{2},\quad &  h(\zeta\sqrt[5]{2})=\zeta^2\sqrt[5]{2}. \cr
 \end{array}
$$
As an abstract group, $\Gal(\KL/\vQ)$ is determined by the relations
$g^4=h^5=1$ and $ghg^{-1}=h^2$, and it is usually called $G_{20}$.

The regular representation of $\Gal(\KL/\vQ)$ decomposes as $\sigma\oplus\rho^4$
where $\sigma\iso\vTRIV\oplus\chi_1\oplus\chi_2\oplus\chi_3$ is the sum of
distinct one-dimensional representations and $\rho$ is irreducible.
Both $\sigma$ and $\rho$ are 4-dimensional and $\sigma$ factors through
$\Gal(K/\Q)$.

The character table of $\Gal(\KL/\vQ)$
(including the character of $\sigma$) is
\begin{equation}
\label{chartable}
\begin{array}{c|rrrrr}
    & 1 & C_h & C_{gh} & C_{g^2h} & C_{g^3h} \cr
\hline
 \vTRIV & 1 & 1 & 1 & 1 & 1 \cr
 \chi_1 & 1 & 1 & i & -1 & -i  \cr
 \chi_2 & 1 & 1 & -1 & 1 & -1 \cr
 \chi_3 & 1 & 1 & -i & -1 & i \cr
 \rho   & 4 &-1 & 0 & 0 & 0  \cr
\hline
 \sigma & 4 & 4 & 0 & 0 & 0 \cr
\end{array}
\end{equation}
We also record the actual characteristic polynomials of the elements of
the Galois group in $\sigma$ and $\rho$,
\begin{equation}\label{charpolys}
\begin{array}{c|lllll}
    & 1 & C_h & C_{gh} & C_{g^2h} & C_{g^3h} \cr
\hline
 \sigma & (T\!-\!1)^4 & (T\!-\!1)^4 & T^4\!-\!1 & (T^2\!-\!1)^2 & T^4\!-\!1 \cr
 \rho   & (T\!-\!1)^4 & T^4\!+\!T^3\!+\!T^2\!+\!T\!+\!1 & T^4\!-\!1 & (T^2\!-\!1)^2 & T^4\!-\!1  \cr
\end{array}
\end{equation}

Note that the characters
and, moreover, the characteristic polynomials of $\rho$ and $\sigma$
are congruent modulo 5.
In fact, the representations themselves are congruent.
(They can be realised over $\Z$ and the semi-simplifications of their
reductions modulo 5 are isomorphic.)

The conductors of $\sigma$ and $\rho$ are given by the
conductor-discriminant formula,
$$
  N(\sigma) = |\Delta_K| = 5^3, \qquad N(\rho) = \frac{|\Delta_L|}{|\Delta_{\Q}|} = 2^4 5^5\>.
$$
Next, we determine the inertia subgroups and Frobenius elements
of $\Gal(\KL/\vQ)$ at all primes.

The primes 2 and 5 ramify in the extension $\KL/\Q$. The prime 5 is totally
ramified (as it ramifies both in $K/\Q$ and in $L/\Q$), so the inertia
subgroup $I_5$ is the whole of $\Gal(\KL/\vQ)$ and the
Frobenius element is trivial. The inertia invariant subspaces
of $\sigma$ and $\rho$ are, respectively, 1 and 0-dimensional.

The prime 2 is inert in $K/\Q$ and ramifies in $\KL/K$. The inertia
subgroup $I_2\subset\Gal(\KL/\vQ)$ is $\Gal(\KL/K)$ and the Frobenius
element is $g$ (modulo inertia). The inertia invariant subspaces
of $\sigma$ and $\rho$ are, respectively, the whole of $\sigma$ and 0.

If $q\ne 2, 5$, then the inertia subgroup $I_q\subset\Gal(\KL/\vQ)$
is trivial, so the representations $\sigma$ and $\rho$ are
unramified at $q$. Primes $q\equiv \pm 2$ mod 5 are inert in $K/\Q$,
so that Frob$_q\in\Gal(\KL/\Q)$ (defined up to conjugation) has order
either 4 or 20. Since $G_{20}$ only has elements of orders 1, 2, 4 and 5,
this implies that $q$ splits in $\KL/K$ and Frob$_q$ has exact order 4.
Similarly, for primes $q\equiv -1$ mod 5, the Frobenius
Frob$_q\in\Gal(\KL/\Q)$ has order 2.

It remains to deal with the primes $q\equiv 1$ mod 5. These split in
$K/\Q$ and may either be totally split in $\KL/\Q$ or remain inert
in $\KL/K$. In the first case Frob$_q=1$, and in the second case Frob$_q$
has order 5. To see which category $q$ falls into, it suffices
to determine whether $q$ is totally split or inert in $L/\Q$. This can be
done by counting roots of $x^5-2$ modulo $q$. For instance,
$x^5-2$ is irreducible modulo 11, so $q=11$ is inert in $L/\Q$. On the
other hand,
$$
  x^5-2=(x-22)(x-25)(x-49)(x-90)(x-116)\vmod 151 \>,
$$
so $q=151$ splits completely in $L/\Q$.

\subsection{Twists by Artin representations}

Now we turn to the conductors and the Euler products
for $L(\E,\sigma,s)$ and $L(\E,\rho,s)$. To begin with, the local
$L$-factors for $\E$ are (see \S\ref{e21red} and \eqref{pvclass}),

$
P_q(\E/\Q,T) = \left\{
   \begin{array}{ll}
     1 \!-\! (q+1-|\tildeE(\F_q)|) T \!+\! q T^2, & q\ne 3, 7\cr
     1 - T, & q=3, \cr
     1 + T, & q=7. \cr
   \end{array}
 \right.
$

The local $L$-factors $\sigma$ and $\rho$ are determined by the action
of Frobenius on the inertia invariant subspace of the representation.
From their classification in \S\ref{e21rep}, we have

$
\qquad\> P_q(\sigma,T) = \left\{
   \begin{array}{lll}
    1-T,       & \qquad\quad & q=5, \cr
    1-T^4,     & & q\equiv\pm 2\mod 5, \cr
    (1-T^2)^2, & & q\equiv -1\mod 5, \cr
    (1-T)^4,   & & q\equiv1\mod 5. \cr
   \end{array}
 \right.
$

\noindent
and

$
\qquad\>\>P_q(\rho,T) = \left\{
   \begin{array}{ll}
    1,               & q=2, \cr
    1,               & q=5, \cr
    1-T^4,           & q\equiv\pm 2\mod 5,\>q\ne 2\cr
    (1-T^2)^2,       & q\equiv -1\mod 5, \cr
    1\!+\!T\!+\!T^2\!+\!T^3\!+\!T^4, & q\equiv1\mod 5, \>\hbox{$x^5-2$ irreducible over $\F_q$}, \cr
    (1-T)^4,         & q\equiv1\mod 5, \>\hbox{$x^5-2$ splits over $\F_q$}. \cr
   \end{array}
 \right.
$

Now we can determine the local polynomials for the twisted $L$-functions
$P_q(\E,\sigma,T)$ and $P_q(\E,\rho,T)$. We illustrate this for $\rho$,
the case of $\sigma$ being identical. By definition,
$$
  P_q(\E,\rho,T) = \det\bigl(1-\vFrob_q^{-1} T\bigm|(M_l(\E)\tensor M_l(\rho))^{I_q}\bigr)
$$
Recall that $M_l(\E)$ is essentially the dual of the $l$-adic
Tate module $T_l(\E)$ and $M_l(\rho)$ is just $\rho$ with coefficients
taken in $\vQlb$ rather than $\vC$. If $q\ne 2,3,5,7$, then the inertia
group $I_q$ acts trivially on both $M_l(\E)$ and $M_l(\rho)$.
Then the eigenvalues of Frobenius on $M_l(\E)\tensor M_l(\rho)$ are
all pairwise products of eigenvalues on $M_l(\E)$ with eigenvalues on
$M_l(\rho)$. In other words, the local polynomial $P_q(\E,\rho,T)$
can be constructed from $P_q(\E,T)$ and $P_q(\rho,T)$ that we already know.

The same holds for $\sigma$ for $q\ne 3,5,7$. For instance, if $q=2$, we have
$$
\begin{array}{lllll}
  P_2(\E,T)&=&1+T+2T^2&=&(1-\alpha_1T)(1-\alpha_2T),\qquad \alpha_{1,2}=\frac{-1\pm\sqrt{-7}}{2},\cr
  P_2(\sigma,T)&=&1-T^4&=&(1-T)(1+T)(1-iT)(1+iT).
\end{array}
$$
Thus,
$$
\begin{array}{l}
  P_2(\E,T)=\cr
\>\> = (1\!-\!\alpha_1 T)(1\!+\!\alpha_1 T)(1\!-\!i\alpha_1 T)(1\!+\!i\alpha_1 T)(1\!-\!\alpha_2 T)(1\!+\!\alpha_2 T)
  (1\!-\!i\alpha_2 T)(1\!+\!i\alpha_2 T)\cr
\>\> =1-T^4+16T^8.
\end{array}
$$

For primes $q=2$,3,5 and 7 it is still true in our case that
\begin{equation}\label{tensinv}
  (M_l(\E)\tensor M_l(\rho))^{I_q} \iso M_l(\E)^{I_q}\tensor M_l(\rho)^{I_q}
\end{equation}
since one of the constituents is always unramified.
So the same process of constructing the local polynomials
$P_q(\E,\rho,T)$ from $P_q(\E,T)$ and $P_q(\rho,T)$ works. Again,
the same holds for $\sigma$ in place of $\rho$.
(In fact, \eqref{tensinv} holds for any semistable curve and any Artin
representation.)

Here is the summary of local polynomials for $q=2,3,5$ and 7:
$$
\begin{array}{|l|l|l|l|l|l|}
\hline
 q & P_q(\E,T) & P_q(\sigma,T) & P_q(\rho,T)  & P_q(\E,\sigma,T) & P_q(\E,\rho,T) \cr
\hline
2 & 1+T+2T^2 & 1-T^4 & 1 & 1-T^4+16T^8 & 1   \cr
3 & 1-T & 1-T^4 & 1-T^4 & 1-T^4 & 1-T^4  \cr
5 & 1+2T+5T^2 & 1-T & 1 & 1+2T+5T^2   & 1      \cr
7 & 1+T &  1-T^4 &  1-T^4  &  1-T^4  & 1-T^4  \cr
\hline
\end{array}
$$
This completes the description of all the local factors.

Finally, the formula \eqref{tenscond} gives us the conductors for the
twisted $L$-functions,
$$
  N(\E,\sigma) = 3^4 5^6 7^4, \qquad N(\E,\rho) = 2^8 3^4 5^{10} 7^4 \>.
$$
Note that, except for $5^{10}$ in case of $N(\E,\rho)$, all other
exponents are 8 (which is $\dim M_l(\E)\tensor M_l(\rho)=\dim M_l(\E)\tensor M_l(\sigma)$)
less the degree of the corresponding local polynomial. In other words,
except for $q=5$ for $\rho$, there is no wild ramification.

\subsection{$L$-functions}

Now that we have the local polynomials for the twisted $L$-functions,
we can compute their Dirichlet expansions,
$$
\begin{array}{lllll}
  L(\E,\sigma,s)&\!\!=\!\!&\prod_q P_q(\E,\sigma,q^{-s})^{-1}&\!\!=\!\!&
    \frac{1}{1^{s}} - \frac{2}{5^{s}} + \frac{16}{11^{s}}
    + \frac{1}{16^{s}} - \frac{1}{25^{s}} + \frac{8}{41^{s}} + \ldots\>, \cr
  L(\E,\rho,s)&\!\!=\!\!&\prod_q P_q(\E,\rho,q^{-s})^{-1}&\!\!=\!\!&
    \frac{1}{1^{s}} - \frac{4}{11^{s}} + \frac{1}{16^{s}} -
    \frac{2}{41^{s}} + \frac{2}{61^{s}} + \frac{1}{81^{s}} + \ldots\>. \cr
\end{array}
$$
We recall that these functions can be extended to entire functions on $\vC$
and satisfy functional equations \eqref{funeq}. We have already computed the
conductors. Theorem \ref{gensign} allows us to compute the signs in
the functional equations,
$$
\begin{array}{l}
  w_{\E,\sigma} = w_{\E}^{\vdim \sigma} (-1)^{\vdim \sigma^{-}} \prod_{q \in \{3,7\}} s_q^{\vdim\sigma-\vdim\sigma^{I_q}} \vdet(\vFrob^{-1}_q|\sigma^{I_q}),\cr
  w_{\E,\rho} = w_{\E}^{\vdim \rho} (-1)^{\vdim \rho^{-}} \prod_{q \in \{3,7\}} s_q^{\vdim\rho-\vdim\rho^{I_q}} \vdet(\vFrob^{-1}_q|\rho^{I_q}).
\end{array}
$$
In the notation of the theorem, $w_{\E}=1$ is the sign for $L(\E/\Q,s)$.
Next, $\dim\sigma^-=\dim\rho^-=2$ since complex conjugation, as an element
of $\Gal(\KL/\Q)$, lies in the conjugacy class of $C_{g^2h}$
(see \eqref{chartable} and \eqref{charpolys}).
The curve has split multiplicative reduction at $q=3$ and non-split
multiplicative reduction at $q=7$, so $s_3=-1$ and $s_7=1$.
These primes are good for $\sigma$ and $\rho$, so $I_3$ and $I_7$ act
trivially. Finally, $\vFrob_3$ and $\vFrob_7$ both have
characteristic polynomial $T^4-1$ on both $\sigma$ and $\rho$
(see \eqref{charpolys}) and, in particular, have determinant 1.
We get
$$
  w_{\E,\sigma} = +1, \qquad w_{\E,\rho} = +1 \>.
$$
The signs +1 are consistent with the Birch--Swinnerton-Dyer conjecture
for $\E/K$ and $\E/L$ as the Mordell-Weil rank is 0 over both fields.

\subsection{$L$-values}
\label{e21lval}

We have now collected all the necessary information to compute the values
$L(\E,\sigma,1)$ and $L(\E,\rho,1)$ as in \S\ref{SSComputeL}.
In fact, for efficiency reasons, we may compute $L(\E,\sigma,1)$
as the product of 1-dimensional twists (see Remark \ref{ldirect}).
We find
$$
  \begin{array}{lll}
    L(\E/\Q,1)     \!\! &\approx\!\! & 0.451115405388, \cr
    L(\E,\sigma,1) \!\! &\approx\!\! & 2.12709564136, \cr
    L(\E,\rho,1)   \!\! &\approx\!\! & 1.70167651313. \cr
  \end{array}
$$
Note that the $L$-values are non-zero as predicted by the Birch--Swinnerton-Dyer
conjecture for $\E/\Q$, $\E/K$ and $\E/L$, given that $\E$
has no points of infinite order over these fields.

The periods of $\E$ are
$$
  \Omega_+ \approx 3.60892324311, \qquad \Omega_- \approx 1.91098978075\>i \>.
$$
Now we can compute the orders of  $\sha(\E/\Q)$, $\sha(\E/K)$ and $\sha(\E/L)$ predicted
by the Birch--Swinnerton-Dyer formula \ref{bsd-l},
\begin{equation}\label{e21shaq}
  \sha_{an}(\E/\Q) = \frac{L(\E,1)}{{\Omega_+(\E)}}
    \frac{|\E(\Q)|^2}{c_3c_7} \approx 1.0000000000 \>,
\end{equation}
\begin{equation}\label{e21shak}
  \sha_{an}(\E/K) = \left|\frac{L(\E,\sigma,1)\Delta_K}{{\Omega_+(\E)^{2}(2\Omega_-(\E))^{2}}}\right|
    \frac{|\E(K)|^2}{c_{v_3}c_{v_7}} \approx 1.0000000000 \>,
\end{equation}
\begin{equation}\label{e21shal}
  \sha_{an}(\E/L) = \left|\frac{L(\E,\rho,1)L(\E/\Q,1)\Delta_L}{{\Omega_+(\E)^{3}(2\Omega_-(\E))^{2}}}\right|
    \frac{|\E(L)|^2}{  c_{z_3^{(1)}}c_{z_3^{(2)}}c_{z_7^{(1)}}c_{z_7^{(2)}}  }
    \approx 1.0000000000 \>,
\end{equation}
in agreement with $\sha(\E/\Q)[2]=1$, $\sha(\E/K)[2]=1$ and $\sha(\E/L)[2]=1$.

\begin{remark}
The $L$-functions of the Artin representations, $L(\sigma,s)$ and
$L(\rho,s)$, also have a meromorphic continuation to $\vC$ and satisfy
a functional equation. In fact, they can be expressed in terms of
Dedekind zeta functions,
$$
   L(\sigma,s) = \zeta_K(s), \qquad L(\rho,s) = \frac{\zeta_L(s)}{\zeta(s)} \>.
$$
Thus their values at $s=1$ have a similar arithmetic interpretation.
Recall that by the class number formula, for a number field $k/\Q$,
$$
  \vRes_{s=1}\zeta_k(s)=\frac{2^{r_{k,1}}(2\pi)^{r_{k,2}}{\text Reg}_k\,|\text{Cl}(\vO_k)|}{|\mu_k|\sqrt{\Delta_k}}\>,
$$
where $r_{k,1}$ (resp. $r_{k,2}$) is the number of real (resp. pairs of complex)
embeddings of~$k$, Reg$_k$ is the regulator, $\text{Cl}(\vO_k)$ is the
class group of $k$ and $\mu_k$ is the set of all roots of unity in $k$.


\end{remark}

\subsection{Iwasawa theory}

To apply results of Iwasawa theory to $\E$, we first need to know that
$\E$ has good ordinary reduction at $p=5$.
This is easily verified, as
$P_5(E/\Q) = 1+2T+5T^2$ has a non-zero linear term modulo 5.

Next, we need the $\mu$-invariant of $\E/K$ to be zero for $p=5$.
This condition is genuinely hard. At present, this requires
knowing $\Selfivei(\E/K)$.
For some curves with rational 5-torsion subgroups, Fisher \cite{FisD,FisS}
carried out 5-descent to determine the structure of their 5-Selmer groups.
Our curve has
$$
  \Gal(\Q(\E[5])/\Q)\iso\vGL_2(\F_5)\qquad\hbox{(order 480)}\>,
$$
so doing 5-descent is probably unrealistic.

We can numerically determine
the conjectural order of $\sha(\E/K)$ from the $L$-value $L(\E,\sigma,1)$
(see \S\ref{e21lval}). We find
$$
  \sha(\E/K)\>\>{\buildrel?\over=}\>\>1\>,
$$
where the equality would hold if we knew the Birch--Swinnerton-Dyer
conjecture for $\E/K$ and precise estimates for our $L$-value computations.
As we already know that the Mordell-Weil rank of $\E/K$ is zero, this would
imply that $\Selfivei(\E/K)$ is trivial.

We can compute the cyclotomic Euler characteristic of $\E/K$, given by \eqref{eq:chicyc},
$$
  \chi_{cyc}(\E/K) = \frac{|\tildeE(\F_5)|^2 c_{v_3}c_{v_7} |\sha(E/K)[5^\infty]|}{|E(K)_{tors}|^2} \vmod\>\Z_5^{*} \>.
$$
As all the terms are coprime to 5, we get
$$
  \chi_{cyc}(\E/K)=1 \>.
$$
Now, by Lemma \ref{chicyc0}, we get
$$
  \mu_{\E/K}=\lambda_{\E/K}=0 \>.
$$
In other words, the characteristic element $f_{\E/K}=1$, so that
$\Selfivei(\E/\Q(\mu_{5^\infty}))$ is finite.

The theorem of Hachimori--Matsuno (Theorem \ref{thHM}) tells us that
$X(\E/\KL)$ is $\Lambda_{\KL}$-torsion and $\mu_{\E/\KL}=0$.
Moreover, as $m=2$ is coprime to the conductor $N(\E)=21$
and $|\E(\F_{16})|=16$ is coprime to 5,
Corollary \ref{corhm} implies that $\lambda_{\E/\KL}=0$.
Again, we get that the characteristic element $f_{\E/\KL}=1$, so that
$\Selfivei(\E/\Q(\mu_{5^\infty},\sqrt[5]2))$ is finite.

The reverse implication of Lemma \ref{chicyc0} now tells us that
$$
  \mu_{\E/\KL}=\lambda_{\E/\KL}=0 \quad\implies\quad \chi_{cyc}(\E/\KL)=1 \>.
$$
In particular, by the Euler characteristic formula \eqref{eq:chicyc},
$$
  \sha(E/\KL)_{5^\infty}=1 \>.
$$

\subsection{Non-abelian Euler characteristics}

To test the main conjecture \ref{maincon} for $\E$, we compute the
non-abelian Euler characteristics $\chi_{na}(\E,\sigma)$ and
$\chi_{na}(\E,\rho)$. First, the formula of Hachimori--Venjakob
(Theorem \ref{thHV}) yields
$$
  \chi_{na}(\E,K) = \chi_{cyc}(\E,K)=1, \qquad \chi_{na}(\E,\KL) = \chi_{cyc}(\E,\KL)=1\>.
$$
This is because $P_1^{(K)}$ and $P_1^{(\KL)}$ are empty since
$m=2$ is coprime to the conductor $N(\E)=21$, and
$P_2^{(K)}$ and $P_2^{(\KL)}$ are empty since
$|\E(\F_{16})|=16$ is coprime to 5.

By Artin formalism for non-abelian Euler characteristics
(Proposition \ref{prArtinChi}),
$$
  \chi_{na}(\E,\sigma)=\chi_{na}(\E,\rho)=1.
$$
This is consistent with Theorem \ref{main}, which states that one
of the Euler characteristics is trivial if and only if the other one is.

The main conjecture \ref{maincon} predicts that the values
of the $p$-adic $L$-function $\cL_{\E}$, given by \eqref{conPL},
satisfy
$$
\vord_5 \cL_{\E}(\sigma) = \vord_5 \chi_{na}(\E,\sigma) = 0
$$
and
$$
\vord_5 \cL_{\E}(\rho) = \vord_5 \chi_{na}(\E,\rho) = 0 \>.
$$
Note that Conjecture \ref{lpred}, that follows from the main conjecture
and Theorem \ref{main}, asserts that these two $5$-adic valuations
are either both zero, or neither of them is.

\subsection{Verifying the main conjecture}

Let us compute $\cL_{\E}(\sigma)$ and $\cL_{\E}(\rho)$. Recall that
$$
  \cL_{\E}(\rho) = \frac{\Lvten ({\E},\rho,1)}{\Omega_+({\E})^{2}\Omega_-({\E})^{2}}
                     \cdot \epsilon_5(\rho) \frac{P_5(\rho,u^{-1})}{P_5(\rho,w^{-1})}\cdot u^{-\n{\rho}}\>,
$$
and similarly for $\sigma$. Here
$$
 \n{\sigma} = \vord_5 N(\sigma) = 3 , \qquad
 \n{\rho} = \vord_5 N(\rho) = 5 \>.
$$
The $5$-adic numbers $u$ and $w$ are determined by
$$
  P_5(E,T) = 1+2T+5T^2 = (1-uT)(1-wT) \qquad u\in\Z_5^{*},
$$
so, using Hensel's Lemma,
$$
\begin{array}{ll}
 u=3 \>+\!\!\!\! & 2\!\cdot\!5 + 4\!\cdot\!5^2 + 2\!\cdot\!5^3 + 5^4 + 4\!\cdot\!5^5 + 2\!\cdot\!5^7 + 5^8 + 5^9 +  O(5^{10}),\cr
 w=    & 2\!\cdot\!5 + 2\!\cdot\!5^3 + 3\!\cdot\!5^4 + 4\!\cdot\!5^6 + 2\!\cdot\!5^7 + 3\!\cdot\!5^8 + 3\!\cdot\!5^9 + O(5^{10})\>.
\end{array}
$$
The local $\epsilon$-factors $\epsilon_5(\sigma)$ and $\epsilon_5(\rho)$
can be computed as in \S\ref{Sepsilon},
$$
   \epsilon_5(\sigma) = -5^{3/2}, \qquad \epsilon_5(\rho) = -5^{5/2} \>.
$$
From \eqref{e21shaq}, \eqref{e21shak} and \eqref{e21shal} we get numerically the values
$$
\begin{array}{l}
  \displaystyle
  L^*(\E,\sigma) = \left|\frac{L(\E,\sigma,1)\sqrt{\Delta_K}}{\Omega_+({\E})^{2}(2\Omega_-({\E}))^{2}}\right| = 1/8 ,
  \\[10pt]
  \displaystyle
  L^*(\E,\rho) = \left|\frac{L(\E,\rho,1)\sqrt{\Delta_L}}{\Omega_+({\E})^{2}(2\Omega_-({\E})^{2}}\right| = 2.
\end{array}
$$
We combine all this information to obtain
$$
\begin{array}{ll}
  \cL_{\E}(\rho) & = \displaystyle
    L^*(\E,\rho)
    \cdot(-2i)^2
    \cdot\prod_{q=2,5} P_q(E,\rho,q^{-1})\cdot
    \frac{\epsilon_5(\rho)}{\sqrt{|\Delta_L|}}
    \frac{P_5(\rho,u^{-1})}{P_5(\rho,w^{-1})}\cdot u^{-\n{\rho}}\cr
  &\displaystyle = 2 \cdot (-4) \cdot (1\cdot 1)\cdot\frac{-5^{5/2}}{2^2 5^{5/2}}\cdot\frac{1}{1}\cdot u^{-5}\cr
  & \displaystyle = 4 + 2\!\cdot\!5 + 2\!\cdot\!5^2 + 4\!\cdot\!5^3 + O(5^4) \>.
  \end{array}
$$
Similarly, for $\sigma$ we have
$$
\begin{array}{ll}
  \cL_{\E}(\sigma) & = \displaystyle
  L^*(\E,\sigma)
    \cdot(-2i)^2
    \cdot\prod_{q=2,5} P_q(E,\sigma,q^{-1})\cdot
    \frac{\epsilon_5(\sigma)}{\sqrt{|\Delta_K|}}
    \frac{P_5(\sigma,u^{-1})}{P_5(\sigma,w^{-1})}\cdot u^{-\n{\sigma}}\cr
  &\displaystyle = \frac18 \cdot (-4) \cdot (1\!-\!2^{-4}\!+\!16\!\cdot\!2^{-8})(1\!+\!2\!\cdot\! 5^{-1}\!+\!5\!\cdot\!5^{-2})
  \cdot\frac{-5^{3/2}}{5^{3/2}}\cdot\frac{1-u^{-1}}{1-w^{-1}}\cdot u^{-3}\cr
  & \displaystyle = 4 + 3\!\cdot\!5^2 + 5^3 + O(5^4) \>.
  \end{array}
$$
As expected, $\cL_{\E}(\sigma)$ and $\cL_{\E}(\rho)$ are 5-adic
units and are congruent modulo 5, in accordance with \eqref{TheCong1}.


\subsection{An example with additive reduction}

The curve $E_{21}$ from the example in \S\ref{Sexample}
is everywhere semistable. If $E/\Q$ has additive reduction
at some prime $l$ of $\Q$, some care has to be taken
when computing the Birch--Swinnerton-Dyer quotient for $E/K$ and $E/\Q(\sqrt[p]{m})$
using a model over $\Q$. Let us illustrate this with an example.

Consider the elliptic curve $E=272C1$ in the notation of Cremona \cite{Cre},
$$
  E: y^2 = x^3-x^2-4x\>.
$$
The curve has bad reduction at $l=2$ and $l=17$ with
$$
  \begin{array}{lll}
  \vord_2 \Delta(E)=8,&\hbox{additive red. at 2 of type $I_0^*$}& c_2=4,\cr
  \vord_{17}\Delta(E)=1,&\hbox{split multiplicative red. at 17 of type $I_1$}& c_{17}=1.
  \end{array}
$$
The given Weierstrass model for $E/\Q$ is minimal at all primes and
$$
  E(\Q)\>\iso\Z/2\Z,
  \qquad  \Omega^+\approx 3.47306346
  \qquad  2\Omega^+\Omega^-\approx 6.75008201i\>.
$$
We can compute the Birch--Swinnerton-Dyer quotient and obtain, numerically,
$$
  \frac{L(E,1)}{\Omega^+}\approx\frac{\prod_l{c_l}}{|E(\Q)|^2}\cdot 1.00000000 \>,
$$
so the analytic order of $\sha(E/\Q)$ is 1.

Now let $p=3$, $m=2$ and $L=\Q(\sqrt[p]m)=\Q(\sqrt[3]2)$.
We use 2-descent to conclude that $E(\Q)=E(L)$, compute the local
Tamagawa numbers $c_v$ for $E/L$ and the value $L(E/L,1)$, leading to
\begin{equation}
\label{sha2}
  \frac{L(E/L,1)}{|\Omega^+\Omega^+2\Omega^-|}\sqrt{|\Delta_L|}\approx\frac{\prod_v{c_v}}{|E(L)|^2}\cdot 2.00000000 \>.
\end{equation}
Clearly the real number on the right is not supposed to be the analytic
order of $\sha$ (being not a square of an integer). The reason is that
the model that we use for $E$ is no longer minimal and the N\'eron
differential of $E/\Q$ is not the N\'eron differential for $E/L$ at $l=2$.
Thus, if we use the same differential as over $\Q$, then the ideal $\vA_{E/L}$
is non-trivial. In fact, it is not hard to see that $\vA_{E/L}=(\sqrt[3]{2})$ 
and its norm $\vNorm_{L/\Q}(\vA_{E/L})=2$
contributes to the Birch--Swinnerton-Dyer quotient,
leading to $\sha_{an}(E/L)=1$ (in agreement with 2-descent.)

Alternatively, note that over $\Q$ the curve $E$ has another minimal
Weierstrass equation,
$$
  y^2-2xy=x^3-2x^2-4x\>,
$$
which can be written over $L$ as
$$
  y^2-\pi^3 xy=x^3-\pi^3x^2-\pi^6x\>,\qquad \pi=\sqrt[3]{2}\>.
$$
Now let $x=\pi^2x'$ and $y=\pi^3y'$, transforming the latter model into
$$
  E':\>\> (y')^2-\pi^2 x'y'=(x')^3-\pi (x')^2-\pi^2x'\>.
$$
From Tate's algorithm, it follows that this is a minimal Weierstrass equation
at all primes, so $\vA'=1$ for the differential of this model.
On the other hand, from the transformation of coordinates it follows that,
with this differential,
$$
  (\Omega^+)'=\pi\,\Omega^+, \qquad (\Omega^-)'=\pi\,\Omega^-.
$$
With these ``corrected'' periods the analytic
order of $\sha$ from \eqref{sha2} becomes 1, as asserted.

\subsection{An example for a $GL_2$-extension}
\label{gl2ex}

One can make $L$-value computations for an elliptic curve $E$
in $p$-adic Lie extensions other than $\Finfty/\Q$.
For example, one may take an extension of the form $\Q(C[p^\infty])/\Q$ for
some elliptic curve $C/\Q$. For non-CM curves $C$, the Galois group of this
extension is an open subgroup of $\vGL_2(\vZ_p)$ and is the whole of $\vGL_2$
for almost all $p$, see \cite{SerP}. These extensions (with $E=C$)
are the main focus of \cite{CFKSV}.
One can again define a non-abelian Euler characteristic $\chi_{na}(E,\tau)$
of $E$ twisted by
an Artin representation $\tau$ that factors through such an extension.
There is also a conjectural $p$-adic $L$-function $\cL_E$, interpolating
the special values of twisted $L$-functions $L(E,\tau,s)$
(\cite{CFKSV}, Conjecture 5.7).
The ``Main Conjecture'' (\cite{CFKSV}, Conjecture 5.8) gives a relation
between the values $\cL_E(\tau)$ and Euler characteristics
$\chi_{na}(E,\tau)$.

\begin{example}
\label{exgl2}
Take $p=5$ and consider the elliptic curve $E=11A3=X_1(11)$ over $\Q$,
$$
  E:\>\> y^2 + y = x^3 - x^2 \>.
$$
The field $\vQGL$ contains the fields $k_1, k_2$ of 5-torsion points
of $E$ and of its isogenous curve $11A2$. Both are extensions of degree 5
of $\Q(\mu_5)$ and can be explicitly described as follows.
In $\Q(\sqrt5)$ one can write $-11=\alpha\beta$ with
$\alpha=\frac{9+5\sqrt{5}}2$ and $\beta=\frac{9-5\sqrt{5}}2$. Then
$$
  k_1 = \Q\Bigl(\mu_5,\sqrt[5]{\frac\alpha\beta}\Bigr),\qquad k_2=\Q(\mu_5,\sqrt[5]{11})\>.
$$
Their Galois groups are the same as in the first layer of the false
Tate curve extension for $p=5$,
$$
  \Gal(k_1/\Q) \iso G_{20} \iso \Gal(k_2/\Q) \>.
$$
In particular, we have two 4-dimensional irreducible Artin representations
$\rho_1, \rho_2$ of $\vGal(\vQGL/\Q)$ that factor through $k_1$
and $k_2$ respectively.

The Main Conjecture implies that (\cite{CFKSV}, 5.7-5.10)
\begin{equation}\label{gl2main}
  \vord_5(\chi_{na}(E,\rho_i)) = \vord_5
    \left(
  \frac{L_R(E,\rho_i,1)}{\Omega_+(E)^2\Omega_-(E)^2}
      \cdot \frac{P_5(\rho_i,u^{-1})}{P_5(\rho_i,w^{-1})}
      \> \epsilon_5(\rho_i)\>u^{-f_{\rho_i}}
    \right), \qquad i=1,2.
\end{equation}
Here $R=\{5,11\}$ is the set of primes that contains $p$ and those primes
that divide the denominator of the $j$-invariant $j(E)=-2^{12}/11$.
The modified $L$-function $L_R(E,\rho_i,s)$ is that of
$E$ of twisted by $\rho_i$ with the local factors at $l\in R$ removed.

We proceed to verify this conjecture, similarly to the false Tate curve case.

First, the curve $E$ has a 5-torsion point over $\Q$, and
$$
  E(\Q(\mu_5))[5^\infty] \iso \Z/5\Z, \qquad
  E(k_1)[5^\infty] \iso (\Z/5\Z)^2, \qquad
  E(k_2)[5^\infty] \iso \Z/5\Z\>.
$$
By results of Fisher \cite{FisD}, $E$ has Mordell-Weil rank 0
over $\Q(\mu_5)$, $k_1$ and $k_2$, and
$$
  \sha(E/\Q(\mu_5))[5^\infty]=0, \qquad \sha(E/k_1)[5^\infty]\iso(\Z/5\Z)^2, \qquad \sha(E/k_2)[5^\infty]=0\>.
$$
Since $E$ has good ordinary reduction at $p=5$, the module $\XE{\Q(\mu_5)}$
is $\Lambda_{\Q(\mu_5)}$-torsion by Theorem \ref{mazcon2}.
Using formula \ref{eq:chicyc}, it is easy to see that
$\chi_{cyc}(E/\Q(\mu_5))=1$.
It follows from Lemma \ref{chicyc0} that
$\lambda_{E/\Q(\mu_5)}=\mu_{E/\Q(\mu_5)}=0$.
By Theorem \ref{thHM}, both $\XE{k_1}$ and $\XE{k_2}$ are torsion and
their $\mu$-invariants are 0.

Next, we compute $\chi_{na}(E,\rho_i)$. By Artin formalism
(\cite{CFKSV} Theorem 3.10),
$$
  \chi_{na}(E,\rho_i) = \sqrt[4]{\frac{\chi_{na}(E/k_i)}{\chi_{na}(E/\Q(\mu_5))}} \>.
$$
A formula for $\chi_{na}$ of $E$ over number fields in $\vQGL$
has been worked out in Coates-Howson \cite{CH}, Theorem 1.1.
In our case, for a number field $k\subset\vQGL$ for which $\Selfivei(E/k)$
is finite, it states
$$
  \chi_{na}(E/k) = \chi_{cyc}(E/k) \prod_{v|11} |L_v(E,1)|_5,
$$
the product taken over primes of $k$ dividing 11.
(That the theorem is applicable here follows from \cite{CH}, Theorem 6.4.)

Now, the above two formulae together with the values of
$|\sha(E/k_i)|$ and $|\sha(E/\Q(\mu_5))|$ yield
\begin{equation}
\label{gl2chnaval}
  \chi_{na}(E/k_1) = 5^{3}, \qquad \chi_{na}(E/k_2) = 5^{1}.
\end{equation}
To test the main conjecture numerically, we evaluate the right-hand side of
\eqref{gl2main}. We compute the $L$-values as described in
\S\S\ref{SSComputeL}-\ref{SSComArith} with mild modifications in the
case of $E/k_1$. We find
$$
  L(E,\rho_1,1) \approx 1.26706100,\qquad L(E,\rho_2,1) \approx 4.05459521\>.
$$
The local $\epsilon$-factors $\epsilon_5(\rho_1)$ and
$\epsilon_5(\rho_2)$ can be determined as in \S\ref{Sepsilon},
$$
  \epsilon_5(\rho_1) = -5^{3/2},\qquad \epsilon_5(\rho_2) = -5^{5/2}\>.
$$
Next,
$$
\begin{array}{ll}
  P_5(\rho_1,T)=1-T, & P_5(\rho_2,T)=1,\cr
  P_5(E,\rho_1,T)=1-T+5T^2, & P_5(E,\rho_2,T)=1,\cr
  P_{11}(E,\rho_1,T)=1, & P_{11}(E,\rho_2,T)=1.
\end{array}
$$
and
$$
  \Omega_1 \approx 6.34604652,\qquad \Omega_2 \approx 1.45881662\>i\>.
$$
Finally,
$$
  u = 1+4\cdot 5+3\cdot 5^2+2\cdot 5^3+\BigO(5^4), \qquad
  w = 5+5^2+2\cdot 5^3+\BigO(5^7)\>.
$$
Combining the above information, we obtain
$$
\begin{array}{ll}
  \displaystyle
  \frac{L_R(E,\rho_1,1)}{\vphantom{\int_{\int}}\Omega_+(E)^2\Omega_-(E)^2}
      \cdot \frac{P_5(\rho_1,u^{-1})}{P_5(\rho_1,w^{-1})}
      \> \epsilon_5(\rho_1)\>u^{-f_{\rho_1}} = 4\cdot5^3+2\cdot 5^4+\BigO(5^5)\cr
  \displaystyle
  \frac{L_R(E,\rho_2,1)}{\Omega_+(E)^2\Omega_-(E)^2}
      \cdot \frac{P_5(\rho_2,u^{-1})}{P_5(\rho_2,w^{-1})}
      \> \epsilon_5(\rho_2)\>u^{-f_{\rho_2}} = 4\cdot5+4\cdot 5^2+\BigO(5^3) \>.
\end{array}
$$
This agrees with \eqref{gl2chnaval}, as predicted by the main conjecture.

\end{example}


\def\theequation{A.\arabic{equation}}
\section*{Appendix A (by J. Coates and R. Sujatha)}
\setcounter{equation}{0}

The aim of this appendix is to prove several theoretical results, which are
related to, and illustrated by, the numerical calculations carried out in
this paper. As earlier, let $p$ be an odd prime number, and put
$K=\vQ(\mu_p)$, $K^{cyc}=\vQ(\mu_{p^{\infty}})$. Let $m>1$ be an integer,
which is $p$-power free, and put
\begin{equation}
\Finfty = K^{cyc}(\sqrt[p^n]{m}: n=1,2..), \qquad G=\vGal(\Finfty/\vQ).
\end{equation}
Let $E$ be an elliptic curve over $\vQ$, which we shall always assume
satisfies:
\begin{hypothesis}
$E$ has good ordinary reduction at $p$, and $m$ is not divisible by
any prime of additive reduction for $E$.
\end{hypothesis}
The hypothesis on $m$ is made to ensure that the principal results
of \cite{HM} and \cite{HV} remain valid for $p=3$ in our case.

We write $X(E/\Finfty)$ for the Pontryagin dual of the $p^{\infty}$-Selmer
group of $E$ over $\Finfty$. If $J$ is any compact $p$-adic Lie group,
$\Lambda(J)$ will denote the Iwasawa algebra of $J$. Then $X(E/\Finfty)$
is in fact $\Lambda(G)$-torsion (see \cite{HV}), thanks to Kato's \cite{KatP}
deep theorem that the dual of the Selmer group of $E$ over $K^{cyc}$
is $\Lambda(\Gamma_K)$-torsion, where $\Gamma_K=\vGal(K^{cyc}/K)$. However,
the ideas of \cite{CFKSV} suggest that something much stronger should always be
true for $X(E/\Finfty)$. Let $H=\vGal(\Finfty/\vQ^{cyc})$, where
$\vQ^{cyc}$ denotes the cyclotomic $\vZp$-extension of $\vQ$, and let
$\fM_H(G)$ denote the category of all finitely generated $\Lambda(G)$-modules
$M$ such that $M/M[p^{\infty}]$ is finitely generated over $\Lambda(H)$;
here $M[p^{\infty}]$ denotes the $p$-primary submodule of $M$.

\begin{conjecture}\label{appCon2}
$X(E/\Finfty)$ belongs to $\fM_H(G)$.
\end{conjecture}

Using similar techniques to those discussed in \S 5 of \cite{CFKSV}, it is easy
to see that Conjecture \ref{appCon2} is equivalent to the assertion that
\begin{equation}
 \mu_{G_K}(X(E/\Finfty)) = \mu_{\Gamma_K}(X(E/K^{cyc})),
\end{equation}
where $G_K$ is the pro-$p$ group given by $G_K=\vGal(\Finfty/K)$, and
$\mu_{G_K}(M)$ denotes the $\mu$-invariant of any finitely generated torsion
$\Lambda(G_K)$-module $M$ as defined in \cite{Hows2,OV}, and similarly
$\mu_{\Gamma_K}(R)$ denotes the $\mu$-invariant of a finitely generated
torsion $\Lambda(\Gamma_K)$-module $R$. Moreover, Conjecture \ref{appCon2}
is true if there exists an elliptic curve $E'$ over $K$, which is isogenous
to $E$ over $K$, such that $X(E'/K^{cyc})$ is a finitely generated
$\vZp$-module.

As always, we shall find that the deepest arithmetic questions arise from the
interplay between the $\Lambda(G)$-module $X(E/\Finfty)$ and the complex
$L$-functions $L(E,\phi,s)$, which are studied earlier in the paper, where
$\phi$ runs over all Artin representations of $G$. We recall that the
$L(E,\phi,s)$ are known to be entire, and to satisfy the standard functional
equation. Moreover, if $k$ denotes any finite extension of $\vQ$ contained
in $\Finfty$, we write $L(E/k,s)$ for the complex $L$-function of
$E$ over $k$, and we define
\begin{equation}
 g_{E/k} = \vrank_{\Z}(E(k)), \qquad\qquad r_{E/k}=\vord_{s=1}(L(E/k,s)) \>.
\end{equation}
Of course, the conjecture of Birch and Swinnerton-Dyer predicts that we
always have $g_{E/k}=r_{E/k}$.

We first discuss the notion of regularity for $E$ over $\Finfty$.

\begin{definition}
We say $E$ is {\em regular} over $\Finfty$ if $X(E/\Finfty)=0$.
\end{definition}

As we shall explain below, the calculations in this paper provide many
examples when $E$ is regular over $\Finfty$.
But first we note the following theorem.

\begin{theorem}
Assume $E$ is regular over $\Finfty$. Then $g_{E/k}=r_{E/k}$
for all finite extensions $k$ of $\vQ$ contained in $\Finfty$ if
and only if $L(E,\phi,1)\neq 0$ for all Artin representations $\phi$ of $G$.
\end{theorem}

\begin{proof}
As earlier, let $\Selpi(E/k)$ \quad (resp. $\Selpi(E/\Finfty)$) denote
the $p^{\infty}$-Selmer group of $E$ over $k$ (resp. $\Finfty$).
Now the kernel of the restriction map
$$
 \Selpi(E/k) \to \Selpi(E/\Finfty)
$$
is always finite (see \cite{HV}). Assuming now that $\Selpi(E/\Finfty)=0$,
it follows that $\Selpi(E/k)$ must be finite, and so $g_{E/k}=0$ for all
finite extensions $k$ of $\vQ$ contained in $\Finfty$. Thus
it remains to show that $r_{E/k}=0$ for all finite extensions $k$ of $\vQ$
contained in $\Finfty$ if and only if $L(E,\phi,1)\neq 0$ for all
irreducible Artin representations $\phi$ of $G$. But on the one hand, for
any finite extension $k$ of $\vQ$ in $\Finfty$, we have
$L(E/k,s)=L(E,\theta,s)$, where $\theta$ is the Artin representation of $G$
induced by the trivial representation of $\vGal(\Finfty/k)$, and
so one direction is clear. The other direction is plain from the holomorphy
of the $L(E,\phi,s)$ at $s=1$, and the fact that, for any finite Galois
extension $k$ of $\vQ$ contained in $\Finfty$ with Galois group $\Delta$,
we have
\begin{equation}\label{eqArtin}
L(E/k,s) = \prod_{\phi\in \tilde\Delta} L(E,\phi,s)^{n_{\phi}} \>,
\end{equation}
where $\tilde\Delta$ denotes the set of all irreducible Artin representations
of $\Delta$, and $n_{\phi}$ denotes the dimension of $\phi$. This
completes the proof.
\end{proof}

The following result leads to a useful numerical criterion for deciding
whether or not $E$ is regular over $\Finfty$. Let $M$
be a compact $\Lambda(G)$-module, and $G'$ an open subgroup of $G$. As
always, we say $M$ has finite $G'$-Euler characteristic if the
$H_i(G',M)$ ($i=0,1,2$) are finite, and we then define
$$
\chi(G',M) = \prod_{i=0}^{2} |H_i(G',M)|^{(-1)^{i}} \>.
$$

\begin{proposition}\label{appProp4}
Assume that $X(E/\Finfty)$ belongs to the category $\fM_H(G)$. Then
$E$ is regular over $\Finfty$ if and only if
\begin{equation}\label{appEul1}
\chi(G_K,X(E/\Finfty))=1,    \qquad\qquad {\rm where} \quad G_K=\vGal(\Finfty/K)\>.
\end{equation}
\end{proposition}

\begin{proof}
Assume \eqref{appEul1} holds. We show that $X(E/\Finfty)=0$ by using
the general remarks made in \cite{CFKSV} in paragraph immediately after
the proof of Lemma 3.9. Of course, $G_K$ is pro-$p$, and entirely
analogous arguments to those given in the proof of Lemma 2.5 of \cite{CSS}
show that $H_1(H_K,X(E/\Finfty))=0$, where
$H_K=\vGal(\Finfty/\vQ(\mu_{p^{\infty}}))$.
Finally, it is proven in \cite{HV} that $X(E/\Finfty)$ has no non-zero
pseudo-null $\Lambda(G)$-submodule. Hence, as explained in \cite{CFKSV},
\eqref{appEul1} implies that $X(E/\Finfty)=0$, as required.
\end{proof}

We remark that Proposition \ref{appProp4} is just a slight generalisation
of Theorem \ref{premain}.

The next result enables us to read off many examples when $E$
is regular over $\Finfty$ from the tables in \refAppB.
As earlier, let $\sigma=\sigma_1:G\to\vGL_r(\vQp)$, where $r=p-1$, be the Artin
representation given by the direct sum of the $p-1$ characters of
$\vGal(K/\vQ)$, and let $\cL_E(\sigma)$ be defined by formula
\eqref{conPL} of \S\ref{sec-nonab}.

\begin{corollary}\label{C5}
Assume that the Mazur--Swinnerton-Dyer $p$-adic $L$-function of $E$ over
$K$ lies in $\Lambda(\Gamma_K)$.
If $\cL_E(\sigma)$ is a $p$-adic unit, then $E$ is regular
over $\Finfty$, i.e. $X(E/\Finfty)=0$.
\end{corollary}

\begin{proof}
Assume that $\cL_E(\sigma)$ is a $p$-adic unit. Define
$$
 A_m = \prod_{v|m, v\notdiv p} L_{v}(E/K,1)^{-1} \>,
$$
where the product is taken over all places $v$ of $K$, which divide $m$,
and which do not divide $p$, and where $L_v(E/K,s)$ denotes the Euler
factor at $v$ of the complex $L$-function of $E$ over $K$.
Since $\cL_E(\sigma)$ is a $p$-adic unit, it follows easily that $A_m$
is a $p$-adic unit and the Mazur-Swinnerton-Dyer $p$-adic $L$-function
of $E$ over $K$ is a unit in $\Lambda(\Gamma_K)$. From this latter assertion
and the results of Kato \cite{KatP} and Matsuno \cite{Mat}, it follows that
$X(E/K^{cyc})=0$, whence $X(E/\Finfty)$ must belong to $\fM_H(G)$.
But it is proven in \cite{HV} that, whenever $X(E/K^{cyc})$ has finite
$\Gamma_K$-Euler characteristic, we have
$$
\chi(G_K,X(E/\Finfty)) = \chi(\Gamma_K,X(E/K^{cyc}))\cdot |A_m|^{-1}_p \>.
$$
In particular, we conclude that the right hand side of this formula, and so
also the left hand side, is equal to $1$. Hence the corollary follows
from Proposition \ref{appProp4}.
\end{proof}

Needless to say, it is conjectured that the Mazur--Swinnerton-Dyer
$p$-adic $L$-function of $E$ over $K$ always lies in $\Lambda(\Gamma_K)$.
This is known when the Galois module $E[p]$ is either irreducible, or
contains either $\vZ/p\vZ$ or $\mu_p$ as a Galois submodule (see \cite{Vat}).

The next result was motivated by the congruences studied in this paper,
together with the ``main conjecture'' of \cite{CFKSV}. If
$\phi:G\to\vGL_m(\vZp)$ is an Artin representation of $G$, which we assume,
for simplicity, can be realised over $\vQp$, and $M$ is a compact
$\Lambda(G)$-module, we recall from \cite{CFKSV} that
\begin{equation}
\vtw_{\phi}(M) = M\vt_{\vZp}\vZp^m
\end{equation}
endowed with the obvious left diagonal action of $G$. For each integer
$n\ge 1$, we let $\rho_n$ be the irreducible Artin representation of
$G$ of dimension $p^{n-1}(p-1)$ defined earlier in the paper
(see \S\ref{Sintro}). Let $\sigma_n$ denote the direct sum of the
$p^{n-1}(p-1)$ characters of the Galois group of $\Q(\mu_{p^n})$ over $\Q$.
Now, assuming the Mazur--Swinnerton-Dyer $p$-adic $L$-function of $E$ over $K$
belongs to $\Lambda(\Gamma_k)$, it follows easily that $\cL_E(\sigma)$ is
a $p$-adic unit if and only if $\cL_E(\sigma_n)$ is a $p$-adic unit for
any $n\ge1$. Thus, assuming both the congruence
$$
  \cL_E(\rho_n) \equiv \cL_E(\sigma_n) \mod p
$$
and Corollary 5.10 of the ``main conjecture'' of \cite{CFKSV}, it would
follow that $E$ is regular over $\Finfty$ if and only if
$\chi(G,\vtw_{\rho_n}(X(E/\Finfty)))=1$ for some $n\ge 1$. We now give an
unconditional proof of this last assertion.

\begin{theorem}\label{Th6}
Let $M$ be a $\Lambda(G)$-module such that
{\rm (i)} $M\in\fM_H(G)$,
{\rm (ii)} $H_1(H',M)$ is finite for all open subgroups $H'$ of $H_K$, and
{\rm (iii)} $M$ has no non-zero pseudo-null $\Lambda(G)$-submodule.
Then $M=0$ if and only if there exists an integer $n\ge 1$ such that
$\chi(G,\vtw_{\rho_n}(M))=1$.
\end{theorem}

\begin{corollary}
\label{CorA19}
Assume that $X(E/\Finfty)$ belongs to $\fM_H(G)$. Then $X(E/\Finfty)=0$
if and only if there exists an integer $n\ge 1$ such that
$\chi(G,\vtw_{\rho_n}(X(E/\Finfty)))=1$.
\end{corollary}

Indeed, (iii)
of Theorem \ref{Th6} is proven for $M=X(E/\Finfty)$
in \cite{HV}, and (ii)
is valid by a similar argument to that used
to prove Lemma 2.5 of \cite{CSS}.

Prior to proving Theorem \ref{Th6}, we state a lemma which will be proven
at the end of the Appendix, after a general discussion Akashi series.
If $G'$ is an open subgroup of $G$, let
$H'=H\cap G'$, $\Gamma'=G'/H'$, and write $Q(\Gamma')$ for the field of
quotients of $\Lambda(\Gamma')$. If $M\in\fM_H(G)$, we recall that
$\vAk_{H'}(M)$ is defined to be the image in
$Q(\Gamma')^{\times}/\Lambda(\Gamma')^{\times}$ of $f_0/f_1$, where
$f_0$ (resp. $f_1$) is a characteristic element in $\Lambda(\Gamma')$
for $H_0(H',M)$ (resp. $H_1(H',M)$). For each integer $n\ge 1$, let
\begin{equation}\label{eqFnFn'}
  F_n = \vQ(\mu_{p^n},\sqrt[p^n]{m}), \qquad
  F_n' = \vQ(\mu_{p^n},\sqrt[p^{n-1}]{m}) \>,
\end{equation}
and let $G_n$ (resp. $G_n'$) be the open subgroup of $G$ fixing $F_n$
(resp. $F_n'$). Put
\begin{equation}
  H_n=G_n\cap H , \qquad H_n' = G_n'\cap H \>,
\end{equation}
and note that
\begin{equation}
  \Gamma_n = G_n/H_n = G_n'/H_n'
\end{equation}
is, in fact, the unique closed subgroup of $\Gamma=G/H$ of index $p^{n-1}$.

\begin{lemma}\label{L8}
Assume $M\in\fM_H(G)$. Then, for all $n\ge1$, we have
\begin{equation}
  \vAk_{H_n}(M) = \vAk_{H_n'}(M)\vAk_H(\vtw_{\rho_n}(M))^{p-1}\>.
\end{equation}
Moreover, if in addition $M$ satisfies condition (ii)
of Theorem
\ref{Th6}, then $\vAk_H(\vtw_{\rho_n}(M))$ belongs to $\Lambda(\Gamma)$.
\end{lemma}

\begin{remark}
\label{integrality}
If one accepts that $X(E/\Finfty)$ belongs to $\fM_H(G)$ and
Case 2 of Conjecture 4.8 of \cite{CFKSV}, the above lemma
implies that the $G$-characteristic element of $X(E/\Finfty)$
is in $\Lambda(G)$.
In view of the ``main conjecture'' of \cite{CFKSV}, the same
should hold for the $p$-adic $L$-function $\cL_E$.
\end{remark}

We now prove Theorem \ref{Th6}. Let us pick a topological generator of
$\Gamma$, and so identify $\Lambda(\Gamma)$ with the formal power series
ring $\vZp[[T]]$. We assume that there exists an integer $n\ge 1$ such
that
\begin{equation}\label{eq11}
  \chi(G,\vtw_{\rho_n}(M))=1\>,
\end{equation}
and we must deduce that $M=0$. Put $f_n=\vAk_H(\vtw_{\rho_n}(M))$.
By Lemma \ref{L8}, $f_n$ belongs to $\vZp[[T]]$. Moreover, by the
connexion between Euler characteristics and Akashi series
(\cite{CFKSV}, Theorem 3.6), we have $f_n(0)\in\vZp^{\times}$ by virtue
of \eqref{eq11}. Hence $f_n$ belongs to $\Lambda(\Gamma)^{\times}$.
Let $g_n$ (resp. $g_n'$) denote a characteristic power series for the
$\Lambda(\Gamma_n)$-module $(M)_{H_n}$ (resp. $(M)_{H_n'}$). Since $f_n$
is a unit in $\Lambda(\Gamma)$, and $M$ satisfies condition (ii)
of Theorem \ref{Th6}, it follows from Lemma \ref{L8} that we can assume
that
\begin{equation}\label{eq12}
g_n = g_n' \>.
\end{equation}
We show that \eqref{eq12} forces $M$ to be zero.

Let $M[p^{\infty}]$ denote the $p$-primary submodule of $M$. We first show
that necessarily $M[p^{\infty}]=0$. If $G'$ is a pro-$p$ open subgroup of $G$,
we recall that $\mu_{G'}(M)$ denotes the $\mu$-invariant of $M$ viewed
as a $\Lambda(G')$-module. Since $G_n \subset G_n' \subset G_1$ are
all pro-$p$ open subgroups of $G$, and $[G_n':G_n]=p$, we have
\begin{equation}\label{eq13}
  \mu_{G_n}(M) = p \cdot \mu_{G_n'}(M)\>.
\end{equation}
On the other hand, we claim that
\begin{equation}\label{eq14}
  \mu_{G_n}(M)  = \mu_{\Gamma_n}((M)_{H_n}), \qquad
  \mu_{G_n'}(M) = \mu_{\Gamma_n}((M)_{H_n'})\>.
\end{equation}
We give the proof of the first equation in \eqref{eq14}, and the proof
of the second is entirely similar. Put $Y(M)=M/M[p^{\infty}]$. Then we have the
exact sequence
\begin{equation}\label{eq15}
 H_1(H_n,Y(M)) \to (M[p^{\infty}])_{H_n} \to (M)_{H_n} \to (Y(M))_{H_n} \to 0\>.
\end{equation}
Since $M$ belongs to $\fM_H(G)$ and $H_n$ is open in $H$, the homology
groups $H_i(H_n,Y(M))$ ($i=0,1$) are finitely generated $\vZp$-modules,
whence \eqref{eq15} implies that
\begin{equation}\label{eq16}
  \mu_{\Gamma_n}((M)_{H_n}) = \mu_{\Gamma_n}((M[p^{\infty}])_{H_n}) \>.
\end{equation}
But, since $G_n$ is pro-$p$, a standard argument with the Hochschild-Serre
spectral sequence (see \cite{CSS}) shows that
\begin{equation}\label{eq17}
 \mu_{G_n}(M) = \mu_{\Gamma_n}((M[p^{\infty}])_{H_n}) - \mu_{\Gamma_n}(H_1(H_n,M[p^{\infty}]))\>.
\end{equation}
But $H_1(H_n,M[p^{\infty}])$ injects into $H_1(H_n,M)$,
and thus $H_1(H_n,M[p^{\infty}])$
is finite because of hypothesis (ii)
on $M$. Thus the term on
the extreme right of \eqref{eq17} is 0, and the first equation in
\eqref{eq14} follows on combining \eqref{eq16} and \eqref{eq17}.
If we now combine \eqref{eq12} and \eqref{eq14}, we conclude that
$\mu_{G_n(M)} = \mu_{G_n'}(M)$. Comparing this last equation with
\eqref{eq13}, it follows that $\mu_{G_n}(M)=0$. Thus $M[p^{\infty}]$ must be
pseudo-null as a $\Lambda(G)$-module, and so $M[p^{\infty}]=0$ because of
hypothesis (iii)
on $M$.

Since $M[p^{\infty}]=0$, $M$ is finitely generated over $\Lambda(H)$, and so
also over $\Lambda(H_n)$ and $\Lambda(H_n')$. If $A$ is a ring with
no zero divisors, and $W$ is a finitely generated $A$-module, we write
$r_A(W)$ for the rank of $W$ over $A$. Since $H_n\subset H_n'\subset H_1$
are all pro-$p$ open subgroups of $H$, and $[H_n':H_n]=p$, we have
\begin{equation}\label{eq18}
  r_{\Lambda(H_n)}(M) = p\cdot r_{\Lambda(H_n')}(M)\>.
\end{equation}
On the other hand, we have the well known formula
\begin{equation}\label{eq19}
r_{\Lambda(H_n)(M)} = \sum_{i=0}^{1} (-1)^i r_{\vZp}(H_i(H_n,M)) \>,
\end{equation}
and similarly for $r_{\Lambda(H_n')(M)}$. But then, as $M$ satisfies
hypothesis (ii),
\eqref{eq12} and \eqref{eq19} together
imply that $r_{\Lambda(H_n)}(M)=r_{\Lambda(H'_n)}(M)$. Hence, in view
of \eqref{eq18}, we have $r_{\Lambda(H_n)}(M)=0$. By a result of
Venjakob \cite{Ven}, we then conclude that $M$ is pseudo-null, whence
$M=0$ by hypothesis (iii).
This completes the proof of Theorem \ref{Th6}.

\bigskip

We next discuss some remarkable arithmetic phenomena which arise from
the interplay between root numbers and Iwasawa theory for our false
Tate curve extension $\Finfty$. For each integer $n\ge 1$, we write
$w(E,\rho_n)=\pm 1$ for the sign in the functional equation of
$L(E,\rho_n,s)$. We also write $w(E/K)$ for the sign in the functional
equation for $L(E/K,s)$. It is an important fact, first proven in
\cite{vd}, that the value $w(E,\rho_n)$ is independent of $n$ (it is
assumed in \cite{vd} that the conductor $N_E$ of $E$ is cube free, but
the assertion remains valid for all $E$ provided $m$ is not divisible
by a prime of additive reduction for $E$, see \eqref{eqT10-22} below.)
Suppose that we have
\begin{equation}\label{eqwn-1}
  w(E,\rho_n)=-1 \qquad \forall n\ge 1\>.
\end{equation}
Since all $L(E,\rho_n,s)$ are holomorphic at $s=1$
and $\rho_n$ has dimension $p^{n-1}(p-1)$, it follows immediately on applying
\eqref{eqArtin} to $F_n/\vQ$ that, if \eqref{eqwn-1} is valid, then
\begin{equation}\label{eqN21}
  r_{E/F_n} \ge p^n - 1 + r_{E/K_n} \qquad \forall n\ge 1\>,
\end{equation}
where $F_n$ is given by \eqref{eqFnFn'}
and $K_n=\Q(\mu_{p^n})$ (note that $r_{E/K_n}$ is bounded as $n\to\infty$
by an important theorem of Rohrlich). We identify
$\Lambda(\Gamma_K)$ with $\vZp[[T]]$ by mapping a fixed topological
generator $\gamma_0$ of $\Gamma_K$ to $1+T$. By Kato's theorem,
$X(E/K^{cyc})$ is a torsion $\Lambda(\Gamma_K)$-module, and we write
$t_{E/K}$ for the multiplicity of the zero at $T=0$ of a characteristic
power series of $X(E/K^{cyc})$. We have $t_{E/K}\ge g_{E/K}$, and
conjecturally there is always equality.

\begin{theorem}\label{Th10}
Assume that $X(E/K^{cyc})$ is a finitely generated $\vZp$-module, and that
$t_{E/K}\equiv r_{E/K} \mod 2$. Then \eqref{eqwn-1} is valid if and only
if $X(E/\Finfty)$ has odd $\Lambda(H_K)$-rank. In particular,
\eqref{eqwn-1} implies that $E$ is not regular over $\Finfty$.
\end{theorem}

\begin{proof}
Let $S$ be the set of primes of multiplicative reduction for $E$ which
divide $m$. Then, as was remarked to us by V. Dokchitser, similar arguments
to those used to prove Propositions 9 and 11 of \cite{vd} show that
\begin{equation}\label{eqT10-22}
  w(E,\rho_n) = w(E/K)\prod_{q\in S} \left(\frac{q}{p}\right), \qquad\qquad n=1,2,..
\end{equation}
On the other hand, let $h(E/\Finfty)$ denote that $\Lambda(H_K)$-rank of
$X(E/\Finfty)$. We can determine $h(E/\Finfty)$ modulo $2$ as follows.
Let $\lambda$ denote the $\vZp$-corank of $X(E/K^{cyc})$. For each prime
$q\neq p$, let $s_q$ denote the number of primes of $K$ above $q$. Since
$K/\vQ$ is a cyclic extension of conductor $p$, one sees easily that $s_q$
is odd if and only if $\left(\frac{q}{p}\right)=-1$. We define $S_1$ to be
the subset of $S$ consisting of all primes $q$
dividing $m$ such that $E$ has split multiplicative reduction at all
primes of $K$ above $q$. Let $U$ denote the set of all primes $q\neq p$
satisfying (i) $q$ divides $m$, (ii) E has good reduction at $q$, and
(iii) $p$ divides the order of $\tilde{E}_q(k_v)$ for each prime $v$ of
$K$ above $q$, where $\tilde{E}_q$ denotes the reduction of $E$ modulo $q$,
and $k_v$ denotes the residue field at $v$. Then a standard argument
applying the formula \eqref{eq19} for $H_K$ and the fact that
$H_1(H_K,X(E/\Finfty))=0$ shows that
\begin{equation}\label{eqT10-23}
 h(E/\Finfty) = \lambda + \sum_{q\in S_1}s_q + 2 \sum_{q\in U}s_q
\end{equation}
(see \cite{HV}, Theorem 3.1 for an alternative proof when $p\ge 5$).
Define $S_2$ to be the set of all primes $q$ of multiplicative reduction
such that $\left(\frac{q}{p}\right)=-1$. Clearly $S_2\subset S_1$ since
a $q$ in $S_2$ is inert in the quadratic subfield of $K$, and so has
split multiplicative reduction at all primes of $K$ above $q$. Hence we have
\begin{equation}\label{eqT10-24}
 \sum_{q\in S_1}s_q \equiv \sum_{q\in S_2}1 \mod 2\>.
\end{equation}
To determine the parity of $\lambda$, we use Greenberg's theorem
\cite{Gre2} which asserts that the characteristic ideal of $X(E/K^{cyc})$
is invariant under the involution of $\Lambda(\Gamma_K)$ which sends
each $\gamma$ in $\Gamma_K$ to $\gamma^{-1}$. Hence, if $f(T)$ denotes
the monic distinguished polynomial which generates the characteristic
ideal of $X(E/K^{cyc})$, and if $\alpha$ is any non-zero root of $f(T)$
lying in the algebraic closure of $\vQp$, then $\frac{1}{1+\alpha}-1$
must also be a root distinct from $\alpha$. Hence the degree $\lambda$
of $f(T)$ must satisfy
\begin{equation}\label{eqT10-25}
 \lambda\equiv t_{E/K} \mod 2\>.
\end{equation}
Combining \eqref{eqT10-22},\eqref{eqT10-23},\eqref{eqT10-24} and
\eqref{eqT10-25}, and recalling the hypothesis that
$t_{E/K}\equiv r_{E/K} \mod 2$, we conclude that \eqref{eqwn-1} holds if
and only if $h(E/\Finfty)$ is odd. This completes the proof of
Theorem \ref{Th10}.
\end{proof}

\begin{proposition}\label{P10}
Assume that $X(E/K^{cyc})$ is a finitely generated $\vZp$-module, and take
$m=q$, where $q$ is a prime of multiplicative reduction for $E$. Then
$X(E/\Finfty)$ has $\Lambda(H_K)$-rank 1 if and only if either
{\rm (i)} $q$ is inert in $K$, and $X(E/K^{cyc})=0$, or
{\rm (ii)} $E$ has non-split multiplicative reduction at the primes of $K$
above $q$, and $X(E/K^{cyc})=\vZp$, with trivial action of $\Gamma_K$.
\end{proposition}

\begin{proof}
The sufficiency follows from \eqref{eqT10-23} above. Conversely,
assume that $X(E/\Finfty)$ has $\Lambda(H_K)$-rank 1. Since
$X(E/K^{cyc})$ is assumed to be a finitely generated $\vZp$-module,
Matsuno's theorem \cite{Mat} shows that $X(E/K^{cyc})$ has no non-zero
finite $\Gamma_K$-submodule. Hence \eqref{eqT10-23} implies that
either
(i) there is a single prime of $K$ above $q$ and $X(E/K^{cyc})=0$, or
(ii) $E$ has non-split multiplicative reduction at the primes of $K$
above $q$, and $X(E/K^{cyc})$ is a free $\vZp$-module of rank 1.
In the second case, the same argument with Greenberg's theorem as
in the proof of Theorem \ref{Th10} shows that $\Gamma_K$ must act trivially
on $X(E/K^{cyc})$.
\end{proof}

We now deal with the two different cases occurring in Proposition \ref{P10}
separately. Recall that $F_n=\vQ(\mu_{p^n},\sqrt[p^n]{q})$.

\begin{theorem}\label{Th11}
Assume that $X(E/K^{cyc})=0$ and $r_{E/K}=0$. Take $m=q$, where $q$ is a
prime of multiplicative reduction for $E$, which is inert in $K^{cyc}$.
Then we have
\begin{equation}\label{eq21}
  g_{E/F_n}\le p^n -1 \le r_{E/F_n}, \qquad \forall n\ge 1\>.
\end{equation}
Moreover, if we assume that $g_{E/F_n}=r_{E/F_n}$ for a given integer
$n\ge 1$, then
\def\theenumi{\roman{enumi}}
\begin{enumerate}
\item\label{Th11i} $g_{E/F_n} = g_{E/F_n^{cyc}} = p^n -1$,
\item\label{Th11ii} $\sha(E/F_n)[p^{\infty}]$ is finite,
\item\label{Th11iii} $\sha(E/F_n^{cyc})[p^{\infty}]=0$,
\item\label{Th11iv} $L(E,\rho_k,s)$ has a simple zero at $s=1$ for $1\le k\le n$, and
\item\label{Th11v} $r_{E/\vQ(\sqrt[p^n]{q})}=n$.
\end{enumerate}
\end{theorem}

In particular, under the hypotheses of Theorem \ref{Th11}, we see that if
we assume the Birch--Swinnerton-Dyer conjecture $g_{E/F_n}=r_{E/F_n}$ for
all $n\ge 1$, then $\sha(E/\Finfty)[p^{\infty}]=0$, and $X(E/\Finfty)$
is dual to $E(\Finfty)\vt\vQp/\vZp$.
Remarkably, if one assumes the hypothesis of Theorem \ref{Th11} and, in
addition, that $E$ has prime conductor, H. Darmon and Y. Tian have informed
us that they can prove that $g_{E/F_n}=r_{E/F_n}$ for all $n\ge 1$.
We also mention that the calculations
carried out in this paper give numerical examples of Theorem \ref{Th11},
for example
\begin{itemize}
\item[] $E=11A3,  p=3,  q=11$,
\item[] $E=38B1,  p=3,  q=2$,
\item[] $E=21A4,  p=5,  q=3$ or $q=7$,
\item[] $E=24A4,  p=5,  q=3$,
\item[] $E=26A1,  p=5,  q=7$,
\item[] $E=84B1,  p=5,  q=3$ or $q=7$,
\item[] $E=17A1,  p=7,  q=17$.
\end{itemize}
The hypothesis that $X(E/K^{cyc})=0$ in each of these cases can be
verified by noting from the tables that $\cL_E(\sigma)$ is a $p$-adic
unit for some choice of $m$ (cf. the proof of Corollary \ref{C5}).

\bigskip

We now prove Theorem \ref{Th11}. By Theorem \ref{Th10} and Proposition
\ref{P10}, we have $w(E,\rho_n)=-1$ for all $n\ge 1$. Hence, by
\eqref{eqN21} and the fact that $r_{E/K}=0$, we have $r_{E/F_n}\ge p^n -1$
for all $n\ge 1$; and if there is equality for a given $n$, $L(E,\rho_k,s)$
must have a simple zero at $s=1$ for $1\le k\le n$. Also we have
\begin{equation}
 L(E/\vQ(\sqrt[p^n]{q}),s) = L(E/\vQ,s)\prod_{k=1}^{n} L(E,\rho_k,s)\>,
\end{equation}
whence assertion \eqref{Th11v} is then clear.

To establish the upper bound for $g_{E/F_n}$ in
\eqref{eq21}, we apply the theorem of Hachimori-Matsuno \cite{HM}
to the Galois extension $F_n^{cyc}/K^{cyc}$ of
degree $p^n$. In this extension, $E$ has split multiplicative reduction
at the unique prime of $K^{cyc}$ above $q$, which is totally ramified,
and no other prime of $K^{cyc}$ not dividing $p$ ramifies. Hence
$X(E/F_n^{cyc})$ is a finitely generated $\vZp$-module of rank $p^n -1$,
and, by Matsuno's theorem \cite{Mat}, it is a free $\vZp$-module. All the
remaining assertions of Theorem \ref{Th11} are now clear if we note
that the restriction map from the $p^{\infty}$-Selmer group of $E$
over $F_n$ to the $p^{\infty}$-Selmer group of $E$ over $F_n^{cyc}$
has finite kernel. This completes the proof.

\begin{theorem}\label{Th12}
Assume that $r_{E/K}=1$, and that $X(E/K^{cyc})=\vZp$, with trivial
action of $\Gamma_K$. Take $m=q$, where $q$ is a prime number such that
$E$ has non-split multiplicative reduction at the primes of $K$ above $q$.
Then we have
\begin{equation}\label{eq25}
 g_{E/F_n} \le p^n \le r_{E/F_n} \qquad \forall n\ge 1\>.
\end{equation}
Moreover, if we assume that $g_{E/F_n}=r_{E/F_n}$ for a given integer
$n\ge 1$, then
\def\theenumi{\roman{enumi}}
\begin{enumerate}
\item\label{Th12i} $g_{E/F_n} = g_{E/F_n^{cyc}} = p^n$,
\item\label{Th12ii} $\sha(E/F_n)[p^{\infty}]$ is finite,
\item\label{Th12iii} $\sha(E/F_n^{cyc})[p^{\infty}]=0$,
\item\label{Th12iv} $L(E,\rho_k,s)$ has a simple zero at $s=1$ for $1\le k\le n$, and
\item\label{Th12v} $r_{E/\vQ(\sqrt[p^n]{q})}=n+r_{E/\vQ}$.
\end{enumerate}
\end{theorem}

Again, under the hypotheses of Theorem \ref{Th12}, we see that if we assume
that $g_{E/F_n}=r_{E/F_n}$ for all $n\ge 1$, then
$\sha(E/\Finfty)[p^{\infty}]=0$, and $X(E/\Finfty)$ is dual to
$E(\Finfty)\vt\vQp/\vZp$. We believe that there must be many numerical
examples of curves $E$ satisfying the hypotheses of Theorem \ref{Th12}.
However, the condition that $X(E/K^{cyc})$ is a free $\vZp$-module of
rank 1 is more delicate to verify. We are very grateful to C. Wuthrich
for providing us with the following example. Let $E$ be the elliptic
curve $79A1$ of Cremona's tables given by
\begin{equation}
  y^2 + xy + y = x^3 + x^2 - 2x \>,
\end{equation}
which has non-split multiplicative reduction at $q=79$. Moreover, $E$
has good ordinary reduction at $p=3$, $r_{E/K}=1$, and Wuthrich has
shown that $E(K^{cyc})$ is a free $\vZp$-module generated by $P=(0,0)$,
and that $X(E/K^{cyc})$ is a free $\vZ_3$-module of rank 1 (we omit
his proof). Finally, 79 splits in $K$, so that $E$ has non-split
multiplicative reduction at both primes of $K$ above 79. Thus all
the hypotheses of Theorem \ref{Th11} are valid in this case with $p=3$
and $q=79$.

\bigskip

The proof of Theorem \ref{Th12} is entirely parallel to that of
Theorem \ref{Th11}. As $r_{E/K}=1$, \eqref{eqN21}
shows that $r_{E/F_n}\ge p^n$, and that if there is equality the
$L(E,\rho_k,s)$ must have a simple zero at $s=1$ for $1\le k\le n$.
On the other hand, as $X(E/K^{cyc})$ is a free $\vZp$-module of rank 1
and $E$ has non-split multiplicative reduction at the primes of $K$
above $q$, the theorems of Hachimori-Matsuno and Matsuno applied to
the extension $F_n^{cyc}/K_n^{cyc}$ of degree $p^n$ show that
$X(E/F_n^{cyc})$ is a free $\vZp$-module of rank $p^n$. The assertions
of Theorem \ref{Th12} now follow as in the proof of Theorem \ref{Th11}.

\medskip

We remark that variants of Theorems \ref{Th11} and \ref{Th12} hold for
certain other choices of the integer $m$ defining the false Tate curve
extension. For example, with the same hypotheses and the same
choice of a prime $q$ of multiplicative reduction of $E$ as before,
the conclusions of Theorems \ref{Th11} and \ref{Th12} remain valid if
we take $m=m'q$, where $m'$ is a $p$-power free integer all of whose
prime factors are either $p$ or primes $l$ such that $E$ has good
ordinary reduction modulo $l$, and the order of $\tilde{E}_l(k_v)$
is prime to $p$; where $\tilde{E}_l$ denotes the reduction of $E$
modulo $l$, $v$ denotes any prime of $K$ above $l$, and $k_v$ is the
residue field of $v$. Similarly, if we assume that $r_{E/K}=1$
and that $X(E/K^{cyc})=\vZp$ with trivial action of $\Gamma_K$, and now
take $m=m'$ as just defined, the conclusions of Theorem \ref{Th12}
again hold. The proofs are entirely analogous to those given before.

We end this appendix by establishing an analogue of the Artin formalism
for Akashi series for an arbitrary compact $p$-adic Lie group $G$ with a
closed normal subgroup $H$ such that $\Gamma=G/H$ is isomorphic to $\Z_p$.
Let $G'$ be an open normal subgroup of $G$, and put
$$
  H'=H\cap G', \qquad \Gamma'=G'/H',
$$
so that there is a natural inclusion of $\Gamma'$ as an open subgroup
in $\Gamma$. Let $\Delta=G/G'$, and write $L$ for some fixed finite extension
of $\Q_p$ such that all absolutely irreducible representations of $\Delta$
can be realised over $L$. Let $\cO$ denote the ring of integers of $L$,
and write $\Lambda_\cO(\Gamma)$ for the Iwasawa algebra of $\Gamma$ with
coefficients in $\cO$. If $M$ is a compact $\Lambda(G)$-module, we define
$M_\cO=M\tensor_{\vZp}\cO$ Assume now that $G$ has no element of order $p$.
For $M\in\fM_H(G)$, we recall (see \S3 of \cite{CFKSV}) that $\vAk_H(M_\cO)$
(resp. $\vAk_{H'}(M_\cO)$) is defined to be the image in
$Q_\cO(\Gamma)^\times/\Lambda_\cO(\Gamma)^\times$
(resp. $Q_\cO(\Gamma')^\times/\Lambda_\cO(\Gamma')^\times$)
of the alternating product of the characteristic elements in
$\Lambda_\cO(\Gamma)$ (resp. in $\Lambda_\cO(\Gamma')$) of the
$H_i(H,M_\cO)$ (resp. $H_i(H',M_\cO)$) for all $i\ge 0$.
Again, $\hat\Delta$ will denote the set of all irreducible representations
of $\Delta$, and we write $n_\rho$ for the dimension of $\rho$ in
$\hat\Delta$.

\begin{theorem}
\label{Th14}
For each $M$ in $\fM_H(G)$, we have
\begin{equation}
\label{eq30}
  \vAk_{H'}(M_\cO)^{[\Gamma:\Gamma']} =
  \prod_{\rho\in\hat\Delta} \vAk_H(\vtw_\rho(M_\cO))^{n_\rho}\>.
\end{equation}
\end{theorem}

Before proving Theorem \ref{Th14}, we note that Lemma \ref{L8}
follows on applying it to the two Galois extensions $F_n/\Q$ and $F_n'/\Q$.
Since the set of all irreducible representations of $\Gal(F_n/\Q)$
consists of $\rho_n$ together with all irreducible representations of
$\Gal(F_n'/\Q)$, we conclude that, for a suitable choice of $L$, we have
$$
  \vAk_{H_n}(M_\cO)^{p^{n-1}} =
  \vAk_{H_n'}(M_\cO)^{p^{n-1}}
  \vAk_{H}(\vtw_{\rho_n}(M_\cO))^{p^{n-1}(p-1)}.
$$
We then obtain Lemma \ref{L8} on extracting $p^{n-1}$-th roots of both sides,
taking the norm map from $\Q_\cO(\Gamma)^\times$ to $\Q(\Gamma)^\times$,
and recalling that $\Lambda(\Gamma)$ is a unique factorisation domain.

We end this Appendix by proving Theorem \ref{Th14}. Let $K_0(\fM_H(G))$
be the Grothendieck group of the category $\fM_H(G)$, and write $[M]$
for the class of a module in this group. As usual, we define
$$
  \vInd_{G'}^G(M_\cO) = \Lambda(G)\tensor_{\Lambda(G')}M_\cO = R\tensor_\cO M_\cO,
$$
where $R=\cO[\Delta]$ is the $\cO$-group ring of $\Delta$.
Let $\cC_R$ be the category of all finitely generated $R$-modules, and
let $K_0(\cC_R)$ be the Grothendieck group in this category.
For each irreducible representation $\rho$ of $\Delta$, let
$L_\rho$ be a free $\cO$-module of rank $n_\rho$ realising $\rho$,
and put
$$
  W = \bigoplus_{\rho\in\hat\Delta} {L_\rho}^{n_\rho}.
$$
Since $W\tensor_\cO L$ is isomorphic to $R\tensor_\cO L$ as
$L[\Delta]$-modules, a theorem of Swan (\cite{SwaK}, Theorem 3) implies that
\begin{equation}
\label{eq31}
[W]=[R] \quad\text{in}\quad K_0(\cC_R).
\end{equation}
But Swan (\cite{SwaK}, Theorem 1.2) has shown that the natural inclusion
of the category $\cD_R$ of all finitely generated $R$-modules which are
$\cO$-free in $\cC_R$ induces an isomorphism from $K_0(\cD_R)$ to $K_0(\cC_R)$.
As $W$ and $R$ are $\cO$-free, it follows from \eqref{eq31} that there
exist a finite number of exact sequences
\begin{equation}
\label{eq32}
0 \lar A_i \lar B_i \lar C_i \lar 0 \qquad (i=1,\ldots,n)
\end{equation}
in $\cD_R$ such that
$$
  R-W = \sum_{i=1}^n (B_i-A_i-C_i)
$$
in the free abelian group on the isomorphism classes of $\cD_R$.
Since the modules in the exact sequence \eqref{eq32} are $\cO$-free,
it remains exact when we tensor it over $\cO$ with any $\cO$-module.
It follows that, for each $M$ in $\fM_H(G)$, we must have
\begin{equation}
\label{eq33}
[R\tensor_\cO M_\cO] = [W\tensor_\cO M_\cO]
\end{equation}
in $K_0(\fM_H(G))$. Now $\vAk_H(\cdot)$ is well defined on
$K_0(\fM_H(G))$. By Shapiro's lemma
$$
  H_i(H,R\tensor M_\cO) = \vTor^{\Lambda(G)}(\Lambda(G/H),
    \vInd_{G'}^G(M_\cO)) = H_i(H',M_\cO)^{[\Gamma:\Gamma']},
$$
and so $\vAk_H(R\tensor_\cO M_\cO)$ is the left hand side of \eqref{eq30}.
But $\vAk_H(W\tensor_\cO M_\cO)$ is, by definition, the right hand
side of \eqref{eq33}. Hence Theorem \ref{Th14} follows from \eqref{eq33}.

Finally, we remark without proof (see Theorem 6.8 of \cite{AW}
for the case of $p$-primary modules in $\fM_H(G)$) that \eqref{eq33}
can also be used to establish the Artin formalism for the
characteristic elements of any module in $\fM_H(G)$.


\def\theequation{B.\arabic{equation}}
\section*{Appendix B. Tables}

Take an odd prime $p$ and an elliptic curve $E/\Q$
with good ordinary reduction at $p$. Define $K=\Q(\mu_p)$ and let $\sigma$
be the regular representation of $\Gal(K/\Q)$. Assume that $E/K$ has Mordell-Weil
rank 0.

In each table we list the prime $p$, the name (as in Cremona's tables \cite{Cre})
and the equation of $E$, the
local polynomial $P_p(E/\Q,T)$ of $E$ at $p$ (see \ref{SSLfEcNf}),
the conductors $N(E)$ and $N(E,\sigma)$ (see \ref{ssCond}) and the orders of
the torsion groups $E(\Q)$ and $E(K)$. Next, we compute
numerically the Birch--Swinnerton-Dyer quotients
$$
  L^*(E/\Q) = \frac{L(E,1)}{\Omega_+({E})}, \qquad
  L^*(E/K) = \left|\frac{L(E,\sigma,1)\sqrt{\Delta_K}}{(2\Omega_+({E})\Omega_-({E}))^{(p-1)/2}}\right|
$$
and the analytic orders of the Tate-Shafarevich groups $\sha(E/\Q)$
and $\sha(E/K)$. We also list the reduction types of $E/\Q_l$ and $E/K_v$
at bad primes $l\in\Q$ and primes $v|l$ of $K$.

\newpage

Finally, for varying $m$ that defines the 
extension $\Finfty=\Q(\mu_{p^\infty},\sqrt[p^\infty]{m})/\Q$,
we tabulate the following data:

\medskip

\begin{tabular}{ll}
\multicolumn{2}{l}{\bf Table columns:}\cr
$m$ & $p$-th power free integer, $m>1$. It specifies the false Tate curve \cr
& extension $\Finfty=\Q(\mu_{p^\infty},\sqrt[p^\infty]{m})/\Q$. (When $m$ and $m'$ define the \cr
& same field, e.g. 20 and 50 for $p=3$, we take the smaller one.)\cr
$N(\rho)$ & conductor of $\rho$
           (the representation of $\Gal(\vQ(\mu_p,\sqrt[p]{m})/\Q)$ in \S\ref{Sintro}).\cr
$N(E,\rho)$ & conductor of the twist of $E$ by $\rho$. \cr
$L^*$       & $\displaystyle\frac{L^*(E/\Q(\sqrt[p]{m}))}{L^*(E/\Q)} =
               \left|\frac{L(E,\rho,1)\sqrt{\Delta_{\Q(\sqrt[p]m)}}}{(2\Omega_+({E})\Omega_-({E}))^{(p-1)/2}}\>\right|.$\cr
$\sha$      & Analytic order of $\sha(E/\Q(\sqrt[p]{m}))$. We write ``--'' if the analytic \cr
            & rank of $E/\Q(\sqrt[p]m)$ is non-zero.\cr
$\cL_E(\sigma)$ & The quantity in \eqref{conPL}. Conjecturally, this is the value of the \cr
   & non-abelian $p$-adic $L$-function for $\Finfty/\Q$ at $\sigma$.\cr
$\cL_E(\rho)$ & The quantity in \eqref{conPL}. Conjecturally, this is the value of the \cr
   & non-abelian $p$-adic $L$-function for $\Finfty/\Q$ at $\rho$.\cr
\end{tabular}

\medskip
\noindent
We do not include $m$ for which $L(E,\rho,s)$ has sign -1 in the functional
equation.


In the table below we give a list of our tables specifying
whether the curve has
one of the contributions in (1), (2) and (3) of \S\ref{Sobs},
and whether the curve is semistable.
We also list what we can say concerning the $\mu$-invariant $\mu_{E/K}$.
We write ``0'' if we can prove that $\Selpi(E/K)$ is finite and $\mu_{E/K}=0$.
We write ``$0^?$'' if we can deduce this from the Birch--Swinnerton-Dyer
conjecture, and write ``?'' otherwise.

\bigskip

{
\smaller[3]
\def\X{$\times$}
\def\O{$\cdot$}
\input{tabinfo.dat}
\medskip
}




\begingroup

\smaller[2]

\parindent 0cm

\def\tabbegin{\pagebreak[2]\begingroup
\baselineskip 8pt
\def\arraystretch{0.8}
\font\txtfnt=cmss8 scaled 857\textfont0=\txtfnt
\font\mthfnt=cmmi7\textfont1=\mthfnt
\par\vskip4pt\noindent$
  \begin{array}{|l|l|l|l|l|l|l|l|}
  \hline
      \vphantom{\int^{X}_X}
      m &
      N(\rho) &
      N(E,\rho) &
      L^* &
      \smallsha &
      \cL_E(\sigma)&
      \cL_E(\rho)
      \cr
  \hline\vphantom{2^{X^{X^a}}}
}
\def\tabend{\end{array}$\endgroup\par\bigskip\vfill\pagebreak[3]}
\def\tabentry#1#2#3#4#5#6#7#8#9{#1&#4&#5&#6&#7&#8&#9}
\def\errentry#1#2{}
\def\noeuler{0}
\def\nosha{-}
\def\extravspace{$\vphantom{\int_{\int}}$}
\def\alsub#1{_{\rlap{$\scriptstyle #1$}}}
\let\sha\smallsha
\input{tables.dat}

\endgroup

\endcomment

\newpage


\end{document}

%% file: tabinfo.dat
\begin{tabular}{|clcccccc|}
\hline
& Table & $E(K)[p]$ & $E(\F_p)[p]$ & $p|\prod_v c_v$ & semistable & $\mu_{E/K}$ 
& page \cr
\hline
p=3 & 3-11A3 & \O & \O & \O & \X & $0^{\phantom{?}}$ & \pageref{3-11A3}\cr
    & 3-20A3 & \X & \X & \X & \O & $?^{\phantom{?}}$ & \pageref{3-20A3}\cr
    & 3-26A1 & \X & \X & \X & \X & $?^{\phantom{?}}$ & \pageref{3-26A1}\cr
    & 3-35A1 & \X & \X & \X & \X & $?^{\phantom{?}}$ & \pageref{3-35A1}\cr
    & 3-38B1 & \O & \O & \O & \X & $0^{\phantom{?}}$ & \pageref{3-38B1}\cr
    & 3-50B1 & \O & \O & \O & \O & $0^?$  & \pageref{3-50B1}\cr
    & 3-56B1 & \O & \O & \O & \O & $0^{\phantom{?}}$ & \pageref{3-56B1}\cr
    & 3-77C1 & \O & \O & \O & \X & $0^{\phantom{?}}$ & \pageref{3-77C1}\cr
    & 3-80B2 & \O & \O & \O & \O & $0^?$  & \pageref{3-80B2}\cr
    & 3-92A1 & \X & \X & \X & \O & $0^?$  & \pageref{3-92A1}\cr
    & 3-116C2 & \O & \O & \X & \O & $0^?$  & \pageref{3-116C2}\cr
    & 3-128B2 & \O & \X & \O & \O & $?^{\phantom{?}}$ & \pageref{3-128B2}\cr
    & 3-152B1 & \O & \X & \O & \O & $?^{\phantom{?}}$ & \pageref{3-152B1}\cr
    & 3-176B1 & \O & \X & \O & \O & $?^{\phantom{?}}$ & \pageref{3-176B1}\cr
    & 3-224B1 & \O & \O & \O & \O & $0^{\phantom{?}}$ & \pageref{3-224B1}\cr
    & 3-260A1 & \O & \O & \X & \O & $?^{\phantom{?}}$ & \pageref{3-260A1}\cr
    & 3-272C1 & \O & \O & \O & \O & $0^{\phantom{?}}$ & \pageref{3-272C1}\cr
    & 3-275B1 & \O & \X & \O & \O & $?^{\phantom{?}}$ & \pageref{3-275B1}\cr
    & 3-395B1 & \O & \O & \X & \X & $?^{\phantom{?}}$ & \pageref{3-395B1}\cr
    & 3-800E1 & \O & \X & \X & \O & $?^{\phantom{?}}$ & \pageref{3-800E1}\cr
\hline
p=5 & 5-11A3 & \X & \X & \O & \X & $0^{\phantom{?}}$ & \pageref{5-11A3}\cr
    & 5-19A3 & \O & \O & \O & \X & $0^?$  & \pageref{5-19A3}\cr
    & 5-21A4 & \O & \O & \O & \X & $0^?$  & \pageref{5-21A4}\cr
    & 5-24A4 & \O & \O & \O & \O & $0^?$  & \pageref{5-24A4}\cr
    & 5-44A1 & \O & \O & \O & \O & $0^?$  & \pageref{5-44A1}\cr
    & 5-56A1 & \O & \O & \O & \O & $0^?$  & \pageref{5-56A1}\cr
    & 5-56B1 & \O & \X & \O & \O & $?^{\phantom{?}}$ & \pageref{5-56B1}\cr
    & 5-84B1 & \O & \O & \O & \O & $0^?$  & \pageref{5-84B1}\cr
\hline
p=7 & 7-17A1 & \O & \O & \O & \X & $0^?$  & \pageref{7-17A1}\cr
    & 7-19A3 & \O & \O & \O & \X & $0^?$  & \pageref{7-19A3}\cr
\hline
\end{tabular}

%% file: tables.dat
\noindent\begin{tabular}{cllll} 
\extravspace
$E=11A3$ & \multicolumn{3}{l}{$y^2+y=x^3-x^2$} & \tabname{3}{11A3}\cr
\extravspace
$p=3$ & \multicolumn{3}{l}{$P_3(E/\Q,T)=1+T+3T^2$}\cr
$/\vQ\alsub{}$   & $N(E)=11$ & $L^*(E/\Q)=5^{-2}$ & $|E(\Q)|=5$ & $|\sha(E/\Q)|=1^2$ \cr
$/K\alsub{}$      & $N(E,\sigma)=3^{2}\!\cdot\!11^{2}$ & $L^*(E/K)=5^{-2}$ & $|E(K)|=5$ & $|\sha(E/K)|=1^2$ \cr
$/\vQ\alsub{11}$ & $\hbox{I}_1 \hbox{ split},\>c_{11}=1$ & $/K\alsub{v}\quad(v|11)$ & $\hbox{I}_1 \hbox{ split},\>c_v=1$\cr
\end{tabular}
\tabbegin
\tabentry{2}{[0, -1, 1, 0, 0]}{x^3-2}{2^{2}3^{3}}{2^{4}3^{6}11^{2}}{1}{1^2}{1\!\cdot\!3^{0}\phantom{\,\,\!+\!2\!\cdot\!3^1}\,\!+\!2\!\cdot\!3^{2}\!+\!O(3^{3})}{1\!\cdot\!3^{0}\phantom{\,\,\!+\!2\!\cdot\!3^1}\,\!+\!1\!\cdot\!3^{2}\!+\!O(3^{3})}\errentry{4\!\cdot\!10^{-17}}{2\!\cdot\!10^{-15}}\cr
\tabentry{3}{[0, -1, 1, 0, 0]}{x^3-3}{3^{5}}{3^{10}11^{2}}{1}{1^2}{2\!\cdot\!3^{0}\!+\!2\!\cdot\!3^{1}\!+\!2\!\cdot\!3^{2}\!+\!O(3^{6})}{2\!\cdot\!3^{0}\!+\!1\!\cdot\!3^{1}\!+\!1\!\cdot\!3^{2}\!+\!O(3^{3})}\errentry{1\!\cdot\!10^{-16}}{2\!\cdot\!10^{-15}}\cr
\tabentry{5}{[0, -1, 1, 0, 0]}{x^3-5}{3^{3}5^{2}}{3^{6}5^{4}11^{2}}{2^{2}}{2^2}{1\!\cdot\!3^{0}\!+\!1\!\cdot\!3^{1}\phantom{\,\,\!+\!2\!\cdot\!3^2}\,\!+\!O(3^{4})}{1\!\cdot\!3^{0}\!+\!2\!\cdot\!3^{1}\!+\!1\!\cdot\!3^{2}\!+\!O(3^{3})}\errentry{2\!\cdot\!10^{-16}}{6\!\cdot\!10^{-16}}\cr
\tabentry{6}{[0, -1, 1, 0, 0]}{x^3-6}{2^{2}3^{5}}{2^{4}3^{10}11^{2}}{1}{1^2}{1\!\cdot\!3^{0}\phantom{\,\,\!+\!2\!\cdot\!3^1}\,\!+\!2\!\cdot\!3^{2}\!+\!O(3^{3})}{1\!\cdot\!3^{0}\!+\!2\!\cdot\!3^{1}\phantom{\,\,\!+\!2\!\cdot\!3^2}\,\!+\!O(3^{3})}\errentry{5\!\cdot\!10^{-17}}{9\!\cdot\!10^{-17}}\cr
\tabentry{7}{[0, -1, 1, 0, 0]}{x^3-7}{3^{3}7^{2}}{3^{6}7^{4}11^{2}}{1}{1^2}{2\!\cdot\!3^{0}\phantom{\,\,\!+\!2\!\cdot\!3^1}\,\!+\!2\!\cdot\!3^{2}\!+\!O(3^{4})}{2\!\cdot\!3^{0}\!+\!2\!\cdot\!3^{1}\!+\!2\!\cdot\!3^{2}\!+\!O(3^{3})}\errentry{2\!\cdot\!10^{-16}}{6\!\cdot\!10^{-16}}\cr
\tabentry{10}{[0, -1, 1, 0, 0]}{x^3-10}{2^{2}3\!\cdot\!5^{2}}{2^{4}3^{2}5^{4}11^{2}}{1}{1^2}{2\!\cdot\!3^{0}\!+\!1\!\cdot\!3^{1}\phantom{\,\,\!+\!2\!\cdot\!3^2}\,\!+\!O(3^{4})}{2\!\cdot\!3^{0}\phantom{\,\,\!+\!2\!\cdot\!3^1}\,\!+\!1\!\cdot\!3^{2}\!+\!O(3^{3})}\errentry{3\!\cdot\!10^{-17}}{1\!\cdot\!10^{-15}}\cr
\tabentry{12}{[0, -1, 1, 0, 0]}{x^3-12}{2^{2}3^{5}}{2^{4}3^{10}11^{2}}{2^{2}}{2^2}{1\!\cdot\!3^{0}\phantom{\,\,\!+\!2\!\cdot\!3^1}\,\!+\!2\!\cdot\!3^{2}\!+\!O(3^{3})}{1\!\cdot\!3^{0}\phantom{\,\,\!+\!2\!\cdot\!3^1}\,\phantom{\,\,\!+\!2\!\cdot\!3^2}\,\!+\!O(3^{3})}\errentry{3\!\cdot\!10^{-16}}{6\!\cdot\!10^{-16}}\cr
\tabentry{13}{[0, -1, 1, 0, 0]}{x^3-13}{3^{3}13^{2}}{3^{6}11^{2}13^{4}}{1}{1^2}{2\!\cdot\!3^{0}\!+\!1\!\cdot\!3^{1}\phantom{\,\,\!+\!2\!\cdot\!3^2}\,\!+\!O(3^{3})}{2\!\cdot\!3^{0}\!+\!1\!\cdot\!3^{1}\!+\!1\!\cdot\!3^{2}\!+\!O(3^{3})}\errentry{8\!\cdot\!10^{-17}}{4\!\cdot\!10^{-17}}\cr
\tabentry{14}{[0, -1, 1, 0, 0]}{x^3-14}{2^{2}3^{3}7^{2}}{2^{4}3^{6}7^{4}11^{2}}{1}{1^2}{1\!\cdot\!3^{0}\!+\!2\!\cdot\!3^{1}\!+\!2\!\cdot\!3^{2}\!+\!O(3^{4})}{1\!\cdot\!3^{0}\!+\!1\!\cdot\!3^{1}\!+\!1\!\cdot\!3^{2}\!+\!O(3^{3})}\errentry{2\!\cdot\!10^{-17}}{1\!\cdot\!10^{-15}}\cr
\tabentry{15}{[0, -1, 1, 0, 0]}{x^3-15}{3^{5}5^{2}}{3^{10}5^{4}11^{2}}{1}{1^2}{1\!\cdot\!3^{0}\!+\!1\!\cdot\!3^{1}\phantom{\,\,\!+\!2\!\cdot\!3^2}\,\!+\!O(3^{4})}{1\!\cdot\!3^{0}\phantom{\,\,\!+\!2\!\cdot\!3^1}\,\!+\!2\!\cdot\!3^{2}\!+\!O(3^{3})}\errentry{2\!\cdot\!10^{-16}}{3\!\cdot\!10^{-16}}\cr
\tabentry{17}{[0, -1, 1, 0, 0]}{x^3-17}{3\!\cdot\!17^{2}}{3^{2}11^{2}17^{4}}{1}{1^2}{1\!\cdot\!3^{0}\!+\!1\!\cdot\!3^{1}\!+\!1\!\cdot\!3^{2}\!+\!O(3^{3})}{1\!\cdot\!3^{0}\!+\!2\!\cdot\!3^{1}\!+\!1\!\cdot\!3^{2}\!+\!O(3^{3})}\errentry{5\!\cdot\!10^{-16}}{2\!\cdot\!10^{-16}}\cr
\tabentry{19}{[0, -1, 1, 0, 0]}{x^3-19}{3\!\cdot\!19^{2}}{3^{2}11^{2}19^{4}}{1}{1^2}{2\!\cdot\!3^{0}\!+\!1\!\cdot\!3^{1}\!+\!1\!\cdot\!3^{2}\!+\!O(3^{5})}{2\!\cdot\!3^{0}\phantom{\,\,\!+\!2\!\cdot\!3^1}\,\!+\!2\!\cdot\!3^{2}\!+\!O(3^{3})}\errentry{9\!\cdot\!10^{-16}}{4\!\cdot\!10^{-16}}\cr
\tabentry{20}{[0, -1, 1, 0, 0]}{x^3-20}{2^{2}3^{3}5^{2}}{2^{4}3^{6}5^{4}11^{2}}{1}{1^2}{2\!\cdot\!3^{0}\!+\!1\!\cdot\!3^{1}\phantom{\,\,\!+\!2\!\cdot\!3^2}\,\!+\!O(3^{4})}{2\!\cdot\!3^{0}\phantom{\,\,\!+\!2\!\cdot\!3^1}\,\phantom{\,\,\!+\!2\!\cdot\!3^2}\,\!+\!O(3^{3})}\errentry{2\!\cdot\!10^{-17}}{1\!\cdot\!10^{-15}}\cr
\tabentry{21}{[0, -1, 1, 0, 0]}{x^3-21}{3^{5}7^{2}}{3^{10}7^{4}11^{2}}{2^{2}}{2^2}{2\!\cdot\!3^{0}\phantom{\,\,\!+\!2\!\cdot\!3^1}\,\!+\!2\!\cdot\!3^{2}\!+\!O(3^{4})}{2\!\cdot\!3^{0}\!+\!2\!\cdot\!3^{1}\phantom{\,\,\!+\!2\!\cdot\!3^2}\,\!+\!O(3^{3})}\errentry{3\!\cdot\!10^{-12}}{1\!\cdot\!10^{-12}}\cr
\tabentry{23}{[0, -1, 1, 0, 0]}{x^3-23}{3^{3}23^{2}}{3^{6}11^{2}23^{4}}{2^{4}}{4^2}{1\!\cdot\!3^{0}\!+\!1\!\cdot\!3^{1}\!+\!1\!\cdot\!3^{2}\!+\!O(3^{3})}{1\!\cdot\!3^{0}\phantom{\,\,\!+\!2\!\cdot\!3^1}\,\phantom{\,\,\!+\!2\!\cdot\!3^2}\,\!+\!O(3^{3})}\errentry{5\!\cdot\!10^{-10}}{5\!\cdot\!10^{-10}}\cr
\tabentry{26}{[0, -1, 1, 0, 0]}{x^3-26}{2^{2}3\!\cdot\!13^{2}}{2^{4}3^{2}11^{2}13^{4}}{2^{2}}{2^2}{1\!\cdot\!3^{0}\!+\!1\!\cdot\!3^{1}\!+\!1\!\cdot\!3^{2}\!+\!O(3^{4})}{1\!\cdot\!3^{0}\phantom{\,\,\!+\!2\!\cdot\!3^1}\,\!+\!2\!\cdot\!3^{2}\!+\!O(3^{4})}\errentry{8\!\cdot\!10^{-17}}{2\!\cdot\!10^{-15}}\cr
\tabentry{28}{[0, -1, 1, 0, 0]}{x^3-28}{2^{2}3\!\cdot\!7^{2}}{2^{4}3^{2}7^{4}11^{2}}{1}{1^2}{1\!\cdot\!3^{0}\!+\!2\!\cdot\!3^{1}\!+\!2\!\cdot\!3^{2}\!+\!O(3^{4})}{1\!\cdot\!3^{0}\!+\!1\!\cdot\!3^{1}\phantom{\,\,\!+\!2\!\cdot\!3^2}\,\!+\!O(3^{3})}\errentry{2\!\cdot\!10^{-16}}{2\!\cdot\!10^{-16}}\cr
\tabentry{29}{[0, -1, 1, 0, 0]}{x^3-29}{3^{3}29^{2}}{3^{6}11^{2}29^{4}}{3^{2}}{3^2}{2\!\cdot\!3^{2}\phantom{\,\,\!+\!2\!\cdot\!3^3}\,\phantom{\,\,\!+\!2\!\cdot\!3^4}\,\!+\!O(3^{5})}{1\!\cdot\!3^{2}\phantom{\,\,\!+\!2\!\cdot\!3^3}\,\!+\!1\!\cdot\!3^{4}\!+\!O(3^{5})}\errentry{1\!\cdot\!10^{-9}}{1\!\cdot\!10^{-9}}\cr
\tabentry{30}{[0, -1, 1, 0, 0]}{x^3-30}{2^{2}3^{5}5^{2}}{2^{4}3^{10}5^{4}11^{2}}{1}{1^2}{2\!\cdot\!3^{0}\!+\!1\!\cdot\!3^{1}\phantom{\,\,\!+\!2\!\cdot\!3^2}\,\!+\!O(3^{4})}{2\!\cdot\!3^{0}\!+\!1\!\cdot\!3^{1}\!+\!2\!\cdot\!3^{2}\!+\!O(3^{4})}\errentry{2\!\cdot\!10^{-9}}{8\!\cdot\!10^{-10}}\cr
\tabentry{31}{[0, -1, 1, 0, 0]}{x^3-31}{3^{3}31^{2}}{3^{6}11^{2}31^{4}}{1}{1^2}{2\!\cdot\!3^{0}\phantom{\,\,\!+\!2\!\cdot\!3^1}\,\!+\!2\!\cdot\!3^{2}\!+\!O(3^{3})}{2\!\cdot\!3^{0}\!+\!1\!\cdot\!3^{1}\phantom{\,\,\!+\!2\!\cdot\!3^2}\,\!+\!O(3^{4})}\errentry{5\!\cdot\!10^{-9}}{1\!\cdot\!10^{-9}}\cr
\tabentry{34}{[0, -1, 1, 0, 0]}{x^3-34}{2^{2}3^{3}17^{2}}{2^{4}3^{6}11^{2}17^{4}}{5^{2}}{5^2}{2\!\cdot\!3^{0}\!+\!1\!\cdot\!3^{1}\!+\!2\!\cdot\!3^{2}\!+\!O(3^{3})}{2\!\cdot\!3^{0}\phantom{\,\,\!+\!2\!\cdot\!3^1}\,\phantom{\,\,\!+\!2\!\cdot\!3^2}\,\!+\!O(3^{6})}\errentry{4\!\cdot\!10^{-8}}{5\!\cdot\!10^{-8}}\cr
\tabentry{35}{[0, -1, 1, 0, 0]}{x^3-35}{3\!\cdot\!5^{2}7^{2}}{3^{2}5^{4}7^{4}11^{2}}{1}{1^2}{1\!\cdot\!3^{0}\phantom{\,\,\!+\!2\!\cdot\!3^1}\,\phantom{\,\,\!+\!2\!\cdot\!3^2}\,\!+\!O(3^{3})}{1\!\cdot\!3^{0}\!+\!2\!\cdot\!3^{1}\phantom{\,\,\!+\!2\!\cdot\!3^2}\,\!+\!O(3^{4})}\errentry{8\!\cdot\!10^{-17}}{1\!\cdot\!10^{-15}}\cr
\tabentry{37}{[0, -1, 1, 0, 0]}{x^3-37}{3\!\cdot\!37^{2}}{3^{2}11^{2}37^{4}}{2^{2}}{2^2}{2\!\cdot\!3^{0}\!+\!2\!\cdot\!3^{1}\phantom{\,\,\!+\!2\!\cdot\!3^2}\,\!+\!O(3^{5})}{2\!\cdot\!3^{0}\!+\!2\!\cdot\!3^{1}\!+\!1\!\cdot\!3^{2}\!+\!O(3^{3})}\errentry{1\!\cdot\!10^{-16}}{5\!\cdot\!10^{-15}}\cr
\tabentry{45}{[0, -1, 1, 0, 0]}{x^3-45}{3^{5}5^{2}}{3^{10}5^{4}11^{2}}{1}{1^2}{1\!\cdot\!3^{0}\!+\!1\!\cdot\!3^{1}\phantom{\,\,\!+\!2\!\cdot\!3^2}\,\!+\!O(3^{4})}{1\!\cdot\!3^{0}\phantom{\,\,\!+\!2\!\cdot\!3^1}\,\!+\!2\!\cdot\!3^{2}\!+\!O(3^{3})}\errentry{2\!\cdot\!10^{-17}}{3\!\cdot\!10^{-16}}\cr
\tabentry{46}{[0, -1, 1, 0, 0]}{x^3-46}{2^{2}3\!\cdot\!23^{2}}{2^{4}3^{2}11^{2}23^{4}}{2^{2}}{2^2}{2\!\cdot\!3^{0}\!+\!1\!\cdot\!3^{1}\!+\!2\!\cdot\!3^{2}\!+\!O(3^{4})}{2\!\cdot\!3^{0}\!+\!2\!\cdot\!3^{1}\!+\!2\!\cdot\!3^{2}\!+\!O(3^{3})}\errentry{9\!\cdot\!10^{-16}}{3\!\cdot\!10^{-16}}\cr
\tabentry{52}{[0, -1, 1, 0, 0]}{x^3-52}{2^{2}3^{3}13^{2}}{2^{4}3^{6}11^{2}13^{4}}{2^{4}}{4^2}{1\!\cdot\!3^{0}\!+\!1\!\cdot\!3^{1}\!+\!1\!\cdot\!3^{2}\!+\!O(3^{4})}{1\!\cdot\!3^{0}\!+\!1\!\cdot\!3^{1}\phantom{\,\,\!+\!2\!\cdot\!3^2}\,\!+\!O(3^{3})}\errentry{5\!\cdot\!10^{-10}}{5\!\cdot\!10^{-10}}\cr
\tabentry{53}{[0, -1, 1, 0, 0]}{x^3-53}{3\!\cdot\!53^{2}}{3^{2}11^{2}53^{4}}{3^{2}}{3^2}{1\!\cdot\!3^{2}\!+\!1\!\cdot\!3^{3}\phantom{\,\,\!+\!2\!\cdot\!3^4}\,\!+\!O(3^{9})}{1\!\cdot\!3^{2}\!+\!2\!\cdot\!3^{3}\!+\!2\!\cdot\!3^{4}\!+\!O(3^{5})}\errentry{3\!\cdot\!10^{-15}}{6\!\cdot\!10^{-15}}\cr
\tabentry{60}{[0, -1, 1, 0, 0]}{x^3-60}{2^{2}3^{5}5^{2}}{2^{4}3^{10}5^{4}11^{2}}{11^{2}}{11^2}{2\!\cdot\!3^{0}\!+\!1\!\cdot\!3^{1}\phantom{\,\,\!+\!2\!\cdot\!3^2}\,\!+\!O(3^{4})}{2\!\cdot\!3^{0}\phantom{\,\,\!+\!2\!\cdot\!3^1}\,\phantom{\,\,\!+\!2\!\cdot\!3^2}\,\!+\!O(3^{3})}\errentry{9\!\cdot\!10^{-10}}{2\!\cdot\!10^{-9}}\cr
\tabentry{62}{[0, -1, 1, 0, 0]}{x^3-62}{2^{2}3\!\cdot\!31^{2}}{2^{4}3^{2}11^{2}31^{4}}{1}{1^2}{1\!\cdot\!3^{0}\!+\!2\!\cdot\!3^{1}\!+\!2\!\cdot\!3^{2}\!+\!O(3^{3})}{1\!\cdot\!3^{0}\!+\!2\!\cdot\!3^{1}\phantom{\,\,\!+\!2\!\cdot\!3^2}\,\!+\!O(3^{3})}\errentry{7\!\cdot\!10^{-13}}{1\!\cdot\!10^{-13}}\cr
\tabentry{63}{[0, -1, 1, 0, 0]}{x^3-63}{3^{5}7^{2}}{3^{10}7^{4}11^{2}}{2^{4}}{4^2}{2\!\cdot\!3^{0}\phantom{\,\,\!+\!2\!\cdot\!3^1}\,\!+\!2\!\cdot\!3^{2}\!+\!O(3^{4})}{2\!\cdot\!3^{0}\!+\!1\!\cdot\!3^{1}\phantom{\,\,\!+\!2\!\cdot\!3^2}\,\!+\!O(3^{4})}\errentry{1\!\cdot\!10^{-14}}{1\!\cdot\!10^{-14}}\cr
\tabentry{68}{[0, -1, 1, 0, 0]}{x^3-68}{2^{2}3^{3}17^{2}}{2^{4}3^{6}11^{2}17^{4}}{2^{4}}{4^2}{2\!\cdot\!3^{0}\!+\!1\!\cdot\!3^{1}\!+\!2\!\cdot\!3^{2}\!+\!O(3^{3})}{2\!\cdot\!3^{0}\phantom{\,\,\!+\!2\!\cdot\!3^1}\,\!+\!1\!\cdot\!3^{2}\!+\!O(3^{5})}\errentry{1\!\cdot\!10^{-7}}{1\!\cdot\!10^{-7}}\cr
\tabentry{71}{[0, -1, 1, 0, 0]}{x^3-71}{3\!\cdot\!71^{2}}{3^{2}11^{2}71^{4}}{0}{\nosha}{1\!\cdot\!3^{2}\phantom{\,\,\!+\!2\!\cdot\!3^3}\,\phantom{\,\,\!+\!2\!\cdot\!3^4}\,\!+\!O(3^{5})}{\noeuler}\errentry{3\!\cdot\!10^{-11}}{7\!\cdot\!10^{-7}}\cr
\tabentry{73}{[0, -1, 1, 0, 0]}{x^3-73}{3\!\cdot\!73^{2}}{3^{2}11^{2}73^{4}}{2^{2}}{2^2}{2\!\cdot\!3^{0}\!+\!1\!\cdot\!3^{1}\!+\!2\!\cdot\!3^{2}\!+\!O(3^{3})}{2\!\cdot\!3^{0}\!+\!2\!\cdot\!3^{1}\!+\!2\!\cdot\!3^{2}\!+\!O(3^{3})}\errentry{6\!\cdot\!10^{-11}}{3\!\cdot\!10^{-11}}\cr
\tabentry{82}{[0, -1, 1, 0, 0]}{x^3-82}{2^{2}3\!\cdot\!41^{2}}{2^{4}3^{2}11^{2}41^{4}}{1}{1^2}{2\!\cdot\!3^{0}\!+\!1\!\cdot\!3^{1}\phantom{\,\,\!+\!2\!\cdot\!3^2}\,\!+\!O(3^{3})}{2\!\cdot\!3^{0}\phantom{\,\,\!+\!2\!\cdot\!3^1}\,\phantom{\,\,\!+\!2\!\cdot\!3^2}\,\!+\!O(3^{5})}\errentry{3\!\cdot\!10^{-10}}{9\!\cdot\!10^{-11}}\cr
\tabentry{89}{[0, -1, 1, 0, 0]}{x^3-89}{3\!\cdot\!89^{2}}{3^{2}11^{2}89^{4}}{3^{2}}{3^2}{1\!\cdot\!3^{2}\!+\!2\!\cdot\!3^{3}\phantom{\,\,\!+\!2\!\cdot\!3^4}\,\!+\!O(3^{5})}{1\!\cdot\!3^{2}\!+\!2\!\cdot\!3^{3}\phantom{\,\,\!+\!2\!\cdot\!3^4}\,\!+\!O(3^{5})}\errentry{2\!\cdot\!10^{-8}}{1\!\cdot\!10^{-8}}\cr
\tabentry{90}{[0, -1, 1, 0, 0]}{x^3-90}{2^{2}3^{5}5^{2}}{2^{4}3^{10}5^{4}11^{2}}{2^{2}}{2^2}{2\!\cdot\!3^{0}\!+\!1\!\cdot\!3^{1}\phantom{\,\,\!+\!2\!\cdot\!3^2}\,\!+\!O(3^{4})}{2\!\cdot\!3^{0}\phantom{\,\,\!+\!2\!\cdot\!3^1}\,\!+\!1\!\cdot\!3^{2}\!+\!O(3^{9})}\errentry{6\!\cdot\!10^{-11}}{3\!\cdot\!10^{-11}}\cr
\tabentry{91}{[0, -1, 1, 0, 0]}{x^3-91}{3\!\cdot\!7^{2}13^{2}}{3^{2}7^{4}11^{2}13^{4}}{1}{1^2}{2\!\cdot\!3^{0}\!+\!2\!\cdot\!3^{1}\phantom{\,\,\!+\!2\!\cdot\!3^2}\,\!+\!O(3^{4})}{2\!\cdot\!3^{0}\phantom{\,\,\!+\!2\!\cdot\!3^1}\,\!+\!1\!\cdot\!3^{2}\!+\!O(3^{3})}\errentry{2\!\cdot\!10^{-9}}{5\!\cdot\!10^{-10}}\cr
\tabentry{116}{[0, -1, 1, 0, 0]}{x^3-116}{2^{2}3\!\cdot\!29^{2}}{2^{4}3^{2}11^{2}29^{4}}{3^{2}}{3^2}{1\!\cdot\!3^{2}\!+\!2\!\cdot\!3^{3}\!+\!1\!\cdot\!3^{4}\!+\!O(3^{8})}{2\!\cdot\!3^{2}\!+\!1\!\cdot\!3^{3}\!+\!1\!\cdot\!3^{4}\!+\!O(3^{6})}\errentry{6\!\cdot\!10^{-14}}{3\!\cdot\!10^{-14}}\cr
\tabentry{150}{[0, -1, 1, 0, 0]}{x^3-150}{2^{2}3^{5}5^{2}}{2^{4}3^{10}5^{4}11^{2}}{1}{1^2}{2\!\cdot\!3^{0}\!+\!1\!\cdot\!3^{1}\phantom{\,\,\!+\!2\!\cdot\!3^2}\,\!+\!O(3^{4})}{2\!\cdot\!3^{0}\!+\!1\!\cdot\!3^{1}\!+\!2\!\cdot\!3^{2}\!+\!O(3^{4})}\errentry{6\!\cdot\!10^{-9}}{1\!\cdot\!10^{-9}}\cr
\tabentry{172}{[0, -1, 1, 0, 0]}{x^3-172}{2^{2}3\!\cdot\!43^{2}}{2^{4}3^{2}11^{2}43^{4}}{1}{1^2}{1\!\cdot\!3^{0}\!+\!1\!\cdot\!3^{1}\!+\!1\!\cdot\!3^{2}\!+\!O(3^{3})}{1\!\cdot\!3^{0}\!+\!1\!\cdot\!3^{1}\!+\!2\!\cdot\!3^{2}\!+\!O(3^{3})}\errentry{1\!\cdot\!10^{-10}}{4\!\cdot\!10^{-11}}\cr
\tabentry{188}{[0, -1, 1, 0, 0]}{x^3-188}{2^{2}3\!\cdot\!47^{2}}{2^{4}3^{2}11^{2}47^{4}}{1}{1^2}{2\!\cdot\!3^{0}\phantom{\,\,\!+\!2\!\cdot\!3^1}\,\!+\!2\!\cdot\!3^{2}\!+\!O(3^{3})}{2\!\cdot\!3^{0}\!+\!1\!\cdot\!3^{1}\!+\!2\!\cdot\!3^{2}\!+\!O(3^{3})}\errentry{1\!\cdot\!10^{-8}}{2\!\cdot\!10^{-9}}\cr
\tabentry{325}{[0, -1, 1, 0, 0]}{x^3-325}{3\!\cdot\!5^{2}13^{2}}{3^{2}5^{4}11^{2}13^{4}}{5^{2}}{5^2}{1\!\cdot\!3^{0}\!+\!2\!\cdot\!3^{1}\phantom{\,\,\!+\!2\!\cdot\!3^2}\,\!+\!O(3^{3})}{1\!\cdot\!3^{0}\!+\!2\!\cdot\!3^{1}\phantom{\,\,\!+\!2\!\cdot\!3^2}\,\!+\!O(3^{5})}\errentry{1\!\cdot\!10^{-12}}{1\!\cdot\!10^{-12}}\cr
\tabentry{350}{[0, -1, 1, 0, 0]}{x^3-350}{2^{2}3\!\cdot\!5^{2}7^{2}}{2^{4}3^{2}5^{4}7^{4}11^{2}}{1}{1^2}{2\!\cdot\!3^{0}\!+\!2\!\cdot\!3^{1}\phantom{\,\,\!+\!2\!\cdot\!3^2}\,\!+\!O(3^{3})}{2\!\cdot\!3^{0}\!+\!2\!\cdot\!3^{1}\!+\!1\!\cdot\!3^{2}\!+\!O(3^{3})}\errentry{1\!\cdot\!10^{-11}}{3\!\cdot\!10^{-12}}\cr
\hline
\tabend

\noindent\begin{tabular}{cllll} 
\extravspace
$E=20A3$ & \multicolumn{3}{l}{$y^2=x^3+x^2-36x-140$} & \tabname{3}{20A3}\cr
\extravspace
$p=3$ & \multicolumn{3}{l}{$P_3(E/\Q,T)=1+2T+3T^2$}\cr
$/\vQ\alsub{}$   & $N(E)=2^{2}\!\cdot\!5$ & $L^*(E/\Q)=2^{-1}$ & $|E(\Q)|=2$ & $|\sha(E/\Q)|=1^2$ \cr
$/K\alsub{}$      & $N(E,\sigma)=2^{4}\!\cdot\!3^{2}\!\cdot\!5^{2}$ & $L^*(E/K)=2^{-1}$ & $|E(K)|=6$ & $|\sha(E/K)|=1^2$ \cr
$/\vQ\alsub{2}$ & $\hbox{IV}^*,\>c_{2}=1$ & $/K\alsub{v}\quad(v|2)$ & $\hbox{IV}^*,\>c_v=3$\cr
$/\vQ\alsub{5}$ & $\hbox{I}_6 \hbox{ non-split},\>c_{5}=2$ & $/K\alsub{v}\quad(v|5)$ & $\hbox{I}_6 \hbox{ split},\>c_v=6$\cr
\end{tabular}
\tabbegin
\tabentry{3}{[0, 1, 0, -36, -140]}{x^3-3}{3^{5}}{2^{4}3^{10}5^{2}}{2\!\cdot\!3^{2}}{1^2}{2\!\cdot\!3^{2}\!+\!2\!\cdot\!3^{3}\!+\!2\!\cdot\!3^{4}\!+\!O(3^{5})}{2\!\cdot\!3^{2}\phantom{\,\,\!+\!2\!\cdot\!3^3}\,\phantom{\,\,\!+\!2\!\cdot\!3^4}\,\!+\!O(3^{7})}\errentry{1\!\cdot\!10^{-17}}{1\!\cdot\!10^{-15}}\cr
\tabentry{7}{[0, 1, 0, -36, -140]}{x^3-7}{3^{3}7^{2}}{2^{4}3^{6}5^{2}7^{4}}{0}{\nosha}{2\!\cdot\!3^{4}\!+\!2\!\cdot\!3^{5}\!+\!1\!\cdot\!3^{6}\!+\!O(3^{7})}{\noeuler}\errentry{5\!\cdot\!10^{-17}}{9\!\cdot\!10^{-10}}\cr
\tabentry{10}{[0, 1, 0, -36, -140]}{x^3-10}{2^{2}3\!\cdot\!5^{2}}{2^{2}3^{2}5^{4}}{2^{2}}{1^2}{1\!\cdot\!3^{3}\!+\!1\!\cdot\!3^{4}\phantom{\,\,\!+\!2\!\cdot\!3^5}\,\!+\!O(3^{6})}{2\!\cdot\!3^{3}\!+\!1\!\cdot\!3^{4}\phantom{\,\,\!+\!2\!\cdot\!3^5}\,\!+\!O(3^{6})}\errentry{1\!\cdot\!10^{-15}}{2\!\cdot\!10^{-16}}\cr
\tabentry{11}{[0, 1, 0, -36, -140]}{x^3-11}{3^{3}11^{2}}{2^{4}3^{6}5^{2}11^{4}}{0}{\nosha}{2\!\cdot\!3^{4}\!+\!1\!\cdot\!3^{5}\phantom{\,\,\!+\!2\!\cdot\!3^6}\,\!+\!O(3^{7})}{\noeuler}\errentry{5\!\cdot\!10^{-17}}{2\!\cdot\!10^{-9}}\cr
\tabentry{13}{[0, 1, 0, -36, -140]}{x^3-13}{3^{3}13^{2}}{2^{4}3^{6}5^{2}13^{4}}{2\!\cdot\!3^{4}}{3^2}{2\!\cdot\!3^{4}\!+\!2\!\cdot\!3^{5}\!+\!1\!\cdot\!3^{6}\!+\!O(3^{7})}{2\!\cdot\!3^{4}\!+\!2\!\cdot\!3^{5}\!+\!1\!\cdot\!3^{6}\!+\!O(3^{7})}\errentry{2\!\cdot\!10^{-17}}{1\!\cdot\!10^{-15}}\cr
\tabentry{17}{[0, 1, 0, -36, -140]}{x^3-17}{3\!\cdot\!17^{2}}{2^{4}3^{2}5^{2}17^{4}}{2\!\cdot\!3^{2}}{1^2}{1\!\cdot\!3^{4}\!+\!1\!\cdot\!3^{5}\phantom{\,\,\!+\!2\!\cdot\!3^6}\,\!+\!O(3^{7})}{1\!\cdot\!3^{4}\!+\!1\!\cdot\!3^{5}\!+\!2\!\cdot\!3^{6}\!+\!O(3^{7})}\errentry{6\!\cdot\!10^{-16}}{8\!\cdot\!10^{-16}}\cr
\tabentry{19}{[0, 1, 0, -36, -140]}{x^3-19}{3\!\cdot\!19^{2}}{2^{4}3^{2}5^{2}19^{4}}{2\!\cdot\!3^{2}}{1^2}{2\!\cdot\!3^{4}\!+\!2\!\cdot\!3^{5}\!+\!2\!\cdot\!3^{6}\!+\!O(3^{7})}{2\!\cdot\!3^{4}\!+\!1\!\cdot\!3^{5}\!+\!1\!\cdot\!3^{6}\!+\!O(3^{7})}\errentry{4\!\cdot\!10^{-16}}{3\!\cdot\!10^{-15}}\cr
\tabentry{21}{[0, 1, 0, -36, -140]}{x^3-21}{3^{5}7^{2}}{2^{4}3^{10}5^{2}7^{4}}{2\!\cdot\!3^{4}}{3^2}{2\!\cdot\!3^{4}\!+\!2\!\cdot\!3^{5}\!+\!1\!\cdot\!3^{6}\!+\!O(3^{7})}{2\!\cdot\!3^{4}\!+\!2\!\cdot\!3^{5}\phantom{\,\,\!+\!2\!\cdot\!3^6}\,\!+\!O(3^{7})}\errentry{2\!\cdot\!10^{-10}}{2\!\cdot\!10^{-10}}\cr
\tabentry{23}{[0, 1, 0, -36, -140]}{x^3-23}{3^{3}23^{2}}{2^{4}3^{6}5^{2}23^{4}}{0}{\nosha}{1\!\cdot\!3^{5}\phantom{\,\,\!+\!2\!\cdot\!3^6}\,\!+\!2\!\cdot\!3^{7}\!+\!O(3^{8})}{\noeuler}\errentry{6\!\cdot\!10^{-9}}{1\!\cdot\!10^{-5}}\cr
\tabentry{30}{[0, 1, 0, -36, -140]}{x^3-30}{2^{2}3^{5}5^{2}}{2^{2}3^{10}5^{4}}{2^{2}3^{2}}{3^2}{1\!\cdot\!3^{3}\!+\!1\!\cdot\!3^{4}\phantom{\,\,\!+\!2\!\cdot\!3^5}\,\!+\!O(3^{6})}{2\!\cdot\!3^{3}\phantom{\,\,\!+\!2\!\cdot\!3^4}\,\!+\!1\!\cdot\!3^{5}\!+\!O(3^{6})}\errentry{7\!\cdot\!10^{-16}}{7\!\cdot\!10^{-16}}\cr
\tabentry{37}{[0, 1, 0, -36, -140]}{x^3-37}{3\!\cdot\!37^{2}}{2^{4}3^{2}5^{2}37^{4}}{0}{\nosha}{2\!\cdot\!3^{6}\phantom{\,\,\!+\!2\!\cdot\!3^7}\,\phantom{\,\,\!+\!2\!\cdot\!3^8}\,\!+\!O(3^{10})}{\noeuler}\errentry{3\!\cdot\!10^{-16}}{2\!\cdot\!10^{-9}}\cr
\tabentry{53}{[0, 1, 0, -36, -140]}{x^3-53}{3\!\cdot\!53^{2}}{2^{4}3^{2}5^{2}53^{4}}{2\!\cdot\!3^{2}}{1^2}{1\!\cdot\!3^{4}\!+\!1\!\cdot\!3^{5}\phantom{\,\,\!+\!2\!\cdot\!3^6}\,\!+\!O(3^{7})}{1\!\cdot\!3^{4}\!+\!1\!\cdot\!3^{5}\phantom{\,\,\!+\!2\!\cdot\!3^6}\,\!+\!O(3^{7})}\errentry{3\!\cdot\!10^{-11}}{8\!\cdot\!10^{-12}}\cr
\tabentry{63}{[0, 1, 0, -36, -140]}{x^3-63}{3^{5}7^{2}}{2^{4}3^{10}5^{2}7^{4}}{0}{\nosha}{2\!\cdot\!3^{4}\!+\!2\!\cdot\!3^{5}\!+\!1\!\cdot\!3^{6}\!+\!O(3^{7})}{\noeuler}\errentry{1\!\cdot\!10^{-9}}{3\!\cdot\!10^{-6}}\cr
\tabentry{70}{[0, 1, 0, -36, -140]}{x^3-70}{2^{2}3^{3}5^{2}7^{2}}{2^{2}3^{6}5^{4}7^{4}}{0}{\nosha}{1\!\cdot\!3^{5}\!+\!1\!\cdot\!3^{6}\!+\!1\!\cdot\!3^{7}\!+\!O(3^{9})}{\noeuler}\errentry{1\!\cdot\!10^{-16}}{7\!\cdot\!10^{-9}}\cr
\tabentry{71}{[0, 1, 0, -36, -140]}{x^3-71}{3\!\cdot\!71^{2}}{2^{4}3^{2}5^{2}71^{4}}{2^{3}3^{2}}{2^2}{1\!\cdot\!3^{4}\!+\!2\!\cdot\!3^{5}\!+\!1\!\cdot\!3^{6}\!+\!O(3^{7})}{1\!\cdot\!3^{4}\!+\!2\!\cdot\!3^{5}\phantom{\,\,\!+\!2\!\cdot\!3^6}\,\!+\!O(3^{7})}\errentry{5\!\cdot\!10^{-8}}{2\!\cdot\!10^{-8}}\cr
\tabentry{73}{[0, 1, 0, -36, -140]}{x^3-73}{3\!\cdot\!73^{2}}{2^{4}3^{2}5^{2}73^{4}}{2\!\cdot\!3^{4}}{3^2}{2\!\cdot\!3^{6}\!+\!2\!\cdot\!3^{7}\!+\!2\!\cdot\!3^{8}\!+\!O(3^{9})}{2\!\cdot\!3^{6}\!+\!1\!\cdot\!3^{7}\!+\!1\!\cdot\!3^{8}\!+\!O(3^{9})}\errentry{3\!\cdot\!10^{-8}}{2\!\cdot\!10^{-8}}\cr
\tabentry{90}{[0, 1, 0, -36, -140]}{x^3-90}{2^{2}3^{5}5^{2}}{2^{2}3^{10}5^{4}}{2^{4}3^{2}}{6^2}{1\!\cdot\!3^{3}\!+\!1\!\cdot\!3^{4}\phantom{\,\,\!+\!2\!\cdot\!3^5}\,\!+\!O(3^{6})}{2\!\cdot\!3^{3}\!+\!2\!\cdot\!3^{4}\!+\!1\!\cdot\!3^{5}\!+\!O(3^{7})}\errentry{2\!\cdot\!10^{-16}}{7\!\cdot\!10^{-15}}\cr
\tabentry{110}{[0, 1, 0, -36, -140]}{x^3-110}{2^{2}3^{3}5^{2}11^{2}}{2^{2}3^{6}5^{4}11^{4}}{2^{2}3^{4}}{9^2}{1\!\cdot\!3^{5}\!+\!2\!\cdot\!3^{6}\phantom{\,\,\!+\!2\!\cdot\!3^7}\,\!+\!O(3^{9})}{1\!\cdot\!3^{5}\!+\!2\!\cdot\!3^{6}\!+\!1\!\cdot\!3^{7}\!+\!O(3^{8})}\errentry{1\!\cdot\!10^{-11}}{2\!\cdot\!10^{-11}}\cr
\tabentry{130}{[0, 1, 0, -36, -140]}{x^3-130}{2^{2}3^{3}5^{2}13^{2}}{2^{2}3^{6}5^{4}13^{4}}{2^{2}3^{4}}{9^2}{1\!\cdot\!3^{5}\!+\!1\!\cdot\!3^{6}\!+\!1\!\cdot\!3^{7}\!+\!O(3^{8})}{2\!\cdot\!3^{5}\!+\!2\!\cdot\!3^{6}\!+\!2\!\cdot\!3^{7}\!+\!O(3^{8})}\errentry{5\!\cdot\!10^{-10}}{1\!\cdot\!10^{-9}}\cr
\tabentry{150}{[0, 1, 0, -36, -140]}{x^3-150}{2^{2}3^{5}5^{2}}{2^{2}3^{10}5^{4}}{2^{2}3^{2}}{3^2}{1\!\cdot\!3^{3}\!+\!1\!\cdot\!3^{4}\phantom{\,\,\!+\!2\!\cdot\!3^5}\,\!+\!O(3^{6})}{2\!\cdot\!3^{3}\phantom{\,\,\!+\!2\!\cdot\!3^4}\,\!+\!1\!\cdot\!3^{5}\!+\!O(3^{6})}\errentry{7\!\cdot\!10^{-17}}{7\!\cdot\!10^{-16}}\cr
\tabentry{170}{[0, 1, 0, -36, -140]}{x^3-170}{2^{2}3\!\cdot\!5^{2}17^{2}}{2^{2}3^{2}5^{4}17^{4}}{2^{2}3^{2}}{3^2}{2\!\cdot\!3^{5}\phantom{\,\,\!+\!2\!\cdot\!3^6}\,\!+\!1\!\cdot\!3^{7}\!+\!O(3^{8})}{1\!\cdot\!3^{5}\!+\!1\!\cdot\!3^{6}\!+\!1\!\cdot\!3^{7}\!+\!O(3^{8})}\errentry{2\!\cdot\!10^{-14}}{1\!\cdot\!10^{-14}}\cr
\tabentry{190}{[0, 1, 0, -36, -140]}{x^3-190}{2^{2}3\!\cdot\!5^{2}19^{2}}{2^{2}3^{2}5^{4}19^{4}}{2^{2}3^{2}}{3^2}{1\!\cdot\!3^{5}\!+\!1\!\cdot\!3^{6}\phantom{\,\,\!+\!2\!\cdot\!3^7}\,\!+\!O(3^{9})}{2\!\cdot\!3^{5}\!+\!1\!\cdot\!3^{6}\!+\!2\!\cdot\!3^{7}\!+\!O(3^{10})}\errentry{5\!\cdot\!10^{-17}}{4\!\cdot\!10^{-15}}\cr
\tabentry{350}{[0, 1, 0, -36, -140]}{x^3-350}{2^{2}3\!\cdot\!5^{2}7^{2}}{2^{2}3^{2}5^{4}7^{4}}{2^{2}3^{2}}{3^2}{1\!\cdot\!3^{5}\!+\!1\!\cdot\!3^{6}\!+\!1\!\cdot\!3^{7}\!+\!O(3^{9})}{2\!\cdot\!3^{5}\phantom{\,\,\!+\!2\!\cdot\!3^6}\,\!+\!2\!\cdot\!3^{7}\!+\!O(3^{8})}\errentry{4\!\cdot\!10^{-16}}{7\!\cdot\!10^{-16}}\cr
\tabentry{370}{[0, 1, 0, -36, -140]}{x^3-370}{2^{2}3\!\cdot\!5^{2}37^{2}}{2^{2}3^{2}5^{4}37^{4}}{2^{2}3^{4}}{9^2}{1\!\cdot\!3^{7}\phantom{\,\,\!+\!2\!\cdot\!3^8}\,\!+\!2\!\cdot\!3^{9}\!+\!O(3^{10})}{2\!\cdot\!3^{7}\!+\!1\!\cdot\!3^{8}\!+\!1\!\cdot\!3^{9}\!+\!O(3^{11})}\errentry{5\!\cdot\!10^{-9}}{1\!\cdot\!10^{-8}}\cr
\tabentry{490}{[0, 1, 0, -36, -140]}{x^3-490}{2^{2}3^{3}5^{2}7^{2}}{2^{2}3^{6}5^{4}7^{4}}{2^{2}3^{4}}{9^2}{1\!\cdot\!3^{5}\!+\!1\!\cdot\!3^{6}\!+\!1\!\cdot\!3^{7}\!+\!O(3^{9})}{2\!\cdot\!3^{5}\phantom{\,\,\!+\!2\!\cdot\!3^6}\,\phantom{\,\,\!+\!2\!\cdot\!3^7}\,\!+\!O(3^{8})}\errentry{3\!\cdot\!10^{-16}}{4\!\cdot\!10^{-15}}\cr
\tabentry{539}{[0, 1, 0, -36, -140]}{x^3-539}{3\!\cdot\!7^{2}11^{2}}{2^{4}3^{2}5^{2}7^{4}11^{4}}{2\!\cdot\!3^{4}}{3^2}{2\!\cdot\!3^{6}\!+\!1\!\cdot\!3^{7}\!+\!2\!\cdot\!3^{8}\!+\!O(3^{9})}{1\!\cdot\!3^{6}\phantom{\,\,\!+\!2\!\cdot\!3^7}\,\phantom{\,\,\!+\!2\!\cdot\!3^8}\,\!+\!O(3^{10})}\errentry{1\!\cdot\!10^{-9}}{8\!\cdot\!10^{-10}}\cr
\tabentry{550}{[0, 1, 0, -36, -140]}{x^3-550}{2^{2}3\!\cdot\!5^{2}11^{2}}{2^{2}3^{2}5^{4}11^{4}}{2^{2}3^{2}}{3^2}{1\!\cdot\!3^{5}\!+\!2\!\cdot\!3^{6}\phantom{\,\,\!+\!2\!\cdot\!3^7}\,\!+\!O(3^{9})}{1\!\cdot\!3^{5}\!+\!2\!\cdot\!3^{6}\!+\!2\!\cdot\!3^{7}\!+\!O(3^{9})}\errentry{1\!\cdot\!10^{-16}}{3\!\cdot\!10^{-15}}\cr
\tabentry{650}{[0, 1, 0, -36, -140]}{x^3-650}{2^{2}3^{3}5^{2}13^{2}}{2^{2}3^{6}5^{4}13^{4}}{2^{4}3^{4}}{18^2}{1\!\cdot\!3^{5}\!+\!1\!\cdot\!3^{6}\!+\!1\!\cdot\!3^{7}\!+\!O(3^{8})}{2\!\cdot\!3^{5}\!+\!1\!\cdot\!3^{6}\!+\!2\!\cdot\!3^{7}\!+\!O(3^{8})}\errentry{8\!\cdot\!10^{-11}}{3\!\cdot\!10^{-10}}\cr
\tabentry{1450}{[0, 1, 0, -36, -140]}{x^3-1450}{2^{2}3\!\cdot\!5^{2}29^{2}}{2^{2}3^{2}5^{4}29^{4}}{2^{4}3^{2}}{6^2}{2\!\cdot\!3^{6}\!+\!1\!\cdot\!3^{7}\phantom{\,\,\!+\!2\!\cdot\!3^8}\,\!+\!O(3^{10})}{1\!\cdot\!3^{5}\phantom{\,\,\!+\!2\!\cdot\!3^6}\,\!+\!1\!\cdot\!3^{7}\!+\!O(3^{10})}\errentry{8\!\cdot\!10^{-13}}{1\!\cdot\!10^{-12}}\cr
\tabentry{2150}{[0, 1, 0, -36, -140]}{x^3-2150}{2^{2}3\!\cdot\!5^{2}43^{2}}{2^{2}3^{2}5^{4}43^{4}}{2^{2}3^{2}5^{2}}{15^2}{1\!\cdot\!3^{9}\!+\!1\!\cdot\!3^{10}\!+\!2\!\cdot\!3^{11}\!+\!O(3^{12})}{2\!\cdot\!3^{5}\!+\!1\!\cdot\!3^{6}\!+\!2\!\cdot\!3^{7}\!+\!O(3^{8})}\errentry{1\!\cdot\!10^{-8}}{4\!\cdot\!10^{-8}}\cr
\tabentry{2350}{[0, 1, 0, -36, -140]}{x^3-2350}{2^{2}3\!\cdot\!5^{2}47^{2}}{2^{2}3^{2}5^{4}47^{4}}{2^{2}3^{4}}{9^2}{1\!\cdot\!3^{7}\phantom{\,\,\!+\!2\!\cdot\!3^8}\,\!+\!1\!\cdot\!3^{9}\!+\!O(3^{10})}{1\!\cdot\!3^{7}\!+\!2\!\cdot\!3^{8}\phantom{\,\,\!+\!2\!\cdot\!3^9}\,\!+\!O(3^{10})}\errentry{3\!\cdot\!10^{-8}}{6\!\cdot\!10^{-8}}\cr
\hline
\tabend

\noindent\begin{tabular}{cllll} 
\extravspace
$E=224B1$ & \multicolumn{3}{l}{$y^2=x^3-x^2+2x$} & \tabname{3}{224B1}\cr
\extravspace
$p=3$ & \multicolumn{3}{l}{$P_3(E/\Q,T)=1-2T+3T^2$}\cr
$/\vQ\alsub{}$   & $N(E)=2^{5}\!\cdot\!7$ & $L^*(E/\Q)=2^{-1}$ & $|E(\Q)|=2$ & $|\sha(E/\Q)|=1^2$ \cr
$/K\alsub{}$      & $N(E,\sigma)=2^{10}\!\cdot\!3^{2}\!\cdot\!7^{2}$ & $L^*(E/K)=2^{-1}$ & $|E(K)|=2$ & $|\sha(E/K)|=1^2$ \cr
$/\vQ\alsub{2}$ & $\hbox{III},\>c_{2}=2$ & $/K\alsub{v}\quad(v|2)$ & $\hbox{III},\>c_v=2$\cr
$/\vQ\alsub{7}$ & $\hbox{I}_1 \hbox{ split},\>c_{7}=1$ & $/K\alsub{v}\quad(v|7)$ & $\hbox{I}_1 \hbox{ split},\>c_v=1$\cr
\end{tabular}
\tabbegin
\tabentry{2}{[0, -1, 0, 2, 0]}{x^3-2}{2^{2}3^{3}}{2^{10}3^{6}7^{2}}{1}{1^2}{1\!\cdot\!3^{0}\!+\!2\!\cdot\!3^{1}\!+\!1\!\cdot\!3^{2}\!+\!O(3^{3})}{1\!\cdot\!3^{0}\phantom{\,\,\!+\!2\!\cdot\!3^1}\,\!+\!2\!\cdot\!3^{2}\!+\!O(3^{4})}\errentry{1\!\cdot\!10^{-16}}{3\!\cdot\!10^{-17}}\cr
\tabentry{6}{[0, -1, 0, 2, 0]}{x^3-6}{2^{2}3^{5}}{2^{10}3^{10}7^{2}}{2^{2}}{2^2}{1\!\cdot\!3^{0}\!+\!2\!\cdot\!3^{1}\!+\!1\!\cdot\!3^{2}\!+\!O(3^{3})}{1\!\cdot\!3^{0}\!+\!2\!\cdot\!3^{1}\!+\!2\!\cdot\!3^{2}\!+\!O(3^{3})}\errentry{2\!\cdot\!10^{-17}}{7\!\cdot\!10^{-15}}\cr
\tabentry{10}{[0, -1, 0, 2, 0]}{x^3-10}{2^{2}3\!\cdot\!5^{2}}{2^{10}3^{2}5^{4}7^{2}}{3^{2}}{3^2}{1\!\cdot\!3^{2}\!+\!1\!\cdot\!3^{3}\!+\!2\!\cdot\!3^{4}\!+\!O(3^{5})}{2\!\cdot\!3^{2}\!+\!1\!\cdot\!3^{3}\!+\!1\!\cdot\!3^{4}\!+\!O(3^{5})}\errentry{1\!\cdot\!10^{-17}}{9\!\cdot\!10^{-16}}\cr
\tabentry{12}{[0, -1, 0, 2, 0]}{x^3-12}{2^{2}3^{5}}{2^{10}3^{10}7^{2}}{1}{1^2}{1\!\cdot\!3^{0}\!+\!2\!\cdot\!3^{1}\!+\!1\!\cdot\!3^{2}\!+\!O(3^{3})}{1\!\cdot\!3^{0}\!+\!1\!\cdot\!3^{1}\!+\!1\!\cdot\!3^{2}\!+\!O(3^{3})}\errentry{1\!\cdot\!10^{-16}}{4\!\cdot\!10^{-16}}\cr
\tabentry{14}{[0, -1, 0, 2, 0]}{x^3-14}{2^{2}3^{3}7^{2}}{2^{10}3^{6}7^{4}}{2^{2}3}{2^2}{1\!\cdot\!3^{2}\!+\!2\!\cdot\!3^{3}\!+\!2\!\cdot\!3^{4}\!+\!O(3^{6})}{1\!\cdot\!3^{1}\!+\!2\!\cdot\!3^{2}\phantom{\,\,\!+\!2\!\cdot\!3^3}\,\!+\!O(3^{4})}\errentry{7\!\cdot\!10^{-17}}{3\!\cdot\!10^{-15}}\cr
\tabentry{20}{[0, -1, 0, 2, 0]}{x^3-20}{2^{2}3^{3}5^{2}}{2^{10}3^{6}5^{4}7^{2}}{2^{2}3^{2}}{6^2}{1\!\cdot\!3^{2}\!+\!1\!\cdot\!3^{3}\!+\!2\!\cdot\!3^{4}\!+\!O(3^{5})}{2\!\cdot\!3^{2}\!+\!2\!\cdot\!3^{3}\!+\!2\!\cdot\!3^{4}\!+\!O(3^{5})}\errentry{7\!\cdot\!10^{-13}}{1\!\cdot\!10^{-12}}\cr
\tabentry{26}{[0, -1, 0, 2, 0]}{x^3-26}{2^{2}3\!\cdot\!13^{2}}{2^{10}3^{2}7^{2}13^{4}}{3^{2}}{3^2}{1\!\cdot\!3^{4}\!+\!1\!\cdot\!3^{5}\!+\!1\!\cdot\!3^{6}\!+\!O(3^{8})}{1\!\cdot\!3^{2}\!+\!1\!\cdot\!3^{3}\!+\!2\!\cdot\!3^{4}\!+\!O(3^{7})}\errentry{7\!\cdot\!10^{-14}}{4\!\cdot\!10^{-14}}\cr
\tabentry{28}{[0, -1, 0, 2, 0]}{x^3-28}{2^{2}3\!\cdot\!7^{2}}{2^{10}3^{2}7^{4}}{3}{1^2}{1\!\cdot\!3^{2}\!+\!2\!\cdot\!3^{3}\!+\!2\!\cdot\!3^{4}\!+\!O(3^{6})}{1\!\cdot\!3^{1}\phantom{\,\,\!+\!2\!\cdot\!3^2}\,\!+\!1\!\cdot\!3^{3}\!+\!O(3^{4})}\errentry{4\!\cdot\!10^{-18}}{1\!\cdot\!10^{-15}}\cr
\tabentry{44}{[0, -1, 0, 2, 0]}{x^3-44}{2^{2}3\!\cdot\!11^{2}}{2^{10}3^{2}7^{2}11^{4}}{1}{1^2}{2\!\cdot\!3^{0}\!+\!2\!\cdot\!3^{1}\phantom{\,\,\!+\!2\!\cdot\!3^2}\,\!+\!O(3^{4})}{2\!\cdot\!3^{0}\!+\!2\!\cdot\!3^{1}\!+\!2\!\cdot\!3^{2}\!+\!O(3^{5})}\errentry{6\!\cdot\!10^{-17}}{3\!\cdot\!10^{-17}}\cr
\tabentry{46}{[0, -1, 0, 2, 0]}{x^3-46}{2^{2}3\!\cdot\!23^{2}}{2^{10}3^{2}7^{2}23^{4}}{2^{6}}{8^2}{2\!\cdot\!3^{0}\!+\!2\!\cdot\!3^{1}\!+\!2\!\cdot\!3^{2}\!+\!O(3^{5})}{2\!\cdot\!3^{0}\!+\!1\!\cdot\!3^{1}\!+\!1\!\cdot\!3^{2}\!+\!O(3^{9})}\errentry{1\!\cdot\!10^{-7}}{2\!\cdot\!10^{-7}}\cr
\tabentry{350}{[0, -1, 0, 2, 0]}{x^3-350}{2^{2}3\!\cdot\!5^{2}7^{2}}{2^{10}3^{2}5^{4}7^{4}}{3^{3}}{3^2}{1\!\cdot\!3^{4}\!+\!1\!\cdot\!3^{5}\phantom{\,\,\!+\!2\!\cdot\!3^6}\,\!+\!O(3^{8})}{2\!\cdot\!3^{3}\phantom{\,\,\!+\!2\!\cdot\!3^4}\,\phantom{\,\,\!+\!2\!\cdot\!3^5}\,\!+\!O(3^{6})}\errentry{9\!\cdot\!10^{-14}}{7\!\cdot\!10^{-14}}\cr
\hline
\tabend

\noindent\begin{tabular}{cllll} 
\extravspace
$E=26A1$ & \multicolumn{3}{l}{$y^2+xy+y=x^3-5x-8$} & \tabname{3}{26A1}\cr
\extravspace
$p=3$ & \multicolumn{3}{l}{$P_3(E/\Q,T)=1-T+3T^2$}\cr
$/\vQ\alsub{}$   & $N(E)=2\!\cdot\!13$ & $L^*(E/\Q)=3^{-1}$ & $|E(\Q)|=3$ & $|\sha(E/\Q)|=1^2$ \cr
$/K\alsub{}$      & $N(E,\sigma)=2^{2}\!\cdot\!3^{2}\!\cdot\!13^{2}$ & $L^*(E/K)=3^{-1}$ & $|E(K)|=9$ & $|\sha(E/K)|=1^2$ \cr
$/\vQ\alsub{2}$ & $\hbox{I}_3 \hbox{ non-split},\>c_{2}=1$ & $/K\alsub{v}\quad(v|2)$ & $\hbox{I}_3 \hbox{ split},\>c_v=3$\cr
$/\vQ\alsub{13}$ & $\hbox{I}_3 \hbox{ split},\>c_{13}=3$ & $/K\alsub{v}\quad(v|13)$ & $\hbox{I}_3 \hbox{ split},\>c_v=3$\cr
\end{tabular}
\tabbegin
\tabentry{3}{[1, 0, 1, -5, -8]}{x^3-3}{3^{5}}{2^{2}3^{10}13^{2}}{2^{2}3}{2^2}{1\!\cdot\!3^{1}\!+\!1\!\cdot\!3^{2}\!+\!1\!\cdot\!3^{3}\!+\!O(3^{5})}{1\!\cdot\!3^{1}\!+\!2\!\cdot\!3^{2}\!+\!2\!\cdot\!3^{3}\!+\!O(3^{4})}\errentry{5\!\cdot\!10^{-17}}{1\!\cdot\!10^{-15}}\cr
\tabentry{5}{[1, 0, 1, -5, -8]}{x^3-5}{3^{3}5^{2}}{2^{2}3^{6}5^{4}13^{2}}{3^{3}}{1^2}{1\!\cdot\!3^{4}\!+\!2\!\cdot\!3^{5}\phantom{\,\,\!+\!2\!\cdot\!3^6}\,\!+\!O(3^{8})}{2\!\cdot\!3^{3}\!+\!1\!\cdot\!3^{4}\!+\!2\!\cdot\!3^{5}\!+\!O(3^{6})}\errentry{4\!\cdot\!10^{-18}}{1\!\cdot\!10^{-16}}\cr
\tabentry{7}{[1, 0, 1, -5, -8]}{x^3-7}{3^{3}7^{2}}{2^{2}3^{6}7^{4}13^{2}}{3^{3}}{3^2}{1\!\cdot\!3^{5}\phantom{\,\,\!+\!2\!\cdot\!3^6}\,\!+\!2\!\cdot\!3^{7}\!+\!O(3^{8})}{1\!\cdot\!3^{3}\phantom{\,\,\!+\!2\!\cdot\!3^4}\,\phantom{\,\,\!+\!2\!\cdot\!3^5}\,\!+\!O(3^{8})}\errentry{4\!\cdot\!10^{-17}}{1\!\cdot\!10^{-15}}\cr
\tabentry{11}{[1, 0, 1, -5, -8]}{x^3-11}{3^{3}11^{2}}{2^{2}3^{6}11^{4}13^{2}}{3^{3}}{3^2}{1\!\cdot\!3^{4}\!+\!1\!\cdot\!3^{5}\phantom{\,\,\!+\!2\!\cdot\!3^6}\,\!+\!O(3^{7})}{2\!\cdot\!3^{3}\!+\!2\!\cdot\!3^{4}\phantom{\,\,\!+\!2\!\cdot\!3^5}\,\!+\!O(3^{7})}\errentry{1\!\cdot\!10^{-15}}{3\!\cdot\!10^{-15}}\cr
\tabentry{13}{[1, 0, 1, -5, -8]}{x^3-13}{3^{3}13^{2}}{2^{2}3^{6}13^{4}}{3^{2}}{1^2}{1\!\cdot\!3^{3}\!+\!1\!\cdot\!3^{4}\!+\!2\!\cdot\!3^{5}\!+\!O(3^{6})}{1\!\cdot\!3^{2}\!+\!1\!\cdot\!3^{3}\!+\!1\!\cdot\!3^{4}\!+\!O(3^{6})}\errentry{9\!\cdot\!10^{-17}}{7\!\cdot\!10^{-17}}\cr
\tabentry{15}{[1, 0, 1, -5, -8]}{x^3-15}{3^{5}5^{2}}{2^{2}3^{10}5^{4}13^{2}}{0}{\nosha}{1\!\cdot\!3^{4}\!+\!2\!\cdot\!3^{5}\phantom{\,\,\!+\!2\!\cdot\!3^6}\,\!+\!O(3^{8})}{\noeuler}\errentry{7\!\cdot\!10^{-12}}{3\!\cdot\!10^{-7}}\cr
\tabentry{17}{[1, 0, 1, -5, -8]}{x^3-17}{3\!\cdot\!17^{2}}{2^{2}3^{2}13^{2}17^{4}}{2^{2}3}{2^2}{2\!\cdot\!3^{3}\!+\!1\!\cdot\!3^{4}\!+\!1\!\cdot\!3^{5}\!+\!O(3^{6})}{2\!\cdot\!3^{3}\phantom{\,\,\!+\!2\!\cdot\!3^4}\,\!+\!1\!\cdot\!3^{5}\!+\!O(3^{6})}\errentry{7\!\cdot\!10^{-18}}{3\!\cdot\!10^{-15}}\cr
\tabentry{19}{[1, 0, 1, -5, -8]}{x^3-19}{3\!\cdot\!19^{2}}{2^{2}3^{2}13^{2}19^{4}}{3^{3}}{3^2}{1\!\cdot\!3^{5}\!+\!2\!\cdot\!3^{6}\!+\!1\!\cdot\!3^{7}\!+\!O(3^{8})}{1\!\cdot\!3^{5}\!+\!1\!\cdot\!3^{6}\!+\!2\!\cdot\!3^{7}\!+\!O(3^{9})}\errentry{2\!\cdot\!10^{-18}}{2\!\cdot\!10^{-15}}\cr
\tabentry{21}{[1, 0, 1, -5, -8]}{x^3-21}{3^{5}7^{2}}{2^{2}3^{10}7^{4}13^{2}}{0}{\nosha}{1\!\cdot\!3^{5}\phantom{\,\,\!+\!2\!\cdot\!3^6}\,\!+\!2\!\cdot\!3^{7}\!+\!O(3^{8})}{\noeuler}\errentry{1\!\cdot\!10^{-8}}{1\!\cdot\!10^{-5}}\cr
\tabentry{35}{[1, 0, 1, -5, -8]}{x^3-35}{3\!\cdot\!5^{2}7^{2}}{2^{2}3^{2}5^{4}7^{4}13^{2}}{3^{3}}{3^2}{1\!\cdot\!3^{8}\!+\!1\!\cdot\!3^{9}\phantom{\,\,\!+\!2\!\cdot\!3^10}\,\!+\!O(3^{11})}{2\!\cdot\!3^{5}\!+\!1\!\cdot\!3^{6}\phantom{\,\,\!+\!2\!\cdot\!3^7}\,\!+\!O(3^{9})}\errentry{4\!\cdot\!10^{-15}}{5\!\cdot\!10^{-15}}\cr
\tabentry{37}{[1, 0, 1, -5, -8]}{x^3-37}{3\!\cdot\!37^{2}}{2^{2}3^{2}13^{2}37^{4}}{3^{3}}{3^2}{1\!\cdot\!3^{5}\phantom{\,\,\!+\!2\!\cdot\!3^6}\,\!+\!1\!\cdot\!3^{7}\!+\!O(3^{9})}{1\!\cdot\!3^{5}\!+\!1\!\cdot\!3^{6}\phantom{\,\,\!+\!2\!\cdot\!3^7}\,\!+\!O(3^{8})}\errentry{5\!\cdot\!10^{-14}}{4\!\cdot\!10^{-14}}\cr
\tabentry{39}{[1, 0, 1, -5, -8]}{x^3-39}{3^{5}13^{2}}{2^{2}3^{10}13^{4}}{3^{4}}{3^2}{1\!\cdot\!3^{3}\!+\!1\!\cdot\!3^{4}\!+\!2\!\cdot\!3^{5}\!+\!O(3^{6})}{1\!\cdot\!3^{4}\phantom{\,\,\!+\!2\!\cdot\!3^5}\,\phantom{\,\,\!+\!2\!\cdot\!3^6}\,\!+\!O(3^{8})}\errentry{1\!\cdot\!10^{-16}}{3\!\cdot\!10^{-15}}\cr
\tabentry{45}{[1, 0, 1, -5, -8]}{x^3-45}{3^{5}5^{2}}{2^{2}3^{10}5^{4}13^{2}}{3^{3}}{3^2}{1\!\cdot\!3^{4}\!+\!2\!\cdot\!3^{5}\phantom{\,\,\!+\!2\!\cdot\!3^6}\,\!+\!O(3^{8})}{2\!\cdot\!3^{3}\!+\!2\!\cdot\!3^{4}\phantom{\,\,\!+\!2\!\cdot\!3^5}\,\!+\!O(3^{6})}\errentry{1\!\cdot\!10^{-11}}{8\!\cdot\!10^{-12}}\cr
\tabentry{53}{[1, 0, 1, -5, -8]}{x^3-53}{3\!\cdot\!53^{2}}{2^{2}3^{2}13^{2}53^{4}}{2^{2}3^{3}}{2^2}{1\!\cdot\!3^{7}\!+\!2\!\cdot\!3^{8}\!+\!2\!\cdot\!3^{9}\!+\!O(3^{10})}{2\!\cdot\!3^{5}\phantom{\,\,\!+\!2\!\cdot\!3^6}\,\phantom{\,\,\!+\!2\!\cdot\!3^7}\,\!+\!O(3^{8})}\errentry{4\!\cdot\!10^{-11}}{2\!\cdot\!10^{-11}}\cr
\tabentry{55}{[1, 0, 1, -5, -8]}{x^3-55}{3\!\cdot\!5^{2}11^{2}}{2^{2}3^{2}5^{4}11^{4}13^{2}}{0}{\nosha}{1\!\cdot\!3^{7}\!+\!2\!\cdot\!3^{8}\!+\!2\!\cdot\!3^{9}\!+\!O(3^{11})}{\noeuler}\errentry{3\!\cdot\!10^{-10}}{2\!\cdot\!10^{-6}}\cr
\tabentry{63}{[1, 0, 1, -5, -8]}{x^3-63}{3^{5}7^{2}}{2^{2}3^{10}7^{4}13^{2}}{3^{3}}{3^2}{1\!\cdot\!3^{5}\phantom{\,\,\!+\!2\!\cdot\!3^6}\,\!+\!2\!\cdot\!3^{7}\!+\!O(3^{8})}{1\!\cdot\!3^{3}\!+\!2\!\cdot\!3^{4}\!+\!2\!\cdot\!3^{5}\!+\!O(3^{7})}\errentry{1\!\cdot\!10^{-8}}{1\!\cdot\!10^{-8}}\cr
\tabentry{65}{[1, 0, 1, -5, -8]}{x^3-65}{3^{3}5^{2}13^{2}}{2^{2}3^{6}5^{4}13^{4}}{2^{2}3^{4}}{6^2}{1\!\cdot\!3^{6}\!+\!2\!\cdot\!3^{7}\!+\!1\!\cdot\!3^{8}\!+\!O(3^{11})}{2\!\cdot\!3^{4}\!+\!1\!\cdot\!3^{5}\phantom{\,\,\!+\!2\!\cdot\!3^6}\,\!+\!O(3^{7})}\errentry{2\!\cdot\!10^{-9}}{3\!\cdot\!10^{-9}}\cr
\tabentry{91}{[1, 0, 1, -5, -8]}{x^3-91}{3\!\cdot\!7^{2}13^{2}}{2^{2}3^{2}7^{4}13^{4}}{3^{4}}{3^2}{1\!\cdot\!3^{7}\phantom{\,\,\!+\!2\!\cdot\!3^8}\,\phantom{\,\,\!+\!2\!\cdot\!3^9}\,\!+\!O(3^{10})}{1\!\cdot\!3^{6}\!+\!1\!\cdot\!3^{7}\phantom{\,\,\!+\!2\!\cdot\!3^8}\,\!+\!O(3^{9})}\errentry{2\!\cdot\!10^{-16}}{6\!\cdot\!10^{-16}}\cr
\tabentry{117}{[1, 0, 1, -5, -8]}{x^3-117}{3^{5}13^{2}}{2^{2}3^{10}13^{4}}{3^{2}}{1^2}{1\!\cdot\!3^{3}\!+\!1\!\cdot\!3^{4}\!+\!2\!\cdot\!3^{5}\!+\!O(3^{6})}{1\!\cdot\!3^{2}\phantom{\,\,\!+\!2\!\cdot\!3^3}\,\phantom{\,\,\!+\!2\!\cdot\!3^4}\,\!+\!O(3^{6})}\errentry{8\!\cdot\!10^{-18}}{7\!\cdot\!10^{-17}}\cr
\tabentry{143}{[1, 0, 1, -5, -8]}{x^3-143}{3\!\cdot\!11^{2}13^{2}}{2^{2}3^{2}11^{4}13^{4}}{2^{2}3^{4}}{6^2}{1\!\cdot\!3^{6}\!+\!1\!\cdot\!3^{7}\!+\!1\!\cdot\!3^{8}\!+\!O(3^{9})}{2\!\cdot\!3^{6}\phantom{\,\,\!+\!2\!\cdot\!3^7}\,\!+\!2\!\cdot\!3^{8}\!+\!O(3^{9})}\errentry{2\!\cdot\!10^{-14}}{4\!\cdot\!10^{-14}}\cr
\tabentry{325}{[1, 0, 1, -5, -8]}{x^3-325}{3\!\cdot\!5^{2}13^{2}}{2^{2}3^{2}5^{4}13^{4}}{0}{\nosha}{1\!\cdot\!3^{6}\!+\!2\!\cdot\!3^{7}\!+\!1\!\cdot\!3^{8}\!+\!O(3^{11})}{\noeuler}\errentry{6\!\cdot\!10^{-17}}{4\!\cdot\!10^{-10}}\cr
\hline
\tabend

\noindent\begin{tabular}{cllll} 
\extravspace
$E=50B1$ & \multicolumn{3}{l}{$y^2+xy+y=x^3+x^2-3x+1$} & \tabname{3}{50B1}\cr
\extravspace
$p=3$ & \multicolumn{3}{l}{$P_3(E/\Q,T)=1+T+3T^2$}\cr
$/\vQ\alsub{}$   & $N(E)=2\!\cdot\!5^{2}$ & $L^*(E/\Q)=5^{-1}$ & $|E(\Q)|=5$ & $|\sha(E/\Q)|=1^2$ \cr
$/K\alsub{}$      & $N(E,\sigma)=2^{2}\!\cdot\!3^{2}\!\cdot\!5^{4}$ & $L^*(E/K)=5^{-1}$ & $|E(K)|=5$ & $|\sha(E/K)|=1^2$ \cr
$/\vQ\alsub{2}$ & $\hbox{I}_5 \hbox{ split},\>c_{2}=5$ & $/K\alsub{v}\quad(v|2)$ & $\hbox{I}_5 \hbox{ split},\>c_v=5$\cr
$/\vQ\alsub{5}$ & $\hbox{II},\>c_{5}=1$ & $/K\alsub{v}\quad(v|5)$ & $\hbox{II},\>c_v=1$\cr
\end{tabular}
\tabbegin
\tabentry{3}{[1, 1, 1, -3, 1]}{x^3-3}{3^{5}}{2^{2}3^{10}5^{4}}{5}{1^2}{1\!\cdot\!3^{0}\!+\!1\!\cdot\!3^{1}\!+\!2\!\cdot\!3^{2}\!+\!O(3^{3})}{1\!\cdot\!3^{0}\!+\!2\!\cdot\!3^{1}\!+\!1\!\cdot\!3^{2}\!+\!O(3^{5})}\errentry{1\!\cdot\!10^{-16}}{1\!\cdot\!10^{-15}}\cr
\tabentry{5}{[1, 1, 1, -3, 1]}{x^3-5}{3^{3}5^{2}}{2^{2}3^{6}5^{4}}{2\!\cdot\!5}{1^2}{1\!\cdot\!3^{0}\!+\!1\!\cdot\!3^{1}\!+\!2\!\cdot\!3^{2}\!+\!O(3^{3})}{1\!\cdot\!3^{0}\!+\!1\!\cdot\!3^{1}\!+\!1\!\cdot\!3^{2}\!+\!O(3^{3})}\errentry{3\!\cdot\!10^{-17}}{9\!\cdot\!10^{-16}}\cr
\tabentry{7}{[1, 1, 1, -3, 1]}{x^3-7}{3^{3}7^{2}}{2^{2}3^{6}5^{4}7^{4}}{5}{1^2}{1\!\cdot\!3^{0}\phantom{\,\,\!+\!2\!\cdot\!3^1}\,\!+\!2\!\cdot\!3^{2}\!+\!O(3^{4})}{1\!\cdot\!3^{0}\!+\!1\!\cdot\!3^{1}\!+\!2\!\cdot\!3^{2}\!+\!O(3^{3})}\errentry{1\!\cdot\!10^{-16}}{2\!\cdot\!10^{-15}}\cr
\tabentry{11}{[1, 1, 1, -3, 1]}{x^3-11}{3^{3}11^{2}}{2^{2}3^{6}5^{4}11^{4}}{3^{2}5}{3^2}{2\!\cdot\!3^{3}\!+\!1\!\cdot\!3^{4}\!+\!2\!\cdot\!3^{5}\!+\!O(3^{6})}{2\!\cdot\!3^{2}\!+\!1\!\cdot\!3^{3}\!+\!1\!\cdot\!3^{4}\!+\!O(3^{5})}\errentry{6\!\cdot\!10^{-12}}{4\!\cdot\!10^{-12}}\cr
\tabentry{13}{[1, 1, 1, -3, 1]}{x^3-13}{3^{3}13^{2}}{2^{2}3^{6}5^{4}13^{4}}{5}{1^2}{1\!\cdot\!3^{0}\!+\!2\!\cdot\!3^{1}\!+\!2\!\cdot\!3^{2}\!+\!O(3^{3})}{1\!\cdot\!3^{0}\!+\!2\!\cdot\!3^{1}\!+\!1\!\cdot\!3^{2}\!+\!O(3^{5})}\errentry{2\!\cdot\!10^{-8}}{6\!\cdot\!10^{-9}}\cr
\tabentry{15}{[1, 1, 1, -3, 1]}{x^3-15}{3^{5}5^{2}}{2^{2}3^{10}5^{4}}{2\!\cdot\!5}{1^2}{1\!\cdot\!3^{0}\!+\!1\!\cdot\!3^{1}\!+\!2\!\cdot\!3^{2}\!+\!O(3^{3})}{1\!\cdot\!3^{0}\phantom{\,\,\!+\!2\!\cdot\!3^1}\,\phantom{\,\,\!+\!2\!\cdot\!3^2}\,\!+\!O(3^{3})}\errentry{6\!\cdot\!10^{-17}}{2\!\cdot\!10^{-15}}\cr
\tabentry{17}{[1, 1, 1, -3, 1]}{x^3-17}{3\!\cdot\!17^{2}}{2^{2}3^{2}5^{4}17^{4}}{0}{\nosha}{2\!\cdot\!3^{2}\!+\!1\!\cdot\!3^{3}\phantom{\,\,\!+\!2\!\cdot\!3^4}\,\!+\!O(3^{5})}{\noeuler}\errentry{6\!\cdot\!10^{-18}}{8\!\cdot\!10^{-10}}\cr
\tabentry{19}{[1, 1, 1, -3, 1]}{x^3-19}{3\!\cdot\!19^{2}}{2^{2}3^{2}5^{4}19^{4}}{0}{\nosha}{1\!\cdot\!3^{2}\phantom{\,\,\!+\!2\!\cdot\!3^3}\,\phantom{\,\,\!+\!2\!\cdot\!3^4}\,\!+\!O(3^{5})}{\noeuler}\errentry{5\!\cdot\!10^{-17}}{4\!\cdot\!10^{-10}}\cr
\tabentry{35}{[1, 1, 1, -3, 1]}{x^3-35}{3\!\cdot\!5^{2}7^{2}}{2^{2}3^{2}5^{4}7^{4}}{2\!\cdot\!5}{1^2}{1\!\cdot\!3^{0}\phantom{\,\,\!+\!2\!\cdot\!3^1}\,\!+\!2\!\cdot\!3^{2}\!+\!O(3^{4})}{1\!\cdot\!3^{0}\!+\!2\!\cdot\!3^{1}\!+\!1\!\cdot\!3^{2}\!+\!O(3^{3})}\errentry{1\!\cdot\!10^{-16}}{1\!\cdot\!10^{-15}}\cr
\tabentry{37}{[1, 1, 1, -3, 1]}{x^3-37}{3\!\cdot\!37^{2}}{2^{2}3^{2}5^{4}37^{4}}{5}{1^2}{1\!\cdot\!3^{0}\phantom{\,\,\!+\!2\!\cdot\!3^1}\,\phantom{\,\,\!+\!2\!\cdot\!3^2}\,\!+\!O(3^{4})}{1\!\cdot\!3^{0}\phantom{\,\,\!+\!2\!\cdot\!3^1}\,\phantom{\,\,\!+\!2\!\cdot\!3^2}\,\!+\!O(3^{3})}\errentry{1\!\cdot\!10^{-10}}{3\!\cdot\!10^{-11}}\cr
\tabentry{45}{[1, 1, 1, -3, 1]}{x^3-45}{3^{5}5^{2}}{2^{2}3^{10}5^{4}}{2^{3}5}{2^2}{1\!\cdot\!3^{0}\!+\!1\!\cdot\!3^{1}\!+\!2\!\cdot\!3^{2}\!+\!O(3^{3})}{1\!\cdot\!3^{0}\!+\!1\!\cdot\!3^{1}\phantom{\,\,\!+\!2\!\cdot\!3^2}\,\!+\!O(3^{3})}\errentry{6\!\cdot\!10^{-17}}{4\!\cdot\!10^{-16}}\cr
\tabentry{55}{[1, 1, 1, -3, 1]}{x^3-55}{3\!\cdot\!5^{2}11^{2}}{2^{2}3^{2}5^{4}11^{4}}{0}{\nosha}{2\!\cdot\!3^{3}\!+\!1\!\cdot\!3^{4}\!+\!2\!\cdot\!3^{5}\!+\!O(3^{6})}{\noeuler}\errentry{1\!\cdot\!10^{-15}}{2\!\cdot\!10^{-9}}\cr
\tabentry{65}{[1, 1, 1, -3, 1]}{x^3-65}{3^{3}5^{2}13^{2}}{2^{2}3^{6}5^{4}13^{4}}{2^{3}5}{2^2}{1\!\cdot\!3^{0}\!+\!2\!\cdot\!3^{1}\!+\!2\!\cdot\!3^{2}\!+\!O(3^{3})}{1\!\cdot\!3^{0}\!+\!1\!\cdot\!3^{1}\phantom{\,\,\!+\!2\!\cdot\!3^2}\,\!+\!O(3^{3})}\errentry{2\!\cdot\!10^{-10}}{1\!\cdot\!10^{-10}}\cr
\tabentry{145}{[1, 1, 1, -3, 1]}{x^3-145}{3\!\cdot\!5^{2}29^{2}}{2^{2}3^{2}5^{4}29^{4}}{2\!\cdot\!3^{2}5}{3^2}{1\!\cdot\!3^{2}\phantom{\,\,\!+\!2\!\cdot\!3^3}\,\!+\!1\!\cdot\!3^{4}\!+\!O(3^{5})}{2\!\cdot\!3^{2}\phantom{\,\,\!+\!2\!\cdot\!3^3}\,\phantom{\,\,\!+\!2\!\cdot\!3^4}\,\!+\!O(3^{5})}\errentry{8\!\cdot\!10^{-14}}{6\!\cdot\!10^{-14}}\cr
\tabentry{175}{[1, 1, 1, -3, 1]}{x^3-175}{3^{3}5^{2}7^{2}}{2^{2}3^{6}5^{4}7^{4}}{2\!\cdot\!5}{1^2}{1\!\cdot\!3^{0}\phantom{\,\,\!+\!2\!\cdot\!3^1}\,\!+\!2\!\cdot\!3^{2}\!+\!O(3^{4})}{1\!\cdot\!3^{0}\!+\!2\!\cdot\!3^{1}\!+\!2\!\cdot\!3^{2}\!+\!O(3^{3})}\errentry{9\!\cdot\!10^{-16}}{6\!\cdot\!10^{-17}}\cr
\tabentry{215}{[1, 1, 1, -3, 1]}{x^3-215}{3\!\cdot\!5^{2}43^{2}}{2^{2}3^{2}5^{4}43^{4}}{2\!\cdot\!5}{1^2}{1\!\cdot\!3^{0}\!+\!2\!\cdot\!3^{1}\!+\!1\!\cdot\!3^{2}\!+\!O(3^{3})}{1\!\cdot\!3^{0}\!+\!2\!\cdot\!3^{1}\phantom{\,\,\!+\!2\!\cdot\!3^2}\,\!+\!O(3^{3})}\errentry{8\!\cdot\!10^{-9}}{2\!\cdot\!10^{-9}}\cr
\tabentry{235}{[1, 1, 1, -3, 1]}{x^3-235}{3\!\cdot\!5^{2}47^{2}}{2^{2}3^{2}5^{4}47^{4}}{2^{3}3^{2}5}{6^2}{2\!\cdot\!3^{3}\!+\!2\!\cdot\!3^{4}\!+\!2\!\cdot\!3^{5}\!+\!O(3^{6})}{2\!\cdot\!3^{2}\!+\!2\!\cdot\!3^{3}\!+\!1\!\cdot\!3^{4}\!+\!O(3^{5})}\errentry{5\!\cdot\!10^{-8}}{7\!\cdot\!10^{-8}}\cr
\tabentry{275}{[1, 1, 1, -3, 1]}{x^3-275}{3^{3}5^{2}11^{2}}{2^{2}3^{6}5^{4}11^{4}}{2\!\cdot\!3^{2}5}{3^2}{2\!\cdot\!3^{3}\!+\!1\!\cdot\!3^{4}\!+\!2\!\cdot\!3^{5}\!+\!O(3^{6})}{2\!\cdot\!3^{2}\phantom{\,\,\!+\!2\!\cdot\!3^3}\,\!+\!1\!\cdot\!3^{4}\!+\!O(3^{5})}\errentry{6\!\cdot\!10^{-12}}{5\!\cdot\!10^{-12}}\cr
\tabentry{325}{[1, 1, 1, -3, 1]}{x^3-325}{3\!\cdot\!5^{2}13^{2}}{2^{2}3^{2}5^{4}13^{4}}{2\!\cdot\!5}{1^2}{1\!\cdot\!3^{0}\!+\!2\!\cdot\!3^{1}\!+\!2\!\cdot\!3^{2}\!+\!O(3^{3})}{1\!\cdot\!3^{0}\phantom{\,\,\!+\!2\!\cdot\!3^1}\,\!+\!2\!\cdot\!3^{2}\!+\!O(3^{4})}\errentry{7\!\cdot\!10^{-17}}{4\!\cdot\!10^{-16}}\cr
\tabentry{575}{[1, 1, 1, -3, 1]}{x^3-575}{3\!\cdot\!5^{2}23^{2}}{2^{2}3^{2}5^{4}23^{4}}{0}{\nosha}{2\!\cdot\!3^{3}\!+\!1\!\cdot\!3^{4}\!+\!1\!\cdot\!3^{5}\!+\!O(3^{6})}{\noeuler}\errentry{5\!\cdot\!10^{-17}}{3\!\cdot\!10^{-9}}\cr
\tabentry{775}{[1, 1, 1, -3, 1]}{x^3-775}{3\!\cdot\!5^{2}31^{2}}{2^{2}3^{2}5^{4}31^{4}}{2\!\cdot\!3^{2}5}{3^2}{1\!\cdot\!3^{2}\!+\!2\!\cdot\!3^{3}\!+\!1\!\cdot\!3^{4}\!+\!O(3^{5})}{1\!\cdot\!3^{2}\phantom{\,\,\!+\!2\!\cdot\!3^3}\,\phantom{\,\,\!+\!2\!\cdot\!3^4}\,\!+\!O(3^{5})}\errentry{9\!\cdot\!10^{-12}}{6\!\cdot\!10^{-12}}\cr
\tabentry{1025}{[1, 1, 1, -3, 1]}{x^3-1025}{3\!\cdot\!5^{2}41^{2}}{2^{2}3^{2}5^{4}41^{4}}{2\!\cdot\!3^{2}5}{3^2}{2\!\cdot\!3^{3}\!+\!1\!\cdot\!3^{4}\!+\!2\!\cdot\!3^{5}\!+\!O(3^{7})}{2\!\cdot\!3^{2}\!+\!2\!\cdot\!3^{3}\phantom{\,\,\!+\!2\!\cdot\!3^4}\,\!+\!O(3^{7})}\errentry{3\!\cdot\!10^{-9}}{2\!\cdot\!10^{-9}}\cr
\hline
\tabend

\noindent\begin{tabular}{cllll} 
\extravspace
$E=35A1$ & \multicolumn{3}{l}{$y^2+y=x^3+x^2+9x+1$} & \tabname{3}{35A1}\cr
\extravspace
$p=3$ & \multicolumn{3}{l}{$P_3(E/\Q,T)=1-T+3T^2$}\cr
$/\vQ\alsub{}$   & $N(E)=5\!\cdot\!7$ & $L^*(E/\Q)=3^{-1}$ & $|E(\Q)|=3$ & $|\sha(E/\Q)|=1^2$ \cr
$/K\alsub{}$      & $N(E,\sigma)=3^{2}\!\cdot\!5^{2}\!\cdot\!7^{2}$ & $L^*(E/K)=3^{-1}$ & $|E(K)|=9$ & $|\sha(E/K)|=1^2$ \cr
$/\vQ\alsub{5}$ & $\hbox{I}_3 \hbox{ non-split},\>c_{5}=1$ & $/K\alsub{v}\quad(v|5)$ & $\hbox{I}_3 \hbox{ split},\>c_v=3$\cr
$/\vQ\alsub{7}$ & $\hbox{I}_3 \hbox{ split},\>c_{7}=3$ & $/K\alsub{v}\quad(v|7)$ & $\hbox{I}_3 \hbox{ split},\>c_v=3$\cr
\end{tabular}
\tabbegin
\tabentry{2}{[0, 1, 1, 9, 1]}{x^3-2}{2^{2}3^{3}}{2^{4}3^{6}5^{2}7^{2}}{0}{\nosha}{1\!\cdot\!3^{3}\phantom{\,\,\!+\!2\!\cdot\!3^4}\,\!+\!1\!\cdot\!3^{5}\!+\!O(3^{6})}{\noeuler}\errentry{2\!\cdot\!10^{-18}}{1\!\cdot\!10^{-9}}\cr
\tabentry{3}{[0, 1, 1, 9, 1]}{x^3-3}{3^{5}}{3^{10}5^{2}7^{2}}{3}{1^2}{1\!\cdot\!3^{1}\!+\!1\!\cdot\!3^{2}\!+\!1\!\cdot\!3^{3}\!+\!O(3^{5})}{1\!\cdot\!3^{1}\!+\!1\!\cdot\!3^{2}\!+\!1\!\cdot\!3^{3}\!+\!O(3^{5})}\errentry{3\!\cdot\!10^{-16}}{4\!\cdot\!10^{-16}}\cr
\tabentry{6}{[0, 1, 1, 9, 1]}{x^3-6}{2^{2}3^{5}}{2^{4}3^{10}5^{2}7^{2}}{3^{3}}{1^2}{1\!\cdot\!3^{3}\phantom{\,\,\!+\!2\!\cdot\!3^4}\,\!+\!1\!\cdot\!3^{5}\!+\!O(3^{6})}{2\!\cdot\!3^{3}\phantom{\,\,\!+\!2\!\cdot\!3^4}\,\!+\!2\!\cdot\!3^{5}\!+\!O(3^{6})}\errentry{3\!\cdot\!10^{-17}}{1\!\cdot\!10^{-16}}\cr
\tabentry{7}{[0, 1, 1, 9, 1]}{x^3-7}{3^{3}7^{2}}{3^{6}5^{2}7^{4}}{3^{2}}{1^2}{1\!\cdot\!3^{3}\!+\!1\!\cdot\!3^{4}\!+\!2\!\cdot\!3^{5}\!+\!O(3^{6})}{1\!\cdot\!3^{2}\phantom{\,\,\!+\!2\!\cdot\!3^3}\,\phantom{\,\,\!+\!2\!\cdot\!3^4}\,\!+\!O(3^{7})}\errentry{1\!\cdot\!10^{-16}}{7\!\cdot\!10^{-16}}\cr
\tabentry{11}{[0, 1, 1, 9, 1]}{x^3-11}{3^{3}11^{2}}{3^{6}5^{2}7^{2}11^{4}}{0}{\nosha}{2\!\cdot\!3^{4}\!+\!1\!\cdot\!3^{5}\phantom{\,\,\!+\!2\!\cdot\!3^6}\,\!+\!O(3^{7})}{\noeuler}\errentry{1\!\cdot\!10^{-13}}{1\!\cdot\!10^{-8}}\cr
\tabentry{12}{[0, 1, 1, 9, 1]}{x^3-12}{2^{2}3^{5}}{2^{4}3^{10}5^{2}7^{2}}{3^{3}}{3^2}{1\!\cdot\!3^{3}\phantom{\,\,\!+\!2\!\cdot\!3^4}\,\!+\!1\!\cdot\!3^{5}\!+\!O(3^{6})}{2\!\cdot\!3^{3}\phantom{\,\,\!+\!2\!\cdot\!3^4}\,\!+\!2\!\cdot\!3^{5}\!+\!O(3^{6})}\errentry{1\!\cdot\!10^{-16}}{4\!\cdot\!10^{-16}}\cr
\tabentry{13}{[0, 1, 1, 9, 1]}{x^3-13}{3^{3}13^{2}}{3^{6}5^{2}7^{2}13^{4}}{0}{\nosha}{1\!\cdot\!3^{5}\!+\!2\!\cdot\!3^{6}\!+\!2\!\cdot\!3^{7}\!+\!O(3^{10})}{\noeuler}\errentry{2\!\cdot\!10^{-12}}{7\!\cdot\!10^{-8}}\cr
\tabentry{14}{[0, 1, 1, 9, 1]}{x^3-14}{2^{2}3^{3}7^{2}}{2^{4}3^{6}5^{2}7^{4}}{3^{4}}{3^2}{1\!\cdot\!3^{5}\phantom{\,\,\!+\!2\!\cdot\!3^6}\,\!+\!2\!\cdot\!3^{7}\!+\!O(3^{8})}{2\!\cdot\!3^{4}\!+\!1\!\cdot\!3^{5}\!+\!1\!\cdot\!3^{6}\!+\!O(3^{7})}\errentry{7\!\cdot\!10^{-18}}{8\!\cdot\!10^{-16}}\cr
\tabentry{17}{[0, 1, 1, 9, 1]}{x^3-17}{3\!\cdot\!17^{2}}{3^{2}5^{2}7^{2}17^{4}}{3}{1^2}{2\!\cdot\!3^{3}\!+\!1\!\cdot\!3^{4}\!+\!1\!\cdot\!3^{5}\!+\!O(3^{6})}{2\!\cdot\!3^{3}\!+\!1\!\cdot\!3^{4}\!+\!2\!\cdot\!3^{5}\!+\!O(3^{6})}\errentry{5\!\cdot\!10^{-16}}{2\!\cdot\!10^{-16}}\cr
\tabentry{19}{[0, 1, 1, 9, 1]}{x^3-19}{3\!\cdot\!19^{2}}{3^{2}5^{2}7^{2}19^{4}}{0}{\nosha}{1\!\cdot\!3^{5}\!+\!2\!\cdot\!3^{6}\!+\!1\!\cdot\!3^{7}\!+\!O(3^{8})}{\noeuler}\errentry{4\!\cdot\!10^{-16}}{1\!\cdot\!10^{-9}}\cr
\tabentry{21}{[0, 1, 1, 9, 1]}{x^3-21}{3^{5}7^{2}}{3^{10}5^{2}7^{4}}{2^{2}3^{2}}{2^2}{1\!\cdot\!3^{3}\!+\!1\!\cdot\!3^{4}\!+\!2\!\cdot\!3^{5}\!+\!O(3^{6})}{1\!\cdot\!3^{2}\phantom{\,\,\!+\!2\!\cdot\!3^3}\,\!+\!2\!\cdot\!3^{4}\!+\!O(3^{7})}\errentry{3\!\cdot\!10^{-16}}{3\!\cdot\!10^{-15}}\cr
\tabentry{26}{[0, 1, 1, 9, 1]}{x^3-26}{2^{2}3\!\cdot\!13^{2}}{2^{4}3^{2}5^{2}7^{2}13^{4}}{3^{3}}{3^2}{1\!\cdot\!3^{7}\!+\!1\!\cdot\!3^{8}\!+\!1\!\cdot\!3^{9}\!+\!O(3^{10})}{2\!\cdot\!3^{5}\!+\!1\!\cdot\!3^{6}\!+\!1\!\cdot\!3^{7}\!+\!O(3^{8})}\errentry{2\!\cdot\!10^{-18}}{3\!\cdot\!10^{-15}}\cr
\tabentry{28}{[0, 1, 1, 9, 1]}{x^3-28}{2^{2}3\!\cdot\!7^{2}}{2^{4}3^{2}5^{2}7^{4}}{3^{2}}{1^2}{1\!\cdot\!3^{5}\phantom{\,\,\!+\!2\!\cdot\!3^6}\,\!+\!2\!\cdot\!3^{7}\!+\!O(3^{8})}{2\!\cdot\!3^{4}\!+\!2\!\cdot\!3^{5}\!+\!2\!\cdot\!3^{6}\!+\!O(3^{7})}\errentry{1\!\cdot\!10^{-16}}{2\!\cdot\!10^{-15}}\cr
\tabentry{37}{[0, 1, 1, 9, 1]}{x^3-37}{3\!\cdot\!37^{2}}{3^{2}5^{2}7^{2}37^{4}}{3^{3}}{3^2}{1\!\cdot\!3^{5}\phantom{\,\,\!+\!2\!\cdot\!3^6}\,\phantom{\,\,\!+\!2\!\cdot\!3^7}\,\!+\!O(3^{8})}{1\!\cdot\!3^{5}\!+\!1\!\cdot\!3^{6}\phantom{\,\,\!+\!2\!\cdot\!3^7}\,\!+\!O(3^{8})}\errentry{4\!\cdot\!10^{-13}}{3\!\cdot\!10^{-13}}\cr
\tabentry{42}{[0, 1, 1, 9, 1]}{x^3-42}{2^{2}3^{5}7^{2}}{2^{4}3^{10}5^{2}7^{4}}{0}{\nosha}{1\!\cdot\!3^{5}\phantom{\,\,\!+\!2\!\cdot\!3^6}\,\!+\!2\!\cdot\!3^{7}\!+\!O(3^{8})}{\noeuler}\errentry{1\!\cdot\!10^{-8}}{9\!\cdot\!10^{-6}}\cr
\tabentry{44}{[0, 1, 1, 9, 1]}{x^3-44}{2^{2}3\!\cdot\!11^{2}}{2^{4}3^{2}5^{2}7^{2}11^{4}}{3^{3}}{3^2}{2\!\cdot\!3^{6}\!+\!2\!\cdot\!3^{7}\phantom{\,\,\!+\!2\!\cdot\!3^8}\,\!+\!O(3^{12})}{1\!\cdot\!3^{5}\phantom{\,\,\!+\!2\!\cdot\!3^6}\,\!+\!2\!\cdot\!3^{7}\!+\!O(3^{9})}\errentry{3\!\cdot\!10^{-17}}{4\!\cdot\!10^{-17}}\cr
\tabentry{46}{[0, 1, 1, 9, 1]}{x^3-46}{2^{2}3\!\cdot\!23^{2}}{2^{4}3^{2}5^{2}7^{2}23^{4}}{3^{3}}{3^2}{2\!\cdot\!3^{6}\!+\!2\!\cdot\!3^{7}\!+\!2\!\cdot\!3^{8}\!+\!O(3^{10})}{1\!\cdot\!3^{5}\!+\!1\!\cdot\!3^{6}\!+\!2\!\cdot\!3^{7}\!+\!O(3^{8})}\errentry{3\!\cdot\!10^{-9}}{2\!\cdot\!10^{-9}}\cr
\tabentry{53}{[0, 1, 1, 9, 1]}{x^3-53}{3\!\cdot\!53^{2}}{3^{2}5^{2}7^{2}53^{4}}{2^{2}3}{2^2}{2\!\cdot\!3^{3}\!+\!2\!\cdot\!3^{4}\phantom{\,\,\!+\!2\!\cdot\!3^5}\,\!+\!O(3^{7})}{2\!\cdot\!3^{3}\phantom{\,\,\!+\!2\!\cdot\!3^4}\,\phantom{\,\,\!+\!2\!\cdot\!3^5}\,\!+\!O(3^{6})}\errentry{1\!\cdot\!10^{-9}}{7\!\cdot\!10^{-10}}\cr
\tabentry{63}{[0, 1, 1, 9, 1]}{x^3-63}{3^{5}7^{2}}{3^{10}5^{2}7^{4}}{0}{\nosha}{1\!\cdot\!3^{3}\!+\!1\!\cdot\!3^{4}\!+\!2\!\cdot\!3^{5}\!+\!O(3^{6})}{\noeuler}\errentry{1\!\cdot\!10^{-16}}{2\!\cdot\!10^{-9}}\cr
\tabentry{84}{[0, 1, 1, 9, 1]}{x^3-84}{2^{2}3^{5}7^{2}}{2^{4}3^{10}5^{2}7^{4}}{3^{4}}{3^2}{1\!\cdot\!3^{5}\phantom{\,\,\!+\!2\!\cdot\!3^6}\,\!+\!2\!\cdot\!3^{7}\!+\!O(3^{8})}{2\!\cdot\!3^{4}\!+\!2\!\cdot\!3^{5}\!+\!2\!\cdot\!3^{6}\!+\!O(3^{7})}\errentry{7\!\cdot\!10^{-8}}{5\!\cdot\!10^{-8}}\cr
\tabentry{91}{[0, 1, 1, 9, 1]}{x^3-91}{3\!\cdot\!7^{2}13^{2}}{3^{2}5^{2}7^{4}13^{4}}{0}{\nosha}{1\!\cdot\!3^{7}\!+\!2\!\cdot\!3^{8}\phantom{\,\,\!+\!2\!\cdot\!3^9}\,\!+\!O(3^{12})}{\noeuler}\errentry{1\!\cdot\!10^{-13}}{2\!\cdot\!10^{-8}}\cr
\tabentry{116}{[0, 1, 1, 9, 1]}{x^3-116}{2^{2}3\!\cdot\!29^{2}}{2^{4}3^{2}5^{2}7^{2}29^{4}}{3^{3}}{3^2}{2\!\cdot\!3^{7}\!+\!1\!\cdot\!3^{8}\!+\!1\!\cdot\!3^{9}\!+\!O(3^{10})}{1\!\cdot\!3^{5}\phantom{\,\,\!+\!2\!\cdot\!3^6}\,\!+\!1\!\cdot\!3^{7}\!+\!O(3^{10})}\errentry{9\!\cdot\!10^{-9}}{6\!\cdot\!10^{-9}}\cr
\tabentry{126}{[0, 1, 1, 9, 1]}{x^3-126}{2^{2}3^{5}7^{2}}{2^{4}3^{10}5^{2}7^{4}}{0}{\nosha}{1\!\cdot\!3^{5}\phantom{\,\,\!+\!2\!\cdot\!3^6}\,\!+\!2\!\cdot\!3^{7}\!+\!O(3^{8})}{\noeuler}\errentry{3\!\cdot\!10^{-9}}{2\!\cdot\!10^{-6}}\cr
\tabentry{154}{[0, 1, 1, 9, 1]}{x^3-154}{2^{2}3\!\cdot\!7^{2}11^{2}}{2^{4}3^{2}5^{2}7^{4}11^{4}}{3^{4}}{3^2}{2\!\cdot\!3^{8}\!+\!2\!\cdot\!3^{9}\!+\!2\!\cdot\!3^{10}\!+\!O(3^{11})}{1\!\cdot\!3^{6}\!+\!1\!\cdot\!3^{7}\!+\!2\!\cdot\!3^{8}\!+\!O(3^{9})}\errentry{1\!\cdot\!10^{-7}}{1\!\cdot\!10^{-7}}\cr
\tabentry{252}{[0, 1, 1, 9, 1]}{x^3-252}{2^{2}3^{5}7^{2}}{2^{4}3^{10}5^{2}7^{4}}{3^{4}}{3^2}{1\!\cdot\!3^{5}\phantom{\,\,\!+\!2\!\cdot\!3^6}\,\!+\!2\!\cdot\!3^{7}\!+\!O(3^{8})}{2\!\cdot\!3^{4}\!+\!2\!\cdot\!3^{5}\!+\!2\!\cdot\!3^{6}\!+\!O(3^{7})}\errentry{2\!\cdot\!10^{-10}}{2\!\cdot\!10^{-10}}\cr
\tabentry{539}{[0, 1, 1, 9, 1]}{x^3-539}{3\!\cdot\!7^{2}11^{2}}{3^{2}5^{2}7^{4}11^{4}}{0}{\nosha}{2\!\cdot\!3^{6}\!+\!1\!\cdot\!3^{7}\!+\!2\!\cdot\!3^{8}\!+\!O(3^{9})}{\noeuler}\errentry{5\!\cdot\!10^{-15}}{3\!\cdot\!10^{-9}}\cr
\hline
\tabend

\noindent\begin{tabular}{cllll} 
\extravspace
$E=38B1$ & \multicolumn{3}{l}{$y^2+xy+y=x^3+x^2+1$} & \tabname{3}{38B1}\cr
\extravspace
$p=3$ & \multicolumn{3}{l}{$P_3(E/\Q,T)=1+T+3T^2$}\cr
$/\vQ\alsub{}$   & $N(E)=2\!\cdot\!19$ & $L^*(E/\Q)=5^{-1}$ & $|E(\Q)|=5$ & $|\sha(E/\Q)|=1^2$ \cr
$/K\alsub{}$      & $N(E,\sigma)=2^{2}\!\cdot\!3^{2}\!\cdot\!19^{2}$ & $L^*(E/K)=5^{-1}$ & $|E(K)|=5$ & $|\sha(E/K)|=1^2$ \cr
$/\vQ\alsub{2}$ & $\hbox{I}_5 \hbox{ split},\>c_{2}=5$ & $/K\alsub{v}\quad(v|2)$ & $\hbox{I}_5 \hbox{ split},\>c_v=5$\cr
$/\vQ\alsub{19}$ & $\hbox{I}_1 \hbox{ non-split},\>c_{19}=1$ & $/K\alsub{v}\quad(v|19)$ & $\hbox{I}_1 \hbox{ non-split},\>c_v=1$\cr
\end{tabular}
\tabbegin
\tabentry{3}{[1, 1, 1, 0, 1]}{x^3-3}{3^{5}}{2^{2}3^{10}19^{2}}{5}{1^2}{1\!\cdot\!3^{0}\!+\!1\!\cdot\!3^{1}\!+\!2\!\cdot\!3^{2}\!+\!O(3^{3})}{1\!\cdot\!3^{0}\!+\!2\!\cdot\!3^{1}\!+\!1\!\cdot\!3^{2}\!+\!O(3^{5})}\errentry{2\!\cdot\!10^{-17}}{1\!\cdot\!10^{-16}}\cr
\tabentry{5}{[1, 1, 1, 0, 1]}{x^3-5}{3^{3}5^{2}}{2^{2}3^{6}5^{4}19^{2}}{5}{1^2}{2\!\cdot\!3^{0}\!+\!1\!\cdot\!3^{1}\!+\!2\!\cdot\!3^{2}\!+\!O(3^{4})}{2\!\cdot\!3^{0}\phantom{\,\,\!+\!2\!\cdot\!3^1}\,\!+\!2\!\cdot\!3^{2}\!+\!O(3^{3})}\errentry{1\!\cdot\!10^{-16}}{1\!\cdot\!10^{-16}}\cr
\tabentry{7}{[1, 1, 1, 0, 1]}{x^3-7}{3^{3}7^{2}}{2^{2}3^{6}7^{4}19^{2}}{2^{2}5}{2^2}{1\!\cdot\!3^{0}\!+\!2\!\cdot\!3^{1}\!+\!2\!\cdot\!3^{2}\!+\!O(3^{4})}{1\!\cdot\!3^{0}\!+\!2\!\cdot\!3^{1}\phantom{\,\,\!+\!2\!\cdot\!3^2}\,\!+\!O(3^{3})}\errentry{7\!\cdot\!10^{-17}}{9\!\cdot\!10^{-16}}\cr
\tabentry{11}{[1, 1, 1, 0, 1]}{x^3-11}{3^{3}11^{2}}{2^{2}3^{6}11^{4}19^{2}}{5}{1^2}{2\!\cdot\!3^{0}\!+\!1\!\cdot\!3^{1}\!+\!2\!\cdot\!3^{2}\!+\!O(3^{3})}{2\!\cdot\!3^{0}\!+\!1\!\cdot\!3^{1}\!+\!1\!\cdot\!3^{2}\!+\!O(3^{3})}\errentry{1\!\cdot\!10^{-13}}{4\!\cdot\!10^{-14}}\cr
\tabentry{13}{[1, 1, 1, 0, 1]}{x^3-13}{3^{3}13^{2}}{2^{2}3^{6}13^{4}19^{2}}{3^{2}5}{3^2}{1\!\cdot\!3^{2}\!+\!1\!\cdot\!3^{3}\!+\!1\!\cdot\!3^{4}\!+\!O(3^{6})}{1\!\cdot\!3^{2}\!+\!2\!\cdot\!3^{3}\!+\!1\!\cdot\!3^{4}\!+\!O(3^{7})}\errentry{3\!\cdot\!10^{-10}}{2\!\cdot\!10^{-10}}\cr
\tabentry{15}{[1, 1, 1, 0, 1]}{x^3-15}{3^{5}5^{2}}{2^{2}3^{10}5^{4}19^{2}}{2^{2}5}{2^2}{2\!\cdot\!3^{0}\!+\!1\!\cdot\!3^{1}\!+\!2\!\cdot\!3^{2}\!+\!O(3^{4})}{2\!\cdot\!3^{0}\phantom{\,\,\!+\!2\!\cdot\!3^1}\,\phantom{\,\,\!+\!2\!\cdot\!3^2}\,\!+\!O(3^{3})}\errentry{1\!\cdot\!10^{-9}}{5\!\cdot\!10^{-10}}\cr
\tabentry{17}{[1, 1, 1, 0, 1]}{x^3-17}{3\!\cdot\!17^{2}}{2^{2}3^{2}17^{4}19^{2}}{0}{\nosha}{2\!\cdot\!3^{2}\!+\!1\!\cdot\!3^{3}\phantom{\,\,\!+\!2\!\cdot\!3^4}\,\!+\!O(3^{5})}{\noeuler}\errentry{2\!\cdot\!10^{-17}}{4\!\cdot\!10^{-10}}\cr
\tabentry{19}{[1, 1, 1, 0, 1]}{x^3-19}{3\!\cdot\!19^{2}}{2^{2}3^{2}19^{4}}{5}{1^2}{1\!\cdot\!3^{0}\!+\!2\!\cdot\!3^{1}\!+\!1\!\cdot\!3^{2}\!+\!O(3^{3})}{1\!\cdot\!3^{0}\phantom{\,\,\!+\!2\!\cdot\!3^1}\,\!+\!2\!\cdot\!3^{2}\!+\!O(3^{3})}\errentry{1\!\cdot\!10^{-17}}{1\!\cdot\!10^{-15}}\cr
\tabentry{35}{[1, 1, 1, 0, 1]}{x^3-35}{3\!\cdot\!5^{2}7^{2}}{2^{2}3^{2}5^{4}7^{4}19^{2}}{5}{1^2}{2\!\cdot\!3^{0}\phantom{\,\,\!+\!2\!\cdot\!3^1}\,\!+\!2\!\cdot\!3^{2}\!+\!O(3^{4})}{2\!\cdot\!3^{0}\!+\!2\!\cdot\!3^{1}\phantom{\,\,\!+\!2\!\cdot\!3^2}\,\!+\!O(3^{3})}\errentry{3\!\cdot\!10^{-13}}{9\!\cdot\!10^{-14}}\cr
\tabentry{37}{[1, 1, 1, 0, 1]}{x^3-37}{3\!\cdot\!37^{2}}{2^{2}3^{2}19^{2}37^{4}}{5}{1^2}{1\!\cdot\!3^{0}\phantom{\,\,\!+\!2\!\cdot\!3^1}\,\phantom{\,\,\!+\!2\!\cdot\!3^2}\,\!+\!O(3^{4})}{1\!\cdot\!3^{0}\phantom{\,\,\!+\!2\!\cdot\!3^1}\,\phantom{\,\,\!+\!2\!\cdot\!3^2}\,\!+\!O(3^{3})}\errentry{9\!\cdot\!10^{-12}}{2\!\cdot\!10^{-12}}\cr
\tabentry{45}{[1, 1, 1, 0, 1]}{x^3-45}{3^{5}5^{2}}{2^{2}3^{10}5^{4}19^{2}}{5}{1^2}{2\!\cdot\!3^{0}\!+\!1\!\cdot\!3^{1}\!+\!2\!\cdot\!3^{2}\!+\!O(3^{4})}{2\!\cdot\!3^{0}\!+\!1\!\cdot\!3^{1}\!+\!1\!\cdot\!3^{2}\!+\!O(3^{3})}\errentry{4\!\cdot\!10^{-10}}{1\!\cdot\!10^{-10}}\cr
\tabentry{53}{[1, 1, 1, 0, 1]}{x^3-53}{3\!\cdot\!53^{2}}{2^{2}3^{2}19^{2}53^{4}}{2^{2}5}{2^2}{2\!\cdot\!3^{0}\!+\!1\!\cdot\!3^{1}\phantom{\,\,\!+\!2\!\cdot\!3^2}\,\!+\!O(3^{4})}{2\!\cdot\!3^{0}\!+\!1\!\cdot\!3^{1}\!+\!1\!\cdot\!3^{2}\!+\!O(3^{3})}\errentry{1\!\cdot\!10^{-8}}{7\!\cdot\!10^{-9}}\cr
\tabentry{55}{[1, 1, 1, 0, 1]}{x^3-55}{3\!\cdot\!5^{2}11^{2}}{2^{2}3^{2}5^{4}11^{4}19^{2}}{5}{1^2}{1\!\cdot\!3^{0}\!+\!1\!\cdot\!3^{1}\!+\!1\!\cdot\!3^{2}\!+\!O(3^{3})}{1\!\cdot\!3^{0}\phantom{\,\,\!+\!2\!\cdot\!3^1}\,\!+\!1\!\cdot\!3^{2}\!+\!O(3^{3})}\errentry{1\!\cdot\!10^{-8}}{4\!\cdot\!10^{-9}}\cr
\tabentry{57}{[1, 1, 1, 0, 1]}{x^3-57}{3^{5}19^{2}}{2^{2}3^{10}19^{4}}{2^{2}5}{2^2}{1\!\cdot\!3^{0}\!+\!2\!\cdot\!3^{1}\!+\!1\!\cdot\!3^{2}\!+\!O(3^{3})}{1\!\cdot\!3^{0}\phantom{\,\,\!+\!2\!\cdot\!3^1}\,\!+\!2\!\cdot\!3^{2}\!+\!O(3^{3})}\errentry{4\!\cdot\!10^{-12}}{2\!\cdot\!10^{-12}}\cr
\tabentry{171}{[1, 1, 1, 0, 1]}{x^3-171}{3^{5}19^{2}}{2^{2}3^{10}19^{4}}{2^{2}5}{2^2}{1\!\cdot\!3^{0}\!+\!2\!\cdot\!3^{1}\!+\!1\!\cdot\!3^{2}\!+\!O(3^{3})}{1\!\cdot\!3^{0}\phantom{\,\,\!+\!2\!\cdot\!3^1}\,\!+\!2\!\cdot\!3^{2}\!+\!O(3^{3})}\errentry{1\!\cdot\!10^{-10}}{6\!\cdot\!10^{-11}}\cr
\hline
\tabend

\noindent\begin{tabular}{cllll} 
\extravspace
$E=56B1$ & \multicolumn{3}{l}{$y^2=x^3-x^2-4$} & \tabname{3}{56B1}\cr
\extravspace
$p=3$ & \multicolumn{3}{l}{$P_3(E/\Q,T)=1-2T+3T^2$}\cr
$/\vQ\alsub{}$   & $N(E)=2^{3}\!\cdot\!7$ & $L^*(E/\Q)=2^{-1}$ & $|E(\Q)|=2$ & $|\sha(E/\Q)|=1^2$ \cr
$/K\alsub{}$      & $N(E,\sigma)=2^{6}\!\cdot\!3^{2}\!\cdot\!7^{2}$ & $L^*(E/K)=2^{-1}$ & $|E(K)|=2$ & $|\sha(E/K)|=1^2$ \cr
$/\vQ\alsub{2}$ & $\hbox{III}^*,\>c_{2}=2$ & $/K\alsub{v}\quad(v|2)$ & $\hbox{III}^*,\>c_v=2$\cr
$/\vQ\alsub{7}$ & $\hbox{I}_1 \hbox{ split},\>c_{7}=1$ & $/K\alsub{v}\quad(v|7)$ & $\hbox{I}_1 \hbox{ split},\>c_v=1$\cr
\end{tabular}
\tabbegin
\tabentry{2}{[0, -1, 0, 0, -4]}{x^3-2}{2^{2}3^{3}}{2^{6}3^{6}7^{2}}{2^{2}}{1^2}{1\!\cdot\!3^{0}\!+\!2\!\cdot\!3^{1}\!+\!1\!\cdot\!3^{2}\!+\!O(3^{3})}{1\!\cdot\!3^{0}\!+\!1\!\cdot\!3^{1}\!+\!2\!\cdot\!3^{2}\!+\!O(3^{3})}\errentry{2\!\cdot\!10^{-17}}{9\!\cdot\!10^{-16}}\cr
\tabentry{6}{[0, -1, 0, 0, -4]}{x^3-6}{2^{2}3^{5}}{2^{6}3^{10}7^{2}}{2^{4}}{2^2}{1\!\cdot\!3^{0}\!+\!2\!\cdot\!3^{1}\!+\!1\!\cdot\!3^{2}\!+\!O(3^{3})}{1\!\cdot\!3^{0}\phantom{\,\,\!+\!2\!\cdot\!3^1}\,\!+\!2\!\cdot\!3^{2}\!+\!O(3^{3})}\errentry{2\!\cdot\!10^{-17}}{1\!\cdot\!10^{-15}}\cr
\tabentry{10}{[0, -1, 0, 0, -4]}{x^3-10}{2^{2}3\!\cdot\!5^{2}}{2^{6}3^{2}5^{4}7^{2}}{2^{2}}{1^2}{2\!\cdot\!3^{0}\phantom{\,\,\!+\!2\!\cdot\!3^1}\,\phantom{\,\,\!+\!2\!\cdot\!3^2}\,\!+\!O(3^{4})}{2\!\cdot\!3^{0}\phantom{\,\,\!+\!2\!\cdot\!3^1}\,\phantom{\,\,\!+\!2\!\cdot\!3^2}\,\!+\!O(3^{4})}\errentry{3\!\cdot\!10^{-17}}{2\!\cdot\!10^{-16}}\cr
\tabentry{12}{[0, -1, 0, 0, -4]}{x^3-12}{2^{2}3^{5}}{2^{6}3^{10}7^{2}}{2^{2}}{1^2}{1\!\cdot\!3^{0}\!+\!2\!\cdot\!3^{1}\!+\!1\!\cdot\!3^{2}\!+\!O(3^{3})}{1\!\cdot\!3^{0}\!+\!2\!\cdot\!3^{1}\!+\!2\!\cdot\!3^{2}\!+\!O(3^{3})}\errentry{1\!\cdot\!10^{-16}}{3\!\cdot\!10^{-16}}\cr
\tabentry{14}{[0, -1, 0, 0, -4]}{x^3-14}{2^{2}3^{3}7^{2}}{2^{6}3^{6}7^{4}}{0}{\nosha}{1\!\cdot\!3^{2}\!+\!2\!\cdot\!3^{3}\!+\!2\!\cdot\!3^{4}\!+\!O(3^{6})}{\noeuler}\errentry{3\!\cdot\!10^{-18}}{9\!\cdot\!10^{-10}}\cr
\tabentry{20}{[0, -1, 0, 0, -4]}{x^3-20}{2^{2}3^{3}5^{2}}{2^{6}3^{6}5^{4}7^{2}}{2^{4}}{2^2}{2\!\cdot\!3^{0}\phantom{\,\,\!+\!2\!\cdot\!3^1}\,\phantom{\,\,\!+\!2\!\cdot\!3^2}\,\!+\!O(3^{4})}{2\!\cdot\!3^{0}\!+\!1\!\cdot\!3^{1}\!+\!2\!\cdot\!3^{2}\!+\!O(3^{3})}\errentry{4\!\cdot\!10^{-17}}{3\!\cdot\!10^{-15}}\cr
\tabentry{22}{[0, -1, 0, 0, -4]}{x^3-22}{2^{2}3^{3}11^{2}}{2^{6}3^{6}7^{2}11^{4}}{0}{\nosha}{1\!\cdot\!3^{2}\phantom{\,\,\!+\!2\!\cdot\!3^3}\,\phantom{\,\,\!+\!2\!\cdot\!3^4}\,\!+\!O(3^{6})}{\noeuler}\errentry{1\!\cdot\!10^{-10}}{1\!\cdot\!10^{-6}}\cr
\tabentry{26}{[0, -1, 0, 0, -4]}{x^3-26}{2^{2}3\!\cdot\!13^{2}}{2^{6}3^{2}7^{2}13^{4}}{2^{2}}{1^2}{1\!\cdot\!3^{0}\!+\!2\!\cdot\!3^{1}\!+\!1\!\cdot\!3^{2}\!+\!O(3^{3})}{1\!\cdot\!3^{0}\!+\!2\!\cdot\!3^{1}\phantom{\,\,\!+\!2\!\cdot\!3^2}\,\!+\!O(3^{4})}\errentry{4\!\cdot\!10^{-17}}{4\!\cdot\!10^{-17}}\cr
\tabentry{28}{[0, -1, 0, 0, -4]}{x^3-28}{2^{2}3\!\cdot\!7^{2}}{2^{6}3^{2}7^{4}}{2^{2}3}{1^2}{1\!\cdot\!3^{2}\!+\!2\!\cdot\!3^{3}\!+\!2\!\cdot\!3^{4}\!+\!O(3^{6})}{1\!\cdot\!3^{1}\!+\!1\!\cdot\!3^{2}\!+\!1\!\cdot\!3^{3}\!+\!O(3^{5})}\errentry{5\!\cdot\!10^{-17}}{1\!\cdot\!10^{-15}}\cr
\tabentry{30}{[0, -1, 0, 0, -4]}{x^3-30}{2^{2}3^{5}5^{2}}{2^{6}3^{10}5^{4}7^{2}}{2^{2}}{1^2}{2\!\cdot\!3^{0}\phantom{\,\,\!+\!2\!\cdot\!3^1}\,\phantom{\,\,\!+\!2\!\cdot\!3^2}\,\!+\!O(3^{4})}{2\!\cdot\!3^{0}\!+\!1\!\cdot\!3^{1}\phantom{\,\,\!+\!2\!\cdot\!3^2}\,\!+\!O(3^{3})}\errentry{4\!\cdot\!10^{-7}}{1\!\cdot\!10^{-7}}\cr
\tabentry{42}{[0, -1, 0, 0, -4]}{x^3-42}{2^{2}3^{5}7^{2}}{2^{6}3^{10}7^{4}}{2^{4}3}{2^2}{1\!\cdot\!3^{2}\!+\!2\!\cdot\!3^{3}\!+\!2\!\cdot\!3^{4}\!+\!O(3^{6})}{1\!\cdot\!3^{1}\!+\!1\!\cdot\!3^{2}\!+\!2\!\cdot\!3^{3}\!+\!O(3^{5})}\errentry{1\!\cdot\!10^{-15}}{6\!\cdot\!10^{-15}}\cr
\tabentry{44}{[0, -1, 0, 0, -4]}{x^3-44}{2^{2}3\!\cdot\!11^{2}}{2^{6}3^{2}7^{2}11^{4}}{2^{2}3^{2}}{3^2}{1\!\cdot\!3^{2}\phantom{\,\,\!+\!2\!\cdot\!3^3}\,\phantom{\,\,\!+\!2\!\cdot\!3^4}\,\!+\!O(3^{6})}{2\!\cdot\!3^{2}\!+\!1\!\cdot\!3^{3}\!+\!2\!\cdot\!3^{4}\!+\!O(3^{6})}\errentry{1\!\cdot\!10^{-17}}{2\!\cdot\!10^{-15}}\cr
\tabentry{46}{[0, -1, 0, 0, -4]}{x^3-46}{2^{2}3\!\cdot\!23^{2}}{2^{6}3^{2}7^{2}23^{4}}{2^{4}}{2^2}{2\!\cdot\!3^{0}\!+\!2\!\cdot\!3^{1}\!+\!2\!\cdot\!3^{2}\!+\!O(3^{5})}{2\!\cdot\!3^{0}\!+\!2\!\cdot\!3^{1}\!+\!1\!\cdot\!3^{2}\!+\!O(3^{3})}\errentry{1\!\cdot\!10^{-15}}{1\!\cdot\!10^{-15}}\cr
\tabentry{52}{[0, -1, 0, 0, -4]}{x^3-52}{2^{2}3^{3}13^{2}}{2^{6}3^{6}7^{2}13^{4}}{2^{4}}{2^2}{1\!\cdot\!3^{0}\!+\!2\!\cdot\!3^{1}\!+\!1\!\cdot\!3^{2}\!+\!O(3^{3})}{1\!\cdot\!3^{0}\!+\!1\!\cdot\!3^{1}\!+\!1\!\cdot\!3^{2}\!+\!O(3^{3})}\errentry{7\!\cdot\!10^{-10}}{3\!\cdot\!10^{-10}}\cr
\tabentry{60}{[0, -1, 0, 0, -4]}{x^3-60}{2^{2}3^{5}5^{2}}{2^{6}3^{10}5^{4}7^{2}}{2^{4}}{2^2}{2\!\cdot\!3^{0}\phantom{\,\,\!+\!2\!\cdot\!3^1}\,\phantom{\,\,\!+\!2\!\cdot\!3^2}\,\!+\!O(3^{4})}{2\!\cdot\!3^{0}\phantom{\,\,\!+\!2\!\cdot\!3^1}\,\!+\!2\!\cdot\!3^{2}\!+\!O(3^{3})}\errentry{8\!\cdot\!10^{-8}}{4\!\cdot\!10^{-8}}\cr
\tabentry{62}{[0, -1, 0, 0, -4]}{x^3-62}{2^{2}3\!\cdot\!31^{2}}{2^{6}3^{2}7^{2}31^{4}}{2^{6}}{4^2}{1\!\cdot\!3^{0}\phantom{\,\,\!+\!2\!\cdot\!3^1}\,\phantom{\,\,\!+\!2\!\cdot\!3^2}\,\!+\!O(3^{3})}{1\!\cdot\!3^{0}\!+\!1\!\cdot\!3^{1}\!+\!1\!\cdot\!3^{2}\!+\!O(3^{4})}\errentry{1\!\cdot\!10^{-11}}{1\!\cdot\!10^{-11}}\cr
\tabentry{70}{[0, -1, 0, 0, -4]}{x^3-70}{2^{2}3^{3}5^{2}7^{2}}{2^{6}3^{6}5^{4}7^{4}}{2^{4}3}{2^2}{2\!\cdot\!3^{2}\phantom{\,\,\!+\!2\!\cdot\!3^3}\,\!+\!2\!\cdot\!3^{4}\!+\!O(3^{6})}{2\!\cdot\!3^{1}\phantom{\,\,\!+\!2\!\cdot\!3^2}\,\!+\!1\!\cdot\!3^{3}\!+\!O(3^{6})}\errentry{7\!\cdot\!10^{-10}}{3\!\cdot\!10^{-10}}\cr
\tabentry{82}{[0, -1, 0, 0, -4]}{x^3-82}{2^{2}3\!\cdot\!41^{2}}{2^{6}3^{2}7^{2}41^{4}}{2^{2}5^{2}}{5^2}{2\!\cdot\!3^{0}\!+\!1\!\cdot\!3^{1}\!+\!2\!\cdot\!3^{2}\!+\!O(3^{4})}{2\!\cdot\!3^{0}\!+\!1\!\cdot\!3^{1}\!+\!1\!\cdot\!3^{2}\!+\!O(3^{4})}\errentry{1\!\cdot\!10^{-9}}{2\!\cdot\!10^{-9}}\cr
\tabentry{84}{[0, -1, 0, 0, -4]}{x^3-84}{2^{2}3^{5}7^{2}}{2^{6}3^{10}7^{4}}{2^{2}3^{3}}{3^2}{1\!\cdot\!3^{2}\!+\!2\!\cdot\!3^{3}\!+\!2\!\cdot\!3^{4}\!+\!O(3^{6})}{1\!\cdot\!3^{3}\phantom{\,\,\!+\!2\!\cdot\!3^4}\,\!+\!2\!\cdot\!3^{5}\!+\!O(3^{6})}\errentry{4\!\cdot\!10^{-16}}{4\!\cdot\!10^{-15}}\cr
\tabentry{90}{[0, -1, 0, 0, -4]}{x^3-90}{2^{2}3^{5}5^{2}}{2^{6}3^{10}5^{4}7^{2}}{2^{6}}{4^2}{2\!\cdot\!3^{0}\phantom{\,\,\!+\!2\!\cdot\!3^1}\,\phantom{\,\,\!+\!2\!\cdot\!3^2}\,\!+\!O(3^{4})}{2\!\cdot\!3^{0}\!+\!2\!\cdot\!3^{1}\!+\!2\!\cdot\!3^{2}\!+\!O(3^{4})}\errentry{9\!\cdot\!10^{-7}}{9\!\cdot\!10^{-7}}\cr
\tabentry{116}{[0, -1, 0, 0, -4]}{x^3-116}{2^{2}3\!\cdot\!29^{2}}{2^{6}3^{2}7^{2}29^{4}}{2^{2}}{1^2}{2\!\cdot\!3^{0}\phantom{\,\,\!+\!2\!\cdot\!3^1}\,\!+\!2\!\cdot\!3^{2}\!+\!O(3^{4})}{2\!\cdot\!3^{0}\!+\!1\!\cdot\!3^{1}\phantom{\,\,\!+\!2\!\cdot\!3^2}\,\!+\!O(3^{3})}\errentry{4\!\cdot\!10^{-13}}{1\!\cdot\!10^{-13}}\cr
\tabentry{126}{[0, -1, 0, 0, -4]}{x^3-126}{2^{2}3^{5}7^{2}}{2^{6}3^{10}7^{4}}{2^{2}3}{1^2}{1\!\cdot\!3^{2}\!+\!2\!\cdot\!3^{3}\!+\!2\!\cdot\!3^{4}\!+\!O(3^{6})}{1\!\cdot\!3^{1}\phantom{\,\,\!+\!2\!\cdot\!3^2}\,\!+\!2\!\cdot\!3^{3}\!+\!O(3^{4})}\errentry{1\!\cdot\!10^{-16}}{1\!\cdot\!10^{-16}}\cr
\tabentry{140}{[0, -1, 0, 0, -4]}{x^3-140}{2^{2}3^{3}5^{2}7^{2}}{2^{6}3^{6}5^{4}7^{4}}{2^{4}3}{2^2}{2\!\cdot\!3^{2}\phantom{\,\,\!+\!2\!\cdot\!3^3}\,\!+\!2\!\cdot\!3^{4}\!+\!O(3^{6})}{2\!\cdot\!3^{1}\phantom{\,\,\!+\!2\!\cdot\!3^2}\,\!+\!1\!\cdot\!3^{3}\!+\!O(3^{6})}\errentry{1\!\cdot\!10^{-8}}{9\!\cdot\!10^{-9}}\cr
\tabentry{150}{[0, -1, 0, 0, -4]}{x^3-150}{2^{2}3^{5}5^{2}}{2^{6}3^{10}5^{4}7^{2}}{2^{2}5^{2}}{5^2}{2\!\cdot\!3^{0}\phantom{\,\,\!+\!2\!\cdot\!3^1}\,\phantom{\,\,\!+\!2\!\cdot\!3^2}\,\!+\!O(3^{4})}{2\!\cdot\!3^{0}\!+\!2\!\cdot\!3^{1}\!+\!1\!\cdot\!3^{2}\!+\!O(3^{3})}\errentry{2\!\cdot\!10^{-7}}{2\!\cdot\!10^{-7}}\cr
\tabentry{154}{[0, -1, 0, 0, -4]}{x^3-154}{2^{2}3\!\cdot\!7^{2}11^{2}}{2^{6}3^{2}7^{4}11^{4}}{2^{2}3^{3}}{3^2}{1\!\cdot\!3^{4}\phantom{\,\,\!+\!2\!\cdot\!3^5}\,\!+\!1\!\cdot\!3^{6}\!+\!O(3^{8})}{2\!\cdot\!3^{3}\phantom{\,\,\!+\!2\!\cdot\!3^4}\,\!+\!1\!\cdot\!3^{5}\!+\!O(3^{6})}\errentry{1\!\cdot\!10^{-13}}{9\!\cdot\!10^{-14}}\cr
\tabentry{172}{[0, -1, 0, 0, -4]}{x^3-172}{2^{2}3\!\cdot\!43^{2}}{2^{6}3^{2}7^{2}43^{4}}{2^{2}3^{2}}{3^2}{1\!\cdot\!3^{4}\phantom{\,\,\!+\!2\!\cdot\!3^5}\,\phantom{\,\,\!+\!2\!\cdot\!3^6}\,\!+\!O(3^{7})}{1\!\cdot\!3^{2}\!+\!1\!\cdot\!3^{3}\phantom{\,\,\!+\!2\!\cdot\!3^4}\,\!+\!O(3^{7})}\errentry{1\!\cdot\!10^{-8}}{1\!\cdot\!10^{-8}}\cr
\tabentry{252}{[0, -1, 0, 0, -4]}{x^3-252}{2^{2}3^{5}7^{2}}{2^{6}3^{10}7^{4}}{2^{6}3}{4^2}{1\!\cdot\!3^{2}\!+\!2\!\cdot\!3^{3}\!+\!2\!\cdot\!3^{4}\!+\!O(3^{6})}{1\!\cdot\!3^{1}\!+\!2\!\cdot\!3^{2}\phantom{\,\,\!+\!2\!\cdot\!3^3}\,\!+\!O(3^{5})}\errentry{2\!\cdot\!10^{-15}}{4\!\cdot\!10^{-15}}\cr
\tabentry{350}{[0, -1, 0, 0, -4]}{x^3-350}{2^{2}3\!\cdot\!5^{2}7^{2}}{2^{6}3^{2}5^{4}7^{4}}{2^{2}3}{1^2}{2\!\cdot\!3^{2}\phantom{\,\,\!+\!2\!\cdot\!3^3}\,\!+\!2\!\cdot\!3^{4}\!+\!O(3^{6})}{2\!\cdot\!3^{1}\!+\!2\!\cdot\!3^{2}\phantom{\,\,\!+\!2\!\cdot\!3^3}\,\!+\!O(3^{4})}\errentry{4\!\cdot\!10^{-17}}{6\!\cdot\!10^{-16}}\cr
\tabentry{476}{[0, -1, 0, 0, -4]}{x^3-476}{2^{2}3\!\cdot\!7^{2}17^{2}}{2^{6}3^{2}7^{4}17^{4}}{2^{2}3^{3}}{3^2}{2\!\cdot\!3^{2}\!+\!2\!\cdot\!3^{3}\!+\!2\!\cdot\!3^{4}\!+\!O(3^{6})}{2\!\cdot\!3^{3}\!+\!1\!\cdot\!3^{4}\!+\!2\!\cdot\!3^{5}\!+\!O(3^{6})}\errentry{5\!\cdot\!10^{-8}}{4\!\cdot\!10^{-8}}\cr
\tabentry{490}{[0, -1, 0, 0, -4]}{x^3-490}{2^{2}3^{3}5^{2}7^{2}}{2^{6}3^{6}5^{4}7^{4}}{2^{2}3^{3}}{3^2}{2\!\cdot\!3^{2}\phantom{\,\,\!+\!2\!\cdot\!3^3}\,\!+\!2\!\cdot\!3^{4}\!+\!O(3^{6})}{2\!\cdot\!3^{3}\!+\!1\!\cdot\!3^{4}\!+\!2\!\cdot\!3^{5}\!+\!O(3^{7})}\errentry{2\!\cdot\!10^{-8}}{1\!\cdot\!10^{-8}}\cr
\tabentry{1666}{[0, -1, 0, 0, -4]}{x^3-1666}{2^{2}3\!\cdot\!7^{2}17^{2}}{2^{6}3^{2}7^{4}17^{4}}{2^{2}3}{1^2}{2\!\cdot\!3^{2}\!+\!2\!\cdot\!3^{3}\!+\!2\!\cdot\!3^{4}\!+\!O(3^{6})}{2\!\cdot\!3^{1}\!+\!1\!\cdot\!3^{2}\!+\!2\!\cdot\!3^{3}\!+\!O(3^{4})}\errentry{4\!\cdot\!10^{-8}}{1\!\cdot\!10^{-8}}\cr
\tabentry{2366}{[0, -1, 0, 0, -4]}{x^3-2366}{2^{2}3\!\cdot\!7^{2}13^{2}}{2^{6}3^{2}7^{4}13^{4}}{2^{2}3^{3}}{3^2}{1\!\cdot\!3^{2}\!+\!2\!\cdot\!3^{3}\!+\!2\!\cdot\!3^{4}\!+\!O(3^{5})}{1\!\cdot\!3^{3}\phantom{\,\,\!+\!2\!\cdot\!3^4}\,\phantom{\,\,\!+\!2\!\cdot\!3^5}\,\!+\!O(3^{7})}\errentry{2\!\cdot\!10^{-10}}{2\!\cdot\!10^{-10}}\cr
\hline
\tabend

\noindent\begin{tabular}{cllll} 
\extravspace
$E=128B2$ & \multicolumn{3}{l}{$y^2=x^3+x^2-2x-2$} & \tabname{3}{128B2}\cr
\extravspace
$p=3$ & \multicolumn{3}{l}{$P_3(E/\Q,T)=1+2T+3T^2$}\cr
$/\vQ\alsub{}$   & $N(E)=2^{7}$ & $L^*(E/\Q)=2^{-2}$ & $|E(\Q)|=2$ & $|\sha(E/\Q)|=1^2$ \cr
$/K\alsub{}$      & $N(E,\sigma)=2^{14}\!\cdot\!3^{2}$ & $L^*(E/K)=2^{-2}$ & $|E(K)|=2$ & $|\sha(E/K)|=1^2$ \cr
$/\vQ\alsub{2}$ & $\hbox{II},\>c_{2}=1$ & $/K\alsub{v}\quad(v|2)$ & $\hbox{II},\>c_v=1$\cr
\end{tabular}
\tabbegin
\tabentry{2}{[0, 1, 0, -2, -2]}{x^3-2}{2^{2}3^{3}}{2^{14}3^{6}}{0}{\nosha}{1\!\cdot\!3^{2}\!+\!1\!\cdot\!3^{3}\!+\!1\!\cdot\!3^{4}\!+\!O(3^{5})}{\noeuler}\errentry{1\!\cdot\!10^{-16}}{1\!\cdot\!10^{-9}}\cr
\tabentry{6}{[0, 1, 0, -2, -2]}{x^3-6}{2^{2}3^{5}}{2^{14}3^{10}}{2\!\cdot\!3^{2}}{3^2}{1\!\cdot\!3^{2}\!+\!1\!\cdot\!3^{3}\!+\!1\!\cdot\!3^{4}\!+\!O(3^{5})}{1\!\cdot\!3^{2}\phantom{\,\,\!+\!2\!\cdot\!3^3}\,\phantom{\,\,\!+\!2\!\cdot\!3^4}\,\!+\!O(3^{7})}\errentry{6\!\cdot\!10^{-17}}{5\!\cdot\!10^{-17}}\cr
\tabentry{10}{[0, 1, 0, -2, -2]}{x^3-10}{2^{2}3\!\cdot\!5^{2}}{2^{14}3^{2}5^{4}}{2}{1^2}{2\!\cdot\!3^{2}\!+\!2\!\cdot\!3^{3}\phantom{\,\,\!+\!2\!\cdot\!3^4}\,\!+\!O(3^{5})}{2\!\cdot\!3^{2}\!+\!1\!\cdot\!3^{3}\phantom{\,\,\!+\!2\!\cdot\!3^4}\,\!+\!O(3^{5})}\errentry{3\!\cdot\!10^{-16}}{2\!\cdot\!10^{-15}}\cr
\tabentry{12}{[0, 1, 0, -2, -2]}{x^3-12}{2^{2}3^{5}}{2^{14}3^{10}}{0}{\nosha}{1\!\cdot\!3^{2}\!+\!1\!\cdot\!3^{3}\!+\!1\!\cdot\!3^{4}\!+\!O(3^{5})}{\noeuler}\errentry{8\!\cdot\!10^{-17}}{2\!\cdot\!10^{-9}}\cr
\tabentry{14}{[0, 1, 0, -2, -2]}{x^3-14}{2^{2}3^{3}7^{2}}{2^{14}3^{6}7^{4}}{2\!\cdot\!3^{2}}{3^2}{1\!\cdot\!3^{2}\!+\!2\!\cdot\!3^{3}\phantom{\,\,\!+\!2\!\cdot\!3^4}\,\!+\!O(3^{5})}{1\!\cdot\!3^{2}\phantom{\,\,\!+\!2\!\cdot\!3^3}\,\!+\!1\!\cdot\!3^{4}\!+\!O(3^{5})}\errentry{3\!\cdot\!10^{-11}}{2\!\cdot\!10^{-11}}\cr
\tabentry{20}{[0, 1, 0, -2, -2]}{x^3-20}{2^{2}3^{3}5^{2}}{2^{14}3^{6}5^{4}}{0}{\nosha}{2\!\cdot\!3^{2}\!+\!2\!\cdot\!3^{3}\phantom{\,\,\!+\!2\!\cdot\!3^4}\,\!+\!O(3^{5})}{\noeuler}\errentry{1\!\cdot\!10^{-17}}{2\!\cdot\!10^{-9}}\cr
\tabentry{26}{[0, 1, 0, -2, -2]}{x^3-26}{2^{2}3\!\cdot\!13^{2}}{2^{14}3^{2}13^{4}}{0}{\nosha}{1\!\cdot\!3^{4}\!+\!1\!\cdot\!3^{5}\!+\!2\!\cdot\!3^{6}\!+\!O(3^{8})}{\noeuler}\errentry{5\!\cdot\!10^{-15}}{4\!\cdot\!10^{-9}}\cr
\tabentry{28}{[0, 1, 0, -2, -2]}{x^3-28}{2^{2}3\!\cdot\!7^{2}}{2^{14}3^{2}7^{4}}{2}{1^2}{1\!\cdot\!3^{2}\!+\!2\!\cdot\!3^{3}\phantom{\,\,\!+\!2\!\cdot\!3^4}\,\!+\!O(3^{5})}{1\!\cdot\!3^{2}\phantom{\,\,\!+\!2\!\cdot\!3^3}\,\!+\!2\!\cdot\!3^{4}\!+\!O(3^{5})}\errentry{5\!\cdot\!10^{-17}}{1\!\cdot\!10^{-17}}\cr
\tabentry{44}{[0, 1, 0, -2, -2]}{x^3-44}{2^{2}3\!\cdot\!11^{2}}{2^{14}3^{2}11^{4}}{2}{1^2}{2\!\cdot\!3^{2}\!+\!1\!\cdot\!3^{3}\phantom{\,\,\!+\!2\!\cdot\!3^4}\,\!+\!O(3^{5})}{2\!\cdot\!3^{2}\!+\!2\!\cdot\!3^{3}\!+\!1\!\cdot\!3^{4}\!+\!O(3^{5})}\errentry{2\!\cdot\!10^{-16}}{1\!\cdot\!10^{-16}}\cr
\tabentry{46}{[0, 1, 0, -2, -2]}{x^3-46}{2^{2}3\!\cdot\!23^{2}}{2^{14}3^{2}23^{4}}{2}{1^2}{2\!\cdot\!3^{2}\!+\!1\!\cdot\!3^{3}\!+\!2\!\cdot\!3^{4}\!+\!O(3^{5})}{2\!\cdot\!3^{2}\!+\!1\!\cdot\!3^{3}\!+\!1\!\cdot\!3^{4}\!+\!O(3^{6})}\errentry{3\!\cdot\!10^{-10}}{7\!\cdot\!10^{-11}}\cr
\tabentry{116}{[0, 1, 0, -2, -2]}{x^3-116}{2^{2}3\!\cdot\!29^{2}}{2^{14}3^{2}29^{4}}{0}{\nosha}{2\!\cdot\!3^{5}\!+\!1\!\cdot\!3^{6}\!+\!2\!\cdot\!3^{7}\!+\!O(3^{9})}{\noeuler}\errentry{1\!\cdot\!10^{-7}}{4\!\cdot\!10^{-5}}\cr
\hline
\tabend

\noindent\begin{tabular}{cllll} 
\extravspace
$E=77C1$ & \multicolumn{3}{l}{$y^2+xy=x^3+x^2+4x+11$} & \tabname{3}{77C1}\cr
\extravspace
$p=3$ & \multicolumn{3}{l}{$P_3(E/\Q,T)=1-2T+3T^2$}\cr
$/\vQ\alsub{}$   & $N(E)=7\!\cdot\!11$ & $L^*(E/\Q)=2^{-1}$ & $|E(\Q)|=2$ & $|\sha(E/\Q)|=1^2$ \cr
$/K\alsub{}$      & $N(E,\sigma)=3^{2}\!\cdot\!7^{2}\!\cdot\!11^{2}$ & $L^*(E/K)=2^{-1}$ & $|E(K)|=2$ & $|\sha(E/K)|=1^2$ \cr
$/\vQ\alsub{7}$ & $\hbox{I}_3 \hbox{ non-split},\>c_{7}=1$ & $/K\alsub{v}\quad(v|7)$ & $\hbox{I}_3 \hbox{ non-split},\>c_v=1$\cr
$/\vQ\alsub{11}$ & $\hbox{I}_2 \hbox{ split},\>c_{11}=2$ & $/K\alsub{v}\quad(v|11)$ & $\hbox{I}_2 \hbox{ split},\>c_v=2$\cr
\end{tabular}
\tabbegin
\tabentry{2}{[1, 1, 0, 4, 11]}{x^3-2}{2^{2}3^{3}}{2^{4}3^{6}7^{2}11^{2}}{2^{3}}{2^2}{2\!\cdot\!3^{0}\!+\!1\!\cdot\!3^{1}\phantom{\,\,\!+\!2\!\cdot\!3^2}\,\!+\!O(3^{4})}{2\!\cdot\!3^{0}\!+\!2\!\cdot\!3^{1}\!+\!1\!\cdot\!3^{2}\!+\!O(3^{3})}\errentry{5\!\cdot\!10^{-17}}{1\!\cdot\!10^{-15}}\cr
\tabentry{3}{[1, 1, 0, 4, 11]}{x^3-3}{3^{5}}{3^{10}7^{2}11^{2}}{2}{1^2}{1\!\cdot\!3^{0}\!+\!2\!\cdot\!3^{1}\!+\!1\!\cdot\!3^{2}\!+\!O(3^{3})}{1\!\cdot\!3^{0}\!+\!2\!\cdot\!3^{1}\!+\!2\!\cdot\!3^{2}\!+\!O(3^{3})}\errentry{3\!\cdot\!10^{-16}}{3\!\cdot\!10^{-16}}\cr
\tabentry{5}{[1, 1, 0, 4, 11]}{x^3-5}{3^{3}5^{2}}{3^{6}5^{4}7^{2}11^{2}}{2^{3}}{2^2}{2\!\cdot\!3^{0}\!+\!1\!\cdot\!3^{1}\!+\!1\!\cdot\!3^{2}\!+\!O(3^{3})}{2\!\cdot\!3^{0}\!+\!1\!\cdot\!3^{1}\!+\!2\!\cdot\!3^{2}\!+\!O(3^{3})}\errentry{6\!\cdot\!10^{-16}}{5\!\cdot\!10^{-15}}\cr
\tabentry{6}{[1, 1, 0, 4, 11]}{x^3-6}{2^{2}3^{5}}{2^{4}3^{10}7^{2}11^{2}}{2}{1^2}{2\!\cdot\!3^{0}\!+\!1\!\cdot\!3^{1}\phantom{\,\,\!+\!2\!\cdot\!3^2}\,\!+\!O(3^{4})}{2\!\cdot\!3^{0}\!+\!2\!\cdot\!3^{1}\!+\!2\!\cdot\!3^{2}\!+\!O(3^{3})}\errentry{9\!\cdot\!10^{-17}}{3\!\cdot\!10^{-16}}\cr
\tabentry{7}{[1, 1, 0, 4, 11]}{x^3-7}{3^{3}7^{2}}{3^{6}7^{4}11^{2}}{2}{1^2}{1\!\cdot\!3^{0}\!+\!1\!\cdot\!3^{1}\!+\!2\!\cdot\!3^{2}\!+\!O(3^{3})}{1\!\cdot\!3^{0}\!+\!2\!\cdot\!3^{1}\phantom{\,\,\!+\!2\!\cdot\!3^2}\,\!+\!O(3^{3})}\errentry{6\!\cdot\!10^{-18}}{1\!\cdot\!10^{-16}}\cr
\tabentry{10}{[1, 1, 0, 4, 11]}{x^3-10}{2^{2}3\!\cdot\!5^{2}}{2^{4}3^{2}5^{4}7^{2}11^{2}}{2}{1^2}{1\!\cdot\!3^{0}\phantom{\,\,\!+\!2\!\cdot\!3^1}\,\phantom{\,\,\!+\!2\!\cdot\!3^2}\,\!+\!O(3^{3})}{1\!\cdot\!3^{0}\phantom{\,\,\!+\!2\!\cdot\!3^1}\,\phantom{\,\,\!+\!2\!\cdot\!3^2}\,\!+\!O(3^{4})}\errentry{5\!\cdot\!10^{-17}}{1\!\cdot\!10^{-17}}\cr
\tabentry{12}{[1, 1, 0, 4, 11]}{x^3-12}{2^{2}3^{5}}{2^{4}3^{10}7^{2}11^{2}}{2}{1^2}{2\!\cdot\!3^{0}\!+\!1\!\cdot\!3^{1}\phantom{\,\,\!+\!2\!\cdot\!3^2}\,\!+\!O(3^{4})}{2\!\cdot\!3^{0}\!+\!2\!\cdot\!3^{1}\!+\!2\!\cdot\!3^{2}\!+\!O(3^{3})}\errentry{4\!\cdot\!10^{-18}}{3\!\cdot\!10^{-16}}\cr
\tabentry{13}{[1, 1, 0, 4, 11]}{x^3-13}{3^{3}13^{2}}{3^{6}7^{2}11^{2}13^{4}}{2}{1^2}{1\!\cdot\!3^{0}\phantom{\,\,\!+\!2\!\cdot\!3^1}\,\phantom{\,\,\!+\!2\!\cdot\!3^2}\,\!+\!O(3^{5})}{1\!\cdot\!3^{0}\phantom{\,\,\!+\!2\!\cdot\!3^1}\,\!+\!1\!\cdot\!3^{2}\!+\!O(3^{3})}\errentry{7\!\cdot\!10^{-8}}{2\!\cdot\!10^{-8}}\cr
\tabentry{14}{[1, 1, 0, 4, 11]}{x^3-14}{2^{2}3^{3}7^{2}}{2^{4}3^{6}7^{4}11^{2}}{2\!\cdot\!5^{2}}{5^2}{2\!\cdot\!3^{0}\!+\!2\!\cdot\!3^{1}\!+\!1\!\cdot\!3^{2}\!+\!O(3^{3})}{2\!\cdot\!3^{0}\phantom{\,\,\!+\!2\!\cdot\!3^1}\,\!+\!2\!\cdot\!3^{2}\!+\!O(3^{5})}\errentry{1\!\cdot\!10^{-16}}{6\!\cdot\!10^{-15}}\cr
\tabentry{17}{[1, 1, 0, 4, 11]}{x^3-17}{3\!\cdot\!17^{2}}{3^{2}7^{2}11^{2}17^{4}}{2^{3}}{2^2}{2\!\cdot\!3^{0}\!+\!1\!\cdot\!3^{1}\phantom{\,\,\!+\!2\!\cdot\!3^2}\,\!+\!O(3^{3})}{2\!\cdot\!3^{0}\!+\!1\!\cdot\!3^{1}\!+\!2\!\cdot\!3^{2}\!+\!O(3^{3})}\errentry{3\!\cdot\!10^{-17}}{6\!\cdot\!10^{-15}}\cr
\tabentry{19}{[1, 1, 0, 4, 11]}{x^3-19}{3\!\cdot\!19^{2}}{3^{2}7^{2}11^{2}19^{4}}{2}{1^2}{1\!\cdot\!3^{0}\phantom{\,\,\!+\!2\!\cdot\!3^1}\,\!+\!2\!\cdot\!3^{2}\!+\!O(3^{3})}{1\!\cdot\!3^{0}\phantom{\,\,\!+\!2\!\cdot\!3^1}\,\!+\!2\!\cdot\!3^{2}\!+\!O(3^{3})}\errentry{4\!\cdot\!10^{-16}}{4\!\cdot\!10^{-16}}\cr
\tabentry{20}{[1, 1, 0, 4, 11]}{x^3-20}{2^{2}3^{3}5^{2}}{2^{4}3^{6}5^{4}7^{2}11^{2}}{2}{1^2}{1\!\cdot\!3^{0}\phantom{\,\,\!+\!2\!\cdot\!3^1}\,\phantom{\,\,\!+\!2\!\cdot\!3^2}\,\!+\!O(3^{3})}{1\!\cdot\!3^{0}\!+\!1\!\cdot\!3^{1}\!+\!1\!\cdot\!3^{2}\!+\!O(3^{3})}\errentry{3\!\cdot\!10^{-10}}{8\!\cdot\!10^{-11}}\cr
\tabentry{21}{[1, 1, 0, 4, 11]}{x^3-21}{3^{5}7^{2}}{3^{10}7^{4}11^{2}}{2\!\cdot\!5^{2}}{5^2}{1\!\cdot\!3^{0}\!+\!1\!\cdot\!3^{1}\!+\!2\!\cdot\!3^{2}\!+\!O(3^{3})}{1\!\cdot\!3^{0}\!+\!2\!\cdot\!3^{1}\!+\!1\!\cdot\!3^{2}\!+\!O(3^{4})}\errentry{3\!\cdot\!10^{-12}}{4\!\cdot\!10^{-12}}\cr
\tabentry{26}{[1, 1, 0, 4, 11]}{x^3-26}{2^{2}3\!\cdot\!13^{2}}{2^{4}3^{2}7^{2}11^{2}13^{4}}{2}{1^2}{2\!\cdot\!3^{0}\phantom{\,\,\!+\!2\!\cdot\!3^1}\,\phantom{\,\,\!+\!2\!\cdot\!3^2}\,\!+\!O(3^{5})}{2\!\cdot\!3^{0}\!+\!2\!\cdot\!3^{1}\!+\!1\!\cdot\!3^{2}\!+\!O(3^{3})}\errentry{1\!\cdot\!10^{-10}}{3\!\cdot\!10^{-11}}\cr
\tabentry{28}{[1, 1, 0, 4, 11]}{x^3-28}{2^{2}3\!\cdot\!7^{2}}{2^{4}3^{2}7^{4}11^{2}}{2}{1^2}{2\!\cdot\!3^{0}\!+\!2\!\cdot\!3^{1}\!+\!1\!\cdot\!3^{2}\!+\!O(3^{3})}{2\!\cdot\!3^{0}\phantom{\,\,\!+\!2\!\cdot\!3^1}\,\!+\!2\!\cdot\!3^{2}\!+\!O(3^{3})}\errentry{4\!\cdot\!10^{-18}}{7\!\cdot\!10^{-16}}\cr
\tabentry{35}{[1, 1, 0, 4, 11]}{x^3-35}{3\!\cdot\!5^{2}7^{2}}{3^{2}5^{4}7^{4}11^{2}}{2\!\cdot\!5^{2}}{5^2}{2\!\cdot\!3^{0}\!+\!2\!\cdot\!3^{1}\!+\!2\!\cdot\!3^{2}\!+\!O(3^{4})}{2\!\cdot\!3^{0}\phantom{\,\,\!+\!2\!\cdot\!3^1}\,\!+\!1\!\cdot\!3^{2}\!+\!O(3^{3})}\errentry{4\!\cdot\!10^{-17}}{3\!\cdot\!10^{-15}}\cr
\tabentry{37}{[1, 1, 0, 4, 11]}{x^3-37}{3\!\cdot\!37^{2}}{3^{2}7^{2}11^{2}37^{4}}{2^{3}}{2^2}{1\!\cdot\!3^{0}\!+\!2\!\cdot\!3^{1}\!+\!1\!\cdot\!3^{2}\!+\!O(3^{3})}{1\!\cdot\!3^{0}\!+\!1\!\cdot\!3^{1}\phantom{\,\,\!+\!2\!\cdot\!3^2}\,\!+\!O(3^{3})}\errentry{5\!\cdot\!10^{-8}}{2\!\cdot\!10^{-8}}\cr
\tabentry{63}{[1, 1, 0, 4, 11]}{x^3-63}{3^{5}7^{2}}{3^{10}7^{4}11^{2}}{2}{1^2}{1\!\cdot\!3^{0}\!+\!1\!\cdot\!3^{1}\!+\!2\!\cdot\!3^{2}\!+\!O(3^{3})}{1\!\cdot\!3^{0}\phantom{\,\,\!+\!2\!\cdot\!3^1}\,\!+\!2\!\cdot\!3^{2}\!+\!O(3^{3})}\errentry{3\!\cdot\!10^{-12}}{9\!\cdot\!10^{-13}}\cr
\tabentry{91}{[1, 1, 0, 4, 11]}{x^3-91}{3\!\cdot\!7^{2}13^{2}}{3^{2}7^{4}11^{2}13^{4}}{2\!\cdot\!7^{2}}{7^2}{1\!\cdot\!3^{0}\!+\!2\!\cdot\!3^{1}\!+\!2\!\cdot\!3^{2}\!+\!O(3^{3})}{1\!\cdot\!3^{0}\!+\!1\!\cdot\!3^{1}\!+\!2\!\cdot\!3^{2}\!+\!O(3^{3})}\errentry{3\!\cdot\!10^{-8}}{6\!\cdot\!10^{-8}}\cr
\tabentry{350}{[1, 1, 0, 4, 11]}{x^3-350}{2^{2}3\!\cdot\!5^{2}7^{2}}{2^{4}3^{2}5^{4}7^{4}11^{2}}{2\!\cdot\!7^{2}}{7^2}{1\!\cdot\!3^{0}\!+\!2\!\cdot\!3^{1}\!+\!2\!\cdot\!3^{2}\!+\!O(3^{3})}{1\!\cdot\!3^{0}\!+\!2\!\cdot\!3^{1}\phantom{\,\,\!+\!2\!\cdot\!3^2}\,\!+\!O(3^{4})}\errentry{9\!\cdot\!10^{-11}}{1\!\cdot\!10^{-10}}\cr
\hline
\tabend

\noindent\begin{tabular}{cllll} 
\extravspace
$E=80B2$ & \multicolumn{3}{l}{$y^2=x^3-x^2-x$} & \tabname{3}{80B2}\cr
\extravspace
$p=3$ & \multicolumn{3}{l}{$P_3(E/\Q,T)=1-2T+3T^2$}\cr
$/\vQ\alsub{}$   & $N(E)=2^{4}\!\cdot\!5$ & $L^*(E/\Q)=2^{-2}$ & $|E(\Q)|=2$ & $|\sha(E/\Q)|=1^2$ \cr
$/K\alsub{}$      & $N(E,\sigma)=2^{8}\!\cdot\!3^{2}\!\cdot\!5^{2}$ & $L^*(E/K)=2^{-2}$ & $|E(K)|=2$ & $|\sha(E/K)|=1^2$ \cr
$/\vQ\alsub{2}$ & $\hbox{II},\>c_{2}=1$ & $/K\alsub{v}\quad(v|2)$ & $\hbox{II},\>c_v=1$\cr
$/\vQ\alsub{5}$ & $\hbox{I}_1 \hbox{ non-split},\>c_{5}=1$ & $/K\alsub{v}\quad(v|5)$ & $\hbox{I}_1 \hbox{ split},\>c_v=1$\cr
\end{tabular}
\tabbegin
\tabentry{2}{[0, -1, 0, -1, 0]}{x^3-2}{2^{2}3^{3}}{2^{8}3^{6}5^{2}}{2}{1^2}{2\!\cdot\!3^{0}\!+\!2\!\cdot\!3^{1}\phantom{\,\,\!+\!2\!\cdot\!3^2}\,\!+\!O(3^{3})}{2\!\cdot\!3^{0}\phantom{\,\,\!+\!2\!\cdot\!3^1}\,\!+\!1\!\cdot\!3^{2}\!+\!O(3^{3})}\errentry{5\!\cdot\!10^{-18}}{1\!\cdot\!10^{-15}}\cr
\tabentry{3}{[0, -1, 0, -1, 0]}{x^3-3}{3^{5}}{2^{8}3^{10}5^{2}}{1}{1^2}{2\!\cdot\!3^{0}\!+\!2\!\cdot\!3^{1}\phantom{\,\,\!+\!2\!\cdot\!3^2}\,\!+\!O(3^{3})}{2\!\cdot\!3^{0}\!+\!2\!\cdot\!3^{1}\!+\!2\!\cdot\!3^{2}\!+\!O(3^{3})}\errentry{5\!\cdot\!10^{-17}}{4\!\cdot\!10^{-16}}\cr
\tabentry{6}{[0, -1, 0, -1, 0]}{x^3-6}{2^{2}3^{5}}{2^{8}3^{10}5^{2}}{2}{1^2}{2\!\cdot\!3^{0}\!+\!2\!\cdot\!3^{1}\phantom{\,\,\!+\!2\!\cdot\!3^2}\,\!+\!O(3^{3})}{2\!\cdot\!3^{0}\!+\!2\!\cdot\!3^{1}\!+\!2\!\cdot\!3^{2}\!+\!O(3^{3})}\errentry{1\!\cdot\!10^{-16}}{4\!\cdot\!10^{-16}}\cr
\tabentry{7}{[0, -1, 0, -1, 0]}{x^3-7}{3^{3}7^{2}}{2^{8}3^{6}5^{2}7^{4}}{2^{2}}{2^2}{2\!\cdot\!3^{0}\phantom{\,\,\!+\!2\!\cdot\!3^1}\,\phantom{\,\,\!+\!2\!\cdot\!3^2}\,\!+\!O(3^{3})}{2\!\cdot\!3^{0}\!+\!1\!\cdot\!3^{1}\!+\!1\!\cdot\!3^{2}\!+\!O(3^{3})}\errentry{6\!\cdot\!10^{-15}}{7\!\cdot\!10^{-15}}\cr
\tabentry{11}{[0, -1, 0, -1, 0]}{x^3-11}{3^{3}11^{2}}{2^{8}3^{6}5^{2}11^{4}}{2^{2}3^{2}}{6^2}{2\!\cdot\!3^{2}\!+\!1\!\cdot\!3^{3}\!+\!1\!\cdot\!3^{4}\!+\!O(3^{5})}{1\!\cdot\!3^{2}\!+\!1\!\cdot\!3^{3}\phantom{\,\,\!+\!2\!\cdot\!3^4}\,\!+\!O(3^{5})}\errentry{1\!\cdot\!10^{-10}}{1\!\cdot\!10^{-10}}\cr
\tabentry{12}{[0, -1, 0, -1, 0]}{x^3-12}{2^{2}3^{5}}{2^{8}3^{10}5^{2}}{2^{3}}{2^2}{2\!\cdot\!3^{0}\!+\!2\!\cdot\!3^{1}\phantom{\,\,\!+\!2\!\cdot\!3^2}\,\!+\!O(3^{3})}{2\!\cdot\!3^{0}\!+\!1\!\cdot\!3^{1}\!+\!2\!\cdot\!3^{2}\!+\!O(3^{3})}\errentry{4\!\cdot\!10^{-17}}{1\!\cdot\!10^{-15}}\cr
\tabentry{14}{[0, -1, 0, -1, 0]}{x^3-14}{2^{2}3^{3}7^{2}}{2^{8}3^{6}5^{2}7^{4}}{2}{1^2}{2\!\cdot\!3^{0}\phantom{\,\,\!+\!2\!\cdot\!3^1}\,\phantom{\,\,\!+\!2\!\cdot\!3^2}\,\!+\!O(3^{3})}{2\!\cdot\!3^{0}\!+\!2\!\cdot\!3^{1}\!+\!1\!\cdot\!3^{2}\!+\!O(3^{4})}\errentry{7\!\cdot\!10^{-14}}{2\!\cdot\!10^{-14}}\cr
\tabentry{17}{[0, -1, 0, -1, 0]}{x^3-17}{3\!\cdot\!17^{2}}{2^{8}3^{2}5^{2}17^{4}}{3^{2}}{3^2}{1\!\cdot\!3^{2}\!+\!1\!\cdot\!3^{3}\!+\!2\!\cdot\!3^{4}\!+\!O(3^{6})}{1\!\cdot\!3^{2}\!+\!1\!\cdot\!3^{3}\!+\!1\!\cdot\!3^{4}\!+\!O(3^{5})}\errentry{4\!\cdot\!10^{-18}}{2\!\cdot\!10^{-15}}\cr
\tabentry{19}{[0, -1, 0, -1, 0]}{x^3-19}{3\!\cdot\!19^{2}}{2^{8}3^{2}5^{2}19^{4}}{1}{1^2}{2\!\cdot\!3^{0}\!+\!1\!\cdot\!3^{1}\phantom{\,\,\!+\!2\!\cdot\!3^2}\,\!+\!O(3^{3})}{2\!\cdot\!3^{0}\!+\!1\!\cdot\!3^{1}\!+\!2\!\cdot\!3^{2}\!+\!O(3^{3})}\errentry{2\!\cdot\!10^{-16}}{5\!\cdot\!10^{-16}}\cr
\tabentry{22}{[0, -1, 0, -1, 0]}{x^3-22}{2^{2}3^{3}11^{2}}{2^{8}3^{6}5^{2}11^{4}}{0}{\nosha}{2\!\cdot\!3^{2}\!+\!1\!\cdot\!3^{3}\!+\!1\!\cdot\!3^{4}\!+\!O(3^{5})}{\noeuler}\errentry{1\!\cdot\!10^{-9}}{6\!\cdot\!10^{-6}}\cr
\tabentry{26}{[0, -1, 0, -1, 0]}{x^3-26}{2^{2}3\!\cdot\!13^{2}}{2^{8}3^{2}5^{2}13^{4}}{2\!\cdot\!3^{2}}{3^2}{2\!\cdot\!3^{2}\!+\!2\!\cdot\!3^{3}\!+\!2\!\cdot\!3^{4}\!+\!O(3^{6})}{2\!\cdot\!3^{2}\!+\!2\!\cdot\!3^{3}\!+\!1\!\cdot\!3^{4}\!+\!O(3^{5})}\errentry{1\!\cdot\!10^{-16}}{1\!\cdot\!10^{-15}}\cr
\tabentry{28}{[0, -1, 0, -1, 0]}{x^3-28}{2^{2}3\!\cdot\!7^{2}}{2^{8}3^{2}5^{2}7^{4}}{2}{1^2}{2\!\cdot\!3^{0}\phantom{\,\,\!+\!2\!\cdot\!3^1}\,\phantom{\,\,\!+\!2\!\cdot\!3^2}\,\!+\!O(3^{3})}{2\!\cdot\!3^{0}\phantom{\,\,\!+\!2\!\cdot\!3^1}\,\!+\!2\!\cdot\!3^{2}\!+\!O(3^{3})}\errentry{7\!\cdot\!10^{-17}}{2\!\cdot\!10^{-15}}\cr
\tabentry{37}{[0, -1, 0, -1, 0]}{x^3-37}{3\!\cdot\!37^{2}}{2^{8}3^{2}5^{2}37^{4}}{2^{2}3^{2}}{6^2}{2\!\cdot\!3^{4}\phantom{\,\,\!+\!2\!\cdot\!3^5}\,\!+\!1\!\cdot\!3^{6}\!+\!O(3^{7})}{2\!\cdot\!3^{2}\phantom{\,\,\!+\!2\!\cdot\!3^3}\,\phantom{\,\,\!+\!2\!\cdot\!3^4}\,\!+\!O(3^{5})}\errentry{1\!\cdot\!10^{-8}}{1\!\cdot\!10^{-8}}\cr
\tabentry{44}{[0, -1, 0, -1, 0]}{x^3-44}{2^{2}3\!\cdot\!11^{2}}{2^{8}3^{2}5^{2}11^{4}}{2\!\cdot\!3^{2}}{3^2}{2\!\cdot\!3^{2}\!+\!1\!\cdot\!3^{3}\!+\!1\!\cdot\!3^{4}\!+\!O(3^{5})}{1\!\cdot\!3^{2}\!+\!2\!\cdot\!3^{3}\!+\!2\!\cdot\!3^{4}\!+\!O(3^{5})}\errentry{2\!\cdot\!10^{-17}}{6\!\cdot\!10^{-16}}\cr
\tabentry{46}{[0, -1, 0, -1, 0]}{x^3-46}{2^{2}3\!\cdot\!23^{2}}{2^{8}3^{2}5^{2}23^{4}}{2\!\cdot\!3^{2}}{3^2}{1\!\cdot\!3^{3}\phantom{\,\,\!+\!2\!\cdot\!3^4}\,\!+\!1\!\cdot\!3^{5}\!+\!O(3^{9})}{1\!\cdot\!3^{2}\phantom{\,\,\!+\!2\!\cdot\!3^3}\,\!+\!2\!\cdot\!3^{4}\!+\!O(3^{5})}\errentry{1\!\cdot\!10^{-13}}{1\!\cdot\!10^{-13}}\cr
\tabentry{62}{[0, -1, 0, -1, 0]}{x^3-62}{2^{2}3\!\cdot\!31^{2}}{2^{8}3^{2}5^{2}31^{4}}{2}{1^2}{2\!\cdot\!3^{0}\!+\!1\!\cdot\!3^{1}\!+\!1\!\cdot\!3^{2}\!+\!O(3^{8})}{2\!\cdot\!3^{0}\!+\!2\!\cdot\!3^{1}\phantom{\,\,\!+\!2\!\cdot\!3^2}\,\!+\!O(3^{4})}\errentry{3\!\cdot\!10^{-9}}{9\!\cdot\!10^{-10}}\cr
\tabentry{116}{[0, -1, 0, -1, 0]}{x^3-116}{2^{2}3\!\cdot\!29^{2}}{2^{8}3^{2}5^{2}29^{4}}{2\!\cdot\!3^{2}}{3^2}{1\!\cdot\!3^{3}\phantom{\,\,\!+\!2\!\cdot\!3^4}\,\!+\!1\!\cdot\!3^{5}\!+\!O(3^{6})}{1\!\cdot\!3^{2}\!+\!2\!\cdot\!3^{3}\!+\!1\!\cdot\!3^{4}\!+\!O(3^{5})}\errentry{7\!\cdot\!10^{-11}}{5\!\cdot\!10^{-11}}\cr
\hline
\tabend

\noindent\begin{tabular}{cllll} 
\extravspace
$E=260A1$ & \multicolumn{3}{l}{$y^2=x^3-x^2-281x+1910$} & \tabname{3}{260A1}\cr
\extravspace
$p=3$ & \multicolumn{3}{l}{$P_3(E/\Q,T)=1-2T+3T^2$}\cr
$/\vQ\alsub{}$   & $N(E)=2^{2}\!\cdot\!5\!\cdot\!13$ & $L^*(E/\Q)=2^{-1}$ & $|E(\Q)|=2$ & $|\sha(E/\Q)|=1^2$ \cr
$/K\alsub{}$      & $N(E,\sigma)=2^{4}\!\cdot\!3^{2}\!\cdot\!5^{2}\!\cdot\!13^{2}$ & $L^*(E/K)=3$ & $|E(K)|=2$ & $|\sha(E/K)|=1^2$ \cr
$/\vQ\alsub{2}$ & $\hbox{IV},\>c_{2}=1$ & $/K\alsub{v}\quad(v|2)$ & $\hbox{IV},\>c_v=3$\cr
$/\vQ\alsub{5}$ & $\hbox{I}_1 \hbox{ non-split},\>c_{5}=1$ & $/K\alsub{v}\quad(v|5)$ & $\hbox{I}_1 \hbox{ split},\>c_v=1$\cr
$/\vQ\alsub{13}$ & $\hbox{I}_2 \hbox{ non-split},\>c_{13}=2$ & $/K\alsub{v}\quad(v|13)$ & $\hbox{I}_2 \hbox{ non-split},\>c_v=2$\cr
\end{tabular}
\tabbegin
\tabentry{3}{[0, -1, 0, -281, 1910]}{x^3-3}{3^{5}}{2^{4}3^{10}5^{2}13^{2}}{3}{1^2}{2\!\cdot\!3^{1}\!+\!1\!\cdot\!3^{2}\phantom{\,\,\!+\!2\!\cdot\!3^3}\,\!+\!O(3^{5})}{2\!\cdot\!3^{1}\!+\!2\!\cdot\!3^{2}\!+\!2\!\cdot\!3^{3}\!+\!O(3^{4})}\errentry{1\!\cdot\!10^{-17}}{1\!\cdot\!10^{-15}}\cr
\tabentry{7}{[0, -1, 0, -281, 1910]}{x^3-7}{3^{3}7^{2}}{2^{4}3^{6}5^{2}7^{4}13^{2}}{0}{\nosha}{2\!\cdot\!3^{3}\!+\!1\!\cdot\!3^{4}\!+\!2\!\cdot\!3^{5}\!+\!O(3^{6})}{\noeuler}\errentry{1\!\cdot\!10^{-8}}{4\!\cdot\!10^{-5}}\cr
\tabentry{10}{[0, -1, 0, -281, 1910]}{x^3-10}{2^{2}3\!\cdot\!5^{2}}{2^{2}3^{2}5^{4}13^{2}}{2\!\cdot\!3^{2}}{3^2}{1\!\cdot\!3^{2}\!+\!2\!\cdot\!3^{3}\phantom{\,\,\!+\!2\!\cdot\!3^4}\,\!+\!O(3^{5})}{2\!\cdot\!3^{3}\!+\!1\!\cdot\!3^{4}\!+\!1\!\cdot\!3^{5}\!+\!O(3^{6})}\errentry{1\!\cdot\!10^{-16}}{5\!\cdot\!10^{-15}}\cr
\tabentry{13}{[0, -1, 0, -281, 1910]}{x^3-13}{3^{3}13^{2}}{2^{4}3^{6}5^{2}13^{4}}{3\!\cdot\!7^{2}}{7^2}{2\!\cdot\!3^{1}\!+\!1\!\cdot\!3^{2}\phantom{\,\,\!+\!2\!\cdot\!3^3}\,\!+\!O(3^{4})}{2\!\cdot\!3^{1}\phantom{\,\,\!+\!2\!\cdot\!3^2}\,\phantom{\,\,\!+\!2\!\cdot\!3^3}\,\!+\!O(3^{4})}\errentry{6\!\cdot\!10^{-17}}{4\!\cdot\!10^{-15}}\cr
\tabentry{17}{[0, -1, 0, -281, 1910]}{x^3-17}{3\!\cdot\!17^{2}}{2^{4}3^{2}5^{2}13^{2}17^{4}}{3}{1^2}{1\!\cdot\!3^{1}\!+\!2\!\cdot\!3^{2}\!+\!1\!\cdot\!3^{3}\!+\!O(3^{4})}{1\!\cdot\!3^{1}\!+\!1\!\cdot\!3^{2}\!+\!1\!\cdot\!3^{3}\!+\!O(3^{4})}\errentry{6\!\cdot\!10^{-10}}{1\!\cdot\!10^{-10}}\cr
\tabentry{19}{[0, -1, 0, -281, 1910]}{x^3-19}{3\!\cdot\!19^{2}}{2^{4}3^{2}5^{2}13^{2}19^{4}}{3\!\cdot\!5^{2}}{5^2}{2\!\cdot\!3^{1}\phantom{\,\,\!+\!2\!\cdot\!3^2}\,\!+\!1\!\cdot\!3^{3}\!+\!O(3^{4})}{2\!\cdot\!3^{1}\!+\!2\!\cdot\!3^{2}\phantom{\,\,\!+\!2\!\cdot\!3^3}\,\!+\!O(3^{4})}\errentry{3\!\cdot\!10^{-9}}{4\!\cdot\!10^{-9}}\cr
\tabentry{20}{[0, -1, 0, -281, 1910]}{x^3-20}{2^{2}3^{3}5^{2}}{2^{2}3^{6}5^{4}13^{2}}{2\!\cdot\!3^{2}}{3^2}{1\!\cdot\!3^{2}\!+\!2\!\cdot\!3^{3}\phantom{\,\,\!+\!2\!\cdot\!3^4}\,\!+\!O(3^{5})}{2\!\cdot\!3^{3}\phantom{\,\,\!+\!2\!\cdot\!3^4}\,\!+\!2\!\cdot\!3^{5}\!+\!O(3^{11})}\errentry{4\!\cdot\!10^{-17}}{1\!\cdot\!10^{-16}}\cr
\tabentry{30}{[0, -1, 0, -281, 1910]}{x^3-30}{2^{2}3^{5}5^{2}}{2^{2}3^{10}5^{4}13^{2}}{2\!\cdot\!3^{4}}{9^2}{1\!\cdot\!3^{2}\!+\!2\!\cdot\!3^{3}\phantom{\,\,\!+\!2\!\cdot\!3^4}\,\!+\!O(3^{5})}{2\!\cdot\!3^{5}\!+\!2\!\cdot\!3^{6}\phantom{\,\,\!+\!2\!\cdot\!3^7}\,\!+\!O(3^{9})}\errentry{1\!\cdot\!10^{-10}}{3\!\cdot\!10^{-10}}\cr
\tabentry{60}{[0, -1, 0, -281, 1910]}{x^3-60}{2^{2}3^{5}5^{2}}{2^{2}3^{10}5^{4}13^{2}}{0}{\nosha}{1\!\cdot\!3^{2}\!+\!2\!\cdot\!3^{3}\phantom{\,\,\!+\!2\!\cdot\!3^4}\,\!+\!O(3^{5})}{\noeuler}\errentry{7\!\cdot\!10^{-12}}{4\!\cdot\!10^{-7}}\cr
\tabentry{90}{[0, -1, 0, -281, 1910]}{x^3-90}{2^{2}3^{5}5^{2}}{2^{2}3^{10}5^{4}13^{2}}{2^{5}3^{2}}{6^2}{1\!\cdot\!3^{2}\!+\!2\!\cdot\!3^{3}\phantom{\,\,\!+\!2\!\cdot\!3^4}\,\!+\!O(3^{5})}{2\!\cdot\!3^{3}\phantom{\,\,\!+\!2\!\cdot\!3^4}\,\!+\!2\!\cdot\!3^{5}\!+\!O(3^{6})}\errentry{6\!\cdot\!10^{-12}}{9\!\cdot\!10^{-12}}\cr
\tabentry{130}{[0, -1, 0, -281, 1910]}{x^3-130}{2^{2}3^{3}5^{2}13^{2}}{2^{2}3^{6}5^{4}13^{4}}{2\!\cdot\!3^{2}7^{2}}{21^2}{1\!\cdot\!3^{2}\!+\!2\!\cdot\!3^{3}\phantom{\,\,\!+\!2\!\cdot\!3^4}\,\!+\!O(3^{5})}{2\!\cdot\!3^{3}\phantom{\,\,\!+\!2\!\cdot\!3^4}\,\!+\!1\!\cdot\!3^{5}\!+\!O(3^{6})}\errentry{1\!\cdot\!10^{-9}}{6\!\cdot\!10^{-9}}\cr
\tabentry{150}{[0, -1, 0, -281, 1910]}{x^3-150}{2^{2}3^{5}5^{2}}{2^{2}3^{10}5^{4}13^{2}}{2\!\cdot\!3^{2}}{3^2}{1\!\cdot\!3^{2}\!+\!2\!\cdot\!3^{3}\phantom{\,\,\!+\!2\!\cdot\!3^4}\,\!+\!O(3^{5})}{2\!\cdot\!3^{3}\!+\!2\!\cdot\!3^{4}\phantom{\,\,\!+\!2\!\cdot\!3^5}\,\!+\!O(3^{7})}\errentry{4\!\cdot\!10^{-12}}{3\!\cdot\!10^{-12}}\cr
\tabentry{350}{[0, -1, 0, -281, 1910]}{x^3-350}{2^{2}3\!\cdot\!5^{2}7^{2}}{2^{2}3^{2}5^{4}7^{4}13^{2}}{2^{3}3^{2}}{3^2}{1\!\cdot\!3^{4}\!+\!2\!\cdot\!3^{5}\!+\!1\!\cdot\!3^{6}\!+\!O(3^{8})}{2\!\cdot\!3^{3}\!+\!2\!\cdot\!3^{4}\phantom{\,\,\!+\!2\!\cdot\!3^5}\,\!+\!O(3^{6})}\errentry{1\!\cdot\!10^{-14}}{6\!\cdot\!10^{-15}}\cr
\tabentry{550}{[0, -1, 0, -281, 1910]}{x^3-550}{2^{2}3\!\cdot\!5^{2}11^{2}}{2^{2}3^{2}5^{4}11^{4}13^{2}}{2\!\cdot\!3^{4}}{9^2}{2\!\cdot\!3^{2}\!+\!2\!\cdot\!3^{3}\!+\!1\!\cdot\!3^{4}\!+\!O(3^{5})}{1\!\cdot\!3^{5}\!+\!2\!\cdot\!3^{6}\!+\!1\!\cdot\!3^{7}\!+\!O(3^{8})}\errentry{3\!\cdot\!10^{-9}}{7\!\cdot\!10^{-9}}\cr
\tabentry{650}{[0, -1, 0, -281, 1910]}{x^3-650}{2^{2}3^{3}5^{2}13^{2}}{2^{2}3^{6}5^{4}13^{4}}{2^{3}3^{4}}{18^2}{1\!\cdot\!3^{2}\!+\!2\!\cdot\!3^{3}\phantom{\,\,\!+\!2\!\cdot\!3^4}\,\!+\!O(3^{5})}{2\!\cdot\!3^{5}\phantom{\,\,\!+\!2\!\cdot\!3^6}\,\phantom{\,\,\!+\!2\!\cdot\!3^7}\,\!+\!O(3^{8})}\errentry{3\!\cdot\!10^{-10}}{1\!\cdot\!10^{-9}}\cr
\tabentry{1300}{[0, -1, 0, -281, 1910]}{x^3-1300}{2^{2}3^{3}5^{2}13^{2}}{2^{2}3^{6}5^{4}13^{4}}{0}{\nosha}{1\!\cdot\!3^{2}\!+\!2\!\cdot\!3^{3}\phantom{\,\,\!+\!2\!\cdot\!3^4}\,\!+\!O(3^{5})}{\noeuler}\errentry{9\!\cdot\!10^{-9}}{2\!\cdot\!10^{-5}}\cr
\hline
\tabend

\noindent\begin{tabular}{cllll} 
\extravspace
$E=275B1$ & \multicolumn{3}{l}{$y^2+y=x^3+x^2-8x+19$} & \tabname{3}{275B1}\cr
\extravspace
$p=3$ & \multicolumn{3}{l}{$P_3(E/\Q,T)=1-T+3T^2$}\cr
$/\vQ\alsub{}$   & $N(E)=5^{2}\!\cdot\!11$ & $L^*(E/\Q)=1$ & $|E(\Q)|=1$ & $|\sha(E/\Q)|=1^2$ \cr
$/K\alsub{}$      & $N(E,\sigma)=3^{2}\!\cdot\!5^{4}\!\cdot\!11^{2}$ & $L^*(E/K)=1$ & $|E(K)|=1$ & $|\sha(E/K)|=1^2$ \cr
$/\vQ\alsub{5}$ & $\hbox{I}^*_0,\>c_{5}=1$ & $/K\alsub{v}\quad(v|5)$ & $\hbox{I}^*_0,\>c_v=1$\cr
$/\vQ\alsub{11}$ & $\hbox{I}_1 \hbox{ split},\>c_{11}=1$ & $/K\alsub{v}\quad(v|11)$ & $\hbox{I}_1 \hbox{ split},\>c_v=1$\cr
\end{tabular}
\tabbegin
\tabentry{2}{[0, 1, 1, -8, 19]}{x^3-2}{2^{2}3^{3}}{2^{4}3^{6}5^{4}11^{2}}{3^{2}}{3^2}{2\!\cdot\!3^{2}\!+\!1\!\cdot\!3^{3}\!+\!2\!\cdot\!3^{4}\!+\!O(3^{5})}{2\!\cdot\!3^{2}\!+\!2\!\cdot\!3^{3}\!+\!1\!\cdot\!3^{4}\!+\!O(3^{5})}\errentry{4\!\cdot\!10^{-16}}{4\!\cdot\!10^{-15}}\cr
\tabentry{3}{[0, 1, 1, -8, 19]}{x^3-3}{3^{5}}{3^{10}5^{4}11^{2}}{3^{2}}{3^2}{1\!\cdot\!3^{2}\!+\!1\!\cdot\!3^{3}\!+\!1\!\cdot\!3^{4}\!+\!O(3^{6})}{1\!\cdot\!3^{2}\!+\!1\!\cdot\!3^{3}\!+\!1\!\cdot\!3^{4}\!+\!O(3^{6})}\errentry{1\!\cdot\!10^{-16}}{4\!\cdot\!10^{-15}}\cr
\tabentry{5}{[0, 1, 1, -8, 19]}{x^3-5}{3^{3}5^{2}}{3^{6}5^{4}11^{2}}{0}{\nosha}{1\!\cdot\!3^{2}\!+\!1\!\cdot\!3^{3}\!+\!1\!\cdot\!3^{4}\!+\!O(3^{6})}{\noeuler}\errentry{6\!\cdot\!10^{-17}}{1\!\cdot\!10^{-9}}\cr
\tabentry{6}{[0, 1, 1, -8, 19]}{x^3-6}{2^{2}3^{5}}{2^{4}3^{10}5^{4}11^{2}}{3^{2}}{3^2}{2\!\cdot\!3^{2}\!+\!1\!\cdot\!3^{3}\!+\!2\!\cdot\!3^{4}\!+\!O(3^{5})}{2\!\cdot\!3^{2}\phantom{\,\,\!+\!2\!\cdot\!3^3}\,\!+\!2\!\cdot\!3^{4}\!+\!O(3^{5})}\errentry{1\!\cdot\!10^{-9}}{8\!\cdot\!10^{-10}}\cr
\tabentry{10}{[0, 1, 1, -8, 19]}{x^3-10}{2^{2}3\!\cdot\!5^{2}}{2^{4}3^{2}5^{4}11^{2}}{5}{1^2}{2\!\cdot\!3^{2}\!+\!1\!\cdot\!3^{3}\!+\!2\!\cdot\!3^{4}\!+\!O(3^{5})}{2\!\cdot\!3^{2}\phantom{\,\,\!+\!2\!\cdot\!3^3}\,\!+\!2\!\cdot\!3^{4}\!+\!O(3^{5})}\errentry{1\!\cdot\!10^{-16}}{1\!\cdot\!10^{-15}}\cr
\tabentry{12}{[0, 1, 1, -8, 19]}{x^3-12}{2^{2}3^{5}}{2^{4}3^{10}5^{4}11^{2}}{0}{\nosha}{2\!\cdot\!3^{2}\!+\!1\!\cdot\!3^{3}\!+\!2\!\cdot\!3^{4}\!+\!O(3^{5})}{\noeuler}\errentry{3\!\cdot\!10^{-8}}{1\!\cdot\!10^{-5}}\cr
\tabentry{15}{[0, 1, 1, -8, 19]}{x^3-15}{3^{5}5^{2}}{3^{10}5^{4}11^{2}}{3^{2}5}{3^2}{1\!\cdot\!3^{2}\!+\!1\!\cdot\!3^{3}\!+\!1\!\cdot\!3^{4}\!+\!O(3^{6})}{1\!\cdot\!3^{2}\!+\!1\!\cdot\!3^{3}\!+\!1\!\cdot\!3^{4}\!+\!O(3^{6})}\errentry{1\!\cdot\!10^{-16}}{4\!\cdot\!10^{-15}}\cr
\tabentry{17}{[0, 1, 1, -8, 19]}{x^3-17}{3\!\cdot\!17^{2}}{3^{2}5^{4}11^{2}17^{4}}{1}{1^2}{2\!\cdot\!3^{2}\phantom{\,\,\!+\!2\!\cdot\!3^3}\,\!+\!2\!\cdot\!3^{4}\!+\!O(3^{8})}{2\!\cdot\!3^{2}\!+\!1\!\cdot\!3^{3}\!+\!2\!\cdot\!3^{4}\!+\!O(3^{5})}\errentry{2\!\cdot\!10^{-8}}{7\!\cdot\!10^{-9}}\cr
\tabentry{19}{[0, 1, 1, -8, 19]}{x^3-19}{3\!\cdot\!19^{2}}{3^{2}5^{4}11^{2}19^{4}}{1}{1^2}{1\!\cdot\!3^{2}\!+\!2\!\cdot\!3^{3}\phantom{\,\,\!+\!2\!\cdot\!3^4}\,\!+\!O(3^{6})}{1\!\cdot\!3^{2}\!+\!1\!\cdot\!3^{3}\!+\!2\!\cdot\!3^{4}\!+\!O(3^{6})}\errentry{8\!\cdot\!10^{-9}}{1\!\cdot\!10^{-9}}\cr
\tabentry{20}{[0, 1, 1, -8, 19]}{x^3-20}{2^{2}3^{3}5^{2}}{2^{4}3^{6}5^{4}11^{2}}{3^{2}5}{3^2}{2\!\cdot\!3^{2}\!+\!1\!\cdot\!3^{3}\!+\!2\!\cdot\!3^{4}\!+\!O(3^{5})}{2\!\cdot\!3^{2}\!+\!2\!\cdot\!3^{3}\!+\!1\!\cdot\!3^{4}\!+\!O(3^{5})}\errentry{1\!\cdot\!10^{-16}}{1\!\cdot\!10^{-16}}\cr
\tabentry{28}{[0, 1, 1, -8, 19]}{x^3-28}{2^{2}3\!\cdot\!7^{2}}{2^{4}3^{2}5^{4}7^{4}11^{2}}{3^{2}}{3^2}{2\!\cdot\!3^{4}\!+\!1\!\cdot\!3^{5}\!+\!1\!\cdot\!3^{6}\!+\!O(3^{7})}{2\!\cdot\!3^{4}\!+\!2\!\cdot\!3^{5}\!+\!2\!\cdot\!3^{6}\!+\!O(3^{7})}\errentry{5\!\cdot\!10^{-12}}{4\!\cdot\!10^{-12}}\cr
\tabentry{30}{[0, 1, 1, -8, 19]}{x^3-30}{2^{2}3^{5}5^{2}}{2^{4}3^{10}5^{4}11^{2}}{3^{2}5}{3^2}{2\!\cdot\!3^{2}\!+\!1\!\cdot\!3^{3}\!+\!2\!\cdot\!3^{4}\!+\!O(3^{5})}{2\!\cdot\!3^{2}\phantom{\,\,\!+\!2\!\cdot\!3^3}\,\!+\!2\!\cdot\!3^{4}\!+\!O(3^{5})}\errentry{2\!\cdot\!10^{-11}}{2\!\cdot\!10^{-11}}\cr
\tabentry{35}{[0, 1, 1, -8, 19]}{x^3-35}{3\!\cdot\!5^{2}7^{2}}{3^{2}5^{4}7^{4}11^{2}}{3^{2}5}{3^2}{1\!\cdot\!3^{4}\!+\!1\!\cdot\!3^{5}\!+\!2\!\cdot\!3^{6}\!+\!O(3^{7})}{1\!\cdot\!3^{4}\!+\!2\!\cdot\!3^{5}\!+\!2\!\cdot\!3^{6}\!+\!O(3^{8})}\errentry{3\!\cdot\!10^{-16}}{2\!\cdot\!10^{-15}}\cr
\tabentry{45}{[0, 1, 1, -8, 19]}{x^3-45}{3^{5}5^{2}}{3^{10}5^{4}11^{2}}{3^{2}5}{3^2}{1\!\cdot\!3^{2}\!+\!1\!\cdot\!3^{3}\!+\!1\!\cdot\!3^{4}\!+\!O(3^{6})}{1\!\cdot\!3^{2}\!+\!1\!\cdot\!3^{3}\!+\!1\!\cdot\!3^{4}\!+\!O(3^{6})}\errentry{1\!\cdot\!10^{-16}}{4\!\cdot\!10^{-15}}\cr
\tabentry{60}{[0, 1, 1, -8, 19]}{x^3-60}{2^{2}3^{5}5^{2}}{2^{4}3^{10}5^{4}11^{2}}{3^{2}5}{3^2}{2\!\cdot\!3^{2}\!+\!1\!\cdot\!3^{3}\!+\!2\!\cdot\!3^{4}\!+\!O(3^{5})}{2\!\cdot\!3^{2}\phantom{\,\,\!+\!2\!\cdot\!3^3}\,\!+\!2\!\cdot\!3^{4}\!+\!O(3^{5})}\errentry{1\!\cdot\!10^{-8}}{1\!\cdot\!10^{-8}}\cr
\tabentry{90}{[0, 1, 1, -8, 19]}{x^3-90}{2^{2}3^{5}5^{2}}{2^{4}3^{10}5^{4}11^{2}}{0}{\nosha}{2\!\cdot\!3^{2}\!+\!1\!\cdot\!3^{3}\!+\!2\!\cdot\!3^{4}\!+\!O(3^{5})}{\noeuler}\errentry{3\!\cdot\!10^{-9}}{1\!\cdot\!10^{-5}}\cr
\tabentry{150}{[0, 1, 1, -8, 19]}{x^3-150}{2^{2}3^{5}5^{2}}{2^{4}3^{10}5^{4}11^{2}}{3^{2}5}{3^2}{2\!\cdot\!3^{2}\!+\!1\!\cdot\!3^{3}\!+\!2\!\cdot\!3^{4}\!+\!O(3^{5})}{2\!\cdot\!3^{2}\phantom{\,\,\!+\!2\!\cdot\!3^3}\,\!+\!2\!\cdot\!3^{4}\!+\!O(3^{5})}\errentry{1\!\cdot\!10^{-8}}{7\!\cdot\!10^{-9}}\cr
\tabentry{325}{[0, 1, 1, -8, 19]}{x^3-325}{3\!\cdot\!5^{2}13^{2}}{3^{2}5^{4}11^{2}13^{4}}{3^{2}5}{3^2}{1\!\cdot\!3^{6}\phantom{\,\,\!+\!2\!\cdot\!3^7}\,\!+\!2\!\cdot\!3^{8}\!+\!O(3^{10})}{1\!\cdot\!3^{4}\phantom{\,\,\!+\!2\!\cdot\!3^5}\,\phantom{\,\,\!+\!2\!\cdot\!3^6}\,\!+\!O(3^{8})}\errentry{1\!\cdot\!10^{-12}}{1\!\cdot\!10^{-12}}\cr
\tabentry{350}{[0, 1, 1, -8, 19]}{x^3-350}{2^{2}3\!\cdot\!5^{2}7^{2}}{2^{4}3^{2}5^{4}7^{4}11^{2}}{3^{2}5}{3^2}{2\!\cdot\!3^{4}\!+\!1\!\cdot\!3^{5}\!+\!1\!\cdot\!3^{6}\!+\!O(3^{7})}{2\!\cdot\!3^{4}\!+\!2\!\cdot\!3^{5}\!+\!2\!\cdot\!3^{6}\!+\!O(3^{7})}\errentry{7\!\cdot\!10^{-12}}{5\!\cdot\!10^{-12}}\cr
\hline
\tabend

\noindent\begin{tabular}{cllll} 
\extravspace
$E=272C1$ & \multicolumn{3}{l}{$y^2=x^3-x^2-4x$} & \tabname{3}{272C1}\cr
\extravspace
$p=3$ & \multicolumn{3}{l}{$P_3(E/\Q,T)=1-2T+3T^2$}\cr
$/\vQ\alsub{}$   & $N(E)=2^{4}\!\cdot\!17$ & $L^*(E/\Q)=2^{-1}$ & $|E(\Q)|=2$ & $|\sha(E/\Q)|=1^2$ \cr
$/K\alsub{}$      & $N(E,\sigma)=2^{8}\!\cdot\!3^{2}\!\cdot\!17^{2}$ & $L^*(E/K)=1$ & $|E(K)|=2$ & $|\sha(E/K)|=1^2$ \cr
$/\vQ\alsub{2}$ & $\hbox{I}^*_0,\>c_{2}=2$ & $/K\alsub{v}\quad(v|2)$ & $\hbox{I}^*_0,\>c_v=4$\cr
$/\vQ\alsub{17}$ & $\hbox{I}_1 \hbox{ split},\>c_{17}=1$ & $/K\alsub{v}\quad(v|17)$ & $\hbox{I}_1 \hbox{ split},\>c_v=1$\cr
\end{tabular}
\tabbegin
\tabentry{2}{[0, -1, 0, -4, 0]}{x^3-2}{2^{2}3^{3}}{2^{8}3^{6}17^{2}}{2}{1^2}{2\!\cdot\!3^{0}\!+\!1\!\cdot\!3^{1}\phantom{\,\,\!+\!2\!\cdot\!3^2}\,\!+\!O(3^{4})}{2\!\cdot\!3^{0}\phantom{\,\,\!+\!2\!\cdot\!3^1}\,\!+\!1\!\cdot\!3^{2}\!+\!O(3^{3})}\errentry{9\!\cdot\!10^{-17}}{3\!\cdot\!10^{-15}}\cr
\tabentry{3}{[0, -1, 0, -4, 0]}{x^3-3}{3^{5}}{2^{8}3^{10}17^{2}}{2^{2}}{1^2}{2\!\cdot\!3^{0}\!+\!1\!\cdot\!3^{1}\phantom{\,\,\!+\!2\!\cdot\!3^2}\,\!+\!O(3^{4})}{2\!\cdot\!3^{0}\!+\!1\!\cdot\!3^{1}\!+\!2\!\cdot\!3^{2}\!+\!O(3^{3})}\errentry{2\!\cdot\!10^{-17}}{1\!\cdot\!10^{-15}}\cr
\tabentry{5}{[0, -1, 0, -4, 0]}{x^3-5}{3^{3}5^{2}}{2^{8}3^{6}5^{4}17^{2}}{2^{2}}{1^2}{1\!\cdot\!3^{0}\phantom{\,\,\!+\!2\!\cdot\!3^1}\,\phantom{\,\,\!+\!2\!\cdot\!3^2}\,\!+\!O(3^{3})}{1\!\cdot\!3^{0}\!+\!2\!\cdot\!3^{1}\!+\!2\!\cdot\!3^{2}\!+\!O(3^{3})}\errentry{6\!\cdot\!10^{-11}}{1\!\cdot\!10^{-11}}\cr
\tabentry{6}{[0, -1, 0, -4, 0]}{x^3-6}{2^{2}3^{5}}{2^{8}3^{10}17^{2}}{2}{1^2}{2\!\cdot\!3^{0}\!+\!1\!\cdot\!3^{1}\phantom{\,\,\!+\!2\!\cdot\!3^2}\,\!+\!O(3^{4})}{2\!\cdot\!3^{0}\!+\!2\!\cdot\!3^{1}\!+\!2\!\cdot\!3^{2}\!+\!O(3^{3})}\errentry{2\!\cdot\!10^{-18}}{2\!\cdot\!10^{-16}}\cr
\tabentry{7}{[0, -1, 0, -4, 0]}{x^3-7}{3^{3}7^{2}}{2^{8}3^{6}7^{4}17^{2}}{0}{\nosha}{2\!\cdot\!3^{2}\!+\!1\!\cdot\!3^{3}\!+\!2\!\cdot\!3^{4}\!+\!O(3^{5})}{\noeuler}\errentry{2\!\cdot\!10^{-7}}{2\!\cdot\!10^{-5}}\cr
\tabentry{10}{[0, -1, 0, -4, 0]}{x^3-10}{2^{2}3\!\cdot\!5^{2}}{2^{8}3^{2}5^{4}17^{2}}{2}{1^2}{1\!\cdot\!3^{0}\phantom{\,\,\!+\!2\!\cdot\!3^1}\,\phantom{\,\,\!+\!2\!\cdot\!3^2}\,\!+\!O(3^{3})}{1\!\cdot\!3^{0}\phantom{\,\,\!+\!2\!\cdot\!3^1}\,\phantom{\,\,\!+\!2\!\cdot\!3^2}\,\!+\!O(3^{4})}\errentry{3\!\cdot\!10^{-17}}{8\!\cdot\!10^{-17}}\cr
\tabentry{12}{[0, -1, 0, -4, 0]}{x^3-12}{2^{2}3^{5}}{2^{8}3^{10}17^{2}}{2^{3}}{2^2}{2\!\cdot\!3^{0}\!+\!1\!\cdot\!3^{1}\phantom{\,\,\!+\!2\!\cdot\!3^2}\,\!+\!O(3^{4})}{2\!\cdot\!3^{0}\!+\!1\!\cdot\!3^{1}\!+\!2\!\cdot\!3^{2}\!+\!O(3^{3})}\errentry{2\!\cdot\!10^{-16}}{2\!\cdot\!10^{-15}}\cr
\tabentry{14}{[0, -1, 0, -4, 0]}{x^3-14}{2^{2}3^{3}7^{2}}{2^{8}3^{6}7^{4}17^{2}}{2\!\cdot\!3^{2}}{3^2}{2\!\cdot\!3^{2}\!+\!1\!\cdot\!3^{3}\!+\!2\!\cdot\!3^{4}\!+\!O(3^{5})}{2\!\cdot\!3^{2}\!+\!2\!\cdot\!3^{3}\!+\!1\!\cdot\!3^{4}\!+\!O(3^{6})}\errentry{6\!\cdot\!10^{-7}}{4\!\cdot\!10^{-7}}\cr
\tabentry{19}{[0, -1, 0, -4, 0]}{x^3-19}{3\!\cdot\!19^{2}}{2^{8}3^{2}17^{2}19^{4}}{2^{4}}{2^2}{2\!\cdot\!3^{0}\phantom{\,\,\!+\!2\!\cdot\!3^1}\,\!+\!1\!\cdot\!3^{2}\!+\!O(3^{3})}{2\!\cdot\!3^{0}\!+\!2\!\cdot\!3^{1}\!+\!1\!\cdot\!3^{2}\!+\!O(3^{4})}\errentry{7\!\cdot\!10^{-8}}{3\!\cdot\!10^{-8}}\cr
\tabentry{20}{[0, -1, 0, -4, 0]}{x^3-20}{2^{2}3^{3}5^{2}}{2^{8}3^{6}5^{4}17^{2}}{2^{5}}{4^2}{1\!\cdot\!3^{0}\phantom{\,\,\!+\!2\!\cdot\!3^1}\,\phantom{\,\,\!+\!2\!\cdot\!3^2}\,\!+\!O(3^{3})}{1\!\cdot\!3^{0}\phantom{\,\,\!+\!2\!\cdot\!3^1}\,\!+\!2\!\cdot\!3^{2}\!+\!O(3^{3})}\errentry{1\!\cdot\!10^{-11}}{1\!\cdot\!10^{-11}}\cr
\tabentry{26}{[0, -1, 0, -4, 0]}{x^3-26}{2^{2}3\!\cdot\!13^{2}}{2^{8}3^{2}13^{4}17^{2}}{2\!\cdot\!3^{2}}{3^2}{2\!\cdot\!3^{2}\!+\!1\!\cdot\!3^{3}\!+\!2\!\cdot\!3^{4}\!+\!O(3^{6})}{2\!\cdot\!3^{2}\!+\!2\!\cdot\!3^{3}\!+\!1\!\cdot\!3^{4}\!+\!O(3^{5})}\errentry{4\!\cdot\!10^{-13}}{3\!\cdot\!10^{-13}}\cr
\tabentry{28}{[0, -1, 0, -4, 0]}{x^3-28}{2^{2}3\!\cdot\!7^{2}}{2^{8}3^{2}7^{4}17^{2}}{2\!\cdot\!3^{2}}{3^2}{2\!\cdot\!3^{2}\!+\!1\!\cdot\!3^{3}\!+\!2\!\cdot\!3^{4}\!+\!O(3^{5})}{2\!\cdot\!3^{2}\phantom{\,\,\!+\!2\!\cdot\!3^3}\,\!+\!2\!\cdot\!3^{4}\!+\!O(3^{5})}\errentry{3\!\cdot\!10^{-17}}{2\!\cdot\!10^{-15}}\cr
\tabentry{44}{[0, -1, 0, -4, 0]}{x^3-44}{2^{2}3\!\cdot\!11^{2}}{2^{8}3^{2}11^{4}17^{2}}{0}{\nosha}{2\!\cdot\!3^{3}\!+\!1\!\cdot\!3^{4}\!+\!1\!\cdot\!3^{5}\!+\!O(3^{6})}{\noeuler}\errentry{2\!\cdot\!10^{-15}}{4\!\cdot\!10^{-9}}\cr
\hline
\tabend

\noindent\begin{tabular}{cllll} 
\extravspace
$E=92A1$ & \multicolumn{3}{l}{$y^2=x^3+x^2+2x+1$} & \tabname{3}{92A1}\cr
\extravspace
$p=3$ & \multicolumn{3}{l}{$P_3(E/\Q,T)=1-T+3T^2$}\cr
$/\vQ\alsub{}$   & $N(E)=2^{2}\!\cdot\!23$ & $L^*(E/\Q)=3^{-1}$ & $|E(\Q)|=3$ & $|\sha(E/\Q)|=1^2$ \cr
$/K\alsub{}$      & $N(E,\sigma)=2^{4}\!\cdot\!3^{2}\!\cdot\!23^{2}$ & $L^*(E/K)=3^{-1}$ & $|E(K)|=3$ & $|\sha(E/K)|=1^2$ \cr
$/\vQ\alsub{2}$ & $\hbox{IV},\>c_{2}=3$ & $/K\alsub{v}\quad(v|2)$ & $\hbox{IV},\>c_v=3$\cr
$/\vQ\alsub{23}$ & $\hbox{I}_1 \hbox{ non-split},\>c_{23}=1$ & $/K\alsub{v}\quad(v|23)$ & $\hbox{I}_1 \hbox{ split},\>c_v=1$\cr
\end{tabular}
\tabbegin
\tabentry{3}{[0, 1, 0, 2, 1]}{x^3-3}{3^{5}}{2^{4}3^{10}23^{2}}{3}{1^2}{1\!\cdot\!3^{1}\!+\!1\!\cdot\!3^{2}\!+\!1\!\cdot\!3^{3}\!+\!O(3^{5})}{1\!\cdot\!3^{1}\!+\!1\!\cdot\!3^{2}\!+\!1\!\cdot\!3^{3}\!+\!O(3^{5})}\errentry{1\!\cdot\!10^{-17}}{1\!\cdot\!10^{-15}}\cr
\tabentry{5}{[0, 1, 0, 2, 1]}{x^3-5}{3^{3}5^{2}}{2^{4}3^{6}5^{4}23^{2}}{0}{\nosha}{1\!\cdot\!3^{3}\phantom{\,\,\!+\!2\!\cdot\!3^4}\,\phantom{\,\,\!+\!2\!\cdot\!3^5}\,\!+\!O(3^{6})}{\noeuler}\errentry{4\!\cdot\!10^{-17}}{1\!\cdot\!10^{-9}}\cr
\tabentry{7}{[0, 1, 0, 2, 1]}{x^3-7}{3^{3}7^{2}}{2^{4}3^{6}7^{4}23^{2}}{0}{\nosha}{1\!\cdot\!3^{3}\!+\!1\!\cdot\!3^{4}\!+\!2\!\cdot\!3^{5}\!+\!O(3^{6})}{\noeuler}\errentry{1\!\cdot\!10^{-12}}{1\!\cdot\!10^{-8}}\cr
\tabentry{11}{[0, 1, 0, 2, 1]}{x^3-11}{3^{3}11^{2}}{2^{4}3^{6}11^{4}23^{2}}{3^{3}}{3^2}{1\!\cdot\!3^{3}\!+\!2\!\cdot\!3^{4}\!+\!1\!\cdot\!3^{5}\!+\!O(3^{6})}{2\!\cdot\!3^{3}\!+\!2\!\cdot\!3^{4}\phantom{\,\,\!+\!2\!\cdot\!3^5}\,\!+\!O(3^{7})}\errentry{4\!\cdot\!10^{-8}}{3\!\cdot\!10^{-8}}\cr
\tabentry{17}{[0, 1, 0, 2, 1]}{x^3-17}{3\!\cdot\!17^{2}}{2^{4}3^{2}17^{4}23^{2}}{3}{1^2}{2\!\cdot\!3^{3}\phantom{\,\,\!+\!2\!\cdot\!3^4}\,\phantom{\,\,\!+\!2\!\cdot\!3^5}\,\!+\!O(3^{9})}{2\!\cdot\!3^{3}\!+\!1\!\cdot\!3^{4}\!+\!2\!\cdot\!3^{5}\!+\!O(3^{6})}\errentry{5\!\cdot\!10^{-17}}{2\!\cdot\!10^{-15}}\cr
\tabentry{19}{[0, 1, 0, 2, 1]}{x^3-19}{3\!\cdot\!19^{2}}{2^{4}3^{2}19^{4}23^{2}}{0}{\nosha}{1\!\cdot\!3^{5}\!+\!2\!\cdot\!3^{6}\!+\!1\!\cdot\!3^{7}\!+\!O(3^{8})}{\noeuler}\errentry{1\!\cdot\!10^{-14}}{3\!\cdot\!10^{-9}}\cr
\tabentry{35}{[0, 1, 0, 2, 1]}{x^3-35}{3\!\cdot\!5^{2}7^{2}}{2^{4}3^{2}5^{4}7^{4}23^{2}}{3^{3}}{3^2}{1\!\cdot\!3^{5}\phantom{\,\,\!+\!2\!\cdot\!3^6}\,\!+\!1\!\cdot\!3^{7}\!+\!O(3^{8})}{2\!\cdot\!3^{5}\!+\!1\!\cdot\!3^{6}\phantom{\,\,\!+\!2\!\cdot\!3^7}\,\!+\!O(3^{9})}\errentry{1\!\cdot\!10^{-8}}{1\!\cdot\!10^{-8}}\cr
\tabentry{46}{[0, 1, 0, 2, 1]}{x^3-46}{2^{2}3\!\cdot\!23^{2}}{2^{2}3^{2}23^{4}}{2\!\cdot\!3^{-1}}{1^2}{2\!\cdot\!3^{2}\!+\!2\!\cdot\!3^{3}\phantom{\,\,\!+\!2\!\cdot\!3^4}\,\!+\!O(3^{6})}{1\!\cdot\!3^{2}\!+\!1\!\cdot\!3^{3}\!+\!2\!\cdot\!3^{4}\!+\!O(3^{5})}\errentry{5\!\cdot\!10^{-17}}{4\!\cdot\!10^{-17}}\cr
\tabentry{92}{[0, 1, 0, 2, 1]}{x^3-92}{2^{2}3^{3}23^{2}}{2^{2}3^{6}23^{4}}{2\!\cdot\!3}{3^2}{2\!\cdot\!3^{2}\!+\!2\!\cdot\!3^{3}\phantom{\,\,\!+\!2\!\cdot\!3^4}\,\!+\!O(3^{6})}{1\!\cdot\!3^{2}\!+\!2\!\cdot\!3^{3}\!+\!1\!\cdot\!3^{4}\!+\!O(3^{7})}\errentry{2\!\cdot\!10^{-16}}{3\!\cdot\!10^{-15}}\cr
\tabentry{138}{[0, 1, 0, 2, 1]}{x^3-138}{2^{2}3^{5}23^{2}}{2^{2}3^{10}23^{4}}{2\!\cdot\!3}{3^2}{2\!\cdot\!3^{2}\!+\!2\!\cdot\!3^{3}\phantom{\,\,\!+\!2\!\cdot\!3^4}\,\!+\!O(3^{6})}{1\!\cdot\!3^{2}\!+\!1\!\cdot\!3^{3}\!+\!2\!\cdot\!3^{4}\!+\!O(3^{5})}\errentry{1\!\cdot\!10^{-8}}{1\!\cdot\!10^{-8}}\cr
\tabentry{276}{[0, 1, 0, 2, 1]}{x^3-276}{2^{2}3^{5}23^{2}}{2^{2}3^{10}23^{4}}{2\!\cdot\!3}{3^2}{2\!\cdot\!3^{2}\!+\!2\!\cdot\!3^{3}\phantom{\,\,\!+\!2\!\cdot\!3^4}\,\!+\!O(3^{6})}{1\!\cdot\!3^{2}\!+\!1\!\cdot\!3^{3}\!+\!2\!\cdot\!3^{4}\!+\!O(3^{5})}\errentry{6\!\cdot\!10^{-9}}{4\!\cdot\!10^{-9}}\cr
\tabentry{414}{[0, 1, 0, 2, 1]}{x^3-414}{2^{2}3^{5}23^{2}}{2^{2}3^{10}23^{4}}{2^{5}3}{12^2}{2\!\cdot\!3^{2}\!+\!2\!\cdot\!3^{3}\phantom{\,\,\!+\!2\!\cdot\!3^4}\,\!+\!O(3^{6})}{1\!\cdot\!3^{2}\phantom{\,\,\!+\!2\!\cdot\!3^3}\,\phantom{\,\,\!+\!2\!\cdot\!3^4}\,\!+\!O(3^{7})}\errentry{1\!\cdot\!10^{-6}}{4\!\cdot\!10^{-6}}\cr
\tabentry{828}{[0, 1, 0, 2, 1]}{x^3-828}{2^{2}3^{5}23^{2}}{2^{2}3^{10}23^{4}}{2^{3}3}{6^2}{2\!\cdot\!3^{2}\!+\!2\!\cdot\!3^{3}\phantom{\,\,\!+\!2\!\cdot\!3^4}\,\!+\!O(3^{6})}{1\!\cdot\!3^{2}\!+\!2\!\cdot\!3^{3}\phantom{\,\,\!+\!2\!\cdot\!3^4}\,\!+\!O(3^{5})}\errentry{8\!\cdot\!10^{-9}}{1\!\cdot\!10^{-8}}\cr
\hline
\tabend

\noindent\begin{tabular}{cllll} 
\extravspace
$E=116C2$ & \multicolumn{3}{l}{$y^2=x^3-x^2-9x+14$} & \tabname{3}{116C2}\cr
\extravspace
$p=3$ & \multicolumn{3}{l}{$P_3(E/\Q,T)=1-2T+3T^2$}\cr
$/\vQ\alsub{}$   & $N(E)=2^{2}\!\cdot\!29$ & $L^*(E/\Q)=2^{-2}$ & $|E(\Q)|=2$ & $|\sha(E/\Q)|=1^2$ \cr
$/K\alsub{}$      & $N(E,\sigma)=2^{4}\!\cdot\!3^{2}\!\cdot\!29^{2}$ & $L^*(E/K)=2^{-2}\!\cdot\!3$ & $|E(K)|=2$ & $|\sha(E/K)|=1^2$ \cr
$/\vQ\alsub{2}$ & $\hbox{IV},\>c_{2}=1$ & $/K\alsub{v}\quad(v|2)$ & $\hbox{IV},\>c_v=3$\cr
$/\vQ\alsub{29}$ & $\hbox{I}_1 \hbox{ non-split},\>c_{29}=1$ & $/K\alsub{v}\quad(v|29)$ & $\hbox{I}_1 \hbox{ split},\>c_v=1$\cr
\end{tabular}
\tabbegin
\tabentry{3}{[0, -1, 0, -9, 14]}{x^3-3}{3^{5}}{2^{4}3^{10}29^{2}}{3}{1^2}{2\!\cdot\!3^{1}\!+\!2\!\cdot\!3^{2}\phantom{\,\,\!+\!2\!\cdot\!3^3}\,\!+\!O(3^{4})}{2\!\cdot\!3^{1}\!+\!2\!\cdot\!3^{2}\!+\!2\!\cdot\!3^{3}\!+\!O(3^{4})}\errentry{1\!\cdot\!10^{-16}}{2\!\cdot\!10^{-15}}\cr
\tabentry{5}{[0, -1, 0, -9, 14]}{x^3-5}{3^{3}5^{2}}{2^{4}3^{6}5^{4}29^{2}}{3}{1^2}{1\!\cdot\!3^{1}\!+\!2\!\cdot\!3^{2}\phantom{\,\,\!+\!2\!\cdot\!3^3}\,\!+\!O(3^{4})}{1\!\cdot\!3^{1}\!+\!1\!\cdot\!3^{2}\!+\!1\!\cdot\!3^{3}\!+\!O(3^{4})}\errentry{7\!\cdot\!10^{-16}}{4\!\cdot\!10^{-16}}\cr
\tabentry{7}{[0, -1, 0, -9, 14]}{x^3-7}{3^{3}7^{2}}{2^{4}3^{6}7^{4}29^{2}}{3}{1^2}{2\!\cdot\!3^{1}\!+\!1\!\cdot\!3^{2}\!+\!2\!\cdot\!3^{3}\!+\!O(3^{4})}{2\!\cdot\!3^{1}\!+\!2\!\cdot\!3^{2}\!+\!1\!\cdot\!3^{3}\!+\!O(3^{5})}\errentry{3\!\cdot\!10^{-12}}{7\!\cdot\!10^{-13}}\cr
\tabentry{17}{[0, -1, 0, -9, 14]}{x^3-17}{3\!\cdot\!17^{2}}{2^{4}3^{2}17^{4}29^{2}}{3}{1^2}{1\!\cdot\!3^{1}\!+\!1\!\cdot\!3^{2}\phantom{\,\,\!+\!2\!\cdot\!3^3}\,\!+\!O(3^{4})}{1\!\cdot\!3^{1}\!+\!1\!\cdot\!3^{2}\!+\!1\!\cdot\!3^{3}\!+\!O(3^{4})}\errentry{4\!\cdot\!10^{-15}}{1\!\cdot\!10^{-15}}\cr
\tabentry{19}{[0, -1, 0, -9, 14]}{x^3-19}{3\!\cdot\!19^{2}}{2^{4}3^{2}19^{4}29^{2}}{2^{2}3}{2^2}{2\!\cdot\!3^{1}\!+\!2\!\cdot\!3^{2}\!+\!1\!\cdot\!3^{3}\!+\!O(3^{4})}{2\!\cdot\!3^{1}\phantom{\,\,\!+\!2\!\cdot\!3^2}\,\!+\!1\!\cdot\!3^{3}\!+\!O(3^{4})}\errentry{5\!\cdot\!10^{-13}}{2\!\cdot\!10^{-13}}\cr
\tabentry{58}{[0, -1, 0, -9, 14]}{x^3-58}{2^{2}3^{3}29^{2}}{2^{2}3^{6}29^{4}}{0}{\nosha}{2\!\cdot\!3^{2}\phantom{\,\,\!+\!2\!\cdot\!3^3}\,\!+\!2\!\cdot\!3^{4}\!+\!O(3^{5})}{\noeuler}\errentry{2\!\cdot\!10^{-17}}{1\!\cdot\!10^{-9}}\cr
\tabentry{406}{[0, -1, 0, -9, 14]}{x^3-406}{2^{2}3\!\cdot\!7^{2}29^{2}}{2^{2}3^{2}7^{4}29^{4}}{2\!\cdot\!3^{4}}{9^2}{2\!\cdot\!3^{2}\!+\!2\!\cdot\!3^{3}\!+\!1\!\cdot\!3^{4}\!+\!O(3^{5})}{2\!\cdot\!3^{5}\!+\!1\!\cdot\!3^{6}\phantom{\,\,\!+\!2\!\cdot\!3^7}\,\!+\!O(3^{10})}\errentry{5\!\cdot\!10^{-9}}{1\!\cdot\!10^{-8}}\cr
\tabentry{1450}{[0, -1, 0, -9, 14]}{x^3-1450}{2^{2}3\!\cdot\!5^{2}29^{2}}{2^{2}3^{2}5^{4}29^{4}}{2^{3}3^{2}}{6^2}{1\!\cdot\!3^{2}\!+\!1\!\cdot\!3^{3}\phantom{\,\,\!+\!2\!\cdot\!3^4}\,\!+\!O(3^{6})}{1\!\cdot\!3^{3}\phantom{\,\,\!+\!2\!\cdot\!3^4}\,\phantom{\,\,\!+\!2\!\cdot\!3^5}\,\!+\!O(3^{6})}\errentry{1\!\cdot\!10^{-14}}{2\!\cdot\!10^{-14}}\cr
\hline
\tabend

\noindent\begin{tabular}{cllll} 
\extravspace
$E=152B1$ & \multicolumn{3}{l}{$y^2=x^3+x^2-8x-16$} & \tabname{3}{152B1}\cr
\extravspace
$p=3$ & \multicolumn{3}{l}{$P_3(E/\Q,T)=1-T+3T^2$}\cr
$/\vQ\alsub{}$   & $N(E)=2^{3}\!\cdot\!19$ & $L^*(E/\Q)=1$ & $|E(\Q)|=1$ & $|\sha(E/\Q)|=1^2$ \cr
$/K\alsub{}$      & $N(E,\sigma)=2^{6}\!\cdot\!3^{2}\!\cdot\!19^{2}$ & $L^*(E/K)=1$ & $|E(K)|=1$ & $|\sha(E/K)|=1^2$ \cr
$/\vQ\alsub{2}$ & $\hbox{II}^*,\>c_{2}=1$ & $/K\alsub{v}\quad(v|2)$ & $\hbox{II}^*,\>c_v=1$\cr
$/\vQ\alsub{19}$ & $\hbox{I}_1 \hbox{ split},\>c_{19}=1$ & $/K\alsub{v}\quad(v|19)$ & $\hbox{I}_1 \hbox{ split},\>c_v=1$\cr
\end{tabular}
\tabbegin
\tabentry{2}{[0, 1, 0, -8, -16]}{x^3-2}{2^{2}3^{3}}{2^{6}3^{6}19^{2}}{0}{\nosha}{1\!\cdot\!3^{2}\!+\!1\!\cdot\!3^{3}\!+\!1\!\cdot\!3^{4}\!+\!O(3^{6})}{\noeuler}\errentry{5\!\cdot\!10^{-17}}{6\!\cdot\!10^{-10}}\cr
\tabentry{6}{[0, 1, 0, -8, -16]}{x^3-6}{2^{2}3^{5}}{2^{6}3^{10}19^{2}}{0}{\nosha}{1\!\cdot\!3^{2}\!+\!1\!\cdot\!3^{3}\!+\!1\!\cdot\!3^{4}\!+\!O(3^{6})}{\noeuler}\errentry{8\!\cdot\!10^{-18}}{1\!\cdot\!10^{-9}}\cr
\tabentry{10}{[0, 1, 0, -8, -16]}{x^3-10}{2^{2}3\!\cdot\!5^{2}}{2^{6}3^{2}5^{4}19^{2}}{0}{\nosha}{1\!\cdot\!3^{4}\phantom{\,\,\!+\!2\!\cdot\!3^5}\,\phantom{\,\,\!+\!2\!\cdot\!3^6}\,\!+\!O(3^{7})}{\noeuler}\errentry{5\!\cdot\!10^{-17}}{9\!\cdot\!10^{-10}}\cr
\tabentry{12}{[0, 1, 0, -8, -16]}{x^3-12}{2^{2}3^{5}}{2^{6}3^{10}19^{2}}{0}{\nosha}{1\!\cdot\!3^{2}\!+\!1\!\cdot\!3^{3}\!+\!1\!\cdot\!3^{4}\!+\!O(3^{6})}{\noeuler}\errentry{1\!\cdot\!10^{-16}}{1\!\cdot\!10^{-9}}\cr
\tabentry{14}{[0, 1, 0, -8, -16]}{x^3-14}{2^{2}3^{3}7^{2}}{2^{6}3^{6}7^{4}19^{2}}{0}{\nosha}{1\!\cdot\!3^{2}\!+\!2\!\cdot\!3^{3}\!+\!1\!\cdot\!3^{4}\!+\!O(3^{5})}{\noeuler}\errentry{1\!\cdot\!10^{-9}}{8\!\cdot\!10^{-7}}\cr
\tabentry{20}{[0, 1, 0, -8, -16]}{x^3-20}{2^{2}3^{3}5^{2}}{2^{6}3^{6}5^{4}19^{2}}{2^{3}3^{2}}{3^2}{1\!\cdot\!3^{4}\phantom{\,\,\!+\!2\!\cdot\!3^5}\,\phantom{\,\,\!+\!2\!\cdot\!3^6}\,\!+\!O(3^{7})}{2\!\cdot\!3^{2}\phantom{\,\,\!+\!2\!\cdot\!3^3}\,\!+\!1\!\cdot\!3^{4}\!+\!O(3^{5})}\errentry{8\!\cdot\!10^{-15}}{4\!\cdot\!10^{-15}}\cr
\tabentry{26}{[0, 1, 0, -8, -16]}{x^3-26}{2^{2}3\!\cdot\!13^{2}}{2^{6}3^{2}13^{4}19^{2}}{2^{3}}{1^2}{1\!\cdot\!3^{2}\!+\!1\!\cdot\!3^{3}\!+\!1\!\cdot\!3^{4}\!+\!O(3^{6})}{1\!\cdot\!3^{2}\!+\!1\!\cdot\!3^{3}\phantom{\,\,\!+\!2\!\cdot\!3^4}\,\!+\!O(3^{6})}\errentry{8\!\cdot\!10^{-17}}{6\!\cdot\!10^{-17}}\cr
\tabentry{28}{[0, 1, 0, -8, -16]}{x^3-28}{2^{2}3\!\cdot\!7^{2}}{2^{6}3^{2}7^{4}19^{2}}{2^{3}}{1^2}{1\!\cdot\!3^{2}\!+\!2\!\cdot\!3^{3}\!+\!1\!\cdot\!3^{4}\!+\!O(3^{5})}{1\!\cdot\!3^{2}\phantom{\,\,\!+\!2\!\cdot\!3^3}\,\!+\!2\!\cdot\!3^{4}\!+\!O(3^{7})}\errentry{1\!\cdot\!10^{-15}}{2\!\cdot\!10^{-15}}\cr
\tabentry{38}{[0, 1, 0, -8, -16]}{x^3-38}{2^{2}3^{3}19^{2}}{2^{6}3^{6}19^{4}}{0}{\nosha}{1\!\cdot\!3^{6}\!+\!2\!\cdot\!3^{7}\!+\!1\!\cdot\!3^{8}\!+\!O(3^{9})}{\noeuler}\errentry{1\!\cdot\!10^{-16}}{1\!\cdot\!10^{-9}}\cr
\tabentry{44}{[0, 1, 0, -8, -16]}{x^3-44}{2^{2}3\!\cdot\!11^{2}}{2^{6}3^{2}11^{4}19^{2}}{2^{3}}{1^2}{2\!\cdot\!3^{2}\!+\!1\!\cdot\!3^{3}\phantom{\,\,\!+\!2\!\cdot\!3^4}\,\!+\!O(3^{6})}{2\!\cdot\!3^{2}\!+\!2\!\cdot\!3^{3}\!+\!1\!\cdot\!3^{4}\!+\!O(3^{5})}\errentry{6\!\cdot\!10^{-17}}{6\!\cdot\!10^{-17}}\cr
\tabentry{46}{[0, 1, 0, -8, -16]}{x^3-46}{2^{2}3\!\cdot\!23^{2}}{2^{6}3^{2}19^{2}23^{4}}{2^{3}}{1^2}{2\!\cdot\!3^{2}\phantom{\,\,\!+\!2\!\cdot\!3^3}\,\!+\!2\!\cdot\!3^{4}\!+\!O(3^{5})}{2\!\cdot\!3^{2}\!+\!1\!\cdot\!3^{3}\!+\!1\!\cdot\!3^{4}\!+\!O(3^{5})}\errentry{2\!\cdot\!10^{-9}}{7\!\cdot\!10^{-10}}\cr
\tabentry{76}{[0, 1, 0, -8, -16]}{x^3-76}{2^{2}3^{3}19^{2}}{2^{6}3^{6}19^{4}}{2^{3}3^{3}}{3^2}{1\!\cdot\!3^{6}\!+\!2\!\cdot\!3^{7}\!+\!1\!\cdot\!3^{8}\!+\!O(3^{9})}{1\!\cdot\!3^{3}\phantom{\,\,\!+\!2\!\cdot\!3^4}\,\!+\!1\!\cdot\!3^{5}\!+\!O(3^{7})}\errentry{1\!\cdot\!10^{-17}}{2\!\cdot\!10^{-16}}\cr
\tabentry{190}{[0, 1, 0, -8, -16]}{x^3-190}{2^{2}3\!\cdot\!5^{2}19^{2}}{2^{6}3^{2}5^{4}19^{4}}{0}{\nosha}{1\!\cdot\!3^{8}\!+\!1\!\cdot\!3^{9}\!+\!2\!\cdot\!3^{10}\!+\!O(3^{11})}{\noeuler}\errentry{1\!\cdot\!10^{-10}}{2\!\cdot\!10^{-6}}\cr
\tabentry{1900}{[0, 1, 0, -8, -16]}{x^3-1900}{2^{2}3\!\cdot\!5^{2}19^{2}}{2^{6}3^{2}5^{4}19^{4}}{0}{\nosha}{1\!\cdot\!3^{8}\!+\!1\!\cdot\!3^{9}\!+\!2\!\cdot\!3^{10}\!+\!O(3^{11})}{\noeuler}\errentry{3\!\cdot\!10^{-9}}{3\!\cdot\!10^{-6}}\cr
\hline
\tabend

\noindent\begin{tabular}{cllll} 
\extravspace
$E=176B1$ & \multicolumn{3}{l}{$y^2=x^3+x^2-5x-13$} & \tabname{3}{176B1}\cr
\extravspace
$p=3$ & \multicolumn{3}{l}{$P_3(E/\Q,T)=1-T+3T^2$}\cr
$/\vQ\alsub{}$   & $N(E)=2^{4}\!\cdot\!11$ & $L^*(E/\Q)=1$ & $|E(\Q)|=1$ & $|\sha(E/\Q)|=1^2$ \cr
$/K\alsub{}$      & $N(E,\sigma)=2^{8}\!\cdot\!3^{2}\!\cdot\!11^{2}$ & $L^*(E/K)=1$ & $|E(K)|=1$ & $|\sha(E/K)|=1^2$ \cr
$/\vQ\alsub{2}$ & $\hbox{II}^*,\>c_{2}=1$ & $/K\alsub{v}\quad(v|2)$ & $\hbox{II}^*,\>c_v=1$\cr
$/\vQ\alsub{11}$ & $\hbox{I}_1 \hbox{ non-split},\>c_{11}=1$ & $/K\alsub{v}\quad(v|11)$ & $\hbox{I}_1 \hbox{ split},\>c_v=1$\cr
\end{tabular}
\tabbegin
\tabentry{2}{[0, 1, 0, -5, -13]}{x^3-2}{2^{2}3^{3}}{2^{8}3^{6}11^{2}}{0}{\nosha}{1\!\cdot\!3^{2}\!+\!1\!\cdot\!3^{3}\!+\!1\!\cdot\!3^{4}\!+\!O(3^{6})}{\noeuler}\errentry{7\!\cdot\!10^{-17}}{1\!\cdot\!10^{-10}}\cr
\tabentry{3}{[0, 1, 0, -5, -13]}{x^3-3}{3^{5}}{2^{8}3^{10}11^{2}}{3^{2}}{3^2}{1\!\cdot\!3^{2}\!+\!1\!\cdot\!3^{3}\!+\!1\!\cdot\!3^{4}\!+\!O(3^{6})}{1\!\cdot\!3^{2}\!+\!1\!\cdot\!3^{3}\!+\!1\!\cdot\!3^{4}\!+\!O(3^{6})}\errentry{3\!\cdot\!10^{-17}}{1\!\cdot\!10^{-15}}\cr
\tabentry{5}{[0, 1, 0, -5, -13]}{x^3-5}{3^{3}5^{2}}{2^{8}3^{6}5^{4}11^{2}}{0}{\nosha}{2\!\cdot\!3^{2}\phantom{\,\,\!+\!2\!\cdot\!3^3}\,\phantom{\,\,\!+\!2\!\cdot\!3^4}\,\!+\!O(3^{6})}{\noeuler}\errentry{6\!\cdot\!10^{-14}}{2\!\cdot\!10^{-8}}\cr
\tabentry{6}{[0, 1, 0, -5, -13]}{x^3-6}{2^{2}3^{5}}{2^{8}3^{10}11^{2}}{0}{\nosha}{1\!\cdot\!3^{2}\!+\!1\!\cdot\!3^{3}\!+\!1\!\cdot\!3^{4}\!+\!O(3^{6})}{\noeuler}\errentry{4\!\cdot\!10^{-17}}{8\!\cdot\!10^{-10}}\cr
\tabentry{7}{[0, 1, 0, -5, -13]}{x^3-7}{3^{3}7^{2}}{2^{8}3^{6}7^{4}11^{2}}{3^{2}}{3^2}{1\!\cdot\!3^{4}\!+\!1\!\cdot\!3^{5}\!+\!2\!\cdot\!3^{6}\!+\!O(3^{7})}{1\!\cdot\!3^{2}\phantom{\,\,\!+\!2\!\cdot\!3^3}\,\phantom{\,\,\!+\!2\!\cdot\!3^4}\,\!+\!O(3^{7})}\errentry{2\!\cdot\!10^{-10}}{1\!\cdot\!10^{-10}}\cr
\tabentry{10}{[0, 1, 0, -5, -13]}{x^3-10}{2^{2}3\!\cdot\!5^{2}}{2^{8}3^{2}5^{4}11^{2}}{2^{3}}{1^2}{2\!\cdot\!3^{2}\phantom{\,\,\!+\!2\!\cdot\!3^3}\,\phantom{\,\,\!+\!2\!\cdot\!3^4}\,\!+\!O(3^{6})}{2\!\cdot\!3^{2}\!+\!1\!\cdot\!3^{3}\phantom{\,\,\!+\!2\!\cdot\!3^4}\,\!+\!O(3^{5})}\errentry{1\!\cdot\!10^{-16}}{1\!\cdot\!10^{-15}}\cr
\tabentry{12}{[0, 1, 0, -5, -13]}{x^3-12}{2^{2}3^{5}}{2^{8}3^{10}11^{2}}{0}{\nosha}{1\!\cdot\!3^{2}\!+\!1\!\cdot\!3^{3}\!+\!1\!\cdot\!3^{4}\!+\!O(3^{6})}{\noeuler}\errentry{4\!\cdot\!10^{-17}}{1\!\cdot\!10^{-9}}\cr
\tabentry{14}{[0, 1, 0, -5, -13]}{x^3-14}{2^{2}3^{3}7^{2}}{2^{8}3^{6}7^{4}11^{2}}{2^{3}3^{2}}{3^2}{1\!\cdot\!3^{4}\!+\!1\!\cdot\!3^{5}\!+\!2\!\cdot\!3^{6}\!+\!O(3^{7})}{1\!\cdot\!3^{2}\!+\!1\!\cdot\!3^{3}\phantom{\,\,\!+\!2\!\cdot\!3^4}\,\!+\!O(3^{7})}\errentry{3\!\cdot\!10^{-10}}{2\!\cdot\!10^{-10}}\cr
\tabentry{17}{[0, 1, 0, -5, -13]}{x^3-17}{3\!\cdot\!17^{2}}{2^{8}3^{2}11^{2}17^{4}}{1}{1^2}{2\!\cdot\!3^{2}\phantom{\,\,\!+\!2\!\cdot\!3^3}\,\!+\!2\!\cdot\!3^{4}\!+\!O(3^{8})}{2\!\cdot\!3^{2}\!+\!1\!\cdot\!3^{3}\!+\!2\!\cdot\!3^{4}\!+\!O(3^{5})}\errentry{8\!\cdot\!10^{-12}}{2\!\cdot\!10^{-12}}\cr
\tabentry{19}{[0, 1, 0, -5, -13]}{x^3-19}{3\!\cdot\!19^{2}}{2^{8}3^{2}11^{2}19^{4}}{1}{1^2}{1\!\cdot\!3^{2}\!+\!2\!\cdot\!3^{3}\phantom{\,\,\!+\!2\!\cdot\!3^4}\,\!+\!O(3^{6})}{1\!\cdot\!3^{2}\!+\!1\!\cdot\!3^{3}\!+\!2\!\cdot\!3^{4}\!+\!O(3^{6})}\errentry{8\!\cdot\!10^{-10}}{2\!\cdot\!10^{-10}}\cr
\tabentry{20}{[0, 1, 0, -5, -13]}{x^3-20}{2^{2}3^{3}5^{2}}{2^{8}3^{6}5^{4}11^{2}}{0}{\nosha}{2\!\cdot\!3^{2}\phantom{\,\,\!+\!2\!\cdot\!3^3}\,\phantom{\,\,\!+\!2\!\cdot\!3^4}\,\!+\!O(3^{6})}{\noeuler}\errentry{2\!\cdot\!10^{-14}}{9\!\cdot\!10^{-9}}\cr
\tabentry{26}{[0, 1, 0, -5, -13]}{x^3-26}{2^{2}3\!\cdot\!13^{2}}{2^{8}3^{2}11^{2}13^{4}}{2^{3}}{1^2}{1\!\cdot\!3^{2}\!+\!2\!\cdot\!3^{3}\!+\!1\!\cdot\!3^{4}\!+\!O(3^{5})}{1\!\cdot\!3^{2}\!+\!1\!\cdot\!3^{3}\phantom{\,\,\!+\!2\!\cdot\!3^4}\,\!+\!O(3^{6})}\errentry{2\!\cdot\!10^{-15}}{5\!\cdot\!10^{-16}}\cr
\tabentry{28}{[0, 1, 0, -5, -13]}{x^3-28}{2^{2}3\!\cdot\!7^{2}}{2^{8}3^{2}7^{4}11^{2}}{0}{\nosha}{1\!\cdot\!3^{4}\!+\!1\!\cdot\!3^{5}\!+\!2\!\cdot\!3^{6}\!+\!O(3^{7})}{\noeuler}\errentry{1\!\cdot\!10^{-16}}{1\!\cdot\!10^{-9}}\cr
\tabentry{46}{[0, 1, 0, -5, -13]}{x^3-46}{2^{2}3\!\cdot\!23^{2}}{2^{8}3^{2}11^{2}23^{4}}{2^{3}}{1^2}{2\!\cdot\!3^{2}\phantom{\,\,\!+\!2\!\cdot\!3^3}\,\!+\!2\!\cdot\!3^{4}\!+\!O(3^{5})}{2\!\cdot\!3^{2}\!+\!1\!\cdot\!3^{3}\!+\!1\!\cdot\!3^{4}\!+\!O(3^{5})}\errentry{2\!\cdot\!10^{-8}}{5\!\cdot\!10^{-9}}\cr
\hline
\tabend

\noindent\begin{tabular}{cllll} 
\extravspace
$E=395B1$ & \multicolumn{3}{l}{$y^2+xy+y=x^3+x^2-40x-128$} & \tabname{3}{395B1}\cr
\extravspace
$p=3$ & \multicolumn{3}{l}{$P_3(E/\Q,T)=1-2T+3T^2$}\cr
$/\vQ\alsub{}$   & $N(E)=5\!\cdot\!79$ & $L^*(E/\Q)=2^{-1}\!\cdot\!3$ & $|E(\Q)|=2$ & $|\sha(E/\Q)|=1^2$ \cr
$/K\alsub{}$      & $N(E,\sigma)=3^{2}\!\cdot\!5^{2}\!\cdot\!79^{2}$ & $L^*(E/K)=2^{-1}\!\cdot\!3$ & $|E(K)|=2$ & $|\sha(E/K)|=1^2$ \cr
$/\vQ\alsub{5}$ & $\hbox{I}_6 \hbox{ split},\>c_{5}=6$ & $/K\alsub{v}\quad(v|5)$ & $\hbox{I}_6 \hbox{ split},\>c_v=6$\cr
$/\vQ\alsub{79}$ & $\hbox{I}_1 \hbox{ non-split},\>c_{79}=1$ & $/K\alsub{v}\quad(v|79)$ & $\hbox{I}_1 \hbox{ non-split},\>c_v=1$\cr
\end{tabular}
\tabbegin
\tabentry{2}{[1, 1, 1, -40, -128]}{x^3-2}{2^{2}3^{3}}{2^{4}3^{6}5^{2}79^{2}}{2^{3}3}{2^2}{2\!\cdot\!3^{1}\!+\!1\!\cdot\!3^{2}\phantom{\,\,\!+\!2\!\cdot\!3^3}\,\!+\!O(3^{5})}{2\!\cdot\!3^{1}\!+\!2\!\cdot\!3^{2}\!+\!1\!\cdot\!3^{3}\!+\!O(3^{4})}\errentry{5\!\cdot\!10^{-17}}{1\!\cdot\!10^{-16}}\cr
\tabentry{3}{[1, 1, 1, -40, -128]}{x^3-3}{3^{5}}{3^{10}5^{2}79^{2}}{2\!\cdot\!3}{1^2}{1\!\cdot\!3^{1}\!+\!2\!\cdot\!3^{2}\!+\!1\!\cdot\!3^{3}\!+\!O(3^{4})}{1\!\cdot\!3^{1}\!+\!2\!\cdot\!3^{2}\!+\!2\!\cdot\!3^{3}\!+\!O(3^{4})}\errentry{3\!\cdot\!10^{-15}}{8\!\cdot\!10^{-16}}\cr
\tabentry{17}{[1, 1, 1, -40, -128]}{x^3-17}{3\!\cdot\!17^{2}}{3^{2}5^{2}17^{4}79^{2}}{0}{\nosha}{1\!\cdot\!3^{5}\phantom{\,\,\!+\!2\!\cdot\!3^6}\,\!+\!2\!\cdot\!3^{7}\!+\!O(3^{8})}{\noeuler}\errentry{3\!\cdot\!10^{-6}}{3\!\cdot\!10^{-5}}\cr
\tabentry{28}{[1, 1, 1, -40, -128]}{x^3-28}{2^{2}3\!\cdot\!7^{2}}{2^{4}3^{2}5^{2}7^{4}79^{2}}{2\!\cdot\!3^{3}}{3^2}{2\!\cdot\!3^{3}\!+\!1\!\cdot\!3^{4}\!+\!2\!\cdot\!3^{5}\!+\!O(3^{6})}{2\!\cdot\!3^{3}\phantom{\,\,\!+\!2\!\cdot\!3^4}\,\!+\!2\!\cdot\!3^{5}\!+\!O(3^{6})}\errentry{8\!\cdot\!10^{-10}}{6\!\cdot\!10^{-10}}\cr
\hline
\tabend

\noindent\begin{tabular}{cllll} 
\extravspace
$E=800E1$ & \multicolumn{3}{l}{$y^2=x^3+x^2-208x-1412$} & \tabname{3}{800E1}\cr
\extravspace
$p=3$ & \multicolumn{3}{l}{$P_3(E/\Q,T)=1-T+3T^2$}\cr
$/\vQ\alsub{}$   & $N(E)=2^{5}\!\cdot\!5^{2}$ & $L^*(E/\Q)=3$ & $|E(\Q)|=1$ & $|\sha(E/\Q)|=1^2$ \cr
$/K\alsub{}$      & $N(E,\sigma)=2^{10}\!\cdot\!3^{2}\!\cdot\!5^{4}$ & $L^*(E/K)=3$ & $|E(K)|=1$ & $|\sha(E/K)|=1^2$ \cr
$/\vQ\alsub{2}$ & $\hbox{I}^*_0,\>c_{2}=1$ & $/K\alsub{v}\quad(v|2)$ & $\hbox{I}^*_0,\>c_v=1$\cr
$/\vQ\alsub{5}$ & $\hbox{IV}^*,\>c_{5}=3$ & $/K\alsub{v}\quad(v|5)$ & $\hbox{IV}^*,\>c_v=3$\cr
\end{tabular}
\tabbegin
\tabentry{2}{[0, 1, 0, -208, -1412]}{x^3-2}{2^{2}3^{3}}{2^{10}3^{6}5^{4}}{2\!\cdot\!3^{3}}{3^2}{1\!\cdot\!3^{3}\!+\!1\!\cdot\!3^{4}\!+\!1\!\cdot\!3^{5}\!+\!O(3^{7})}{1\!\cdot\!3^{3}\!+\!2\!\cdot\!3^{4}\phantom{\,\,\!+\!2\!\cdot\!3^5}\,\!+\!O(3^{8})}\errentry{2\!\cdot\!10^{-17}}{3\!\cdot\!10^{-16}}\cr
\tabentry{6}{[0, 1, 0, -208, -1412]}{x^3-6}{2^{2}3^{5}}{2^{10}3^{10}5^{4}}{0}{\nosha}{1\!\cdot\!3^{3}\!+\!1\!\cdot\!3^{4}\!+\!1\!\cdot\!3^{5}\!+\!O(3^{7})}{\noeuler}\errentry{3\!\cdot\!10^{-10}}{2\!\cdot\!10^{-6}}\cr
\tabentry{12}{[0, 1, 0, -208, -1412]}{x^3-12}{2^{2}3^{5}}{2^{10}3^{10}5^{4}}{2^{3}3^{3}}{6^2}{1\!\cdot\!3^{3}\!+\!1\!\cdot\!3^{4}\!+\!1\!\cdot\!3^{5}\!+\!O(3^{7})}{1\!\cdot\!3^{3}\!+\!2\!\cdot\!3^{4}\!+\!2\!\cdot\!3^{5}\!+\!O(3^{6})}\errentry{1\!\cdot\!10^{-11}}{1\!\cdot\!10^{-11}}\cr
\tabentry{28}{[0, 1, 0, -208, -1412]}{x^3-28}{2^{2}3\!\cdot\!7^{2}}{2^{10}3^{2}5^{4}7^{4}}{2^{3}3}{2^2}{1\!\cdot\!3^{3}\phantom{\,\,\!+\!2\!\cdot\!3^4}\,\!+\!1\!\cdot\!3^{5}\!+\!O(3^{7})}{1\!\cdot\!3^{3}\phantom{\,\,\!+\!2\!\cdot\!3^4}\,\!+\!2\!\cdot\!3^{5}\!+\!O(3^{8})}\errentry{1\!\cdot\!10^{-14}}{1\!\cdot\!10^{-15}}\cr
\tabentry{44}{[0, 1, 0, -208, -1412]}{x^3-44}{2^{2}3\!\cdot\!11^{2}}{2^{10}3^{2}5^{4}11^{4}}{2\!\cdot\!3}{1^2}{2\!\cdot\!3^{3}\phantom{\,\,\!+\!2\!\cdot\!3^4}\,\!+\!1\!\cdot\!3^{5}\!+\!O(3^{6})}{2\!\cdot\!3^{3}\phantom{\,\,\!+\!2\!\cdot\!3^4}\,\!+\!1\!\cdot\!3^{5}\!+\!O(3^{6})}\errentry{5\!\cdot\!10^{-9}}{1\!\cdot\!10^{-9}}\cr
\hline
\tabend

\noindent\begin{tabular}{cllll} 
\extravspace
$E=11A3$ & \multicolumn{3}{l}{$y^2+y=x^3-x^2$} & \tabname{5}{11A3}\cr
\extravspace
$p=5$ & \multicolumn{3}{l}{$P_5(E/\Q,T)=1-T+5T^2$}\cr
$/\vQ\alsub{}$   & $N(E)=11$ & $L^*(E/\Q)=5^{-2}$ & $|E(\Q)|=5$ & $|\sha(E/\Q)|=1^2$ \cr
$/K\alsub{}$      & $N(E,\sigma)=5^{6}\!\cdot\!11^{4}$ & $L^*(E/K)=5^{-2}$ & $|E(K)|=5$ & $|\sha(E/K)|=1^2$ \cr
$/\vQ\alsub{11}$ & $\hbox{I}_1 \hbox{ split},\>c_{11}=1$ & $/K\alsub{v}\quad(v|11)$ & $\hbox{I}_1 \hbox{ split},\>c_v=1$\cr
\end{tabular}
\tabbegin
\tabentry{2}{[0, -1, 1, 0, 0]}{x^5-2}{2^{4}5^{5}}{2^{8}5^{10}11^{4}}{0}{\nosha}{4\!\cdot\!5^{2}\!+\!1\!\cdot\!5^{3}\phantom{\,\,\!+\!2\!\cdot\!5^4}\,\!+\!O(5^{5})}{\noeuler}\errentry{3\!\cdot\!10^{-15}}{5\!\cdot\!10^{-9}}\cr
\tabentry{3}{[0, -1, 1, 0, 0]}{x^5-3}{3^{4}5^{5}}{3^{8}5^{10}11^{4}}{0}{\nosha}{2\!\cdot\!5^{2}\!+\!4\!\cdot\!5^{3}\!+\!1\!\cdot\!5^{4}\!+\!O(5^{7})}{\noeuler}\errentry{4\!\cdot\!10^{-9}}{4\!\cdot\!10^{-6}}\cr
\tabentry{5}{[0, -1, 1, 0, 0]}{x^5-5}{5^{9}}{5^{18}11^{4}}{1}{1^2}{4\!\cdot\!5^{0}\!+\!3\!\cdot\!5^{1}\phantom{\,\,\!+\!2\!\cdot\!5^2}\,\!+\!O(5^{4})}{4\!\cdot\!5^{0}\!+\!1\!\cdot\!5^{1}\!+\!3\!\cdot\!5^{2}\!+\!O(5^{3})}\errentry{2\!\cdot\!10^{-4}}{1\!\cdot\!10^{-3}}\cr
\tabentry{6}{[0, -1, 1, 0, 0]}{x^5-6}{2^{4}3^{4}5^{5}}{2^{8}3^{8}5^{10}11^{4}}{5^{2}}{5^2}{2\!\cdot\!5^{4}\!+\!3\!\cdot\!5^{5}\!+\!1\!\cdot\!5^{6}\!+\!O(5^{8})}{4\!\cdot\!5^{2}\!+\!2\!\cdot\!5^{3}\!+\!4\!\cdot\!5^{4}\!+\!O(5^{5})}\errentry{4\!\cdot\!10^{-4}}{5\!\cdot\!10^{-4}}\cr
\tabentry{7}{[0, -1, 1, 0, 0]}{x^5-7}{5^{3}7^{4}}{5^{6}7^{8}11^{4}}{2^{2}}{2^2}{4\!\cdot\!5^{2}\!+\!4\!\cdot\!5^{3}\!+\!3\!\cdot\!5^{4}\!+\!O(5^{5})}{4\!\cdot\!5^{2}\!+\!4\!\cdot\!5^{3}\!+\!4\!\cdot\!5^{4}\!+\!O(5^{5})}\errentry{9\!\cdot\!10^{-10}}{4\!\cdot\!10^{-10}}\cr
\tabentry{11}{[0, -1, 1, 0, 0]}{x^5-11}{5^{5}11^{4}}{5^{10}11^{8}}{2^{4}5}{4^2}{4\!\cdot\!5^{4}\!+\!3\!\cdot\!5^{5}\phantom{\,\,\!+\!2\!\cdot\!5^6}\,\!+\!O(5^{7})}{4\!\cdot\!5^{1}\!+\!1\!\cdot\!5^{2}\!+\!3\!\cdot\!5^{3}\!+\!O(5^{4})}\errentry{4\!\cdot\!10^{-9}}{4\!\cdot\!10^{-9}}\cr
\tabentry{12}{[0, -1, 1, 0, 0]}{x^5-12}{2^{4}3^{4}5^{5}}{2^{8}3^{8}5^{10}11^{4}}{0}{\nosha}{2\!\cdot\!5^{4}\!+\!3\!\cdot\!5^{5}\!+\!1\!\cdot\!5^{6}\!+\!O(5^{8})}{\noeuler}\errentry{8\!\cdot\!10^{-5}}{6\!\cdot\!10^{-3}}\cr
\tabentry{18}{[0, -1, 1, 0, 0]}{x^5-18}{2^{4}3^{4}5^{3}}{2^{8}3^{8}5^{6}11^{4}}{0}{\nosha}{2\!\cdot\!5^{4}\!+\!3\!\cdot\!5^{5}\!+\!1\!\cdot\!5^{6}\!+\!O(5^{8})}{\noeuler}\errentry{1\!\cdot\!10^{-11}}{4\!\cdot\!10^{-7}}\cr
\tabentry{22}{[0, -1, 1, 0, 0]}{x^5-22}{2^{4}5^{5}11^{4}}{2^{8}5^{10}11^{8}}{5^{3}}{5^2}{4\!\cdot\!5^{6}\!+\!1\!\cdot\!5^{7}\phantom{\,\,\!+\!2\!\cdot\!5^8}\,\!+\!O(5^{9})}{1\!\cdot\!5^{3}\!+\!1\!\cdot\!5^{4}\!+\!2\!\cdot\!5^{5}\!+\!O(5^{8})}\errentry{3\!\cdot\!10^{-4}}{3\!\cdot\!10^{-4}}\cr
\tabentry{44}{[0, -1, 1, 0, 0]}{x^5-44}{2^{4}5^{5}11^{4}}{2^{8}5^{10}11^{8}}{0}{\nosha}{4\!\cdot\!5^{6}\!+\!1\!\cdot\!5^{7}\phantom{\,\,\!+\!2\!\cdot\!5^8}\,\!+\!O(5^{9})}{\noeuler}\errentry{1\!\cdot\!10^{-4}}{1\!\cdot\!10^{-2}}\cr
\tabentry{48}{[0, -1, 1, 0, 0]}{x^5-48}{2^{4}3^{4}5^{5}}{2^{8}3^{8}5^{10}11^{4}}{5^{2}}{5^2}{2\!\cdot\!5^{4}\!+\!3\!\cdot\!5^{5}\!+\!1\!\cdot\!5^{6}\!+\!O(5^{8})}{4\!\cdot\!5^{2}\!+\!2\!\cdot\!5^{3}\!+\!4\!\cdot\!5^{4}\!+\!O(5^{5})}\errentry{1\!\cdot\!10^{-5}}{1\!\cdot\!10^{-5}}\cr
\tabentry{88}{[0, -1, 1, 0, 0]}{x^5-88}{2^{4}5^{5}11^{4}}{2^{8}5^{10}11^{8}}{5^{3}}{5^2}{4\!\cdot\!5^{6}\!+\!1\!\cdot\!5^{7}\phantom{\,\,\!+\!2\!\cdot\!5^8}\,\!+\!O(5^{9})}{1\!\cdot\!5^{3}\!+\!1\!\cdot\!5^{4}\!+\!2\!\cdot\!5^{5}\!+\!O(5^{8})}\errentry{6\!\cdot\!10^{-5}}{6\!\cdot\!10^{-5}}\cr
\tabentry{99}{[0, -1, 1, 0, 0]}{x^5-99}{3^{4}5^{3}11^{4}}{3^{8}5^{6}11^{8}}{3^{2}5}{3^2}{2\!\cdot\!5^{6}\!+\!4\!\cdot\!5^{7}\!+\!1\!\cdot\!5^{8}\!+\!O(5^{9})}{4\!\cdot\!5^{3}\!+\!2\!\cdot\!5^{4}\!+\!3\!\cdot\!5^{5}\!+\!O(5^{6})}\errentry{1\!\cdot\!10^{-5}}{9\!\cdot\!10^{-6}}\cr
\tabentry{176}{[0, -1, 1, 0, 0]}{x^5-176}{2^{4}5^{3}11^{4}}{2^{8}5^{6}11^{8}}{5}{1^2}{4\!\cdot\!5^{6}\!+\!1\!\cdot\!5^{7}\phantom{\,\,\!+\!2\!\cdot\!5^8}\,\!+\!O(5^{9})}{1\!\cdot\!5^{3}\!+\!3\!\cdot\!5^{4}\!+\!4\!\cdot\!5^{5}\!+\!O(5^{6})}\errentry{7\!\cdot\!10^{-11}}{1\!\cdot\!10^{-11}}\cr
\hline
\tabend

\noindent\begin{tabular}{cllll} 
\extravspace
$E=19A3$ & \multicolumn{3}{l}{$y^2+y=x^3+x^2+x$} & \tabname{5}{19A3}\cr
\extravspace
$p=5$ & \multicolumn{3}{l}{$P_5(E/\Q,T)=1-3T+5T^2$}\cr
$/\vQ\alsub{}$   & $N(E)=19$ & $L^*(E/\Q)=3^{-2}$ & $|E(\Q)|=3$ & $|\sha(E/\Q)|=1^2$ \cr
$/K\alsub{}$      & $N(E,\sigma)=5^{6}\!\cdot\!19^{4}$ & $L^*(E/K)=3^{-2}$ & $|E(K)|=3$ & $|\sha(E/K)|=1^2$ \cr
$/\vQ\alsub{19}$ & $\hbox{I}_1 \hbox{ split},\>c_{19}=1$ & $/K\alsub{v}\quad(v|19)$ & $\hbox{I}_1 \hbox{ split},\>c_v=1$\cr
\end{tabular}
\tabbegin
\tabentry{2}{[0, 1, 1, 1, 0]}{x^5-2}{2^{4}5^{5}}{2^{8}5^{10}19^{4}}{1}{1^2}{2\!\cdot\!5^{0}\!+\!2\!\cdot\!5^{1}\!+\!4\!\cdot\!5^{2}\!+\!O(5^{3})}{2\!\cdot\!5^{0}\!+\!1\!\cdot\!5^{1}\!+\!2\!\cdot\!5^{2}\!+\!O(5^{3})}\errentry{9\!\cdot\!10^{-11}}{2\!\cdot\!10^{-11}}\cr
\tabentry{3}{[0, 1, 1, 1, 0]}{x^5-3}{3^{4}5^{5}}{3^{8}5^{10}19^{4}}{2^{2}}{2^2}{3\!\cdot\!5^{0}\!+\!1\!\cdot\!5^{1}\!+\!2\!\cdot\!5^{2}\!+\!O(5^{3})}{3\!\cdot\!5^{0}\!+\!3\!\cdot\!5^{1}\phantom{\,\,\!+\!2\!\cdot\!5^2}\,\!+\!O(5^{3})}\errentry{2\!\cdot\!10^{-8}}{1\!\cdot\!10^{-8}}\cr
\tabentry{5}{[0, 1, 1, 1, 0]}{x^5-5}{5^{9}}{5^{18}19^{4}}{2^{4}}{4^2}{3\!\cdot\!5^{0}\!+\!2\!\cdot\!5^{1}\!+\!3\!\cdot\!5^{2}\!+\!O(5^{3})}{3\!\cdot\!5^{0}\!+\!4\!\cdot\!5^{1}\!+\!1\!\cdot\!5^{2}\!+\!O(5^{3})}\errentry{3\!\cdot\!10^{-4}}{3\!\cdot\!10^{-4}}\cr
\tabentry{6}{[0, 1, 1, 1, 0]}{x^5-6}{2^{4}3^{4}5^{5}}{2^{8}3^{8}5^{10}19^{4}}{2^{2}}{2^2}{2\!\cdot\!5^{0}\!+\!3\!\cdot\!5^{1}\phantom{\,\,\!+\!2\!\cdot\!5^2}\,\!+\!O(5^{3})}{2\!\cdot\!5^{0}\!+\!3\!\cdot\!5^{1}\!+\!2\!\cdot\!5^{2}\!+\!O(5^{5})}\errentry{4\!\cdot\!10^{-4}}{3\!\cdot\!10^{-4}}\cr
\tabentry{7}{[0, 1, 1, 1, 0]}{x^5-7}{5^{3}7^{4}}{5^{6}7^{8}19^{4}}{2^{4}}{4^2}{3\!\cdot\!5^{0}\phantom{\,\,\!+\!2\!\cdot\!5^1}\,\!+\!2\!\cdot\!5^{2}\!+\!O(5^{3})}{3\!\cdot\!5^{0}\phantom{\,\,\!+\!2\!\cdot\!5^1}\,\!+\!4\!\cdot\!5^{2}\!+\!O(5^{3})}\errentry{2\!\cdot\!10^{-8}}{2\!\cdot\!10^{-8}}\cr
\tabentry{12}{[0, 1, 1, 1, 0]}{x^5-12}{2^{4}3^{4}5^{5}}{2^{8}3^{8}5^{10}19^{4}}{7^{2}}{7^2}{2\!\cdot\!5^{0}\!+\!3\!\cdot\!5^{1}\phantom{\,\,\!+\!2\!\cdot\!5^2}\,\!+\!O(5^{3})}{2\!\cdot\!5^{0}\phantom{\,\,\!+\!2\!\cdot\!5^1}\,\!+\!4\!\cdot\!5^{2}\!+\!O(5^{3})}\errentry{1\!\cdot\!10^{-3}}{2\!\cdot\!10^{-3}}\cr
\tabentry{18}{[0, 1, 1, 1, 0]}{x^5-18}{2^{4}3^{4}5^{3}}{2^{8}3^{8}5^{6}19^{4}}{1}{1^2}{2\!\cdot\!5^{0}\!+\!3\!\cdot\!5^{1}\phantom{\,\,\!+\!2\!\cdot\!5^2}\,\!+\!O(5^{3})}{2\!\cdot\!5^{0}\!+\!4\!\cdot\!5^{1}\phantom{\,\,\!+\!2\!\cdot\!5^2}\,\!+\!O(5^{3})}\errentry{4\!\cdot\!10^{-8}}{1\!\cdot\!10^{-8}}\cr
\tabentry{19}{[0, 1, 1, 1, 0]}{x^5-19}{5^{5}19^{4}}{5^{10}19^{8}}{2^{4}5}{4^2}{2\!\cdot\!5^{2}\phantom{\,\,\!+\!2\!\cdot\!5^3}\,\!+\!2\!\cdot\!5^{4}\!+\!O(5^{5})}{3\!\cdot\!5^{1}\!+\!3\!\cdot\!5^{2}\!+\!1\!\cdot\!5^{3}\!+\!O(5^{4})}\errentry{2\!\cdot\!10^{-4}}{1\!\cdot\!10^{-4}}\cr
\tabentry{48}{[0, 1, 1, 1, 0]}{x^5-48}{2^{4}3^{4}5^{5}}{2^{8}3^{8}5^{10}19^{4}}{2^{2}}{2^2}{2\!\cdot\!5^{0}\!+\!3\!\cdot\!5^{1}\phantom{\,\,\!+\!2\!\cdot\!5^2}\,\!+\!O(5^{3})}{2\!\cdot\!5^{0}\!+\!3\!\cdot\!5^{1}\!+\!2\!\cdot\!5^{2}\!+\!O(5^{5})}\errentry{2\!\cdot\!10^{-3}}{1\!\cdot\!10^{-3}}\cr
\tabentry{57}{[0, 1, 1, 1, 0]}{x^5-57}{3^{4}5^{3}19^{4}}{3^{8}5^{6}19^{8}}{3^{4}5}{9^2}{2\!\cdot\!5^{2}\!+\!1\!\cdot\!5^{3}\!+\!2\!\cdot\!5^{4}\!+\!O(5^{6})}{3\!\cdot\!5^{1}\!+\!4\!\cdot\!5^{2}\!+\!3\!\cdot\!5^{3}\!+\!O(5^{4})}\errentry{1\!\cdot\!10^{-4}}{3\!\cdot\!10^{-4}}\cr
\tabentry{76}{[0, 1, 1, 1, 0]}{x^5-76}{2^{4}5^{3}19^{4}}{2^{8}5^{6}19^{8}}{2^{6}5}{8^2}{3\!\cdot\!5^{2}\!+\!4\!\cdot\!5^{3}\!+\!4\!\cdot\!5^{4}\!+\!O(5^{5})}{2\!\cdot\!5^{1}\phantom{\,\,\!+\!2\!\cdot\!5^2}\,\!+\!3\!\cdot\!5^{3}\!+\!O(5^{5})}\errentry{3\!\cdot\!10^{-5}}{6\!\cdot\!10^{-5}}\cr
\hline
\tabend

\noindent\begin{tabular}{cllll} 
\extravspace
$E=21A4$ & \multicolumn{3}{l}{$y^2+xy=x^3+x$} & \tabname{5}{21A4}\cr
\extravspace
$p=5$ & \multicolumn{3}{l}{$P_5(E/\Q,T)=1+2T+5T^2$}\cr
$/\vQ\alsub{}$   & $N(E)=3\!\cdot\!7$ & $L^*(E/\Q)=2^{-3}$ & $|E(\Q)|=4$ & $|\sha(E/\Q)|=1^2$ \cr
$/K\alsub{}$      & $N(E,\sigma)=3^{4}\!\cdot\!5^{6}\!\cdot\!7^{4}$ & $L^*(E/K)=2^{-3}$ & $|E(K)|=4$ & $|\sha(E/K)|=1^2$ \cr
$/\vQ\alsub{3}$ & $\hbox{I}_2 \hbox{ split},\>c_{3}=2$ & $/K\alsub{v}\quad(v|3)$ & $\hbox{I}_2 \hbox{ split},\>c_v=2$\cr
$/\vQ\alsub{7}$ & $\hbox{I}_1 \hbox{ non-split},\>c_{7}=1$ & $/K\alsub{v}\quad(v|7)$ & $\hbox{I}_1 \hbox{ split},\>c_v=1$\cr
\end{tabular}
\tabbegin
\tabentry{2}{[1, 0, 0, 1, 0]}{x^5-2}{2^{4}5^{5}}{2^{8}3^{4}5^{10}7^{4}}{2}{1^2}{4\!\cdot\!5^{0}\phantom{\,\,\!+\!2\!\cdot\!5^1}\,\!+\!3\!\cdot\!5^{2}\!+\!O(5^{3})}{4\!\cdot\!5^{0}\!+\!2\!\cdot\!5^{1}\!+\!2\!\cdot\!5^{2}\!+\!O(5^{3})}\errentry{4\!\cdot\!10^{-11}}{1\!\cdot\!10^{-11}}\cr
\tabentry{5}{[1, 0, 0, 1, 0]}{x^5-5}{5^{9}}{3^{4}5^{18}7^{4}}{2^{3}}{2^2}{4\!\cdot\!5^{0}\phantom{\,\,\!+\!2\!\cdot\!5^1}\,\!+\!3\!\cdot\!5^{2}\!+\!O(5^{3})}{4\!\cdot\!5^{0}\!+\!1\!\cdot\!5^{1}\!+\!1\!\cdot\!5^{2}\!+\!O(5^{3})}\errentry{1\!\cdot\!10^{-2}}{5\!\cdot\!10^{-3}}\cr
\tabentry{21}{[1, 0, 0, 1, 0]}{x^5-21}{3^{4}5^{5}7^{4}}{3^{8}5^{10}7^{8}}{5\!\cdot\!17^{2}}{17^2}{4\!\cdot\!5^{3}\!+\!4\!\cdot\!5^{4}\!+\!2\!\cdot\!5^{5}\!+\!O(5^{7})}{2\!\cdot\!5^{1}\!+\!2\!\cdot\!5^{2}\!+\!3\!\cdot\!5^{3}\!+\!O(5^{4})}\errentry{4\!\cdot\!10^{-5}}{1\!\cdot\!10^{-4}}\cr
\tabentry{63}{[1, 0, 0, 1, 0]}{x^5-63}{3^{4}5^{5}7^{4}}{3^{8}5^{10}7^{8}}{2^{4}5}{4^2}{4\!\cdot\!5^{3}\!+\!4\!\cdot\!5^{4}\!+\!2\!\cdot\!5^{5}\!+\!O(5^{7})}{3\!\cdot\!5^{1}\phantom{\,\,\!+\!2\!\cdot\!5^2}\,\!+\!4\!\cdot\!5^{3}\!+\!O(5^{4})}\errentry{1\!\cdot\!10^{-5}}{2\!\cdot\!10^{-5}}\cr
\tabentry{126}{[1, 0, 0, 1, 0]}{x^5-126}{2^{4}3^{4}5^{3}7^{4}}{2^{8}3^{8}5^{6}7^{8}}{2^{4}5}{4^2}{4\!\cdot\!5^{3}\!+\!4\!\cdot\!5^{4}\!+\!2\!\cdot\!5^{5}\!+\!O(5^{7})}{3\!\cdot\!5^{1}\!+\!1\!\cdot\!5^{2}\phantom{\,\,\!+\!2\!\cdot\!5^3}\,\!+\!O(5^{4})}\errentry{4\!\cdot\!10^{-4}}{4\!\cdot\!10^{-4}}\cr
\tabentry{147}{[1, 0, 0, 1, 0]}{x^5-147}{3^{4}5^{5}7^{4}}{3^{8}5^{10}7^{8}}{5}{1^2}{4\!\cdot\!5^{3}\!+\!4\!\cdot\!5^{4}\!+\!2\!\cdot\!5^{5}\!+\!O(5^{7})}{3\!\cdot\!5^{1}\!+\!1\!\cdot\!5^{2}\!+\!4\!\cdot\!5^{3}\!+\!O(5^{4})}\errentry{3\!\cdot\!10^{-4}}{6\!\cdot\!10^{-6}}\cr
\tabentry{168}{[1, 0, 0, 1, 0]}{x^5-168}{2^{4}3^{4}5^{3}7^{4}}{2^{8}3^{8}5^{6}7^{8}}{2^{2}5}{2^2}{4\!\cdot\!5^{3}\!+\!4\!\cdot\!5^{4}\!+\!2\!\cdot\!5^{5}\!+\!O(5^{7})}{2\!\cdot\!5^{1}\phantom{\,\,\!+\!2\!\cdot\!5^2}\,\phantom{\,\,\!+\!2\!\cdot\!5^3}\,\!+\!O(5^{4})}\errentry{1\!\cdot\!10^{-4}}{2\!\cdot\!10^{-5}}\cr
\tabentry{567}{[1, 0, 0, 1, 0]}{x^5-567}{3^{4}5^{5}7^{4}}{3^{8}5^{10}7^{8}}{2^{6}5}{8^2}{4\!\cdot\!5^{3}\!+\!4\!\cdot\!5^{4}\!+\!2\!\cdot\!5^{5}\!+\!O(5^{7})}{2\!\cdot\!5^{1}\!+\!2\!\cdot\!5^{2}\!+\!1\!\cdot\!5^{3}\!+\!O(5^{4})}\errentry{1\!\cdot\!10^{-3}}{3\!\cdot\!10^{-3}}\cr
\hline
\tabend

\noindent\begin{tabular}{cllll} 
\extravspace
$E=24A4$ & \multicolumn{3}{l}{$y^2=x^3-x^2+x$} & \tabname{5}{24A4}\cr
\extravspace
$p=5$ & \multicolumn{3}{l}{$P_5(E/\Q,T)=1+2T+5T^2$}\cr
$/\vQ\alsub{}$   & $N(E)=2^{3}\!\cdot\!3$ & $L^*(E/\Q)=2^{-3}$ & $|E(\Q)|=4$ & $|\sha(E/\Q)|=1^2$ \cr
$/K\alsub{}$      & $N(E,\sigma)=2^{12}\!\cdot\!3^{4}\!\cdot\!5^{6}$ & $L^*(E/K)=2^{-3}$ & $|E(K)|=4$ & $|\sha(E/K)|=1^2$ \cr
$/\vQ\alsub{2}$ & $\hbox{III},\>c_{2}=2$ & $/K\alsub{v}\quad(v|2)$ & $\hbox{III},\>c_v=2$\cr
$/\vQ\alsub{3}$ & $\hbox{I}_1 \hbox{ non-split},\>c_{3}=1$ & $/K\alsub{v}\quad(v|3)$ & $\hbox{I}_1 \hbox{ split},\>c_v=1$\cr
\end{tabular}
\tabbegin
\tabentry{2}{[0, -1, 0, 1, 0]}{x^5-2}{2^{4}5^{5}}{2^{12}3^{4}5^{10}}{2}{1^2}{4\!\cdot\!5^{0}\phantom{\,\,\!+\!2\!\cdot\!5^1}\,\!+\!3\!\cdot\!5^{2}\!+\!O(5^{3})}{4\!\cdot\!5^{0}\!+\!2\!\cdot\!5^{1}\!+\!2\!\cdot\!5^{2}\!+\!O(5^{3})}\errentry{2\!\cdot\!10^{-16}}{1\!\cdot\!10^{-15}}\cr
\tabentry{10}{[0, -1, 0, 1, 0]}{x^5-10}{2^{4}5^{9}}{2^{12}3^{4}5^{18}}{2^{5}}{4^2}{4\!\cdot\!5^{0}\phantom{\,\,\!+\!2\!\cdot\!5^1}\,\!+\!3\!\cdot\!5^{2}\!+\!O(5^{3})}{4\!\cdot\!5^{0}\!+\!1\!\cdot\!5^{1}\!+\!1\!\cdot\!5^{2}\!+\!O(5^{3})}\errentry{3\!\cdot\!10^{-5}}{3\!\cdot\!10^{-5}}\cr
\tabentry{20}{[0, -1, 0, 1, 0]}{x^5-20}{2^{4}5^{9}}{2^{12}3^{4}5^{18}}{2^{5}}{4^2}{4\!\cdot\!5^{0}\phantom{\,\,\!+\!2\!\cdot\!5^1}\,\!+\!3\!\cdot\!5^{2}\!+\!O(5^{3})}{4\!\cdot\!5^{0}\!+\!1\!\cdot\!5^{1}\!+\!1\!\cdot\!5^{2}\!+\!O(5^{3})}\errentry{8\!\cdot\!10^{-4}}{7\!\cdot\!10^{-4}}\cr
\tabentry{40}{[0, -1, 0, 1, 0]}{x^5-40}{2^{4}5^{9}}{2^{12}3^{4}5^{18}}{2^{5}}{4^2}{4\!\cdot\!5^{0}\phantom{\,\,\!+\!2\!\cdot\!5^1}\,\!+\!3\!\cdot\!5^{2}\!+\!O(5^{3})}{4\!\cdot\!5^{0}\!+\!1\!\cdot\!5^{1}\!+\!1\!\cdot\!5^{2}\!+\!O(5^{3})}\errentry{3\!\cdot\!10^{-4}}{3\!\cdot\!10^{-4}}\cr
\tabentry{80}{[0, -1, 0, 1, 0]}{x^5-80}{2^{4}5^{9}}{2^{12}3^{4}5^{18}}{2}{1^2}{4\!\cdot\!5^{0}\phantom{\,\,\!+\!2\!\cdot\!5^1}\,\!+\!3\!\cdot\!5^{2}\!+\!O(5^{3})}{4\!\cdot\!5^{0}\!+\!4\!\cdot\!5^{1}\!+\!1\!\cdot\!5^{2}\!+\!O(5^{4})}\errentry{6\!\cdot\!10^{-5}}{3\!\cdot\!10^{-5}}\cr
\tabentry{176}{[0, -1, 0, 1, 0]}{x^5-176}{2^{4}5^{3}11^{4}}{2^{12}3^{4}5^{6}11^{8}}{2^{3}}{2^2}{4\!\cdot\!5^{0}\!+\!4\!\cdot\!5^{1}\!+\!3\!\cdot\!5^{2}\!+\!O(5^{4})}{4\!\cdot\!5^{0}\!+\!1\!\cdot\!5^{1}\!+\!2\!\cdot\!5^{2}\!+\!O(5^{4})}\errentry{4\!\cdot\!10^{-4}}{2\!\cdot\!10^{-4}}\cr
\hline
\tabend

\noindent\begin{tabular}{cllll} 
\extravspace
$E=56A1$ & \multicolumn{3}{l}{$y^2=x^3+x+2$} & \tabname{5}{56A1}\cr
\extravspace
$p=5$ & \multicolumn{3}{l}{$P_5(E/\Q,T)=1-2T+5T^2$}\cr
$/\vQ\alsub{}$   & $N(E)=2^{3}\!\cdot\!7$ & $L^*(E/\Q)=2^{-2}$ & $|E(\Q)|=4$ & $|\sha(E/\Q)|=1^2$ \cr
$/K\alsub{}$      & $N(E,\sigma)=2^{12}\!\cdot\!5^{6}\!\cdot\!7^{4}$ & $L^*(E/K)=2^{-2}$ & $|E(K)|=4$ & $|\sha(E/K)|=1^2$ \cr
$/\vQ\alsub{2}$ & $\hbox{I}^*_1,\>c_{2}=4$ & $/K\alsub{v}\quad(v|2)$ & $\hbox{I}^*_1,\>c_v=4$\cr
$/\vQ\alsub{7}$ & $\hbox{I}_1 \hbox{ non-split},\>c_{7}=1$ & $/K\alsub{v}\quad(v|7)$ & $\hbox{I}_1 \hbox{ split},\>c_v=1$\cr
\end{tabular}
\tabbegin
\tabentry{2}{[0, 0, 0, 1, 2]}{x^5-2}{2^{4}5^{5}}{2^{12}5^{10}7^{4}}{2^{4}}{2^2}{3\!\cdot\!5^{0}\phantom{\,\,\!+\!2\!\cdot\!5^1}\,\phantom{\,\,\!+\!2\!\cdot\!5^2}\,\!+\!O(5^{3})}{3\!\cdot\!5^{0}\!+\!2\!\cdot\!5^{1}\!+\!4\!\cdot\!5^{2}\!+\!O(5^{3})}\errentry{4\!\cdot\!10^{-14}}{2\!\cdot\!10^{-14}}\cr
\tabentry{6}{[0, 0, 0, 1, 2]}{x^5-6}{2^{4}3^{4}5^{5}}{2^{12}3^{8}5^{10}7^{4}}{2^{2}3^{2}}{3^2}{2\!\cdot\!5^{0}\!+\!1\!\cdot\!5^{1}\phantom{\,\,\!+\!2\!\cdot\!5^2}\,\!+\!O(5^{3})}{2\!\cdot\!5^{0}\!+\!4\!\cdot\!5^{1}\!+\!4\!\cdot\!5^{2}\!+\!O(5^{4})}\errentry{1\!\cdot\!10^{-4}}{1\!\cdot\!10^{-4}}\cr
\tabentry{12}{[0, 0, 0, 1, 2]}{x^5-12}{2^{4}3^{4}5^{5}}{2^{12}3^{8}5^{10}7^{4}}{2^{4}}{2^2}{2\!\cdot\!5^{0}\!+\!1\!\cdot\!5^{1}\phantom{\,\,\!+\!2\!\cdot\!5^2}\,\!+\!O(5^{3})}{2\!\cdot\!5^{0}\!+\!1\!\cdot\!5^{1}\!+\!3\!\cdot\!5^{2}\!+\!O(5^{3})}\errentry{1\!\cdot\!10^{-3}}{4\!\cdot\!10^{-4}}\cr
\tabentry{18}{[0, 0, 0, 1, 2]}{x^5-18}{2^{4}3^{4}5^{3}}{2^{12}3^{8}5^{6}7^{4}}{2^{2}3^{2}}{3^2}{2\!\cdot\!5^{0}\!+\!1\!\cdot\!5^{1}\phantom{\,\,\!+\!2\!\cdot\!5^2}\,\!+\!O(5^{3})}{2\!\cdot\!5^{0}\!+\!2\!\cdot\!5^{1}\phantom{\,\,\!+\!2\!\cdot\!5^2}\,\!+\!O(5^{3})}\errentry{1\!\cdot\!10^{-9}}{9\!\cdot\!10^{-10}}\cr
\tabentry{48}{[0, 0, 0, 1, 2]}{x^5-48}{2^{4}3^{4}5^{5}}{2^{12}3^{8}5^{10}7^{4}}{2^{2}3^{2}}{3^2}{2\!\cdot\!5^{0}\!+\!1\!\cdot\!5^{1}\phantom{\,\,\!+\!2\!\cdot\!5^2}\,\!+\!O(5^{3})}{2\!\cdot\!5^{0}\!+\!4\!\cdot\!5^{1}\!+\!4\!\cdot\!5^{2}\!+\!O(5^{4})}\errentry{5\!\cdot\!10^{-4}}{3\!\cdot\!10^{-4}}\cr
\hline
\tabend

\noindent\begin{tabular}{cllll} 
\extravspace
$E=56B1$ & \multicolumn{3}{l}{$y^2=x^3-x^2-4$} & \tabname{5}{56B1}\cr
\extravspace
$p=5$ & \multicolumn{3}{l}{$P_5(E/\Q,T)=1+4T+5T^2$}\cr
$/\vQ\alsub{}$   & $N(E)=2^{3}\!\cdot\!7$ & $L^*(E/\Q)=2^{-1}$ & $|E(\Q)|=2$ & $|\sha(E/\Q)|=1^2$ \cr
$/K\alsub{}$      & $N(E,\sigma)=2^{12}\!\cdot\!5^{6}\!\cdot\!7^{4}$ & $L^*(E/K)=2^{-1}$ & $|E(K)|=2$ & $|\sha(E/K)|=1^2$ \cr
$/\vQ\alsub{2}$ & $\hbox{III}^*,\>c_{2}=2$ & $/K\alsub{v}\quad(v|2)$ & $\hbox{III}^*,\>c_v=2$\cr
$/\vQ\alsub{7}$ & $\hbox{I}_1 \hbox{ split},\>c_{7}=1$ & $/K\alsub{v}\quad(v|7)$ & $\hbox{I}_1 \hbox{ split},\>c_v=1$\cr
\end{tabular}
\tabbegin
\tabentry{2}{[0, -1, 0, 0, -4]}{x^5-2}{2^{4}5^{5}}{2^{12}5^{10}7^{4}}{0}{\nosha}{3\!\cdot\!5^{2}\!+\!2\!\cdot\!5^{3}\phantom{\,\,\!+\!2\!\cdot\!5^4}\,\!+\!O(5^{5})}{\noeuler}\errentry{7\!\cdot\!10^{-14}}{2\!\cdot\!10^{-8}}\cr
\tabentry{6}{[0, -1, 0, 0, -4]}{x^5-6}{2^{4}3^{4}5^{5}}{2^{12}3^{8}5^{10}7^{4}}{2^{3}5^{2}}{5^2}{3\!\cdot\!5^{2}\!+\!1\!\cdot\!5^{3}\!+\!4\!\cdot\!5^{4}\!+\!O(5^{5})}{2\!\cdot\!5^{2}\!+\!2\!\cdot\!5^{3}\!+\!3\!\cdot\!5^{4}\!+\!O(5^{5})}\errentry{1\!\cdot\!10^{-3}}{1\!\cdot\!10^{-3}}\cr
\tabentry{12}{[0, -1, 0, 0, -4]}{x^5-12}{2^{4}3^{4}5^{5}}{2^{12}3^{8}5^{10}7^{4}}{0}{\nosha}{3\!\cdot\!5^{2}\!+\!1\!\cdot\!5^{3}\!+\!4\!\cdot\!5^{4}\!+\!O(5^{5})}{\noeuler}\errentry{4\!\cdot\!10^{-4}}{1\!\cdot\!10^{-2}}\cr
\tabentry{18}{[0, -1, 0, 0, -4]}{x^5-18}{2^{4}3^{4}5^{3}}{2^{12}3^{8}5^{6}7^{4}}{2^{3}}{1^2}{3\!\cdot\!5^{2}\!+\!1\!\cdot\!5^{3}\!+\!4\!\cdot\!5^{4}\!+\!O(5^{5})}{3\!\cdot\!5^{2}\phantom{\,\,\!+\!2\!\cdot\!5^3}\,\!+\!4\!\cdot\!5^{4}\!+\!O(5^{6})}\errentry{1\!\cdot\!10^{-8}}{2\!\cdot\!10^{-9}}\cr
\tabentry{48}{[0, -1, 0, 0, -4]}{x^5-48}{2^{4}3^{4}5^{5}}{2^{12}3^{8}5^{10}7^{4}}{2^{3}5^{2}}{5^2}{3\!\cdot\!5^{2}\!+\!1\!\cdot\!5^{3}\!+\!4\!\cdot\!5^{4}\!+\!O(5^{5})}{2\!\cdot\!5^{2}\!+\!2\!\cdot\!5^{3}\!+\!3\!\cdot\!5^{4}\!+\!O(5^{5})}\errentry{1\!\cdot\!10^{-4}}{1\!\cdot\!10^{-4}}\cr
\hline
\tabend

\noindent\begin{tabular}{cllll} 
\extravspace
$E=44A1$ & \multicolumn{3}{l}{$y^2=x^3+x^2+3x-1$} & \tabname{5}{44A1}\cr
\extravspace
$p=5$ & \multicolumn{3}{l}{$P_5(E/\Q,T)=1+3T+5T^2$}\cr
$/\vQ\alsub{}$   & $N(E)=2^{2}\!\cdot\!11$ & $L^*(E/\Q)=3^{-1}$ & $|E(\Q)|=3$ & $|\sha(E/\Q)|=1^2$ \cr
$/K\alsub{}$      & $N(E,\sigma)=2^{8}\!\cdot\!5^{6}\!\cdot\!11^{4}$ & $L^*(E/K)=3^{-1}$ & $|E(K)|=3$ & $|\sha(E/K)|=1^2$ \cr
$/\vQ\alsub{2}$ & $\hbox{IV}^*,\>c_{2}=3$ & $/K\alsub{v}\quad(v|2)$ & $\hbox{IV}^*,\>c_v=3$\cr
$/\vQ\alsub{11}$ & $\hbox{I}_1 \hbox{ non-split},\>c_{11}=1$ & $/K\alsub{v}\quad(v|11)$ & $\hbox{I}_1 \hbox{ non-split},\>c_v=1$\cr
\end{tabular}
\tabbegin
\tabentry{2}{[0, 1, 0, 3, -1]}{x^5-2}{2^{4}5^{5}}{2^{8}5^{10}11^{4}}{2^{3}}{1^2}{4\!\cdot\!5^{0}\!+\!2\!\cdot\!5^{1}\!+\!3\!\cdot\!5^{2}\!+\!O(5^{3})}{4\!\cdot\!5^{0}\!+\!3\!\cdot\!5^{1}\!+\!1\!\cdot\!5^{2}\!+\!O(5^{4})}\errentry{7\!\cdot\!10^{-15}}{2\!\cdot\!10^{-15}}\cr
\tabentry{3}{[0, 1, 0, 3, -1]}{x^5-3}{3^{4}5^{5}}{2^{8}3^{8}5^{10}11^{4}}{0}{\nosha}{2\!\cdot\!5^{2}\!+\!1\!\cdot\!5^{3}\!+\!3\!\cdot\!5^{4}\!+\!O(5^{5})}{\noeuler}\errentry{3\!\cdot\!10^{-5}}{3\!\cdot\!10^{-3}}\cr
\tabentry{6}{[0, 1, 0, 3, -1]}{x^5-6}{2^{4}3^{4}5^{5}}{2^{8}3^{8}5^{10}11^{4}}{2^{3}5^{2}}{5^2}{2\!\cdot\!5^{2}\!+\!1\!\cdot\!5^{3}\!+\!3\!\cdot\!5^{4}\!+\!O(5^{5})}{1\!\cdot\!5^{2}\!+\!3\!\cdot\!5^{3}\!+\!4\!\cdot\!5^{4}\!+\!O(5^{5})}\errentry{1\!\cdot\!10^{-5}}{1\!\cdot\!10^{-5}}\cr
\tabentry{7}{[0, 1, 0, 3, -1]}{x^5-7}{5^{3}7^{4}}{2^{8}5^{6}7^{8}11^{4}}{0}{\nosha}{4\!\cdot\!5^{2}\!+\!3\!\cdot\!5^{3}\!+\!2\!\cdot\!5^{4}\!+\!O(5^{5})}{\noeuler}\errentry{2\!\cdot\!10^{-4}}{9\!\cdot\!10^{-3}}\cr
\tabentry{11}{[0, 1, 0, 3, -1]}{x^5-11}{5^{5}11^{4}}{2^{8}5^{10}11^{8}}{2^{6}3}{8^2}{4\!\cdot\!5^{0}\!+\!3\!\cdot\!5^{1}\!+\!3\!\cdot\!5^{2}\!+\!O(5^{3})}{4\!\cdot\!5^{0}\!+\!3\!\cdot\!5^{1}\!+\!2\!\cdot\!5^{2}\!+\!O(5^{3})}\errentry{1\!\cdot\!10^{-4}}{2\!\cdot\!10^{-4}}\cr
\tabentry{12}{[0, 1, 0, 3, -1]}{x^5-12}{2^{4}3^{4}5^{5}}{2^{8}3^{8}5^{10}11^{4}}{0}{\nosha}{2\!\cdot\!5^{2}\!+\!1\!\cdot\!5^{3}\!+\!3\!\cdot\!5^{4}\!+\!O(5^{5})}{\noeuler}\errentry{4\!\cdot\!10^{-4}}{1\!\cdot\!10^{-2}}\cr
\tabentry{18}{[0, 1, 0, 3, -1]}{x^5-18}{2^{4}3^{4}5^{3}}{2^{8}3^{8}5^{6}11^{4}}{0}{\nosha}{2\!\cdot\!5^{2}\!+\!1\!\cdot\!5^{3}\!+\!3\!\cdot\!5^{4}\!+\!O(5^{5})}{\noeuler}\errentry{6\!\cdot\!10^{-11}}{7\!\cdot\!10^{-7}}\cr
\tabentry{22}{[0, 1, 0, 3, -1]}{x^5-22}{2^{4}5^{5}11^{4}}{2^{8}5^{10}11^{8}}{2^{3}3^{4}}{9^2}{4\!\cdot\!5^{0}\!+\!3\!\cdot\!5^{1}\!+\!3\!\cdot\!5^{2}\!+\!O(5^{3})}{4\!\cdot\!5^{0}\!+\!1\!\cdot\!5^{1}\!+\!4\!\cdot\!5^{2}\!+\!O(5^{3})}\errentry{9\!\cdot\!10^{-5}}{2\!\cdot\!10^{-4}}\cr
\tabentry{44}{[0, 1, 0, 3, -1]}{x^5-44}{2^{4}5^{5}11^{4}}{2^{8}5^{10}11^{8}}{2^{3}3^{4}}{9^2}{4\!\cdot\!5^{0}\!+\!3\!\cdot\!5^{1}\!+\!3\!\cdot\!5^{2}\!+\!O(5^{3})}{4\!\cdot\!5^{0}\!+\!1\!\cdot\!5^{1}\!+\!4\!\cdot\!5^{2}\!+\!O(5^{3})}\errentry{4\!\cdot\!10^{-4}}{1\!\cdot\!10^{-3}}\cr
\tabentry{48}{[0, 1, 0, 3, -1]}{x^5-48}{2^{4}3^{4}5^{5}}{2^{8}3^{8}5^{10}11^{4}}{2^{3}5^{2}}{5^2}{2\!\cdot\!5^{2}\!+\!1\!\cdot\!5^{3}\!+\!3\!\cdot\!5^{4}\!+\!O(5^{5})}{1\!\cdot\!5^{2}\!+\!3\!\cdot\!5^{3}\!+\!4\!\cdot\!5^{4}\!+\!O(5^{5})}\errentry{9\!\cdot\!10^{-5}}{1\!\cdot\!10^{-4}}\cr
\tabentry{88}{[0, 1, 0, 3, -1]}{x^5-88}{2^{4}5^{5}11^{4}}{2^{8}5^{10}11^{8}}{2^{5}7^{2}}{14^2}{4\!\cdot\!5^{0}\!+\!3\!\cdot\!5^{1}\!+\!3\!\cdot\!5^{2}\!+\!O(5^{3})}{4\!\cdot\!5^{0}\!+\!3\!\cdot\!5^{1}\!+\!3\!\cdot\!5^{2}\!+\!O(5^{3})}\errentry{1\!\cdot\!10^{-3}}{3\!\cdot\!10^{-3}}\cr
\tabentry{176}{[0, 1, 0, 3, -1]}{x^5-176}{2^{4}5^{3}11^{4}}{2^{8}5^{6}11^{8}}{2^{5}3^{2}}{6^2}{4\!\cdot\!5^{0}\!+\!3\!\cdot\!5^{1}\!+\!3\!\cdot\!5^{2}\!+\!O(5^{3})}{4\!\cdot\!5^{0}\!+\!3\!\cdot\!5^{1}\!+\!2\!\cdot\!5^{2}\!+\!O(5^{4})}\errentry{4\!\cdot\!10^{-10}}{6\!\cdot\!10^{-10}}\cr
\hline
\tabend

\noindent\begin{tabular}{cllll} 
\extravspace
$E=84B1$ & \multicolumn{3}{l}{$y^2=x^3-x^2-x-2$} & \tabname{5}{84B1}\cr
\extravspace
$p=5$ & \multicolumn{3}{l}{$P_5(E/\Q,T)=1-4T+5T^2$}\cr
$/\vQ\alsub{}$   & $N(E)=2^{2}\!\cdot\!3\!\cdot\!7$ & $L^*(E/\Q)=2^{-1}$ & $|E(\Q)|=2$ & $|\sha(E/\Q)|=1^2$ \cr
$/K\alsub{}$      & $N(E,\sigma)=2^{8}\!\cdot\!3^{4}\!\cdot\!5^{6}\!\cdot\!7^{4}$ & $L^*(E/K)=2^{-1}\!\cdot\!3$ & $|E(K)|=2$ & $|\sha(E/K)|=1^2$ \cr
$/\vQ\alsub{2}$ & $\hbox{IV},\>c_{2}=1$ & $/K\alsub{v}\quad(v|2)$ & $\hbox{IV},\>c_v=3$\cr
$/\vQ\alsub{3}$ & $\hbox{I}_1 \hbox{ non-split},\>c_{3}=1$ & $/K\alsub{v}\quad(v|3)$ & $\hbox{I}_1 \hbox{ split},\>c_v=1$\cr
$/\vQ\alsub{7}$ & $\hbox{I}_2 \hbox{ non-split},\>c_{7}=2$ & $/K\alsub{v}\quad(v|7)$ & $\hbox{I}_2 \hbox{ split},\>c_v=2$\cr
\end{tabular}
\tabbegin
\tabentry{2}{[0, -1, 0, -1, -2]}{x^5-2}{2^{4}5^{5}}{2^{8}3^{4}5^{10}7^{4}}{2^{2}}{1^2}{1\!\cdot\!5^{0}\!+\!3\!\cdot\!5^{1}\!+\!4\!\cdot\!5^{2}\!+\!O(5^{3})}{1\!\cdot\!5^{0}\!+\!4\!\cdot\!5^{1}\!+\!1\!\cdot\!5^{2}\!+\!O(5^{3})}\errentry{5\!\cdot\!10^{-10}}{1\!\cdot\!10^{-10}}\cr
\tabentry{126}{[0, -1, 0, -1, -2]}{x^5-126}{2^{4}3^{4}5^{3}7^{4}}{2^{8}3^{8}5^{6}7^{8}}{2\!\cdot\!3^{2}5^{2}}{15^2}{1\!\cdot\!5^{3}\!+\!4\!\cdot\!5^{4}\!+\!3\!\cdot\!5^{5}\!+\!O(5^{6})}{3\!\cdot\!5^{2}\phantom{\,\,\!+\!2\!\cdot\!5^3}\,\!+\!1\!\cdot\!5^{4}\!+\!O(5^{8})}\errentry{1\!\cdot\!10^{-5}}{7\!\cdot\!10^{-5}}\cr
\tabentry{168}{[0, -1, 0, -1, -2]}{x^5-168}{2^{4}3^{4}5^{3}7^{4}}{2^{8}3^{8}5^{6}7^{8}}{2\!\cdot\!3^{2}5^{2}}{15^2}{1\!\cdot\!5^{3}\!+\!4\!\cdot\!5^{4}\!+\!3\!\cdot\!5^{5}\!+\!O(5^{6})}{3\!\cdot\!5^{2}\phantom{\,\,\!+\!2\!\cdot\!5^3}\,\!+\!1\!\cdot\!5^{4}\!+\!O(5^{8})}\errentry{6\!\cdot\!10^{-5}}{3\!\cdot\!10^{-4}}\cr
\hline
\tabend

\noindent\begin{tabular}{cllll} 
\extravspace
$E=17A1$ & \multicolumn{3}{l}{$y^2+xy+y=x^3-x^2-x-14$} & \tabname{7}{17A1}\cr
\extravspace
$p=7$ & \multicolumn{3}{l}{$P_7(E/\Q,T)=1-4T+7T^2$}\cr
$/\vQ\alsub{}$   & $N(E)=17$ & $L^*(E/\Q)=2^{-2}$ & $|E(\Q)|=4$ & $|\sha(E/\Q)|=1^2$ \cr
$/K\alsub{}$      & $N(E,\sigma)=7^{10}\!\cdot\!17^{6}$ & $L^*(E/K)=2^{4}$ & $|E(K)|=4$ & $|\sha(E/K)|=8^2$ \cr
$/\vQ\alsub{17}$ & $\hbox{I}_4 \hbox{ split},\>c_{17}=4$ & $/K\alsub{v}\quad(v|17)$ & $\hbox{I}_4 \hbox{ split},\>c_v=4$\cr
\end{tabular}
\tabbegin
\tabentry{2}{[1, -1, 1, -1, -14]}{x^7-2}{2^{6}7^{7}}{2^{12}7^{14}17^{6}}{2^{6}}{4^2}{5\!\cdot\!7^{0}\!+\!2\!\cdot\!7^{1}\!+\!2\!\cdot\!7^{2}\!+\!O(7^{3})}{5\!\cdot\!7^{0}\!+\!5\!\cdot\!7^{1}\!+\!1\!\cdot\!7^{2}\!+\!O(7^{3})}\errentry{4\!\cdot\!10^{-2}}{2\!\cdot\!10^{-2}}\cr
\tabentry{3}{[1, -1, 1, -1, -14]}{x^7-3}{3^{6}7^{7}}{3^{12}7^{14}17^{6}}{2^{6}7^{2}}{28^2}{5\!\cdot\!7^{2}\!+\!3\!\cdot\!7^{3}\!+\!4\!\cdot\!7^{4}\!+\!O(7^{5})}{2\!\cdot\!7^{2}\!+\!4\!\cdot\!7^{3}\!+\!2\!\cdot\!7^{4}\!+\!O(7^{5})}\errentry{3\!\cdot\!10^{-2}}{2\!\cdot\!10^{-3}}\cr
\hline
\tabend

\noindent\begin{tabular}{cllll} 
\extravspace
$E=19A3$ & \multicolumn{3}{l}{$y^2+y=x^3+x^2+x$} & \tabname{7}{19A3}\cr
\extravspace
$p=7$ & \multicolumn{3}{l}{$P_7(E/\Q,T)=1+T+7T^2$}\cr
$/\vQ\alsub{}$   & $N(E)=19$ & $L^*(E/\Q)=3^{-2}$ & $|E(\Q)|=3$ & $|\sha(E/\Q)|=1^2$ \cr
$/K\alsub{}$      & $N(E,\sigma)=7^{10}\!\cdot\!19^{6}$ & $L^*(E/K)=3^{-2}$ & $|E(K)|=3$ & $|\sha(E/K)|=1^2$ \cr
$/\vQ\alsub{19}$ & $\hbox{I}_1 \hbox{ split},\>c_{19}=1$ & $/K\alsub{v}\quad(v|19)$ & $\hbox{I}_1 \hbox{ split},\>c_v=1$\cr
\end{tabular}
\tabbegin
\tabentry{2}{[0, 1, 1, 1, 0]}{x^7-2}{2^{6}7^{7}}{2^{12}7^{14}19^{6}}{1}{1^2}{1\!\cdot\!7^{0}\!+\!3\!\cdot\!7^{1}\!+\!1\!\cdot\!7^{2}\!+\!O(7^{3})}{1\!\cdot\!7^{0}\phantom{\,\,\!+\!2\!\cdot\!7^1}\,\!+\!1\!\cdot\!7^{2}\!+\!O(7^{3})}\errentry{1\!\cdot\!10^{-2}}{5\!\cdot\!10^{-4}}\cr
\tabentry{3}{[0, 1, 1, 1, 0]}{x^7-3}{3^{6}7^{7}}{3^{12}7^{14}19^{6}}{2^{2}}{2^2}{3\!\cdot\!7^{0}\!+\!4\!\cdot\!7^{1}\!+\!3\!\cdot\!7^{2}\!+\!O(7^{3})}{3\!\cdot\!7^{0}\phantom{\,\,\!+\!2\!\cdot\!7^1}\,\!+\!4\!\cdot\!7^{2}\!+\!O(7^{3})}\errentry{5\!\cdot\!10^{-2}}{3\!\cdot\!10^{-2}}\cr
\hline
\tabend